\documentclass[10pt,notitlepage]{amsart} 
\usepackage{amsmath}
\usepackage{amsthm} 
\usepackage{amsfonts}
\usepackage{scalerel}
\usepackage{amssymb,amscd,epsf,verbatim,mathtools,framed}
\usepackage{mathrsfs}
\usepackage{graphicx} 
\usepackage{latexsym}
\usepackage{lscape}
\usepackage[colorlinks=true]{hyperref}
\hypersetup{colorlinks, citecolor=blue, filecolor=black, linkcolor=purple, urlcolor=violet}
\usepackage{epstopdf}
\usepackage{pdflscape}
 
\usepackage{tikz}
\usepackage[all]{xy}
\usetikzlibrary{calc}
\usetikzlibrary{matrix,arrows,decorations.pathmorphing}
\usepackage{tabularx}
\usepackage{tikz-cd}
\tikzcdset{scale cd/.style={every label/.append style={scale=#1},
    cells={nodes={scale=#1}}}}
\usepackage{color}
\usepackage{multirow}  
\usepackage{scalefnt}
\usepackage{fancyhdr}
\usepackage[margin=1in]{geometry}
\usepackage[T1]{fontenc}
\newcommand{\hooklongrightarrow}{\lhook\joinrel\longrightarrow}

\usepackage{relsize} 
\usepackage[bbgreekl]{mathbbol} 
\DeclareSymbolFontAlphabet{\mathbb}{AMSb} %to ensure that the meaning of \mathbb does not change
\DeclareSymbolFontAlphabet{\mathbbl}{bbold}

\newcommand{\Z}{\mathbb{Z}}
\newcommand{\D}{\mathbb{D}}

\newcommand{\R}{\mathbb{R}}
\newcommand{\Q}{\mathbb{Q}}
\newcommand{\A}{\mathbb{A}}
\renewcommand{\L}{\mathbb{L}}
\renewcommand{\P}{\mathbb{P}}
\newcommand{\C}{\mathbb{C}}

\newcommand{\mD}{\mathcal{D}}
\newcommand{\mE}{\mathcal{E}}

\newcommand{\mO}{\mathcal{O}}

\newcommand{\sH}{\mathscr{H}}
\newcommand{\sK}{\mathscr{K}}
\newcommand{\fX}{\mathfrak{X}}

\newcommand{\fm}{\mathfrak{m}}

\newcommand{\fA}{\mathbf{A}}

\newcommand{\fC}{\mathbf{C}}
\newcommand{\fN}{\mathbf{N}}

\newcommand{\cl}{\overline}

\renewcommand{\le}{\leqslant}
\renewcommand{\ge}{\geqslant}
\newcommand{\ra}{\rightarrow}

\newcommand{\sq}{\widetilde}
\newcommand{\minus}{\backslash}

\DeclareMathOperator{\End}{End}

\DeclareMathOperator{\GL}{GL}

\DeclareMathOperator{\lra}{\: \longrightarrow \:} 
\DeclareMathOperator{\isom}{\;\xrightarrow{\: {}_{\sim} \:} \;} 
\DeclareMathOperator{\Gal}{Gal}

\newcommand{\gr}[1]{\langle {#1} \rangle} 
\DeclareMathOperator{\spec}{Spec}
\DeclareMathOperator{\spf}{Spf}
\DeclareMathOperator{\spa}{Spa}

\DeclareMathOperator{\ett}{\mathrm{\acute{e}t}}

\DeclareMathOperator{\proet}{\mathrm{pro\acute{e}t}}

\theoremstyle{theorem}
\newtheorem{theorem}{Theorem}[section] 
\newtheorem{mainthm}{Theorem}

\newcommand{\tu}{\mathrm} 

\theoremstyle{definition}
\newtheorem{mainquestion}[mainthm]{Question}
\newtheorem{mainrmk}[mainthm]{Remark}
\newtheorem{maindef}[mainthm]{Definition}
\newtheorem{example}[theorem]{Example}
\newtheorem{definition}[theorem]{Definition}
\newtheorem{proposition}[theorem]{Proposition}

\newtheorem{lemma}[theorem]{Lemma}
\newtheorem{corollary}[theorem]{Corollary}
\newtheorem{construction}[theorem]{Construction} 

\newtheorem{notation}[theorem]{Notation}
\newtheorem{remark}[theorem]{Remark}

\makeatother
\usepackage[T1]{fontenc}
\usepackage{lmodern}
\usetikzlibrary{arrows.meta,calc,patterns}

\definecolor{sectorink}{HTML}{B6D0E2}
\definecolor{sectorboundary}{HTML}{4169E1}%{0818A8}%{4169E1}
\definecolor{sectoraxis}{HTML}{120802}
\definecolor{sectorguide}{HTML}{8B642C}
\definecolor{sectorred}{HTML}{8B445C}
\definecolor{orange}{HTML}{FFDBBB}

\definecolor{ink}{HTML}{233240}
\definecolor{ann}{HTML}{16816F}
\definecolor{disk}{HTML}{120802}
\definecolor{point}{HTML}{B86120}
\definecolor{excluded}{HTML}{120802}
\definecolor{level}{HTML}{8B642C}
\tikzset{
 every node/.style={font=\large,text=ink},
 tree/.style={line width=1.2pt,draw=excluded},
 inside/.style={line width=3pt,draw=ann},
 branch/.style={line width=2.7pt,draw=disk},
 open/.style={circle,draw=ann,line width=1.2pt,fill=white,
              inner sep=0pt,minimum size=6pt},
 dot/.style={circle,fill=ink,inner sep=0pt,minimum size=5pt},
 xdot/.style={circle,fill=point,inner sep=0pt,minimum size=7pt}
}

\title{A $p$-adic monodromy theorem for curves}

\author{Hansheng Diao, Yong Suk Moon, and Zijian Yao}

\date{}

\numberwithin{equation}{section}
\begin{document}

\begin{abstract}
We prove that every de Rham $p$-adic local system on a smooth projective curve over a $p$-adic field is potentially semistable; that is, it becomes semistable after pulling back along a finite cover of the curve. This establishes a relative version of the classical $p$-adic monodromy theorem of Berger and Andr\'e--Kedlaya--Mebkhout. Along the way, we show that every $p$-adic differential equation near a type I\!V point on a curve becomes trivial after a finite \'etale extension. 
\end{abstract}

\maketitle
\thispagestyle{empty}
\tableofcontents

\section{Introduction}\label{section: introduction}
% 1.1
\subsection{Main result} 
\noindent
\vspace{0.1in}

\noindent Let $p$ be a prime number and let $K$ be a $p$-adic field, that is, a complete discrete valuation field of mixed characteristic $(0,p)$ with perfect residue field. Let $\Gal_K=\Gal(\overline{K}/K)$ denote the absolute Galois group.
The main goal of this article is to establish the following relative version of the $p$-adic monodromy theorem. 

\begin{mainthm}\label{thm:main_intro}
Let $X$ be a smooth projective curve over $K$ and let $\L$ be an \'etale $\Z_p$-local system on $X$. If $\L$ is \emph{de Rham}, then $\L$ is \emph{potentially semistable}. That is, after replacing the base field $K$ by a finite extension, there exists a finite cover $f\colon X' \ra X$ such that $X'$ is a smooth projective curve with semistable reduction over $\mathcal{O}_K$ and $f^* \L$ is \emph{semistable}. 
\end{mainthm}
 
Before explaining the terminology, let us place our results in context. First consider a smooth projective curve $\cl S$ over $\C$ and let $S = \cl S - \{t_1, \ldots, t_r\}$ be a Zariski open subset. Let $\pi\colon Y \ra S$ be a smooth proper map and denote by $Y_s$ the fiber over a point $s \in S$. For a fixed integer $i \ge 0$, the singular cohomology groups $H^i(Y_s, \C)$ form a local system on $S$. The ``local monodromy theorem'' of Borel--Griffiths--Grothendieck--Katz--Landman \cite{Griffiths_period_Rice, Griffiths_period_Bulletin, Landman, SGA7, Katz_differential_equations} asserts that, if we restrict this local system over a small punctured disk $\Delta^*$ around one of the ``missing points'' $t_i$ and pick a base point $s_0 \in \Delta^*$, then the monodromy action of $\pi_1 (\Delta^*, s_0)$ on  $H^i(Y_{s_0}, \C)$ is \emph{quasi-unipotent}. More precisely, the automorphism $T$ on $H^i(Y_{s_0}, \C)$ induced by a generator of $\pi_1(\Delta^*, s_0)$ admits a Jordan decomposition $T = U\cdot D = D \cdot U $, where $U$ is unipotent and $D$ is semisimple of finite order. This result plays a fundamental role in the study of degeneration of Hodge structures (for example, see \cite{Schmid, Steenbrink, Clemens, CKS}). Moreover, it  implies, rather strikingly, that all the eigenvalues of $T$ are roots of unity.  

Grothendieck's proof \cite{SGA7} is best explained from the arithmetic perspective. In fact, his argument is robust enough to imply the following. Let $X$ be a smooth proper variety over a finite extension $K/\Q_p$ and let $\L$ be an $l$-adic \'etale local system on $X$ for a prime $l \ne p$. Then there exists a finite \'etale cover $f\colon X' \ra X$, such that for any finite extension $K'/K$ and any $K'$-point $s \in X'(K')$, the $l$-adic representation of $\Gal_{K'}$ given by the pullback $f^* \L |_{s}$ to the point $s$ is unipotent when restricted to the inertia subgroup $I_{K'}$. Since the property of $I_K$ acting unipotently on $H^i_{\ett} (Y_{\cl K}, \Q_l)$ for some smooth projective variety $Y$ over $K$ is intimately related to $Y$ acquiring \emph{semistable} reduction (for example, by Grothendieck's semistable reduction criterion for abelian varieties; also see \cite{Rapoport_Zink}), Fontaine conjectured a $p$-adic analogue of the local monodromy theorem, which states that a $p$-adic representation $V$ of $\Gal_K$ is de Rham if and only if it is potentially semistable (see \cite{Fontaine_94}, and also Colmez's Bourbaki survey \cite{Colmez_Bourbaki}); that is,  $V|_{\Gal_{K'}}$ is semistable for some finite extension $K'/K$. The conjecture of Fontaine was proved by Berger and Andr\'e--Kedlaya--Mebkhout  \cite{Berger-differential, Andre_p-adic_monodormy, Kedlaya_p-adic_monodormy, Mebkhout_p-adic_monodormy} and later revisited by Colmez and others (see \cite{Colmez08, Fargues-Fontaine, ABLBRCS}). 

Since then, it has remained an important open problem whether the $p$-adic local monodromy theorem admits a generalization to the relative geometric setting (i.e., the \emph{relative $p$-adic monodromy conjecture}). Our main result provides such a generalization for smooth projective curves.  

\begin{mainrmk}
The relative $p$-adic monodromy conjecture predicts that ``de Rham local systems are potentially semistable''. But what should ``potential'' really mean in the relative setting? Lawrence and Li \cite{Lawrence_Li} construct an example which shows that simply extending the base field is not enough. On the other hand, in \cite{Liu_Zhu}, Liu and Zhu make the following conjecture: a de Rham local system on $X$ should become semistable at every classical point over a \emph{finite \'etale} cover of $X$. From this perspective, our main result proves a slightly weaker version of the ultimate hope of Liu and Zhu.
\end{mainrmk}
 
From a more geometric perspective, our result suggests the following question. 
\begin{mainquestion}
 For a smooth proper map $Y \ra X$ of projective varieties over $K$, up to replacing the base field $K$ by a finite extension, does there exist a projective alteration $X' \ra X$, that is, a proper generically finite \'etale cover with $X'$ projective, such that $Y \times_{X} X' \ra X'$ has (fiber-wise) semistable reduction? 
\end{mainquestion} The affirmative answer to this question would provide a mixed characteristic analogue of some well-known conjectures on semistable reductions in algebraic geometry, notably the semistable reduction conjecture of Abramovich--Karu  \cite[Conjecture~0.2]{Abramovich_Karu} (also see the introduction of \cite{KKMS} and \cite{ALT}). This mixed characteristic 
analogue has a positive answer when $Y$ is a family of abelian varieties over $X$, by considering the relative $l$-adic Tate module on $X$ for a prime $l \ne p$. 
When $X$ is a point, a related question is  
\begin{mainquestion}
Does a smooth projective variety $Y$ over $K$ always admit semistable reduction after a finite \'etale extension?  
\end{mainquestion}
It is a folklore conjecture that this latter question has a positive answer.

\vspace{0.1in}
%% 1.2
\subsection{$p$-adic local systems}\label{ss: p-adic local systems}
\noindent
\vspace{0.1in}

\noindent Now we recall the notions of de Rham and semistable $p$-adic local systems. Let $X$ be a smooth rigid analytic variety over $K$, viewed as an adic space over $\spa(K, \mathcal{O}_K)$. In this article, by a \emph{$p$-adic local system} on $X$, we shall exclusively mean a lisse $\Z_p$-sheaf on the \'etale site $X_{\ett}$. Intuitively, one may view a $p$-adic local system on $X$ as a geometric family of $p$-adic Galois representations parameterized by $X$. Indeed, for any classical point $x\in X$ with residue field $\kappa(x)$, the restriction of the local system at $x$ gives rise to (a $\Z_p$-lattice inside) a $p$-adic representation of $\Gal_{\kappa(x)}$.

The notion of \emph{de Rham} local systems is defined in \cite{Scholze_p_adic_Hodge} (also see \cite{Faltings_almost, Brinon, Liu_Zhu}).\footnote{In fact, Scholze \cite{Scholze_p_adic_Hodge} defines the notion of de Rham $p$-adic local systems in a slightly more general context. In particular, one can make sense of this notion even for $\Q_p$-local systems.} 

\begin{maindef}[{\cite[Definition~8.3]{Scholze_p_adic_Hodge}}]
Let $X$ be a smooth rigid analytic variety over $K$. A $p$-adic local system $\L$ on $X$ is called \emph{de Rham} if there exists a filtered vector bundle with integrable connection $(\mD, \nabla, \mathrm{Fil}^\bullet)$ on $X$ in the sense of \cite[Definition~7.4]{Scholze_p_adic_Hodge}, together with an isomorphism   
\[
\widehat{\L} \otimes_{\widehat{\Z}_p} \mO \mathbb B_{\mathrm{dR}, X} \cong \mD \otimes_{\mO_X}  \mO \mathbb B_{\mathrm{dR}, X} 
\]
of pro-\'etale sheaves on $X$, compatible with filtrations and connections. Here $\widehat{\L}$ denotes the lisse $\widehat{\Z}_p$-sheaf on $X_{\proet}$ associated with $\L$. The readers are referred to \S 6 of \emph{op.\ cit.} for the construction of the de Rham period sheaf $\mO \mathbb B_{\mathrm{dR}, X}$ as well as the filtration and connection on it. 
\end{maindef} 
  
On the other hand, for the notion of \emph{semistable} local systems, developed in various generalities in  \cite{Faltings_almost, Brinon, Tsuji_notes, Andreatta_Iovita_st,  DLMS2, GY}, etc., one typically requires a semistable integral model of $X$ to start with. The following definition is borrowed from \cite{GY}. 

\begin{maindef}[{\cite[Definition~4.4]{GY}}] \label{def:st_local_sys} Let $\fX$ be a $p$-adic formal scheme over $\mO_K$ with semistable reduction and let $X$ be its adic generic fiber. Let $\fX_s$ denote the special fiber of $\fX$ and let $\fX_s^{\log}$ denote the log scheme $(\fX_s, \alpha_s)$ where the log structure $\alpha_s$ is obtained from the pullback of the divisorial log structure  
\[\alpha\colon \mO_{\fX, \ett} \cap (\mO_{\fX, \ett}[1/p])^{\times} \ra \mO_{\fX, \ett} \]  on $\fX$. Let $\L$ be a $p$-adic local system on $X$. We say that $\L$ is \emph{semistable} if there exists a log $F$-isocrystal\footnote{By \emph{log $F$-isocrystals}, we mean \emph{finite locally free} isocrystals over the log crystalline site, equipped with a Frobenius structure.} $(\mE, \varphi)$ over $\fX_s^{\log}$, together with a Frobenius equivariant isomorphism of $\mathbb{B}_{\mathrm{crys}}$-vector bundles
\[ 
\widehat{\L}\otimes_{\widehat{\Z}_p} \mathbb{B}_{\mathrm{crys}}\cong\mathbb{B}_{\mathrm{crys}} (\mE).
\] 
Here $\mathbb{B}_{\mathrm{crys}}$ is the horizontal crystalline period sheaf and  $\mathbb{B}_{\mathrm{crys}}(\mE)$ is the sheaf of $\mathbb{B}_{\mathrm{crys}}$-modules on the pro-\'etale site $X_{\proet}$ associated with $(\mE, \varphi)$. For the constructions of $\mathbb{B}_{\mathrm{crys}}$ and $\mathbb{B}_{\mathrm{crys}} (\mE)$, see \cite[Definition~3.8, Construction~4.1]{GY}.
\end{maindef} 

\begin{mainrmk}
There is also a notion of \emph{crystalline} local systems. For this, we further require the $p$-adic formal scheme $\fX$ to be smooth. Then crystalline local systems can be defined in the same way as in Definition~\ref{def:st_local_sys}, simply omitting all ``log'' therein (see \cite[Definition~4.4]{GY}).
\end{mainrmk}

\begin{mainrmk}[Rigidity] 
A remarkable theorem of Liu--Zhu \cite{Liu_Zhu} asserts that de Rham local systems are rigid. That is, if a $p$-adic local system on a smooth geometrically connected rigid analytic variety is de Rham at one classical point, then it is already a de Rham local system, and hence is de Rham at every classical point. Although the most naive analogues of this rigidity phenomenon fail for semistable and crystalline local systems, it is natural to ask to what extent analogous rigidity results remain valid in these settings. In the crystalline case, some of the authors (H.D. and Z.Y.) recently established the following rigidity result.
\end{mainrmk}

\begin{mainthm}[{\cite[Theorem~1.1]{DY-rigidity}}]\label{thm:rigidity_intro}
Let $\fX$ be a geometrically connected smooth projective scheme over $\mO_K$ and let $X$ be its generic fiber. Let $\L$ be a $p$-adic local system on $X$. If there exists one classical point $x_0$ on $X$ such that $\L|_{x_0}$ is potentially crystalline, then $\L$ is potentially crystalline at every classical point on $X$. 
\end{mainthm}

It is worth pointing out that the proof of Theorem~\ref{thm:rigidity_intro}   crucially relies on the main result of this article. 
We conclude this review section with two recent results on the semistability of $p$-adic local systems.

\begin{mainthm}[{Pointwise criterion \cite[Theorem~1.1]{GY}}] \label{thm:pointwise_criterion}
Let $\fX$ be a $p$-adic formal scheme over $\mO_K$ with semistable reduction and let $X$ be its adic generic fiber. Let $\L$ be a $p$-adic local system on $X$. If the restriction of $\L$ to \emph{every} classical point on $X$ is semistable, then $\L$ is semistable as a local system. 
\end{mainthm}

\begin{mainthm}[{Purity theorem \cite[Theorem~1.1]{DLMS2}}]\label{thm: purity DLMS}
Let $\fX$ be a $p$-adic formal scheme over $\mO_K$ with semistable reduction and let $X$ be its adic generic fiber. Let $\L$ be a $p$-adic local system on $X$. Then $\L$ is semistable if and only if its restriction to each Shilov point\footnote{We refer the reader to \cite{DLMS2} for the definition of Shilov points.} is semistable.
\end{mainthm}

\noindent  Both Theorem~\ref{thm:pointwise_criterion} and Theorem~\ref{thm: purity DLMS} play important roles in our proof of Theorem~\ref{thm:main_intro}.

\vspace{0.1in}
%%% 1.3
\subsection{Related results}
\noindent 
\vspace{0.1in}

\noindent 
We make a few more remarks on some relevant works in the literature. There have been numerous attempts to generalize the $p$-adic monodromy theorem to the relative settings; for example, see \cite{Shimizu_2, kedlaya-monodromy-in-families, Du}\footnote{The proof of \cite[Theorem~3.2.6]{kedlaya-monodromy-in-families} relies on working with a ``relative connection'' in order to apply the techniques developed in \cite{kedlaya_book} and \cite[\S 1--2]{kedlaya-monodromy-in-families}. More precisely, the argument requires working with a relative annulus of the form $A \times_K X'$ over a base $X'$ together with a vector bundle $\mathcal{E}$ equipped with an \emph{$\mathcal{O}_{X'}$-linear} connection. However, one cannot expect such an $\mathcal{O}_{X'}$-linear connection to exist, because the arithmetic and geometric directions are inseparably intertwined. This appears to leave a serious gap in the argument in \cite[Theorem~3.2.6]{kedlaya-monodromy-in-families}.} 
and the introductions thereof. Most of these works are of local nature. Notably, Shimizu  \cite{Shimizu_2} shows that for any smooth rigid analytic variety $X$ over $K$, any de Rham local system $\L$ on $X$, and any classical point $x_0 \in X$, there exists an affinoid open neighborhood $U$ of $x_0$ and a finite extension $K'/K$ such that $\L$ becomes semistable at every classical point in $U_{K'}$. In a more recent preprint, Du \cite{Du} shows that there exist a dense open subset $U \subset X$ and a (generally non-finite)  \'etale cover $V \ra U$ on which $\L$ becomes semistable at all classical points by studying  Newton polygons. In contrast, the  main result of this article provides the first general $p$-adic monodromy theorem in the global setting. It is certainly worth mentioning that, in the global setting, Shimizu \cite{Shimizu_1} proves a similar result to our main theorem when the local system $\L$ has a single Hodge--Tate weight for a general class of analytic varieties.  Let us also mention the work of Ohkubo \cite{Ohkubo-p-adic-monodromy} and Morita \cite{Morita-p-adic-monodromy}, where they extend the classical $p$-adic monodromy theorem to more general $p$-adic fields; namely, they treat complete discretely valued fields with imperfect residue fields. For type I\!I\!I points, Kedlaya \cite{kedlaya_fake,kedlaya_st3} obtains several closely related results on $p$-adic differential equations over certain ``fake Robba rings'' (although a reasonable theory of relative $p$-adic differential equation was missing at that time, see \S \ref{ss:relative_p_adic_Hodge}).  In fact, the main results of \cite{Ohkubo-p-adic-monodromy, Morita-p-adic-monodromy, kedlaya_st3} are used as ingredients in our approach when we analyze $p$-adic local systems near  points of type I\!I and I\!I\!I.  
In a different direction, Berger and Colmez \cite{Berger_Colmez_S_p_adic_monodromy} prove the following result for \emph{arithmetic families} of $p$-adic Galois representations: consider a reduced affinoid algebra $S$ over $\Q_p$ and an $S$-family of representations of $\Gal_K$, that is, a finite locally free $S$-module $V$ equipped with a continuous $S$-linear action of $\Gal_K$. If for each point $\fm_x \in \mathrm{Spm} (S)$ in the maximal spectrum of $S$, the representation $V_x = V \otimes_{S} S/\fm_x$ of $\Gal_K$ is de Rham with Hodge--Tate weights contained in a fixed interval $[a,b]$, then there exists a finite extension $K'/K$ such that each $V_x$ is a semistable representation of $\Gal_{K'}.$ 
It seems an interesting problem to combine the approach in this article with tools developed in \cite{Berger_Colmez_S_p_adic_monodromy} to establish a $p$-adic monodromy result for $S$-valued local systems.

\vspace{0.1in}
%%% 1.4
\subsection{Relative $p$-adic Hodge theory} \label{ss:relative_p_adic_Hodge}
\noindent 
\vspace{0.1in}

\noindent 
Before explaining the proof of our main result, it may be instructive to first recall the proof of Berger \cite{Berger-differential} in the classical setting. By the works of Fontaine \cite{Fontaine_90} and Cherbonnier--Colmez \cite{Cherbonnier_Colmez}, one can attach  a $(\varphi, \Gamma_K)$-module $D^\dagger_K(V)$ over the \emph{overconvergent} imperfect period ring $\mathbf{B}^\dagger_K$ to a $p$-adic $\Gal_K$-representation $V$. Assuming that $V$ is in addition de Rham, Berger \cite{Berger-differential} then constructs a $p$-adic differential equation $(\mathbf{N}_{\mathrm{dR}}(V), \partial_V)$ over the Robba ring $\mathbf{C}^\dagger_{K}$,\footnote{The period ring $\mathbf{B}^\dagger_{\mathrm{rig}, K}$ from Berger's article \cite{Berger-differential} will be denoted as $\mathbf{C}^\dagger_{K}$ in this article. Also see Appendix \ref{appendix:notation}.} equipped with a suitable Frobenius structure.  This is achieved by systematically removing the singularities of the differential operator $\partial_V$ on the base change $D^{\dagger}_{\mathrm{rig}, K}(V)=D^\dagger_{K}(V)\otimes_{\mathbf{B}^{\dagger}_K}\mathbf{C}^{\dagger}_K$. By the result of Andr\'e--Kedlaya--Mebkhout (Theorem~\ref{thm:all_pDE_with_Frob_are_quasiunipotent}), such a $p$-adic differential equation is  \emph{quasi-unipotent}, meaning that  there exists a finite extension $L/K$ such that $\partial_V$ admits a full basis of solutions over $\mathbf{C}^\dagger_{L} [\log \pi]$ --- this is the Robba ring $\mathbf{C}^{\dagger}_L$ adjoining the variable $\log \pi$ (see \S \ref{sss:C log}). Berger then shows that, up to a further finite extension of the base field, the semistable periods of $V$ contain these horizontal sections, thus settling $p$-adic monodromy theorem in the classical setting. We refer the reader to \S \ref{sec:classical_phi_Gamma} for more details. 

The first difficulty we encounter in extending this strategy to the relative setting concerns the foundations of relative $p$-adic Hodge theory. We would need a relative version of $(\varphi, \Gamma)$-modules, where $\Gamma$ is a $p$-adic Lie group of the form $\Gamma\cong \Gamma_{\mathrm{geom}}\rtimes \Gamma_K$ (here the ``geometric part'' $\Gamma_{\mathrm{geom}}$ is isomorphic to $\Z_p^{d}$ where $d$ is the dimension of $X$), as well as a relative analogue of the main result of \cite{Cherbonnier_Colmez}. In order to achieve this, one first needs to construct suitable \emph{imperfect relative period rings} $\mathbf{A}^{\dagger}_{X}, \mathbf{B}^{\dagger}_{X}, \mathbf{C}_{X}^{\dagger}$ inside their perfect counterparts $\widetilde{\mathbf{A}}^{\dagger}_{X}, \widetilde{\mathbf{B}}^{\dagger}_{X}, \widetilde{\mathbf{C}}_{X}^{\dagger},$\footnote{These \emph{perfect relative period rings} are constructed by Kedlaya and Liu \cite{Kedlaya-LiuI} using the theory of perfectoid spaces.} and then study decompletion of relative $(\varphi, \Gamma)$-modules. %from  perfect period rings to imperfect ones. 
When this article is written, the only reference in the literature that treats such imperfect period rings (and the decompletion theory) in the generality we need is due to Kedlaya--Liu \cite{Kedlaya-LiuII}. In \emph{loc.\ cit.}, they initiate the study of imperfect period rings relative to a \emph{(weakly) decompleting tower} of finite covers of the base $X$. However, to our best knowledge, some of the assertions and proofs in \cite[in particular, \S 5]{Kedlaya-LiuII} concerning (weakly) decompleting towers appear to be incomplete and require further justification.  
In particular, the general descent of relative $(\varphi, \Gamma)$-modules claimed in \cite{Kedlaya-LiuII} remains difficult to verify, as the proposed argument therein seems to omit several essential details.  

In order to have a complete and coherent theory to work with, we develop a foundation of relative $p$-adic Hodge theory in a separate article \cite{relative}, where we systematically study imperfect relative period rings and $(\varphi, \Gamma)$-modules over them. The treatment in \cite{relative} is inspired by \cite{Kedlaya-LiuII} but differs substantially from \emph{loc.\ cit.}, especially in how imperfect period rings are treated. A central notion we introduce in \cite{relative} is that of \emph{decompleting data} (Definition~\ref{defn: weakly decompleting datum}, \ref{definition:decompleting_datum}), which is engineered to capture the idea of decompleting towers from \cite[Definition~5.2.3]{Kedlaya-LiuII}.  This allows us to obtain a relative analogue of the Fontaine--Cherbonnier--Colmez equivalence. 

\begin{mainthm}[{\cite[Theorem~1]{relative}}] \label{intro:theorem_equiv_of_categories}
   Let $\mathbf{B}_{\Lambda}^{\dagger}$ and $\mathbf{C}_{\Lambda}^{\dagger}$ be the imperfect relative period rings for $X$ associated to a decompleting data $\Lambda$. Then there are natural equivalences 
\[
D^{\dagger}\colon \mathrm{Loc}_{\Z_p} (X)_{\Q} \isom \mathrm{Mod}^{(\varphi, \Gamma), \ett}_{/\mathbf{B}_{\Lambda}^{\dagger}} 
\]
\[
D^{\dagger}_{\mathrm{rig}}\colon \mathrm{Loc}_{\Z_p} (X)_{\Q} \isom \mathrm{Mod}^{(\varphi, \Gamma), \ett}_{/\mathbf{C}_{\Lambda}^{\dagger}} 
\]
from the isogeny category of \'etale $\Z_p$-local systems on $X$ to the category of finite projective \'etale $(\varphi, \Gamma)$-modules over $\mathbf{B}_{\Lambda}^{\dagger}$ (resp. $\mathbf{C}_{\Lambda}^{\dagger}$). 
\end{mainthm}

To generalize the strategy of Berger to the relative setting, another missing ingredient is a relative analogue of $p$-adic differential equations. With this objective in mind, we introduce the notion of \emph{strongly decompleting data}   (Definition~\ref{definition:strongly_decompleting_datum}) in \cite{relative}, which axiomatizes the process of ``removing singularities'' in \cite{Berger-differential}. This eventually allows us to construct \emph{relative $p$-adic differential equations}. 

\begin{mainthm}[{\cite[Theorem~2]{relative}}, crude form]  \label{intro:theorem_relative_Berger}
Assume that $\Lambda$ is in addition strongly decompleting, and let $\mathbf{C}_{\Lambda}^{\dagger}$ be the associated imperfect relative period ring. Then there is a natural functor  
\[
\mathbf{N}_{\mathrm{dR}}\colon \mathrm{Loc}^{\mathrm{dR}}_{\Z_p}(X)_{\Q} \lra \mathrm{Mod}^{(\varphi, \Gamma)}_{/\mathbf{C}^{\dagger}_{\Lambda}}
\]
from the isogeny category of de Rham $\Z_p$-local systems to the category of finite projective $(\varphi, \Gamma)$-modules over $\mathbf{C}^{\dagger}_{\Lambda}$, where $\mathbf{N}_{\mathrm{dR}} (\L)$ carries not only the arithmetic differential operator arising from the cyclotomic direction, but also geometric differential operators describing the variation of the local system $\L$ over $X$. 
\end{mainthm} 

For the reader's convenience, we review the precise definitions of decompleting and strongly decompleting data  in \S \ref{section:relative_period_rings}, as well as some explicit examples of such decompleting and strongly decompleting data (see \S \ref{ss: examples of strongly decompleting data} and \S \ref{ss: annulus--pDE}). These explicit examples turn out to be sufficient for the application towards the $p$-adic monodromy theorem.

\vspace{0.1in}
%%% 1.5
\subsection{Reduction of Theorem~\ref{thm:main_intro} to a local statement}\label{subsection: outline Thm 1.1}
\noindent 
\vspace{0.1in}

\noindent Finally, we explain the proof of our main result, which relies on careful analysis of (various forms of) $p$-adic differential equations, Theorem~\ref{thm:pointwise_criterion} and \ref{thm: purity DLMS}, together with several delicate arguments in $p$-adic geometry. A key step is a spread-out argument (Theorem~\ref{thm: spread-out thm}), which allows us to reduce Theorem~\ref{thm:main_intro} to the following local statement.

\begin{mainthm}[{Theorem~\ref{thm: local main thm}}, crude form]\label{thm:local theorem_intro} Let $X$ and $\L$ be as in Theorem~\ref{thm:main_intro} and let $x \in X$ be \emph{any} adic point. Then, up to replacing $K$ by a finite extension, there exists an affinoid neighborhood $U$ of $x$ in $X$, and a finite \'etale cover $g\colon V \ra U$, such that the pullback $g^* (\L|_U)$ is semistable at every classical point on $V$. 
\end{mainthm}

Granting Theorem~\ref{thm:local theorem_intro}, we complete the proof of Theorem~\ref{thm:main_intro} as follows. By compactness, $X$ admits a cover by finitely many $U$'s as in Theorem~\ref{thm:local theorem_intro}.\footnote{
Let us emphasize that it is crucial to consider \emph{all} adic points on $X$. In particular, it is not enough to only consider \emph{classical points} arising from algebraic geometry (which are certain type I points on the curve), because a union of affinoid neighborhoods taken over all classical points need not cover the entire $X$.} 
Our task is to ``glue'' these local pieces into a global cover of $X$. To achieve this, we spread-out each finite \'etale cover $V \ra U$ to a  global finite cover over $X$ and take the fiber product of all of them to form a single finite cover over $X$. This spread-out step is achieved using certain $p$-adic approximation techniques and is closely related to a result of Raynaud \cite{Raynaud}, often viewed as the $p$-adic analogue of Runge's theorem on holomorphic functions on complex Riemann surfaces. 

For Theorem~\ref{thm:local theorem_intro}, our strategy is to perform a case-by-case study based on the type of the point $x$.\footnote{See \S \ref{section: prelim} for a review of the classification of points of type I--V on adic curves.} As mentioned earlier, the case of type I points is settled by \cite{Shimizu_2}; while the case of type I\!I points follows from \cite{Ohkubo-p-adic-monodromy, Morita-p-adic-monodromy} by an argument similar to that of \cite[\S 6--7]{GY}. The remaining cases are points of type  I\!I\!I, I\!V, and V. For each of these types, we will need a different argument, relying on a systematic study of various forms of $p$-adic differential equations, which we will briefly summarize in the remainder of this introduction.

\vspace{0.1in}
%%% 1.6
\subsection{Type I\!I\!I points}\label{subsection: outline type3}
\noindent 
\vspace{0.1in}

\noindent 
By enlarging the base field $K$ if necessary, a type I\!I\!I point $x \in X$ admits a neighborhood isomorphic to a closed annulus inside the standard unit disk, where $x$ corresponds to a circle of radius $\sigma$ with $\sigma\notin p^\Q$. In this situation, a system of open neighborhoods of $x$ can be chosen to consist of closed annuli $A_K ([\sigma_1, \sigma_2]) = \{\sigma_1 \le |T| \le \sigma_2\}$ with $\sigma_1 < \sigma < \sigma_2$ and $\sigma_1, \sigma_2 \in p^\Q$. Our strategy is to consider the imperfect relative period rings $\mathbf{C}^\dagger_{[\sigma_1, \sigma_2]}$ associated with these annuli $A_K ([\sigma_1, \sigma_2])$ and study the $p$-adic differential equations $\mathbf{N}_{\mathrm{dR}}(\L|_{A_K ([\sigma_1, \sigma_2])})$ over $\mathbf{C}^\dagger_{[\sigma_1, \sigma_2]}$ (cf. Theorem~\ref{intro:theorem_relative_Berger}). Thanks to the explicit description of $\mathbf{C}^\dagger_{[\sigma_1, \sigma_2]}$ in \cite{relative}, a suitable completed colimit of these period rings as $\sigma_1, \sigma_2 \ra \sigma$ is related to the \emph{fake Robba rings} studied by Kedlaya in \cite{kedlaya_fake} and \cite{kedlaya_st3}. Using a form of $p$-adic local monodromy theorem proved in \emph{loc.\ cit.} for these fake Robba rings, we are able to achieve unipotence of the differential operators ``at the point $x$'' after pulling back to a finite \'etale cover of a sufficiently ``thin'' closed annulus. The idea is then to ``propagate'' this unipotence from the point $x$ to nearby points. This propagation step relies crucially on the explicit descriptions for various relative period rings associated with closed annuli. In particular, via a \emph{monomial change of coordinates}, we show that $\mathbf{C}^\dagger_{[\sigma_1, \sigma_2]}$ behaves like the ring of analytic functions on a product of two closed annuli, after passing to suitable colimits.\footnote{In fact, we further show that this behavior persists when the annulus is replaced by a certain class of finite \'etale covers, which suffices for our applications.} Given this, we can apply the multivariable $p$-adic Fuchs theorem and conclude the desired unipotence. Finally, as in the classical case (cf. \cite{Berger-differential}; also see \S \ref{ss:classical_pDE}), the unipotence of the differential operators allows us to produce enough semistable periods in a neighborhood of $x$ after passing to a finite \'etale cover if necessary.

\vspace{0.1in}
%%% 1.7
\subsection{Type I\!V points}\label{subsection: outline type3}
\noindent 
\vspace{0.1in}

\noindent 
Type I\!V points  are perhaps the most mysterious points  on curves. They correspond to  nested sequences of closed disks $\D_1 \supset \D_2 \supset \cdots $ with empty intersection, that is, $\cap_{i \ge 0} \D_i$ contains no type I point. As the disks do not share the same center, it is more difficult to analyze the relative $p$-adic differential equations $\mathbf{N}_{\tu{dR}} (\L |_{\D_i})$ systematically in a compatible way. The idea is to prove the following geometric result on vector bundles with connection near a type I\!V point.

\begin{mainthm} [Theorem~\ref{thm:pDE_near_type_4_trivial}, crude form]
 \label{mainthm:intro_pde}
    Let $x \in X$ be a  type I\!V point and let $(M, \nabla)$ be a vector bundle with connection on a neighborhood of $x$. Then after shrinking the neighborhood and passing to a finite \'etale cover, $M$ becomes  a trivial connection near $x$.   
\end{mainthm} 

In other words, Theorem~\ref{mainthm:intro_pde} says that all $p$-adic differential equations near a type I\!V point can be solved locally after passing to a finite \'etale cover. A key observation is that, if a vector bundle with connection near a type I\!V point
with trivial determinant is irreducible of rank $n$, then either $p|n$ or $n=1$ (cf. Proposition~\ref{prop:type_4_key}). Consequently, we deduce that if a connection near a type I\!V point remains irreducible after passing to all possible finite \'etale covers, then it must be trivial of rank 1. This eventually implies  Theorem~\ref{mainthm:intro_pde}. Applying this result directly to $D_{\mathrm{dR}}(\L)$, we deduce that $\L$ is \emph{horizontal de Rham} in the sense of \cite{Shimizu_2}, after pulling back along a finite \'etale cover. Then the main result of \emph{loc.\ cit.} yields the following local version of the $p$-adic monodromy theorem around type I\!V points. 

\begin{mainthm}[Theorem~\ref{thm:horizontal_st_type4}, crude form]
Let $X$ and $\L$ be as in Theorem~\ref{thm:main_intro} and let $x \in X$ be a type I\!V point. Then, up to replacing $K$ by a finite extension, there exists an affinoid neighborhood $U$ of $x$ in $X$, and a finite \'etale cover $g\colon V \ra U$, such that the pullback $g^* (\L|_U)$ is horizontal semistable in the sense of  \cite{Shimizu_2}. 
\end{mainthm}

\vspace{0.1in}
%%% 1.8
\subsection{Type V points}\label{subsection: outline type3}
\noindent 
\vspace{0.1in}

\noindent 
The only remaining case is that of type V points. Let $x$ be a type V point and let $x_{\mathrm{I\!I}}$ denote the associated type I\!I point. By the discussion above, we may assume that $\L$ is already semistable in an open neighborhood of $x_{\mathrm{I\!I}}$ (but this neighborhood might not contain $x$). With some careful analysis on integral models (see \S \ref{ss: type 5}), we may further reduce to the special situation where $x$ admits an open neighborhood $U$ isomorphic to a closed annulus $A_K([\sigma, 1])$ with finitely many open unit disks removed from the boundary circle $\{|T|=1\}$, such that $x_{\mathrm{I\!I}}$ is precisely the Shilov point corresponding to the boundary circle and $x$ is the type V point corresponding to the open unit disk centered at 0. We refer to such $U$ as a ``Swiss-cheese neighborhood'' of $x$. In fact, Swiss-cheese neighborhoods form a basis of open neighborhoods of $x$.

Our approach is to study the $p$-adic differential equations $\mathbf{N}_{\mathrm{dR}} (\L)$ over these neighborhoods. Similar to the case of type I\!I\!I points, we want to achieve unipotence of the differential operators after passing to finite \'etale covers. Since $\L$ is assumed to be semistable at the type I\!I point $x_{\mathrm{I\!I}}$, we already have unipotence there. Then the idea is to propagate this unipotence to a Swiss-cheese neighborhood. This time, however, only information about the differential operator in the arithmetic direction can be extracted from $x_{\mathrm{I\!I}}$. Therefore, we deal with the arithmetic direction first, before turning to the geometric direction (this is one essential difference between the case of type I\!I\!I points and that of type V points). More precisely, we show that $\mathbf{N}_{\mathrm{dR}} (\L)$ is arithmetically unipotent on a small enough neighborhood $U$. Then, we make full use of the explicit description of the relative period rings over such neighborhoods from \cite{relative} to prove that, after further enlarging $K$ and shrinking $U$, the arithmetic horizontal sections of $\mathbf{N}_{\mathrm{dR}} (\L)[\log \pi_K]$ looks like a $p$-adic differential equation with Frobenius structure over an annulus. Finally, we apply the classical $p$-adic monodromy theorem of Andr\'e--Kedlaya--Mebkhout  to conclude.

\vspace{0.1in}
%%%
\subsection{Notation and convention} \label{sss:notation_cyclo}
\noindent 
\vspace{0.1in}

\noindent Throughout this article, we follow the conventions on adic geometry in \cite{Huber_etale} and \cite{perfectoid}, and we adopt the following notation unless otherwise specified.

\begin{itemize}
    \item Let $K$ be a $p$-adic field with ring of integers $\mO_K$ and residue field $k=k_K$. Let $\varpi$ be a uniformizer of $K$. Let $K_0=W(k)[1/p]$. Let $\Gal_K=\Gal(\overline{K}/K)$ be the absolute Galois group of $K$.
    \item Fix an algebraic closure $\overline{K}$ of $K$ and let $\C_K$ be its completion. Note that $\C_K$ is a perfectoid field; let $\C_K^{\flat}$ denote its tilt. Let $\mO_{\C_K}$ (resp. $\mO_{\C_K}^\flat$) denote the ring of integers in $\C_K$ (resp. $\C_K^\flat$). 
    \item Fix a compatible system of $p$-power roots  $\{\zeta_{p^n} \}_{n \ge 0}$  of unity and a compatible system of $p$-power roots $\{p^{1/p^n}\}_{n \ge 0}$ of $p$ in $\mO_{\C_K}$. In particular, this gives rise to elements $\varepsilon = (1, \zeta_p, \zeta_{p^2}, \ldots)$ and $p^{\flat}=(p, p^{1/p}, p^{1/p^2}, \ldots)$ in $\mO_{\C_K}^\flat$. Let $\cl \pi = \varepsilon - 1$, and let $\pi = [\varepsilon] -1$ viewed as an element in $W(\mO_{\C_K}^{\flat})$ lifting $\cl \pi$. Let $t=\log[\varepsilon]=\log(\pi+1)$ be ``Fontaine's period'' in the de Rham period ring $\mathbf{B}_{\mathrm{dR}}$.
    \item For each $n\ge 0$, let $K_n=K(\zeta_{p^n})$ and let $K_{\infty}=\bigcup_{n\ge 0} K_n$ be the cyclotomic extension. Let $\widehat{K}_{\infty}$ be the $p$-adic completion of $K_{\infty}$ and let $\widehat{K}_{\infty}^\flat$ be its tilt. Let $\mO_{\widehat{K}_{\infty}}$ (resp.  $\mO^\flat_{\widehat{K}_{\infty}}$) denote the ring of integers in $\widehat{K}_{\infty}$ (resp.  $\widehat{K}_{\infty}^\flat$). 
    \item Let $e_K = [K_{\infty}: (K_0)_{\infty}]$.
    \item Let $K'_0$ denote the maximal unramified extension of $K_0$ inside ${K}_{\infty}$. Let $\kappa=\kappa_K$ denote the residue field of $K'_0$. 
    \item The standard norm $|\cdot|$ on $\mathbb{C}_K$, normalized by $|p|=\frac{1}{p}$, induces a norm (still denoted by $|\cdot|$) on $\mathbb{C}_K^{\flat}$. In particular, $|\cdot|$ restricts to a norm on $\widehat{K}_{\infty}^{\flat}$. Note that $|\overline{\pi}| = |\varepsilon-1|=p^{-\frac{p}{p-1}}$.  
    \item Let $v\colon \mathbb{C}_K\rightarrow \mathbb{R}\cup \{+\infty\}$ (resp. $v\colon \mathbb{C}_K^{\flat}\rightarrow \mathbb{R}\cup \{+\infty\}$) denote the valuation given by $v(y) \coloneqq - \log_p|y|$. We have $v(\overline{\pi}) = v(\varepsilon - 1)=\frac{p}{p-1}$.\footnote{In \cite{Berger-differential} and \cite{Colmez08}, the valuation $v$ is denoted by $v_{\mathbf{E}}$, the ring $\mO^\flat_{\C_K}$ (resp. $ \C_K^\flat$) is denoted by $\widetilde{\mathbf{E}}^+$ (resp. $\widetilde{\mathbf{E}}$), and the ring $\mO^\flat_{\widehat{K}_{\infty}}$ (resp. $ {\widehat{K}_{\infty}}^\flat$) is denoted by $\widetilde{\mathbf{E}}^+_K$ (resp. $\widetilde{\mathbf{E}}_K$).}
    \item Let $H_K = \Gal(\overline{K}/K_{\infty})$ be the Galois group of $K_\infty$, which equals the kernel of the cyclotomic character $\chi\colon\Gal_K\rightarrow \mathbb{Z}_p^{\times}$. 
    \item Let $\Gamma_K = \mathrm{Gal}(K_{\infty}/K)$. Note that  $\Gamma_{K} \cong \Gal_K /H_K$. 
    \item Let $\D_K(\alpha,\sigma)$ (resp., $\D_K(\alpha,<\sigma)$) denote the closed disk (resp., open disk) defined over $K$ centered at $\alpha$ with radius $\sigma$. (See Definition~\ref{defn: closed and open disks} for precise definitions.)
    \item  
    For an (open, closed, or half-open) interval $I \subset  \R_{> 0}$ with end-points in $p^{\mathbb{Q}}$, we denote by $A_K (I)$ the annulus over $K$ (viewed as an adic space) centered at $0$ with radius $|T| \in I$. When $I=[\sigma, \sigma]$ consists of a singleton, the annulus $A_K[\sigma, \sigma]$ is referred to as a \emph{circle}. We write $\Gamma_K (I) = \mO(A_K (I))$ for the ring of analytic functions on $A_K (I)$. 
    Sometimes we emphasize the coordinate $T$ for the annulus by writing $A_K (I; T)$ and $\Gamma_K (I; T)$ instead of $A_K(I)$ and $\Gamma_K (I)$. 
    \item Throughout the article, for all the period rings $\widetilde{\mathbf{A}}^r_K, \widetilde{\mathbf{B}}^r_K, \widetilde{\mathbf{C}}^r_K, \widetilde{\mathbf{C}}^{[s,r]}_K$ and $\mathbf{A}^r_K, \mathbf{B}^r_K,\mathbf{C}^r_K, \mathbf{C}^{[s,r]}_K$, etc., as well as their relative counterparts, we always assume $r,s\in \Q$ in the superscripts.

    \item Let $c_K\coloneqq \frac{p}{p-1}\cdot \frac{1}{e_K}$. For any Banach $K_0$-algebra $A$ equipped with a norm $|\cdot|$ extending the $p$-adic norm on $K_0'$, and any $r, s \in \Q$ with $0 < s \le r$, consider the ring 
    \[
    \mathscr{R}_A^{[s, r]}\coloneqq \left\{ \sum_{i\in \Z} a_i x^i ~\Big|~ a_i \in A, ~~ |a_i|\cdot p^{-it\cdot c_K} \rightarrow 0 \textrm{ as } i\rightarrow \pm\infty, \textrm{ for all } t\in [s,r] \right\}.
    \]
    In practice, the variable $x$ will mostly be taken to be $\pi_K$, except in \S \ref{ss:arithmetic_hs} where we take the variable to be $[T^{\flat}]$. For $r > 0$, we write $\displaystyle \mathscr{R}_A^r\coloneqq \varprojlim_{s \rightarrow 0^+} \mathscr{R}_A^{[s, r]}$ and we define the \emph{Robba ring} over $A$ to be $\displaystyle \mathscr{R}_A \coloneqq \varinjlim_{r \rightarrow 0^+} \mathscr{R}^r_A$. Note that when $A=K'_0$ for example, we have $\mathscr{R}_{K_0'}^{[s, r]} = \Gamma_{K_0'}([p^{-r\cdot c_K}, p^{-s\cdot c_K}]; x)$ and $\mathscr{R}_{K_0'}^r = \Gamma_{K_0'}([p^{-r\cdot c_K}, 1); x)$.
\end{itemize}

\subsection*{Acknowledgement}   The authors would like to thank Tomoyuki Abe, Tinhinane Azzouz, Alexander Beilinson, Dori Bejleri, Heng Du, Matthew Emerton, Lei Fu, Haoyang Guo, Mark Kisin, Kevin Lin, Alexander Petrov, Lie Qian, Peter Scholze, Koji Shimizu, Junyi Xie, and Ziquan Yang for helpful discussions during the preparation of this article, during which, the proof of our main result has gone through several iterations. In an earlier draft of this manuscript, the authors claimed a result for annuli that was too strong to be correct. We thank Peter Scholze for pointing this out. We are grateful to Junyi Xie for discussions around Proposition~\ref{prop: p_adic_Runge} and also on type I\!V points on Berkovich curves. We thank Tinhinane Azzouz for discussions around Corollary \ref{cor:1_dim_differential_equation_all_solvable}. During the project, H.D. was partially supported by the National Natural Science Foundation of China (Grant No. 12422101) and the National Key R\&D Program of China (Grant Nos. 2023YFA1009703 and 2021YFA1000704). Y.S.M. was partially supported by the Beijing Natural Science Foundation (Grant No. IS25030). Finally, Z.Y. wishes to thank Hausdorff Institute of Mathematics and Yau Mathematical Sciences Center at Tsinghua University for hospitality during his stays while completing this project.

\vspace{0.1in}
%%% Section 2
\section{Preliminaries on $p$-adic geometry}\label{section: prelim}

In this section, we recall some basic results from $p$-adic geometry that will be used throughout the article, including a classification of points on a smooth rigid analytic curve (regarded as an adic space). We adopt the notation introduced in \S \ref{sss:notation_cyclo}. Although this article only concerns the case where $K$ is a $p$-adic field, all the results in \S \ref{section: prelim} remain valid for an arbitrary discretely valued non-archimedean field equipped with a nontrivial valuation.

%\vspace{0.1in}
%% 2.1
\subsection{Classification of points on closed unit disks}\label{subsection: points on closed unit disks}
\noindent 
\vspace{0.1in}

\noindent
We start with the classification of points on a closed unit disk over $\mathbb{C}_K$. We refer the reader to \cite[Example 7.57]{Wedhorn} (also see \cite[Example 2.20]{perfectoid}) for a detailed exposition on this subject. Consider the closed unit disk centered at 0
\[ \D_{\mathbb{C}_K}  \coloneqq \spa \left(\mathbb{C}_K \gr{T}, \mathcal{O}_{\mathbb{C}_K} \gr{T}\right)=\D_{\C_K}(0,1)\]  
over $\mathbb{C}_K$, with coordinate $T$. For any $\alpha\in \mO_{\mathbb{C}_K}$ and any $0<\sigma\le 1$, let $\D_{\mathbb{C}_K}  (\alpha, \sigma) = \left\{|T - \alpha| \le \sigma\right\} \subset \D_{\mathbb{C}_K}$ denote the closed disk of radius $r$ centered at $\alpha$. Similarly, let $\D_{\mathbb{C}_K}  (\alpha, <\sigma) = \left\{|T - \alpha| < \sigma\right\}$ denote the open disk of radius $\sigma$.
The points on $\D_{\mathbb{C}_K}$ fall into the following five types. \\ 

\noindent \emph{Type I}. \: Each $\alpha \in \mathcal{O}_{\mathbb{C}_K}$ gives rise to a continuous valuation $x_{\alpha}\colon \mathbb{C}_K \gr{T} \ra \R_{\ge 0}$ sending $f \mapsto |f(\alpha)|$. Such a valuation is referred to as a \emph{type I} point. There is a one-to-one correspondence between $\mathcal{O}_{\mathbb{C}_K}$ and the type I points on $\D_{\mathbb{C}_K}$.

\vspace*{0.2cm}
\noindent \emph{Type I\!I and I\!I\!I}. \: For each closed disk $\D_{\mathbb{C}_K}  (\alpha, \sigma)$ with $\alpha \in \mO_{\mathbb{C}_K}$ and $0 < \sigma \le 1$, consider the continuous valuation 
\[ 
x_{\alpha, \sigma}\colon \mathbb{C}_K \gr{T} \ra \R_{\ge 0}, \quad f \mapsto \sup_{\substack{\beta \in \mO_{\mathbb{C}_K} \\ |\beta-\alpha|\le \sigma}} |f(\beta)|. 
\] 
Such a point $x_{\alpha, \sigma}$ is \emph{type I\!I} if $\sigma \in \left|\mathbb{C}_K^\times\right|=p^{\mathbb{Q}}$; otherwise $x_{\alpha, \sigma}$ is \emph{type I\!I\!I}. In the case of type I\!I, the point $x_{\alpha, \sigma}$ is called the \emph{Gaussian point} of the disk $\D_{\mathbb{C}_K}  (\alpha, \sigma)$.

\vspace*{0.2cm}
\noindent \emph{Type I\!V}.  Consider a sequence $\eta$ of closed disks $\D_1 \supset \D_2 \supset \cdots$ in $\D_{\mathbb{C}_K}$ such that their intersection $\cap_{i \in \Z_{>0}} \D_i$ contains no type I point. Then $\eta$ gives rise to a continuous valuation 
\[x_{\eta}\colon \mathbb{C}_K \gr{T} \ra \R_{\ge 0}, \quad f \mapsto \inf_{i \in \Z_{>0}} \sup_{\substack{\beta\in \mO_{\mathbb{C}_K} \\ \beta \in \D_i}} |f(\beta)|. \]
Such a point is called \emph{type I\!V}. We often write $x_{\eta} =$ ``$\cap_{i \in \Z_{>0}} \D_i$'' and view the intersection of these nested disks as the type I\!V point itself.  

\vspace*{0.2cm}
\noindent \emph{Type V}.  These are the points of rank-$2$ valuations (and thus do not appear in the theory of Berkovich spaces). Let us start with the type I\!I point $x_{\alpha, r}$ corresponding to the closed disk $\D_{\mathbb{C}_K} (\alpha, \sigma)$. For a sign $? \in \{>, <\}$, consider the totally ordered abelian group 
\[
\Gamma_{? \sigma}= \begin{cases}
    \Gamma_{<\sigma} = \R_{> 0} \times \gamma^{\Z} \quad \textrm{where } \gamma = \sigma^{-}, \textrm{ i.e., } \sigma'< \gamma < \sigma \textrm{ for all } \sigma' < \sigma  \\ 
    \Gamma_{>\sigma} = \R_{> 0} \times \gamma^{\Z} \quad \textrm{where } \gamma = \sigma^{+},  \textrm{ i.e., } \sigma < \gamma < \sigma' \textrm{ for all } \sigma' > \sigma
\end{cases}.
\]
(When $\sigma=1$, we only allow $?=$ ``$<$''.) Let $\beta \in \mO_{\mathbb{C}_K}$ such that $|\beta-\alpha|\le \sigma$. In particular, $\beta$ corresponds to a type I point $x_{\beta}$ in $\D_{\mathbb{C}_K} (\alpha, \sigma)$. Notice that every element $f\in \mathbb{C}_K \gr{T} $ can be uniquely written as $f = \sum_{i=0}^{\infty} a_i (T - \beta)^i$ with $a_i\in \mathbb{C}_K$. Consider the continuous rank-2 valuation $\mathbb{C}_K \gr{T} \ra \Gamma_{? \sigma}\cup \{0\}$ sending 
\[ 
f = \sum_{i=0}^{\infty} a_i (T - \beta)^i  \:  \longmapsto \: \sup_{i} |a_i| \gamma^i. 
\]
Such a valuation is referred to as a \emph{type V} point.\footnote{We have excluded the case $\sigma=1$ and $?=$ ``$>$''. In this case, the valuation does not live in $\D_{\mathbb{C}_K}$, but rather lives in the closure of $\D_{\mathbb{C}_K}$ in $\mathbb{A}^{1,\mathrm{an}}_{\mathbb{C}_K}$.} When $?$ is ``$>$'', the valuation does not depend on the choice of $x_{\beta}$; it only depends on the disk $\D_{\mathbb{C}_K}  (\alpha, \sigma)$. We denote this type V point by $x_{\alpha, > \sigma}$. On the other hand, when $?$ is ``$<$'', we denote the corresponding point by $x_{\beta, < \sigma}$. This point $x_{\beta, < \sigma}$ only depends on the open disk $\D_{\mathbb{C}_K}  (\beta, < \sigma) \subset \D_{\mathbb{C}_K}  (\alpha, r)$; it does not depend on the specific choice of $x_{\beta}$ in this open disk.

\begin{remark}
Points of type I, I\!I\!I, I\!V, and V are closed. The closure of a type I\!I point $x_{\alpha, \sigma}$ precisely consists of $x_{\alpha, \sigma}$ and the type V points $x_{\beta, < \sigma}$ and $x_{\alpha, >\sigma}$ described above; this closure is homeomorphic to $\P^1_{\bar{k}}$ when $0<\sigma<1$ (resp., $\A^1_{\bar{k}}$ when $\sigma=1$), with $x_{\alpha, \sigma}$ being the generic point.  
\end{remark} 

Next, we describe points on the closed unit disk (centered at 0)
\[\D_K \coloneqq \spa (K \gr{T}, \mathcal{O}_K \gr{T})=\D_K(0,1)\] defined over $K$. The following lemma and corollary allow us to identify $\D_K$ as a quotient of $\D_{\mathbb{C}_K}$.

\begin{lemma}\label{lemma: Galois orbit}
Let $A$ be a Banach $K$-algebra. Let $L/K$ be a finite normal extension and consider the Banach $L$-algebra $A_L \coloneqq A\otimes_KL$. Let $A^{\circ}$ and $A_L^{\circ}$ be the subring of power-bounded elements in $A$ and $A_L$, respectively. Let $\pi\colon\spa (A_L, A_L^{\circ})\rightarrow \spa (A, A^{\circ})$ be the natural map (of sets) induced by the inclusion $A\rightarrow A_L$ sending $a\mapsto a\otimes 1$.
\begin{enumerate}
\item[(1)] If $L/K$ is separable, there is a natural action of $\Gal(L/K)$ on $\spa (A_L, A_L^{\circ})$ and $\pi$ induces a homeomorphism
$\spa (A_L, A_L^{\circ})/\Gal(L/K)\cong \spa (A, A^{\circ})$,  
where the left-hand side is equipped with the quotient topology.
\item[(2)] If $L/K$ is purely inseparable, $\pi$ is a homeomorphism; namely, $\spa (A_L, A_L^{\circ})\cong \spa (A, A^{\circ}).$ 
\end{enumerate}
\end{lemma}

\begin{proof}
These are likely well-known to experts. 
For completeness, we include a proof (also see \cite[Proposition 1.3.5]{Berkovich_spectral} for a proof in the setup of Berkovich spectra). For part (1), the natural action of $\Gal(L/K)$ on $A_L$ preserves $A_L^{\circ}$ and hence induces an action on $\spa (A_L, A_L^{\circ})$. Notice that \[\left(A_L^{\Gal(L/K)}, (A_L^{\circ})^{\Gal(L/K)}\right)=(A, A^{\circ}).\] The desired result follows from \cite[Theorem 1.2]{Hansen}. For part (2), we construct an inverse of $\pi$. By assumption, we must have $\mathrm{char} (K)=p>0$ for some prime number $p$, and there exists an integer $m\ge 1$ such that $x^{p^m}\in A$ for all $x\in A_L$. Let $|\cdot|_A\in \spa (A, A^{\circ})$; namely, $|\cdot|_A\colon A\rightarrow \Gamma\cup \{0\}$ is a continuous valuation such that $|x|_A\le 1$ for all $x\in A^{\circ}$. By enlarging $\Gamma$ if necessary, we may assume that $\Gamma\cup\{0\}$ contains $\sqrt[p^m]{|A|_A}$. Consider a continuous valuation $|\cdot|_{A_L}\colon A_L\rightarrow \Gamma\cup\{0\}$ given by
$
|x|_{A_L} \coloneqq (|x^{p^m}  |_A)^{\frac{1}{p^m}}.$ 
Then $|\cdot|_{A_L}$ is a continuous valuation on $A_L$ satisfying $|x|_{A_L}\le 1$ for all $x\in A_L^{\circ}$. Using the multiplicativity of valuations, one sees that $|\cdot|_A\mapsto |\cdot|_{A_L}$ is indeed the inverse of $\pi$. 
\end{proof}

\begin{corollary}\label{cor: Galois orbit}
Let $A$ be a Banach $K$-algebra and let $A_{\mathbb{C}_K}$ be the completed tensor $A\widehat{\otimes}_K \mathbb{C}_K$. Let $A^{\circ}$ and $A_{\mathbb{C}_K}^{\circ}$ be the subring of power-bounded elements in $A$ and $A_{\mathbb{C}_K}$, respectively. Then the natural map 
$\pi\colon\spa (A_{\mathbb{C}_K}, A_{\mathbb{C}_K}^{\circ})\rightarrow \spa (A, A^{\circ})$ induces a homeomorphism
\[\spa (A_{\mathbb{C}_K}, A_{\mathbb{C}_K}^{\circ})/\Gal_K\cong \spa (A, A^{\circ}).\]
\end{corollary}

\begin{proof}
Note that the natural map
\[\spa(A_{\mathbb{C}_K}, A_{\mathbb{C}_K}^{\circ})\rightarrow \varprojlim_{L/K} \spa (A_L, A_L^{\circ})\]
is a homeomorphism, where the limit runs through all finite normal extensions $L/K$ inside $\overline{K}$. Then the desired statement follows from the previous lemma.
\end{proof}

Applying Corollary~\ref{cor: Galois orbit} to our situation, the natural map $\pi\colon\D_{\mathbb{C}_K}\rightarrow \D_K$ is surjective and induces a homeomorphism $\D_{\mathbb{C}_K}/\Gal_K\cong \D_K$; in particular, the points on $\D_K$ can be identified as the $\Gal_K$-orbits of the points on $\D_{\mathbb{C}_K}$. 

\begin{definition}
For $\star \in \{\mathrm{I, I\!I,  I\!I\!I,  I\!V, V}\}$, we say that a point on $\D_K$ has \emph{type} $\star$ if it is the image of a point on $\D_{\mathbb{C}_K}$ of type $\star$. 
\end{definition}

For later use, we also introduce the following notation for closed/open disks over $K$. 

\begin{definition}\label{defn: closed and open disks}
For $\alpha\in \mathcal{O}_{\mathbb{C}_K}$ and $0<\sigma\le 1$, let $\D_K(\alpha,\sigma)$ (resp. $\D_K(\alpha,<\sigma)$) denote the image of $\D_{\mathbb{C}_K} (\alpha,\sigma)$ (resp. $\D_{\mathbb{C}_K} (\alpha,<\sigma)$) under the projection $\pi\colon\D_{\mathbb{C}_K}\rightarrow \D_K$. We still refer to them as \emph{closed} (resp. \emph{open}) \emph{disks} of radius $\sigma$ centered at $\alpha$.
\end{definition}

We warn the reader that, even when $\sigma\in |K^{\times}|$, the disk $\D_K(\alpha,\sigma)$ is not necessarily isomorphic to $\D_K(0,\sigma)$ as rigid analytic spaces over $K$; in fact, $\D_K(\alpha,\sigma)\cong \D_K(0,\sigma)$ if and only if $\D_{\mathbb{C}_K} (\alpha,\sigma)$ contains some $x_{\beta}$ with $\beta\in \mathcal{O}_K$.

\vspace{0.1in}
%% 2.2
\subsection{Points on smooth curves}\label{subsection: points on smooth curves}
\noindent 
\vspace{0.1in}

\noindent
Now, we describe the points on a general smooth curve. Let $Y$ be a smooth rigid analytic curve over $K$, viewed as an adic space over $\spa (K, \mathcal{O}_K)$. By Corollary~\ref{cor: Galois orbit}, there is a natural surjection $\pi\colon Y_{\mathbb{C}_K}\rightarrow Y$ and the points on $Y$ can be identified as the $\Gal_K$-orbits of the points on $Y_{\mathbb{C}_K}$.

\begin{definition}\label{defn: points on a smooth curve}
For $\star \in \{\mathrm{I, I\!I,  I\!I\!I,  I\!V, V}\}$, we say that a point $x \in Y_{\mathbb{C}_K}$ has \emph{type} $\star$ if for an/any \'etale map $f\colon U \ra \D_{\mathbb{C}_K}$ from an open neighborhood $U$ of $x$ in $Y_{\mathbb{C}_K}$, the image $f(x) \in \D_{\mathbb{C}_K}$ has type $\star$. We say that a point on $Y$ has \emph{type} $\star$ if it is the image of a point on $Y_{\mathbb{C}_K}$ of type $\star$ under the surjection $\pi\colon Y_{\mathbb{C}_K}\rightarrow Y$.
\end{definition}

Note that an \'etale morphism between smooth rigid analytic curves over $K$ sends a point to a point of the same type.

\begin{definition}[{$K$-rational points}]\label{defn: rational points} 
Let $Y$ be a smooth rigid analytic curve over $K$, viewed as an adic space over $\spa(K, \mathcal{O}_K)$. Recall that a point on $Y$ corresponds to a $\Gal_K$-orbit of points on $Y_{\mathbb{C}_K}$. 
\begin{enumerate}
    \item[(1)] We say that a point on $Y_{\mathbb{C}_K}$ is \emph{$K$-rational} if its $\Gal_K$-orbit consists of a single element.
    \item[(2)] We say that a point on $Y$ is \emph{$K$-rational} if it is the image of a $K$-rational point on $Y_{\mathbb{C}_K}$. More generally, we say that a point on $Y$ is \emph{$K'$-rational} for some finite extension $K'/K$ if it is the image of a $K'$-rational point on $Y_{\mathbb{C}_K}$.
    \item[(3)] We say that a type I point on $Y$ is \emph{classical} if it is $K'$-rational for some finite extension $K'/K$. Classical points are precisely those points on $Y$ that belong to the underlying rigid analytic space.
\end{enumerate}
\end{definition}
 
\begin{remark}\label{remark:K_rational_point_under_field_extension}
If $x$ is a $K$-rational point on $Y$, and $K'/K$ is a finite extension, then the pre-image of $x$ under the projection $Y_{K'} \ra Y$ consists of a single element; thus, it can also be viewed as a $K'$-rational point on $Y_{K'}$. 
\end{remark} 

\begin{example} 
Consider the closed unit disk $\D_{K}=\D_K(0,1)$ over $K$. A type I point on $\D_K$ is a $\Gal_K$-orbit of a type I point $x_{\alpha}$ on $\D_{\mathbb{C}_K}$, with $\alpha \in \mO_{\mathbb{C}_K}$. This point is $K$-rational precisely when $\alpha \in \mO_K$ (in which case the $\Gal_K$-orbit consists of just $x_{\alpha}$ itself), and it is $K'$-rational if it is the $\Gal_K$-orbit of some $x_{\alpha}$ where $\alpha \in \mO_{K'}$. A point of type I\!I or I\!I\!I, given by the $\Gal_K$-orbit of $x_{\alpha, \sigma} \in \D_{\mathbb{C}_K}$, is $K'$-rational if one can choose $\alpha \in \mO_{K'}$. A similar remark holds for type V points given by the $\Gal_K$-orbit of $x_{\alpha, >\sigma}$ or $x_{\beta, <\sigma}$. In particular, every point of type I\!I, I\!I\!I, or V is $K'$-rational for some finite extension $K'/K$.
\end{example}

\begin{example}
More generally, for a smooth rigid analytic curve $Y$ over $K$, every point $x \in Y$ of type I\!I, I\!I\!I, or V is $K'$-rational for some finite extension $K'/K$. In this case, all points in the pre-image of this $x$ under the projection $Y_{K'} \ra Y$ are $K'$-rational. 
\end{example}

\vspace{0.1in}
%% 2.3
\subsection{Stalks and residue fields} \label{subsection: stalks}
\noindent 
\vspace{0.1in}

\noindent
For any point $x$ on a smooth rigid analytic curve $Y$ over $K$ (viewed as an adic space), we write $\mO_{Y, x}$ for the stalk of the structure sheaf at $x$, and denote the residue field of $\mO_{Y, x}$ by $\mathbf{k}(x)$. The continuous (multiplicative) valuation corresponding to $x$ extends to a valuation $v_x$ on $\mathbf{k}(x)$. The completion of $\mathbf{k}(x)$ with respect to $v_x$ is denoted by $\sH (x)$. The field $\mathbf{k}(x)$ (resp. $\sH (x)$) is referred to as the \emph{residue field} (resp. \emph{completed residue field}) at $x$. 

\begin{example}[Stalks of non-type I points]\label{example: stalk_on_curve}
For a non-type I point $\xi$ on a smooth rigid analytic curve $Y$ over $K$, its stalk $\mO_{Y, \xi}$ is a field and is equal to $\mathbf k(\xi)$. Hence, the completed stalk agrees with $\sH(\xi)$. Moreover, suppose $\xi$ is a type V point and $\xi_{\mathrm{I\!I}}$ is the unique type I\!I point specializing to $\xi$, then we have an injection $\iota\colon\mathbf k(\xi)\hookrightarrow \mathbf k(\xi_{\mathrm{I\!I}})$ by \cite[Lemma 1.1.10]{Huber_etale}. In fact, $\iota$ is a homeomorphism onto its image and $\iota(\mathbf k(\xi))$ is dense in $\mathbf k(\xi_{\mathrm{I\!I}})$ (also see \cite[Theorem I.1.5.4(b)]{Morel_adic}). Passing to completions, we obtain an isomorphism $\sH (\xi) \cong \sH (\xi_{\mathrm{I\!I}})$. 
\end{example}  

The following lemma will be used in the proof of the main theorem.

\begin{lemma} \label{lemma: base_change_from_type2_to_type5}
Let $Y$ be a smooth rigid analytic curve over $K$, viewed as an adic space. Let $\xi \in Y$ be a type V point and let $\xi_{\mathrm{I\!I}}$ be the unique type I\!I generalization of $\xi$. Then the inclusion 
\[ \iota\colon \mO_{Y, \xi} = \mathbf{k}(\xi) \hookrightarrow \mO_{Y, \xi_{\mathrm{I\!I}}}  = \mathbf{k}(\xi_{\mathrm{I\!I}})\] induces an equivalence between the category of finite extensions of $\mO_{Y, \xi}$ and that of $\mO_{Y, \xi_{\mathrm{I\!I}}}$.
\end{lemma}

\begin{proof}
For simplicity, we identify $\mathbf{k}(\xi)$ as a subfield of $\mathbf{k}(\xi_{\mathrm{I\!I}})$ via $\iota$. We first fix the following notation. 
\begin{itemize}
    \item Let $v$ denote the rank-$2$ (multiplicative) valuation on $\mathbf{k}(\xi)$ and let $w$ denote the rank-$1$ (multiplicative) valuation on $\mathbf{k}(\xi_{\mathrm{I\!I}})$.
    \item Let $u  \coloneqq  w|_{\mathbf{k}(\xi)}$ denote the restriction of $w$ to $\mathbf{k}(\xi)$, which is a nontrivial rank-$1$ valuation on $\mathbf{k}(\xi)$.
    \item  Let $\mathbf{k}^+(\xi) = \mO_{Y, \xi}^+  \coloneqq  \{a \in \mO_{Y, \xi} \,|\, v (a) \le 1\}$, $\mathbf{k}^+(\xi_{\mathrm{I\!I}}) = \mO_{Y, \xi_{\mathrm{I\!I}}}^+  \coloneqq  \{a \in \mO_{Y, \xi_{\mathrm{I\!I}}} \,|\, w (a) \le 1\}$, and $\mathbf{k}_u^+(\xi)  \coloneqq  \{a \in \mO_{Y, \xi}\,|\, u (a) \le 1\}$ be the corresponding valuation subrings. For simplicity, we write $V = \mathbf{k}^+(\xi)$ and $W = \mathbf{k}_u^+(\xi)$.  
\end{itemize} 
Note that $W =  \mathbf{k}(\xi) \cap \mathbf{k}^+(\xi_{\mathrm{I\!I}})$ and $V \subset W$. Also note that $(W, u)$ is a rank-$1$ valuation ring, while $(V, v)$ is a rank-$2$ valuation ring. Let $\mathfrak p$ (resp. $\mathfrak q$) denote the unique height-$1$ prime in $V$ (resp. in $W$, that is, the maximal ideal in $W$).  
 
We claim that the non-archimedean valued field $(\mathbf{k}(\xi), u)$ is \emph{quasi-complete} in the sense of \cite[Definition 2.3.1]{Berkovich_etale}. By \cite[Theorem 1.9.13]{Temkin_notes}, this is equivalent to claiming that $(W, \mathfrak q)$ is a Henselian pair. It follows from \cite[Corollary I.1.4.4]{Morel_adic} that $W = V_{\mathfrak p}$ and $\mathfrak q = \mathfrak p V_{\mathfrak p}$. In fact, since $V$ is a valuation ring, we must have $\mathfrak p = \mathfrak p V_{\mathfrak p}$ as subsets of $\mathbf{k}(\xi)$. Therefore, by Gabber's criterion (cf. \cite[Tag 09XI (5)]{stacks-project}), it suffices to show that $(V, \mathfrak p)$ is Henselian. Next, we observe that $\mathfrak p = \mathrm{rad} (\varpi)$. To see this, note that for every prime ideal $\mathfrak q' \subsetneq \sqrt{\mathfrak p}$ properly contained in the radical of $\mathfrak p$, we have $\varpi^n\notin \mathfrak{q}$ for all $n$, and hence $\mathfrak q' \subset (\varpi^n)$ for all $n$. This forces $\mathfrak q' = (0)$ by the topological nilpotence of $\varpi$. Since being a Henselian pair is insensitive to passing to radicals of the relevant ideal, it suffices to show that $(V, (\varpi))$ is Henselian, but this is \cite[Proposition III.6.3.7]{Morel_adic} (also see \cite[Proposition 7.5.5]{Bhatt_adic_notes}). This finishes the proof of the claim. In fact, the same argument implies that $(\mathbf{k}(\xi_{\mathrm{I\!I}}), w)$ is also quasi-complete. 

Let $\mathbf{k}(\xi)^\wedge_u$ denote the completion of $\mathbf{k}(\xi)$ with respect to $u$, and similarly define $\mathbf{k}(\xi)^\wedge_v$ and $\mathbf{k}(\xi_{\mathrm{I\!I}})^\wedge_w$. Note that $\mathbf{k}(\xi)^\wedge_v$ (resp. $\mathbf{k}(\xi_{\mathrm{I\!I}})^\wedge_w$) is precisely the completed residue field $\mathscr{H} (\xi)$ (resp. $\mathscr{H} (\xi_{\mathrm{I\!I}})$).  We have
\begin{equation}\label{eq: completed residue fields are the same}
\mathbf{k}(\xi)^{\wedge}_{v} =  \mathbf{k}(\xi)^{\wedge}_{u} = \mathbf{k}(\xi_{\mathrm{I\!I}})^\wedge_w = \mathscr{H}(\xi)=\mathscr{H}(\xi_{\mathrm{I\!I}})
\end{equation}
by  Example \ref{example: stalk_on_curve}. 
Now, by \cite[Proposition 2.4.1]{Berkovich_etale}, taking base change along the inclusion $\mathbf{k}(\xi) \hookrightarrow \mathbf{k}(\xi)^\wedge_u $ (resp.  $\mathbf{k}(\xi_{\mathrm{I\!I}}) \hookrightarrow  \mathbf{k}(\xi_{\mathrm{I\!I}})^\wedge_w$) induces an equivalence between finite extensions of $\mathbf{k}(\xi)$ and that of $\mathbf{k}(\xi)^\wedge_u $ (resp. finite extensions of $\mathbf{k}(\xi_{\mathrm{I\!I}})$ and that of $\mathbf{k}(\xi_{\mathrm{I\!I}})^\wedge_w$). Combining this with (\ref{eq: completed residue fields are the same}), we are done.
\end{proof}

\vspace{0.1in}
%% 2.4
\subsection{Comparison with Berkovich curves}\label{subsection: comparison_between_stalks}
\noindent 
\vspace{0.1in}

\noindent
Let $Y$ be a smooth affinoid rigid analytic curve over $K$ associated with an affinoid $K$-algebra $\mathscr A$. Let $Y^{\mathrm{Ber}}$ (resp. $Y^{\mathrm{ad}}$) be its associated Berkovich space (resp. adic space). Note that $Y^{\mathrm{Ber}}$ is identified as the maximal Hausdorff quotient of $Y^{\mathrm{ad}}$ as a topological space. Let $\mO_{Y^{\mathrm{Ber}}}$ denote the structure sheaf of the Berkovich space. For an open subset\footnote{This is an open subset in the Hausdorff topology of $Y^{\mathrm{Ber}}$, also referred to as the \emph{canonical} topology.
}  
$U \subset Y^{\mathrm{Ber}}$, the global sections are given by 
$
\Gamma(U, \mO_{Y^\mathrm{Ber}}) = \lim_{V \subset U} \mathscr A_V,
$ 
where the limit is taken over all $V \subset U$ that are finite unions of affinoid subdomains (that is, \emph{special subsets}) in $Y$, and $\mathscr A_V$ is the ring of analytic functions on $V$ (see \cite[2.2.6]{Berkovich_spectral}). Let $x$ be a point on $Y^{\mathrm{Ber}}$, which we also view as a point on $Y^{\mathrm{ad}}$ (it has type I, I\!I, I\!I\!I, or I\!V). By the Gerritzen--Grauert Theorem (see \cite[7.3.5/3]{BGR} and \cite[Proposition 2.2.3]{Berkovich_spectral}), the Berkovich stalk $\mO_{Y^{\mathrm{Ber}}, x}$ can be computed as 
$
\mO_{Y^{\mathrm{Ber}}, x} = \underset{x \in U}{\mathrm{colim}} \lim_{V \subset U} \mathscr A_V,
$  
where the colimit is taken over all open neighborhoods $U$ of $x$ (in the canonical topology), and for each such $U$ the limit is taken over all \emph{rational subsets} $V$ contained in $U$. On the other hand, the stalk of the structure sheaf $\mO_{Y^{\mathrm{ad}}}$ of the adic space $Y^{\mathrm{ad}}$ is given by 
$
\mO_{Y^{\mathrm{ad}}, x} =  \underset{x \in V}{\mathrm{colim} } \: \mathscr A_V
$  
where the colimit is taken over all rational subsets $V \subset Y$ that contain the point $x$. Thus, we have a natural inclusion 
\[
\mO_{Y^{\mathrm{Ber}}, x}  \hookrightarrow \mO_{Y^{\mathrm{ad}}, x}\] of local rings. This inclusion is an isomorphism when $x$ is a type I point but \emph{not} necessarily an isomorphism in general. When $x$ has type I\!I, I\!I\!I, or I\!V, both the Berkovich stalk $\mO_{Y^{\mathrm{Ber}}, x}$ and adic stalk $\mO_{Y^{\mathrm{ad}}, x}$ are fields (cf. Example \ref{example: stalk_on_curve}). The valuation corresponding to $x$ extends to a valuation $v_x$ on both stalks, making them topological fields. We claim that, in this case, the inclusion $\mO_{Y^{\mathrm{Ber}}, x}  \hookrightarrow\mO_{Y^{\mathrm{ad}}, x}$ has a dense image and thus becomes an isomorphism upon taking completions, which we identify and simply denote by $\sH (x)$ as in Example \ref{example: stalk_on_curve}. To see the claim, let $h \in  \Gamma (V, \mO_{Y^{\mathrm{ad}}})$ be an analytic function on a rational subset $V$ containing $x$, we want to show that for any $\epsilon \in \R_{> 0}$, there exists an open subset $U \subset Y^{\mathrm{Ber}}$ (in the canonical topology) containing $V$ and an analytic function $h' \in \Gamma (U, \mO_{Y^{\mathrm{Ber}}})$, such that $v_x(h - h') < \epsilon$. Without loss of generality, we may write 
\[
V = Y \left(\frac{f_1, \ldots, f_n}{g}\right)  =   \left\{ y \in Y \:  \Big\vert \: v_y(f_i) \le v_y(g), \: \forall i    \right\}
\]
where $f_1, \ldots, f_n \in A$ generates the unit ideal, $g \in A$, and $f_n = \varpi^N$ for some sufficiently large integer $N$ (see \cite[Remark 2.8]{perfectoid}). By the construction of rational subsets, there exists some $h'$ of the form 
\[ 
h' = \sum_{\substack{I = (i_1, \ldots, i_n) \\ \mathrm{finite \,\, sum} }} a_I \cdot (f_1/g)^{i_1} (f_2/g)^{i_2} \cdots (f_n/g)^{i_n}, \qquad a_I \in A,
\]
such that $|h - h'|_{V} < \epsilon$, where $|\cdot|_{V}$ denotes the spectral norm on $V$. In particular, $v_x(h - h') < \epsilon$. Since $h'$ is defined on the Zariski open subset $Y \minus V(g)$, it is defined on some open subset $U \subset Y^{\mathrm{Ber}}$ (in the canonical topology) that contains $V$. This justifies our claim. 

Finally, let us remark that, when $x$ has type I\!I\!I or I\!V, the inclusion $\mO_{Y^{\mathrm{Ber}}, x}  \hookrightarrow \mO_{Y^{\mathrm{ad}}, x}$ is in fact an isomorphism since Berkovich neighborhoods are cofinal among adic neighborhoods at a point of  type I\!I\!I  or  I\!V (cf. \cite[Proposition 2.2.1] {Berkovich_differential}).

\vspace{0.1in}
%%% Section 3
\section{$(\varphi, \Gamma_K)$-modules and Galois representations} \label{sec:classical_phi_Gamma}

In this section, we review the classical theory on how to associate overconvergent  $(\varphi, \Gamma_K)$-modules and $p$-adic differential equations to (de Rham) $p$-adic Galois representations, following \cite{Cherbonnier_Colmez, Berger-differential, Colmez08}. We largely follow the notation of Kedlaya--Liu \cite{Kedlaya-LiuI} on period rings although the content in this section predates \emph{op.\ cit.} and is mostly developed in \cite{Cherbonnier_Colmez} and \cite{Berger-differential}.
Experts familiar with  \cite{Cherbonnier_Colmez} and \cite{Berger-differential} should feel free to skip most of this section. 

With much regret, we remind the reader that each standard reference in the literature uses a slightly different notation for period rings. To orient the reader, we compile a list of different notations used in various references in Appendix \ref{appendix:notation}.   

\vspace{0.1in}
%%% 3.1
\subsection{Perfect period rings} \label{ss:point perfect period ring}

\subsubsection{Perfect period rings of type $\mathbf{A}, \mathbf{B}$}
Let $\widetilde{\mathbf{A}} \coloneqq W(\C_K^{\flat})$ and $\widetilde{\mathbf{A}}^+  \coloneqq W(\mO_{\C_K}^{\flat})$, equipped with natural $\Gal_K$- and Frobenius actions. Taking $H_K$-invariants, we obtain  
$\widetilde{\mathbf{A}}_K  \coloneqq  (\widetilde{\mathbf{A}})^{H_K} = W(\widehat{K}_{\infty}^{\flat}) $ and $  \widetilde{\mathbf{A}}^+_K  \coloneqq  (\widetilde{\mathbf{A}}^+)^{H_K}=W(\mO_{\widehat{K}_{\infty}^{\flat}}).
$ 
Every element $x\in \widetilde{\mathbf{A}}$ (resp. $x\in \widetilde{\mathbf{A}}_K$) can be uniquely expressed as $x=\sum_{n=0}^{\infty}p^n[y_n]$ with $y_n\in \C_K^{\flat}$ (resp. $y_n\in \widehat{K}_{\infty}^{\flat}$). The ring $\widetilde{\mathbf{A}}^+$ is denoted by $\mathbf{A}_{\mathrm{inf}}$ in \cite{Fontaine_94}. We refer the reader to \emph{op.\ cit.} for the construction of period rings $\mathbf{B}_{\mathrm{dR}}, \mathbf{B}_{\mathrm{crys}}, \mathbf{B}_{\mathrm{st}}$ and their variants.

\begin{definition}[{\cite[\S 2.1]{Berger-differential}, \cite[\S 5.2]{Colmez08}}]
\label{def:perfect_period_rings_type_A}
\begin{enumerate}
\item[(1)]  
For every $r>0$, define 
\begin{align*}
\widetilde{\mathbf{A}}^{r}  \coloneqq \left\{x=\sum_{n=0}^{\infty}p^n[y_n] \in \widetilde{\mathbf{A}}\,\Big|\,\lim_{n\rightarrow\infty}v(y_n)+\frac{n}{r}=+\infty\right\}
\end{align*}
and define $\widetilde{\mathbf{A}}_K^{r}  \coloneqq  (\widetilde{\mathbf{A}}^{r})^{H_K}$. 
\item[(2)] 
For any $s\in (0,r]$, we define a norm $\lambda_s$ on $\widetilde{\mathbf{A}}^{r}$ by
\[\lambda_{s}(x) \coloneqq \sup_{n\ge 0} \,\{p^{-n}\cdot |y_n|^s\}.\]
This restricts to a norm $\lambda_s$ on $\widetilde{\mathbf{A}}_K^{r}$. Note that $\widetilde{\mathbf{A}}^{r}$ and $\widetilde{\mathbf{A}}_{K}^{r}$ are complete with respect to $\lambda_{r}$.
\item[(3)] Define $\widetilde{\mathbf{A}}^{\dagger}  \coloneqq\varinjlim_{r>0} \widetilde{\mathbf{A}}^{r}$ and  $ \widetilde{\mathbf{A}}^{\dagger}_K \coloneqq\varinjlim_{r>0} \widetilde{\mathbf{A}}^{r}_K,$
where the transition maps are natural inclusions. 
\item[(4)] The period rings $\widetilde{\mathbf{B}}, \widetilde{\mathbf{B}}^+, \widetilde{\mathbf{B}}^{r}, \widetilde{\mathbf{B}}^{\dagger},
    \widetilde{\mathbf{B}}_K$, $\widetilde{\mathbf{B}}^+_K$, $\widetilde{\mathbf{B}}^{r}_K$, and $\widetilde{\mathbf{B}}^{\dagger}_K$ are obtained from the corresponding period rings of type $\mathbf{A}$ with $p$ inverted.  
     \end{enumerate}
\end{definition}

One checks that $\widetilde{\mathbf{A}}$ (resp. $\widetilde{\mathbf{A}}_K$) is the $p$-adic completion of $\widetilde{\mathbf{A}}^{\dagger}$ (resp. $\widetilde{\mathbf{A}}_K^{\dagger}$).

\vspace{0.1in}
\subsubsection{Perfect period rings of type $\mathbf{C}$} 

\begin{definition}[{\cite[\S 2.3, Definition 2.16]{Berger-differential}, \cite[\S 5.4]{Colmez08}%\cite[\S 5.1]{Kedlaya-LiuI}
}]
\label{def:perfect_robba_ring_type_C_over_K} \noindent 
\begin{enumerate}
\item[(1)] For $0 < s \le r$,  we define $\widetilde{\mathbf{C}}^{[s, r]}$  (resp. $\widetilde{\mathbf{C}}^{[s, r]}_K$ )  to be the completion of $\widetilde{\mathbf{B}}^r$  (resp. of $\widetilde{\mathbf{B}}^r_K$) with respect to \(\max\{\lambda_{s}, \lambda_{r}\}\) (or equivalently, completion with respect to the family of norms $\{\lambda_t\}_{t\in [s,r]}$).  
\item[(2)] 
For every $r > 0$, we define $\widetilde{\mathbf{C}}^{r}$ (resp.  $\widetilde{\mathbf{C}}^{r}_K$) to be the Fr\'echet completion of  $\widetilde{\mathbf{B}}^r$ (resp.  $\widetilde{\mathbf{B}}^r_K$) with respect to the family of norms $\{\lambda_{s}\}_{s \in  (0, r]}$.  
\item[(3)] 
Define $\sq{\mathbf C}^{\dagger}$ and $\sq{\mathbf C}_K^{\dagger}$ to be the union $
\sq{\mathbf{C}}^{\dagger}\coloneqq \varinjlim_{r>0} \widetilde{\mathbf{C}}^{r}$ and $\sq{\mathbf{C}}_K^{\dagger}\coloneqq \varinjlim_{r>0} \widetilde{\mathbf{C}}^{r}_K.
$
\end{enumerate} 
\end{definition}

\begin{remark} \label{remark:subscript_K_as_invariant_perfect}
\begin{enumerate}
    \item[(i)] 
For each $r>0$, we have  $\widetilde{\mathbf{A}}_K^{r} =   (\widetilde{\mathbf{A}}^{r})^{H_K} = \widetilde{\mathbf{A}}^{r} \cap \widetilde{\mathbf{A}}_K$ as subrings of $\widetilde{\mathbf{A}}$. Similarly we have 
$\widetilde{\mathbf{B}}_K^{r} =   (\widetilde{\mathbf{B}}^{r})^{H_K}$ and $ \widetilde{\mathbf{C}}_K^{r} = (\widetilde{\mathbf{C}}^{r})^{H_K}$, and therefore %filtered colimit preserves finite limits
$\widetilde{\mathbf{A}}_K^{\dagger} =   (\widetilde{\mathbf{A}}^{\dagger})^{H_K}, \widetilde{\mathbf{B}}_K^{\dagger} = (\widetilde{\mathbf{B}}^{\dagger})^{H_K},$ and $  \widetilde{\mathbf{C}}_K^{\dagger} = (\widetilde{\mathbf{C}}^{\dagger})^{H_K}$. In particular, our definitions agree with the ones in 
 \cite{Berger-differential} and \cite{Colmez08}.
 \item[(ii)] 
For each $r > 0$, the natural map $\widetilde{\mathbf{C}}^{r} \rightarrow \varprojlim_{0<s\le r} \widetilde{\mathbf{C}}^{[s,r]}$ is an isomorphism. A similar isomorphism holds for $\widetilde{\mathbf{C}}_K^{r}$ and $\widetilde{\mathbf{C}}_K^{[s,r]}$.
 \end{enumerate}
\end{remark}

\begin{remark}[Frobenius structures] \label{remark:Frobenius_classical}
The Frobenius maps on $\widetilde{\mathbf{A}}^+$ and $\widetilde{\mathbf{A}}$ extend to Frobenius maps on $\widetilde{\mathbf{B}}^+$ and $\widetilde{\mathbf{B}}$, which we denote by $\varphi$. Moreover, for $0 < s \le r$, we have Frobenius isomorphisms 
$
\varphi\colon \widetilde{\mathbf{A}}^{r}\isom  \widetilde{\mathbf{A}}^{r/p}, $ $ \varphi\colon \widetilde{\mathbf{B}}^{r}\isom  \widetilde{\mathbf{B}}^{r/p}, $ $ \varphi\colon \widetilde{\mathbf{C}}^{r}\isom  \widetilde{\mathbf{C}}^{r/p},  $ $ \varphi\colon \widetilde{\mathbf{C}}^{[s, r]}\isom  \widetilde{\mathbf{C}}^{[s/p,\,r/p]}
$ 
extending the Frobenius on $\widetilde{\mathbf{A}}^+$. Passing to colimits, we obtain Frobenius maps $\varphi$ on the period rings $\widetilde{\mathbf{A}}^\dagger, \widetilde{\mathbf{B}}^\dagger$, and $\widetilde{\mathbf{C}}^\dagger$. 
Similarly, there are natural Frobenius structures on $\widetilde{\mathbf{A}}_K^r, \widetilde{\mathbf{B}}_K^r$,  $\widetilde{\mathbf{C}}_K^r$, $\widetilde{\mathbf{C}}_K^{[s,r]}$, $\widetilde{\mathbf{A}}_K^\dagger, \widetilde{\mathbf{B}}_K^\dagger$, and $\widetilde{\mathbf{C}}_K^\dagger$.
\end{remark}

\begin{remark}\label{remark:notation_iota_n}
For later use, we recall that there is a natural inclusion $\iota_0\colon \widetilde{\mathbf{A}}^r\hookrightarrow \mathbf{B}_{\mathrm{dR}}^+$ (see \cite[Proposition 2.11]{Berger-differential}) for any $r\ge 1$. In general, for any $r>0$ and any non-negative integer $n\ge -\log_pr$, there is an injection
\[\iota_n \coloneqq \iota_0\circ \varphi^{-n}\colon \widetilde{\mathbf{A}}^r\xrightarrow[]{\varphi^{-n}} \widetilde{\mathbf{A}}^{p^nr}\xrightarrow[]{\iota_0}\mathbf{B}_{\mathrm{dR}}^+.\]
By continuity, for any $0<s\le r$ and any non-negative integer $n$ such that $p^{-n}\in [s,r]$, the map $\iota_n$ above extends to injections
$
\iota_n\colon \sq{\mathbf{C}}^{[s,r]} \ra \mathbf{B}_{\mathrm{dR}}^+$ and $ \iota_n\colon \sq{\mathbf{C}}^r \ra \mathbf{B}_{\mathrm{dR}}^+.
$ 
\end{remark}

%%% 3.2
\vspace{0.1in}
\subsection{Imperfect period rings} \label{ss:point imperfect period ring}
\noindent
\vspace{0.1in}

\noindent Let $\mathbf{E} \coloneqq k(\!(\cl \pi)\!)^{\mathrm{sep}}$ be the separable closure of $k(\!(\cl \pi)\!)$ inside $\C_K^\flat$ and let $\mathbf{E}_K  \coloneqq  (\mathbf{E})^{H_K}$. It turns out $\mathbf{E}_K$ is a discretely-valued field with residue field $\kappa_K$. Let $\mathbf{E}_K^+$ be its ring of integers. By choosing a uniformizer $\cl \pi_K$ of $\mathbf E_K^+$, we have $\mathbf{E}_K\cong \kappa_K(\!(\overline{\pi}_K)\!)$ and $\mathbf E_K^+ \cong \kappa_K [\![\cl \pi_K]\!]$. Moreover, the theory of fields of norms yields an isomorphism 
\[
\varprojlim_{n} \mO_{K_n}/\mathfrak a \isom \mathbf{E}_K^+ 
\]
where $\mathfrak{a}$ is the ideal consisting of elements $x$ satisfying $|x| \le p^{-\frac{1}{p}}$ and the inverse limit is taken over the Frobenius maps  $ \mO_{K_{n+1}}/\mathfrak a \ra \mO_{K_n}/\mathfrak a$ sending $x \mapsto x^p$.

\begin{remark}
Notice that $\mathbf{E}_K$ is a totally ramified extension of $\mathbf{E}_{K'_0}\cong \kappa_K(\!(\overline{\pi})\!)$ of degree $e_K=[K_{\infty}:(\Q_p)_{\infty}]$ as in \S \ref{sss:notation_cyclo}. Thus, we have $|\overline{\pi}_K|=p^{-\frac{p}{p-1}\cdot\frac{1}{e_K}}$. %Consequently, for any $r>0$, we have $\bar{\lambda}_r(\overline{\pi}_K)=p^{-\frac{p}{p-1}\cdot\frac{r}{e_K}}$.
\end{remark}

\subsubsection{Imperfect period rings} 
There are subrings $\mathbf{A} \subset \widetilde{\mathbf{A}}$ and  $\mathbf{A}_K\subset \widetilde{\mathbf{A}}_K$, which are both $p$-adically complete and satisfy 
\[ 
\mathbf{A}/(p)=\mathbf{E}, \qquad \mathbf{A}_K/(p)=\mathbf{E}_K\cong \kappa_K(\!(\overline{\pi}_K)\!).
\] 
Moreover, $\mathbf{A}$ is stable under the $\varphi$- and $\Gal_K$-actions on $\widetilde{\mathbf{A}}$, and we have $\mathbf{A}_K=(\mathbf{A})^{H_K}$, which thus carries an action of $\Gamma_K$.  We omit the constructions here and refer the reader to \cite[\S 6.1]{Colmez08} for details, noting only that, for a finite unramified extension $L_0 = W(k_0)[1/p]$ over $\Q_p$, the ring $\mathbf{A}_{L_0}$ is precisely the $p$-adic completion of $W(k_0) [\![\pi]\!] [1/\pi]$ in $\widetilde{\mathbf{A}}$. We write $\mathbf A_{L_0}^+ \subset \mathbf A_{L_0}$ for the subring $W(k_0) [\![\pi]\!].$ 

\begin{construction}[{\cite[\S 6.2]{Colmez08}}] \label{construction:pi_K}
We construct a lift $\pi_K\in \mathbf{A}_K $ of $\overline{\pi}_K$ as follows. Let $\overline{P}_K\in \mathbf{E}_{K'_0}^+[X]=\kappa_K [\![\cl \pi]\!][X]$ be the minimal polynomial of $\overline{\pi}_K$ over $\mathbf{E}_{K'_0}$. % $\mathbf{E}_{K'_0} = \kappa_K (\!(\cl \pi)\!)$.
Choose a monic polynomial $P_K\in \mathbf{A}^+_{K'_0}[X]$ of the same degree that reduces to $\overline{P}_K$ modulo $p$. By Hensel's Lemma, $P_K$ has a unique root $\pi_K$ in $\mathbf{A}_K$ lifting $\overline{\pi}_K$. Moreover, $\mathbf{A}_K$ is a free $\mathbf{A}_{K'_0}$-module with basis $\{1, \pi_K, \ldots, \pi_K^{e_K-1}\}$. %2In fact, we can (and will) pick $\pi_K$ such that $\pi_K^{e_K} = \pi$. If $L/K$ is another finite extension, we can further pick $\pi_K \in \mathbf{A}_K$ such that $\pi_K^{e_{L/K}} = \pi_K$ where $e_{L/K} = e_L/e_K$ is the ramification index of $L/K$. 
\end{construction} 

Using \cite[Proposition 7.5]{Colmez08}, one can identify $\mathbf{A}_K$ with the $p$-adic completion of $\mathcal{O}_{K'_0}[\![\pi_K]\!][\frac{1}{\pi_K}]$. More explicitly, we have  
\[\mathbf{A}_K=\left\{\sum_{i\in \mathbb{Z}}c_i\pi_K^i\,\Big|\,c_i\in \mathcal{O}_{K'_0}\; \textrm{ and }c_i\rightarrow 0 \textrm{ $p$-adically as }i\rightarrow -\infty\right\}.\]
%Next, let us recall the definition of various imperfect period rings that will appear in this article (for the reader familiar with other notations, see Remark~\ref{remark:notation1} and Appendix \ref{appendix:notation}). %It is worth pointing out that these rings are defined as subrings of the corresponding perfect period rings.  %$\widetilde{\mathbf{A}}_K^{r}$. 

\begin{definition}[Classical imperfect period rings] \label{definition:A^r_as_intersection}
\begin{enumerate}
\item[(1)] For $r>0$, define 
$\mathbf{A}^{r} \coloneqq \mathbf{A} \cap \widetilde{\mathbf{A}}^{r}$ and $\mathbf{A}_K^{r} \coloneqq \mathbf{A}_K\cap \widetilde{\mathbf{A}}^{r}_K$, where the intersection is taken inside $\widetilde{\mathbf{A}}$ and $\widetilde{\mathbf{A}}_K$ respectively.
\item[(2)] Define ${\mathbf{A}}^{\dagger}$ and ${\mathbf{A}}^{\dagger}_K$ as the unions ${\mathbf{A}}^{\dagger} \coloneqq\varinjlim_{r>0}  {\mathbf{A}}^{r}$ and $  {\mathbf{A}}^{\dagger}_K \coloneqq\varinjlim_{r>0}  {\mathbf{A}}^{r}_K$. 
\item[(3)] The period rings ${\mathbf{B}},  {\mathbf{B}}_K$, $\mathbf{B}^r$,  $ {\mathbf{B}}^{r}_K$, $\mathbf{B}^{\dagger}$, and $ {\mathbf{B}}^{\dagger}_K$ are obtained from the corresponding period rings of type $\mathbf{A}$ with $p$ inverted.    %Note that $\widetilde{\mathbf{B}}^{\dagger}_K$ is a field.  
\item[(4)] For every $0 < s \le r$,  we define ${\mathbf{C}}^{[s, r]}$ (resp.  ${\mathbf{C}}^{[s, r]}_K$) to be the completion of ${\mathbf{B}}^r$ (resp. ${\mathbf{B}}^r_K$) with respect to \(\max\{\lambda_{s}, \lambda_{r}\}\) (or equivalently, completion with respect to the family of norms $\{\lambda_t\}_{t\in [s,r]}$) inside $\widetilde{\mathbf{C}}^{[s, r]}$ (resp. $\widetilde{\mathbf{C}}^{[s, r]}_K$). Note that, for any $r'\in [s, r]$, the norm $\lambda_{r'}$ extends to $\mathbf{C}^{[s,r]}$ and $\mathbf{C}^{[s,r]}_K$.
\item[(5)] 
For every $r > 0$, we define ${\mathbf{C}}^{r}$ (resp. ${\mathbf{C}}^{r}_K$) to be the Fr\'echet completion of  ${\mathbf{B}}^r$ (resp. ${\mathbf{B}}^r_K$) inside $\widetilde{\mathbf{C}}^{r}$ (resp. $\widetilde{\mathbf{C}}^{r}_K$) with respect to the family of norms $\{\lambda_s\}_{s \in (0, r]}$.  
\item[(6)] Finally, define ${\mathbf C}^{\dagger}$ and ${\mathbf C}_K^{\dagger}$ as the unions ${\mathbf{C}}^{\dagger}\coloneqq \varinjlim_{r>0} {\mathbf{C}}^{r}$ and $
{\mathbf{C}}_K^{\dagger}\coloneqq \varinjlim_{r>0} {\mathbf{C}}^{r}_K$.
\end{enumerate}
\end{definition}

\begin{remark}\label{remark:completion_of_A_dagger}
\begin{enumerate}
\item[(i)] Since $\widetilde{\mathbf{A}}_K$ is the $p$-adic completion of $\widetilde{\mathbf{A}}_K^{\dagger}$, we deduce that $\mathbf{A}_K$ is the $p$-adic completion of $\mathbf{A}_K^{\dagger}$. 
\item[(ii)] For each $r > 0$, the natural map $\mathbf{C}^{r} \rightarrow \varprojlim_{0<s\le r} \mathbf{C}^{[s,r]}$ is an isomorphism. A similar isomorphism holds for $\mathbf{C}^{r}_K$ and $\mathbf{C}_K^{[s,r]}$. % (and thus the map $\mathbf{C}^{r}_K \isom \bigcap_{0<s<r}  {\mathbf{C}}^{[s,r]}_K$)
\item[(iii)] Remark~\ref{remark:subscript_K_as_invariant_perfect} continues to hold for imperfect period rings defined above. In particular, we have  $ {\mathbf{A}}_K^{\dagger} =   ( {\mathbf{A}}^{\dagger})^{H_K}$, $  {\mathbf{B}}_K^{\dagger} =   ( {\mathbf{B}}^{\dagger})^{H_K}$, and $ {\mathbf{C}}_K^{\dagger} =   ( {\mathbf{C}}^{\dagger})^{H_K}$. 
 \item[(iv)] Finally, 
 we remind the reader again that the rings ${\mathbf C}^{\dagger}$ and ${\mathbf C}_K^{\dagger}$ are denoted by ${\mathbf B}_{\mathrm{rig}}^{\dagger}$ and ${\mathbf B}_{\mathrm{rig}, K}^{\dagger}$ in \cite{Berger-differential} and \cite{Colmez08} (see Appendix \ref{appendix:notation}). 
\end{enumerate}
\end{remark}

\begin{remark}[Topology]\label{remark: topology--classical}
We endow $\widetilde{\mathbf{A}}^r_K$ with the topology induced by the $\lambda_r$-norm and endow $\widetilde{\mathbf{A}}^{\dagger}_K$ with the inductive limit topology (i.e., the finest topology making $\widetilde{\mathbf{A}}^r_K\hookrightarrow \widetilde{\mathbf{A}}^{\dagger}_K$ continuous for all $r>0$). Inverting $p$, we obtain natural topologies on $\widetilde{\mathbf{B}}^r_K$ and $\widetilde{\mathbf{B}}^{\dagger}_K$. The ring $\widetilde{\mathbf{C}}^{[s,r]}_K$ is endowed with the topology induced by $\max\{\lambda_r, \lambda_s\}$, which then induces a Fr\'echet topology on $\widetilde{\mathbf{C}}^r_K$. Finally, we equip $\widetilde{\mathbf{C}}^{\dagger}_K$ with the inductive limit topology (i.e., the LF topology). The topologies on the imperfect period rings are defined similarly.
\end{remark}

\begin{remark}[$(\varphi, \Gamma_K)$-structures]\label{remark: Frobenius--classical}
\begin{enumerate}
    \item[(i)]
The Frobenius $\varphi$ on $\widetilde{\mathbf{A}}_K$ (resp. on $\widetilde{\mathbf{A}}^{\dagger}_K$) preserves ${\mathbf{A}}_K$ (resp. ${\mathbf{A}}^{\dagger}_K$). For $0<s\le r$, we also have Frobenius maps
$\varphi\colon \mathbf{A}^r_K\rightarrow \mathbf{A}^{r/p}_K, \varphi\colon \mathbf{B}^r_K\rightarrow \mathbf{B}^{r/p}_K,  \varphi\colon \mathbf{C}^r_K\rightarrow \mathbf{C}^{r/p}_K, \varphi\colon \mathbf{C}^{[s,r]}_K\rightarrow \mathbf{C}^{[s/p, \,r/p]}_K$
as well as Frobenius maps on $\mathbf{B}^{\dagger}_K$ and $\mathbf{C}^{\dagger}_K$. We caution the reader that $\varphi\colon \mathbf{A}^r_K\rightarrow \mathbf{A}^{r/p}_K$ is no longer bijective.
\item[(ii)] %\label{remark: Gamma action--classical}
All period rings $\widetilde{\mathbf{A}}^r_K$, $\widetilde{\mathbf{A}}^{\dagger}_K$, $\widetilde{\mathbf{B}}^r_K$, $\widetilde{\mathbf{B}}^{\dagger}_K$, $\widetilde{\mathbf{C}}^{[s,r]}_K$, $\widetilde{\mathbf{C}}^r_K$, $\widetilde{\mathbf{C}}^{\dagger}_K$, as well as their imperfect counterparts, are equipped with continuous actions of $\Gamma_K$, where $\Gamma_K$ is endowed with the profinite topology. The $\Gamma_K$-actions commute with the Frobenius actions.
\end{enumerate}
\end{remark}

\begin{notation}
Let $t = \log (1 + \pi)$. This is an element that lies in $\mathbf{C}^{\dagger}_K$ but not in $\mathbf{B}^{\dagger}_K$.
\end{notation}

%For later purpose, let us record the following result from \cite{Colmez08}. 

\begin{lemma}[{\cite[Lemme 6.5]{Colmez08}}]\label{lemma:pi vs bar-pi}
Let
\[r_K=\begin{cases}
        (2v(\mathfrak{d}_{\mathbf{E}_K/\mathbf{E}_{\mathbb{Q}_p}}))^{-1} & \text{if } \mathbf{E}_K/\mathbf{E}_{\mathbb{Q}_p} \text{ is ramified,}\\
        1 & \text{if } \mathbf{E}_K/\mathbf{E}_{\mathbb{Q}_p} \text{ is unramified.}
    \end{cases}\]
Then for all $0<r<r_K$, we have $\pi_K \in \mathbf{A}_K ^r$ and the element $u=\frac{\pi_K}{[\overline{\pi}_K]}$ is a unit in the ring of integers of $\widetilde{\mathbf{A}}_K^r$ (with respect to the norm $\lambda_r$). More precisely, if we write $u=\sum_{n=0}^{\infty}p^n[u_n]$ with $u_n\in \widehat{K}_{\infty}^{\flat}$, then $u_0=1$ and $rv(u_n)+n>0$ for all $n\ge 1$.
\end{lemma}

In particular, for $0 < r < r_K$, we know that $\lambda_r (\pi_K) =  |\overline \pi_K|^r = p^{-\frac{p}{p-1} \cdot \frac{r}{e_K}}$.

\subsubsection{Explicit (geometric) descriptions}

There are explicit descriptions of imperfect period rings in terms of double-ended power series subject to certain convergence conditions.

\begin{proposition} \label{prop:identificaiton_for_A^r}
If $0<r<r_K$, there is an isomorphism
\begin{equation*}
\left\{\sum_{i\in \mathbb{Z}}c_iT^i\,\Big|\,c_i\in \mathcal{O}_{K'_0}, \,\, |c_i|\cdot p^{-\frac{p}{p-1} \cdot \frac{ir}{e_K}}\rightarrow 0 \textrm{ as }i\rightarrow -\infty\right\}  \isom   \mathbf{A}_K^r
\end{equation*}
sending $f\mapsto f(\pi_K)$. Moreover, the isomorphism preserves the norms on both sides, where the norm on the left hand side is $\sup_{i\in \mathbb{Z}} \left\{|c_i|\cdot p^{-\frac{p}{p-1} \cdot \frac{ir}{e_K}}\right\}$ and the one on the right hand side is $\lambda_r$. 
\end{proposition}

\begin{proof}
This is \cite[Proposition II.2.1]{Cherbonnier_Colmez} (also see \cite[Proposition 1.4]{Berger-differential} and \cite[Proposition 7.5(i)]{Colmez08}).
\end{proof}

\begin{remark} \label{remark:completion_of_A_dagger_2} An immediate consequence is that, for $0 < r < r_K$, the isomorphism in Proposition~\ref{prop:identificaiton_for_A^r} induces isomorphisms
\[ 
\mathbf{A}^r_K/(p) \cong \mathbf{A}^\dagger_K/(p) \cong \mathbf{A}_K/(p) \cong \kappa_K (\!(\overline \pi_K)\!) \cong \mathbf{E}_K  
\]
In particular, $\mathbf{A}_K$ is the $p$-adic completion of $\mathbf{A}^r_K$ for any $r \in (0, r_K)$. This gives another way to recover $\mathbf{A}_K$ as the $p$-adic completion of $\mathbf{A}^\dagger_K$ (cf. Remark~\ref{remark:completion_of_A_dagger}). 
\end{remark}

We also have explicit descriptions of $\mathbf{C}_K^r$ and $\mathbf{C}_K^{[s,r]}$. 
 
\begin{proposition}\label{prop:identification_C_with_annulus}
\begin{enumerate}
\item[(1)]
For $0<r<r_K$, the isomorphism in Proposition~\ref{prop:identificaiton_for_A^r} induces an isomorphism \[\left\{\sum_{i\in \mathbb{Z}}c_iT^i\,\Big|\,c_i\in K'_0, \,\,|c_i|\cdot p^{-\frac{p}{p-1} \cdot \frac{ir'}{e_K}}\rightarrow 0 \textrm{ as }i\rightarrow \pm\infty \textrm{ for any } r'\in (0, r] \right\}\isom  \mathbf{C}_K^r\]
sending $f\mapsto f(\pi_K)$.
In other words, we can identify $\mathbf{C}_K^r$ with the ring of analytic functions $\Gamma_{K'_0}([{\delta, 1});\pi_K)$ on the half-open annulus $A_{K'_0}([{\delta, 1}); \pi_K)$ with open outer boundary at radius $1$ and closed inner boundary at radius $  \delta = p^{- \frac{p}{p-1} \cdot \frac{r}{e_K}}$ (see the notation in \S \ref{sss:notation_cyclo}). 
\item[(2)] Similarly, for any $0<s\le r<r_K$, we have an identification $\Gamma_{K'_0}([\delta_1, \delta_2]; \pi_K) \isom \mathbf{C}_K^{[s, r]}$  
where $\delta_1 = p^{- \frac{p}{p-1} \cdot \frac{r}{e_K}}$ and $\delta_2 = p^{-\frac{p}{p-1} \cdot \frac{s}{e_K}}$. Under this isomorphism, for any $r' \in [s, r]$, the $\rho$-Gauss norm on $\Gamma_{K'_0}([\delta_1, \delta_2]; \pi_K)$ with $\rho = p^{- \frac{p}{p-1} \cdot \frac{r'}{e_K}}$ coincides with the norm $\lambda_{r'}$ on $\mathbf{C}_K^{[s, r]}$.
\item[(3)] Consequently, the ring $\mathbf{C}_K^{\dagger}$ can be identified with the Robba ring $\mathscr{R}_{K'_0}$ with coefficients in $K'_0$ (cf. \S \ref{sss:notation_cyclo}).
\end{enumerate}
\end{proposition}

\begin{proof}
These follow from Proposition~\ref{prop:identificaiton_for_A^r}, as well as \cite[Proposition 7.5 \& 7.6]{Colmez08} and their proofs. 
\end{proof}

\subsubsection{The ring $\mathbf{C}_{\log}^{\dagger}$ and $\mathbf{C}_{\log, K}^{\dagger}$}\label{sss:C log}

To study semistable representations, we also need the ``logarithmic variants'' of the imperfect period rings (cf. \cite[\S 2.6]{Berger-differential}).\footnote{The ring $\mathbf{C}^{\dagger}_{\log}$ (resp. $\mathbf{C}^{\dagger}_{\log, K}$) is denoted by $\mathbf{B}^{\dagger}_{\log}$ (resp. $\mathbf{B}^{\dagger}_{\log, K}$) in \cite{Berger-differential}.}

\begin{definition}
For $0 < s \le r$, we define 
\[
 {\mathbf{C}}^{r}_{\log, K} \coloneqq  \mathbf{C}^{r}_{K}[\log \pi], \quad 
  {\mathbf{C}}^{\dagger}_{\log, K} \coloneqq  \mathbf{C}^{\dagger}_{K}[\log \pi] = \varinjlim_{r>0}    {\mathbf{C}}^{r}_{\log, K} 
\]
equipped with commuting $\varphi$ and $\Gamma_{K}$-actions extending the ones on %$\mathbf{C}^{\dagger}$ 
$\mathbf{C}^{\dagger}_K$. The Frobenius map $\varphi$ on $\log \pi$ is induced from $\varphi(\pi) = (\pi+1)^p - 1$ and can be written as 
$\varphi (\log \pi) = p \log \pi + \log (\varphi(\pi)/\pi^p).$  The $\Gamma_K$-action on $\log \pi$ is given by $\gamma (\log \pi) = \log \pi + \log (\gamma (\pi)/\pi)$  for $\gamma \in \Gamma_K$.   
\end{definition}

%%% 3.3
\vspace{0.1in}
\subsection{$(\varphi, \Gamma_K)$-modules and $p$-adic differential equations} \label{ss:classical_pDE}
\noindent
\vspace{0.1in}

\noindent Now we review the theory of overconvergent $(\varphi, \Gamma_K)$-modules associated with $p$-adic Galois representations, as well as Berger's construction of $p$-adic differential equations associated with \emph{de Rham} representations.  

\subsubsection{($\varphi, \Gamma_K$)-modules}\label{sss:overconvergent classical}
For every $\Q_p$-representation $V$ of $\Gal_K$, 
one can associate a $(\varphi, \Gamma_K)$-module $D^{\dagger}_{K}(V) \coloneqq   (V\otimes_{\mathbb{Q}_p}\mathbf{B}^{\dagger})^{H_K}$ over $\mathbf{B}^{\dagger}_K$, as well as $(\varphi, \Gamma_K)$-modules over $\mathbf{C}^{\dagger}_K$ and $\mathbf{C}^{\dagger}_{\log, K}$ given by 
\begin{align*}
D^{\dagger}_{\mathrm{rig},K}(V)   \coloneqq  D^{\dagger}_{K}(V)\otimes_{\mathbf{B}^{\dagger}_{K}} \mathbf{C}^{\dagger}_{K},   \qquad  D^{\dagger}_{\log, K}(V) \coloneqq  D^{\dagger}_{K}(V)\otimes_{\mathbf{B}^{\dagger}_{K}} \mathbf{C}^{\dagger}_{\log, K}.
\end{align*}

\begin{theorem}[{\cite{Cherbonnier_Colmez}, \cite[Theorem 9.5.6]{Kedlaya-LiuI}}]\label{thm: Galois representations and phi Gamma modules}
The functors $D^{\dagger}_K$ and  $D^{\dagger}_{\mathrm{rig}, K}$ induce equivalences of categories
\begin{align*}
\: \: \:  \mathrm{Rep}_{\mathbb{Q}_p}(\Gal_{K}) & \isom \big\{\textup{\'etale }(\varphi, \Gamma_{K})\textup{-modules over }\mathbf{B}^{\dagger}_{K}\big\} 
 \\
& \isom \big\{\textup{\'etale }(\varphi, \Gamma_{K})\textup{-modules over }\mathbf{C}^{\dagger}_{K}\big\}
\end{align*} 
between the category of finite dimensional continuous $\Gal_K$-representations over $\mathbb{Q}_p$-vector spaces and the category of (finite free) \'etale  $(\varphi, \Gamma_K)$-modules over  $\mathbf{B}_K^{\dagger}$ (resp. $\mathbf{C}_K^{\dagger}$).\footnote{
Recall that a $(\varphi, \Gamma_K)$-module over  $\mathbf{B}_K^{\dagger}$ is a finite free $\mathbf{B}_K^{\dagger}$-module $M$  with a semilinear  Frobenius  action $\varphi$ inducing an isomorphism  
$\varphi \otimes \mathrm{id}\colon M \otimes_{\mathbf{B}_K^{\dagger}, \varphi} \mathbf{B}_K^{\dagger} \isom M$ and a continuous semilinear $\Gamma_K$-action that commutes with $\varphi$. There is a similar notion over $\mathbf{C}_K^{\dagger}$. Such a $(\varphi, \Gamma_K)$-module over (either $\mathbf{B}_K^{\dagger}$ or $\mathbf{C}_K^{\dagger}$) is \emph{\'etale} if there is a finite free module $M_0$ over $\mathbf{A}_K^{\dagger}$ equipped with a continuous semilinear $\Gamma_K$-action and a Frobenius semilinear endomorphism $\varphi\colon M_0 \ra M_0$ commuting with $\Gamma_{K}$, such that the induced linear map $\varphi \otimes \mathrm{id}\colon M_0 \otimes_{\mathbf{A}_K^{\dagger}, \varphi} \mathbf{A}_K^{\dagger} \isom M_0 $ is an isomorphism, and that $M$ is the base change of $M_0$ as a $(\varphi, \Gamma_K)$-module.
}
\end{theorem}

\begin{remark}\label{rmk: descent to alpha}
Let $M^{\dagger} = D^{\dagger}_{\mathrm{rig},K} (V)$ be the \'etale $(\varphi, \Gamma_K)$-modules over $\mathbf{C}_K^{\dagger}$ associated with some $V\in \mathrm{Rep}_{\mathbb{Q}_p}(\Gal_{K})$. Then there exists $r > 0$ (and hence for all $r$ sufficiently close to $0$) such that $M^{\dagger}\cong M^{r}\otimes_{\mathbf{C}_K^{r}}\mathbf{C}_K^{\dagger}$ for a unique finite free \'etale $(\varphi, \Gamma_K)$-module $M^{r}$ over $\mathbf{C}_K^r$.\footnote{
This means that there exists a finite free module $M^r_0$ over $\mathbf{A}_K^r$ equipped with a continuous semilinear $\Gamma_K$-action and an isomorphism
$\varphi\colon M^r_0 \otimes_{\mathbf{A}_K^r, \varphi} \mathbf{A}_K^{r/p}\isom M^r_0 \otimes_{\mathbf{A}_K^r} \mathbf{A}_K^{r/p}
$ 
commuting with the $\Gamma_K$-action, where the second tensor is along the natural inclusion $\mathbf{A}_K^r\subset \mathbf{A}_K^{r/p}$, such that $M^r$ is the base change of $M_0^r$.}
\end{remark}

\begin{remark}\label{rmk: phi Gamma modules over perfect period rings--classical}
One can also consider $(\varphi, \Gamma_K)$-modules over perfect period rings. Indeed, we have equivalences of categories between the category $ \mathrm{Rep}_{\mathbb{Z}_p}(\Gal_{K}) $ of continuous $\Gal_K$-representations over finite free $\Z_p$-modules, that of $(\varphi, \Gamma_{K})$-modules over $\widetilde{\mathbf{A}}^{\dagger}_{K}$, and that of $(\varphi, \Gamma_{K})$-modules over $\widetilde{\mathbf{A}}_{K}$. Similarly, we have equivalences of categories between $\mathrm{Rep}_{\mathbb{Q}_p}(\Gal_{K})$, the category of \'etale $(\varphi, \Gamma_{K})$-modules over $\widetilde{\mathbf{B}}^{\dagger}_{K}$, and that of \'etale $(\varphi, \Gamma_{K})$-modules over $\widetilde{\mathbf{C}}^{\dagger}_{K}$. 
For $V\in \mathrm{Rep}_{\mathbb{Z}_p}(\Gal_{K})$ or $\mathrm{Rep}_{\mathbb{Q}_p}(\Gal_{K})$, let $\widetilde{D}^{\mathrm{int}}_K(V)$, $\widetilde{D}^{\mathrm{int}, \dagger}_K(V)$, $\widetilde{D}^{\dagger}_K(V)$, $\widetilde{D}^{\dagger}_{\mathrm{rig}, K}(V)$ denote the associated $(\varphi, \Gamma_K)$-modules over $\widetilde{\mathbf{A}}_K$, $\widetilde{\mathbf{A}}^{\dagger}_K$, $\widetilde{\mathbf{B}}^{\dagger}_K$, $\widetilde{\mathbf{C}}^{\dagger}_K$, respectively. In particular, $\widetilde{D}^{\dagger}_K(V)$ (resp. $\widetilde{D}^{\dagger}_{\mathrm{rig}, K}(V)$) is the base change of $D^{\dagger}_K(V)$ (resp. $D^{\dagger}_{\mathrm{rig}, K}(V)$) along $\mathbf{B}^{\dagger}_K\rightarrow \widetilde{\mathbf{B}}^{\dagger}_K$ (resp. $\mathbf{C}^{\dagger}_K\rightarrow \widetilde{\mathbf{C}}^{\dagger}_K$). We refer the reader to \cite{Cherbonnier_Colmez} and \cite{Kedlaya-LiuI} for more details.
\end{remark}

\subsubsection{Differential modules} \label{sss:pDE_over_K}

Before proceeding further, let us  recall some basic notions on $p$-adic differential equations. Firstly, a \emph{differential ring} is a commutative ring $R$ equipped with a derivation $d\colon R \ra R$ (cf. \cite[\S 5]{kedlaya_book}), that is, an additive map satisfying the Leibniz rule $d (ab) = a d(b) + b d(a)$ for all $a, b \in R$. 
A \emph{differential module} over a differential ring $(R, d)$ is an $R$-module $M$ equipped with a differential operator $D\colon M \ra M$ compatible with $d$, that is, an additive map satisfying $D(ax) = a D(x) + d(a) x$ for all $a \in R, x \in M$. When $R$ is equipped with more than one differential operator $d_1, \ldots, d_m$, we abuse the terminology and say that $(M, D_1, \ldots, D_m)$ is a \emph{differential module} over $(R, d_1, \ldots, d_m)$ if $(M, D_i)$ is a differential module over $(R, d_i)$ for all $i=1,\ldots m$.

We are primarily interested in the Robba ring $\mathbf{C}^{\dagger}_K$ on which the action of the Lie algebra of $\Gamma_K$ induces a differential operator $\nabla$ defined by $\nabla = \log (\gamma)/\log (\chi (\gamma))$ for any $\gamma$ sufficiently close to $1$, where $\chi\colon\Gamma_K\rightarrow \Z_p^{\times}$ is the usual cyclotomic character. By \cite[Lemme 4.2]{Berger-differential}, one has 
\begin{equation} \label{eq:nabla_vs_partial}
\nabla = t \cdot \partial = t \cdot (1+\pi) \cdot \frac{d}{d\pi}  
\end{equation} 
where $t=\log (1+\pi)$ and $\partial =  (1+\pi) \cdot \frac{d}{d\pi}$ is a well-defined differential operator on $\mathbf{C}^{\dagger}_K$. Then a differential module $(M, \partial_M)$ over $(\mathbf{C}^{\dagger}_K, \partial)$ can be viewed as a connection 
\[
\nabla_{M}\colon M \ra M \otimes_{\mathbf{C}^{\dagger}_K} \Omega^1_{\mathbf{C}^{\dagger}_K}
\]
where $\nabla_M (x) = \partial_M (x) \cdot d \log (1+ \pi) =  \frac{\partial_M (x)}{1 + \pi} \cdot d \pi$.\footnote{Note that this is a connection extending $d x = \frac{\partial x}{1+ \pi} d\pi$.} A finite free differential module over $(\mathbf{C}^{\dagger}_K, \partial)$ is known as a \emph{$p$-adic differential equation} in \cite{Berger-differential}. 

\begin{remark}\label{remark:partial vs d/d pi_K}
In practice, one may replace $(\mathbf{C}^{\dagger}_K, \partial)$ with the differential ring $(\mathbf{C}^{\dagger}_K, \frac{d}{d\pi_K})$ without changing the underlying theory of $p$-adic differential equations, as $\partial=u_0\cdot \frac{d}{d\pi_K}$ for some nonzero element $u_0\in \mathbf{B}^{\dagger}_K$ (hence, a unit in $\mathbf{C}^{\dagger}_K$).
\end{remark}

\begin{definition} \label{def:p-adic_diff_eqn_with_Frob}
 Let $(M, \partial_M)$ be a finite free differential module over $(\mathbf{C}^{\dagger}_K, \partial)$.  A \emph{Frobenius structure} on $M$ is a $\varphi$-semilinear map $\varPhi\colon M \ra  M$ such that 
\begin{enumerate}
    \item[(1)] $\varPhi$ induces an isomorphism $\varphi^*(M) \isom M$, i.e.,  $\varPhi(M)$ generates $M$ as a $\mathbf{C}_{K}^{\dagger}$-module; 
    \item[(2)] $\varPhi$ is compatible with the differential structures; namely, the diagram
     \[ \begin{tikzcd}
         M \arrow[r, "\nabla_M"] \arrow[d, swap, "\varPhi"] & M \otimes_{\mathbf{C}^{\dagger}_K} \Omega^1_{\mathbf{C}^{\dagger}_K} \arrow[d, "\varPhi \otimes d \varphi"]  \\
         M \arrow[r, "\nabla_M"]  & M  \otimes_{\mathbf{C}^{\dagger}_K} \Omega^1_{\mathbf{C}^{\dagger}_K}
     \end{tikzcd} \]
    commutes. Equivalently, $\varPhi$ satisfies $\partial_M \circ \varPhi = p \,\varPhi \circ \partial_M$.    
\end{enumerate}
\end{definition}

\subsubsection{$p$-adic differential equations attached to de Rham representations} \label{sss:Berger_N_dR}

Let $D^{\dagger}_{\mathrm{rig}, K} (V)$ be the \'etale $(\varphi, \Gamma_K)$-module over  $\mathbf{C}^{\dagger}_K$ associated with some $V\in \mathrm{Rep}_{\mathbb{Q}_p}(\Gal_{K})$.
The action of the Lie algebra of $\Gamma_K$ induces a differential operator $\nabla_V =   \log (\gamma)/\log (\chi (\gamma))$ on $D^{\dagger}_{\mathrm{rig}, K} (V)$. To obtain a differential module over $(\mathbf{C}_{K}^{\dagger}, \partial)$, we consider 
\[ 
\partial_V = \frac{1}{t} \nabla_V = \frac{1}{\log (1+\pi)} \nabla_V.
\] 
In general, this operator does not preserve $D^{\dagger}_{\mathrm{rig}, K} (V)$, as it might introduce poles at the vanishing locus of $t$. The main technical result of \cite{Berger-differential} states that, when $V$ is de Rham, these poles can be removed systematically to produce an actual $p$-adic differential equation over $(\mathbf{C}^{\dagger}_K, \partial)$.

\begin{theorem}[{\cite[Th\'eor\`eme 5.10]{Berger-differential}}] \label{thm:Berger_N_dR}
    Let $V$ be a de Rham representation of $\Gal_K$ of dimension $d$ with non-negative Hodge--Tate weights.\footnote{For general de Rham representations, one can first twist the representations by powers of cyclotomic character to make the Hodge-Tate weights non-negative.} Then there exists a unique free $\mathbf{C}^{\dagger}_K$-submodule $\mathbf{N}_{\mathrm{dR}} (V) \subset D^{\dagger}_{\mathrm{rig}, K} (V)$ of rank $d$ that is stable under $\partial_{V}$.  Moreover, $\mathbf{N}_{\mathrm{dR}}(V)$ is stable under $\varphi$ and $ \Gamma_K$-actions,\footnote{Note that, as a $(\varphi, \Gamma_K)$-module, $\mathbf{N}_{\mathrm{dR}}(V)$ is not necessarily \'etale.} and, consequently, it is a finite free differential module over $(\mathbf C^{\dagger}_K, \partial)$ with a Frobenius structure. Furthermore, the inclusion $\mathbf{N}_{\mathrm{dR}} (V) \subset D^{\dagger}_{\mathrm{rig}, K} (V)$ induces an isomorphism after inverting $t$.
\end{theorem}

\begin{remark}
The proof of uniqueness in Theorem~\ref{thm:Berger_N_dR} crucially uses a result of Lazard \cite{Lazard} which states that $\mathbf{C}^{\dagger}_K$ is a B\'ezout domain. 
\end{remark}

\begin{remark}\label{remark: classical NdR general HT wt}
The construction of $\mathbf{N}_{\mathrm{dR}}(V)$ can be generalized to a functor
\[
\mathbf{N}_{\mathrm{dR}}:\mathrm{Rep}^{\mathrm{dR}}_{\mathbb{Q}_p}(\Gal_{K})  \lra \big\{(\varphi, \Gamma_{K})\textrm{-modules over }\mathbf{C}^{\dagger}_{K}\big\},
\]
from the category of all de Rham representations of $\Gal_K$. We refer the reader to \cite[\S 3.4.2]{relative} for details.
\end{remark}

%%%
\subsubsection{The effect of finite extensions of $K$}
For later use, let us record the following compatibility result. 
Let $L/K$ be a finite extension. Then there are natural $(\varphi, \Gamma_L)$-equivariant injections $\mathbf{C}^{r}_K\rightarrow \mathbf{C}^{r}_L$ and $\mathbf{C}^{\dagger}_K\rightarrow \mathbf{C}^{\dagger}_L$.

\begin{proposition}[{\cite[Proposition 3.44]{relative}}]\label{prop: NdR finite extenion of K--classical}  
 Let $V\in \mathrm{Rep}^{\mathrm{dR}}_{\Q_p}(\Gal_K)$ and let $L/K$ be a finite extension. Then there is a natural isomorphism
\[\mathbf{N}_{\mathrm{dR}}(V)\otimes_{\mathbf{C}^{\dagger}_K}\mathbf{C}^{\dagger}_L\isom  \mathbf{N}_{\mathrm{dR}}(V|_{\Gal_L})\]
of $(\varphi, \Gamma_L)$-modules; hence it is also an isomorphism of differential modules over $(\mathbf{C}^{\dagger}_L, \partial)$.
\end{proposition}

%%% 3.4
\vspace{0.1in}
\subsection{The classical $p$-adic monodromy theorem}\label{ss: proof of classical p-adic monodromy theorem}
\noindent
\vspace{0.1in}

\noindent To conclude this section, we recall the proof of Berger’s classical $p$-adic monodromy theorem (\cite[Th\'eor\`eme 0.7]{Berger-differential}), which asserts that every de Rham representation is potentially semistable. The first step is the following result of Berger, which recognizes the $F$-isocrystal $D_{\mathrm{st}}(V)$ from the $(\varphi, \Gamma_K)$-module $D^{\dagger}_{\log, K}(V)$. 
\begin{theorem}[{\cite[Th\'eor\`eme 3.6]{Berger-differential}}]\label{theorem: Berger 3.6}
For $V\in \mathrm{Rep}_{\mathbb{Q}_p}(\Gal_{K})$, we have
\[
D_{\mathrm{crys}}(V) =\Big(D^{\dagger}_{\mathrm{rig}, K}(V)[\frac{1}{t}]\Big)^{\Gamma_{K}}\qquad\textrm{ and}\qquad
D_{\mathrm{st}}(V)  =\Big(D^{\dagger}_{\log, K}(V)[\frac{1}{t}]\Big)^{\Gamma_{K}}
\]
as subspaces of $(V\otimes_{\Q_p}\widetilde{\mathbf{C}}^{\dagger} [\frac{1}{t}])^{\Gal_K}$ and  $(V\otimes_{\Q_p}\widetilde{\mathbf{C}}_{\log}^{\dagger} [\frac{1}{t}])^{\Gal_K}$ respectively. 
\end{theorem}
 
The second ingredient (perhaps the most important one) is a key input from Andr\'e--Kedlaya--Mebkhout. It follows from the so-called ``$p$-adic local monodromy theorem'', independently proved by Andr\'e, Kedlaya, and Mebkhout (see \cite{Andre_p-adic_monodormy, Kedlaya_p-adic_monodormy, Mebkhout_p-adic_monodormy}).

\begin{theorem}[Andr\'e, Kedlaya, Mebkhout] \label{thm:all_pDE_with_Frob_are_quasiunipotent}
Let  $(M, \partial_M)$ be a differential module over  $(\mathbf{C}^{\dagger}_K, \partial)$  equipped with a Frobenius structure. Then there exists a finite extension $L/K$ such that $M\otimes_{\mathbf{C}^{\dagger}_K}\mathbf{C}^{\dagger}_L$ becomes unipotent with respect to $\partial_M$. In particular, $(M, \partial_M)$ admits a full basis of solutions in $\mathbf{C}^{\dagger}_{\log, L}=\mathbf{C}^{\dagger}_{L} [\log \pi]$.  
\end{theorem}

Now, we are ready to finish the proof.

\begin{theorem}[The classical $p$-adic monodromy theorem]
Let $V\in \mathrm{Rep}^{\mathrm{dR}}_{\Q_p}(\Gal_K)$. Then there exists a finite extension $L/K$ such that $V|_{\Gal_L}$ is a semistable representation of $\Gal_L$.
\end{theorem}

\begin{proof}
Without loss of generality, assume that $V$ is a de Rham representation of dimension $d$ with non-negative Hodge--Tate weights. Let $(\mathbf{N}_{\mathrm{dR}}(V), \partial_V)$ be the rank $d$ differential module over $(\mathbf C^{\dagger}_K, \partial)$ with a Frobenius structure as in Theorem~\ref{thm:Berger_N_dR} and Remark~\ref{remark: classical NdR general HT wt}. By Theorem~\ref{thm:all_pDE_with_Frob_are_quasiunipotent}, there exists a finite extension $L/K$ such that 
\[
    (\mathbf{N}_{\mathrm{dR}} (V) \otimes_{\mathbf{C}^{\dagger}_{K}} \mathbf{C}^{\dagger}_{\log, L} )^{\partial_V = 0}
\]
has full rank. Therefore, there exists a large enough $m$ such that, after passing to the finite extension $L_m = L(\zeta_{p^m})$ over $L$,  the $\Gamma_{L_m}$-invariant subspace
\[(\mathbf{N}_{\mathrm{dR}} (V) \otimes_{\mathbf{C}^{\dagger}_{K}} \mathbf{C}^{\dagger}_{\log, L_m} )^{\Gamma_{L_m}} \] has full rank. By Theorem~\ref{theorem: Berger 3.6} and the fact that $\mathbf{N}_{\mathrm{dR}} (V)$ coincides with $D^{\dagger}_{\mathrm{rig}, K} (V)$ after inverting $t$, we conclude that $D_{\mathrm{st}} (V|_{\Gal_{L_m}})$ has full rank and thus $V|_{\Gal_{L_m}}$ is semistable. 
\end{proof}

\vspace*{0.1in}

%%% Section 4
\section{Relative period rings and relative $(\varphi, \Gamma)$-modules}  \label{section:relative_period_rings}
 
The goal of this section is to summarize the results from \cite{relative} on imperfect relative period rings and relative $(\varphi, \Gamma)$-modules over such rings.\footnote{Although this article only requires the theory of relative $(\varphi, \Gamma)$-modules over a circle, we present the theory in a more general setting for the sake of completeness.} These generalize the $(\varphi, \Gamma_K)$-module theory in the classical setting. To maintain the flow of the exposition, we present only those definitions and theorems required for this article and refer the reader to \emph{loc.\ cit.} for complete proofs. As this section is somewhat technical, the reader may wish to skip the technical details on a first reading and return to them later if necessary. This section is organized as follows.

\begin{itemize}
\item In \S \ref{ss:toric_tower} and \S \ref{ss:relative_period_rings:perfect}, we recall the notion of toric tower and perfect period rings following \cite{Kedlaya-LiuI}.
\item In \S \ref{sec:weakly_decompleting_towers} and \S \ref{ss:relative_phi_gamma_mod}, we introduce imperfect period rings associated with \emph{weakly decompleting data}, as well as the notion of $(\varphi, \Gamma)$-modules over such period rings. We also collect some descent results under the assumption of weakly decompleting data. (cf. \cite[\S 5]{relative}.)
\item In \S \ref{sec:decompleting_towers}, we work with a more restrictive condition called \emph{decompleting data} and present an important descent theorem for $(\varphi, \Gamma)$-modules over $\mathbf{C}$-type period rings. (cf. \cite[\S 6 -- 7]{relative}.)
\item In \S \ref{sec:strongly decompleting data} and \S \ref{ss:relative NdR}, we further upgrade the condition to \emph{strongly decompleting data}. This notion is an entirely new development in \cite{relative}, which provides the key ingredient for extending Berger's  $\mathbf{N}_{\mathrm{dR}}$ functor to the relative setting. (cf. \cite[\S 8 -- 10]{relative}.)\footnote{We refer the reader to \S \ref{ss: annulus--pDE} for an explicit description in the case of circles; in this case, the relative $\mathbf{N}_{\mathrm{dR}}$ functor gives rise to $p$-adic differential equations equipped with the ``correct'' differential operators.}
\item Finally, in \S \ref{ss: examples of strongly decompleting data}, we provide an explicit description of imperfect period rings in the case of annuli. 
\end{itemize}

%%% 4.1
\vspace{0.1in}
\subsection{Toric chart and toric tower}  \label{ss:toric_tower}
\noindent
\vspace{0.1in}

\noindent
Throughout \S \ref{section:relative_period_rings}, let $A_X$ be a reduced affinoid $K$-algebra (in the sense of Tate's rigid analytic geometry) and let $X=\spa (A_X, A_X^+)$ be the associated adic space over $\spa(K,\mO_K)$. In particular, $A_X^+=A_X^{\circ}$. 
For any integer $d\ge 1$, let
$\D_K^d \coloneqq \spa (K \gr{T_1, \ldots, T_d}, \mO_K \gr{T_1, \ldots, T_d})$ 
denote the $d$-dimensional closed unit polydisk over $K$, equipped with the log structure associated with the normal crossing divisor $\{T_1\cdots T_d=0\}$.

\begin{definition}\label{defn: toric chart}
A \emph{toric chart} of $X$ is an \'etale morphism $h\colon X \ra \D_K^d$ of adic spaces (for some $d\ge 1$) which is a composition of rational localizations and finite \'etale maps.
\end{definition}

For each $n \ge 0$, let
\[
\D^d_{K, n} \coloneqq  \spa \left(K_n \gr{T_1^{\frac{1}{p^n}}, \ldots, T_d^{\frac{1}{p^n}}}, \mO_{K_n} \gr{T_1^{\frac{1}{p^n}}, \ldots, T_d^{\frac{1}{p^n}}}\right). 
\]
This is the $d$-dimensional closed unit polydisk over $K_n$, with coordinates $T_1^{\frac{1}{p^n}}, \ldots, T_d^{\frac{1}{p^n}}$; we equip it with the log structure associated with the normal crossing divisor $\{T_1^{\frac{1}{p^n}}\cdots T_d^{\frac{1}{p^n}}=0\}$. When the base field $K$ and dimension $d$ are understood, we often write $\D_{n} = \D^d_{K, n}$ to simplify the notation. The natural maps $\psi_n\colon \D_n \ra \D^d_K$ (induced by the natural inclusion of rings) put together a tower 
$\psi\colon \cdots \ra \D_n\ra \cdots \ra\D_1\ra \D^d_K$  of Kummer \'etale maps. 
This tower is a Galois \emph{pro-Kummer-\'etale tower} in the sense of \cite[Definition 6.1.2]{DLLZ1}, with Galois group $\Gamma \cong \Gamma_{\mathrm{geom}} \rtimes \Gamma_K$ where $\Gamma_{\mathrm{geom}}$ is the Galois group of 
\[\D_{\infty}=\varprojlim_{n\ge 0} \D_n\rightarrow \varprojlim_{n \ge 0}\D^d_K\times_{\spa(K, \mO_K)} \spa(K_n, \mO_{K_n}).\]
In particular, $\Gamma_{\mathrm{geom}}$ acts on $K_{\infty}$ trivially and permutes the $p$-power roots of $T_i$'s. On the other hand, $\Gamma_K = \Gal(K_\infty/K)$ acts on the $p$-power roots of unity via the usual cyclotomic character $\chi$ and acts trivially on the prescribed $p$-power roots $T_i^{1/p^n}$'s of $T_i$.\footnote{In fact, a choice of compatible system of $p$-power roots of $T_i$'s determines a splitting $\Gamma\cong \Gamma_{\mathrm{geom}}\rtimes \Gamma_K$.} Note that $\Gamma_{\mathrm{geom}} \cong \Z_p^{\oplus d}$ and $\Gamma_K \subset \Z_p^\times$ is an open subgroup and for any $\sigma \in \Gamma_K$ and $\tau \in \Gamma_{\mathrm{geom}}$, we have  $\sigma \tau \sigma^{-1} = \tau^{\chi(\sigma)}.$ 
Finally, let us note that the tower is ``log perfectoid'' with associated perfectoid space 
\[ 
\widehat{\D}_{\infty} = \spa (\widehat{K}_\infty \gr{T_1^{\frac{1}{p^\infty}}, \ldots, T_d^{\frac{1}{p^\infty}}}, \mO_{\widehat{K}_\infty} \gr{T_1^{\frac{1}{p^\infty}}, \ldots, T_d^{\frac{1}{p^\infty}}}).
\]

\begin{definition}\label{defn: toric tower X} Let $h\colon X\ra \D^d_K$ be a toric chart. For each $n\ge 0$, let $X_n\coloneqq X \times_{\D^d_K} \D_n.$ The tower $\psi=\{\psi_n\colon\D_n\ra\D^d_K\}$ then pulls back to a tower $\psi = \{\psi_n\colon X_n \ra X\}$ over $X$, which we refer to as the \emph{toric tower} over $X$ relative to the toric chart $h$. 
\end{definition} 

The pro-Kummer \'etale cover $\psi\colon X_{\infty}=\varprojlim X_n \ra X$ also admits an action of $\Gamma$. We write $\Gamma_X$ whenever it is necessary to emphasize the role of $X$. The tower $\psi$ is \emph{log perfectoid} with associated perfectoid space $\widehat{X}_\infty = X \times_{\D^d_K} \widehat{\D}_\infty$. If the image of $h$ does not intersect with the locus $\{T_1\cdots T_d=0\}$, then each $\psi_n$ is finite \'etale, thus the tower $\psi\colon X_{\infty} \ra X$ is an affinoid perfectoid object in the pro-\'etale site $X_{\proet}$ of $X$.

\begin{definition} \label{construction:A_psi_X_infty}
\begin{enumerate}
\item[(1)] For every $n\ge 0$, write $X_n = \spa (A_{X, n}, A_{X, n}^+)$. In particular, we have $ A_{X,0} =A_X$ and $A_{X,0}^+ =A_X^+$. Each $A_{X, n}$ is equipped with the spectral norm $|\cdot|$ (also known as the supremum norm), which extends to a norm $|\cdot|$ on 
\[ A_{X, \infty} \coloneqq \bigcup_{n\ge 0} A_{X, n}.\] Note that, when $h = \mathrm{id}$, the norm $|x|$ of an element $x = \sum_{j_1, \ldots, j_d \in \Z[1/p]} a_{j_1, \ldots, j_d} \cdot T_1^{j_1}\cdots T_d^{j_d}$ with $a_{j_1, \ldots, j_d} \in K_\infty$ is simply given by $|x| = \max_{j_1, \ldots, j_d \in \Z[1/p]} |a_{j_1, \ldots, j_d}|$. 
\item[(2)] Let $\widehat A_{X, \infty}$ denote the completion of $A_{X, \infty}$ with respect to $|\cdot|$ and define 
$\widehat A_{X, \infty}^+  \coloneqq  \left\{x \in \widehat A_{X, \infty} \: \big| \: |x| \le 1 \right\}.$  Then $\widehat A_{X, \infty}$ is a perfectoid affinoid $\Q_p$-Banach algebra,  $(\widehat A_{X, \infty}, \widehat A_{X, \infty}^+)$  is a perfectoid Huber pair, and we have 
$\widehat{X}_\infty = \spa (\widehat A_{X, \infty}, \widehat A_{X, \infty}^+)$.  
\item[(3)] Denote the tilt of  $(\widehat A_{X, \infty}, \widehat A_{X, \infty}^+)$ by \[(\sq R_X, \sq R_X^+)  \coloneqq  (\widehat A_{X, \infty}^\flat, \widehat A_{X, \infty}^{\flat, +}).\] 
We again write $|\cdot|$ for the norm on $\sq R_X$ induced from the norm $|\cdot|$ on $\widehat A_{X, \infty}$.   
\item[(4)] Define a valuation $v\colon \sq R_X \ra \R \cup \{+\infty\}$ by setting $v(y)  \coloneqq  - \log_p |y|$,  extending the valuation $v$ on $\widehat K_\infty^\flat \subset \sq R_X$ defined at the beginning of  \S \ref{sec:classical_phi_Gamma}.    
\item[(5)] Note that $\varpi$ is a pseudouniformizer in $A_{X, \infty}^+$. 
We fix a tilt $\varpi^\flat \in  \widetilde{R}_X^+$ of $\varpi$. 
\end{enumerate} 
\end{definition} 

\begin{remark}\label{remark: extends the base field}
Throughout the paper, we often need to ``extend the base field'', that is, replacing $K$ by a finite extension $L$. Whenever we do so, we implicitly replace the toric chart $h\colon X\rightarrow \D_K$ by its base change $h\colon X_L\rightarrow \D_L$. In particular, the group $\Gamma_X=\Gamma_{\mathrm{geom}}\rtimes \Gamma_K$ is replaced by $\Gamma_{X_L} = \Gamma_{\mathrm{geom}}\rtimes \Gamma_L$.
\end{remark}

%%% 4.2
\vspace{0.1in}
\subsection{Perfect period rings}\label{ss:relative_period_rings:perfect}
\noindent
\vspace{0.1in}

\noindent We recall the constructions of perfect period rings in the relative setting following \cite{Kedlaya-LiuI}. These are natural extensions of the classical perfect period rings reviewed in \S \ref{ss:point perfect period ring}. 

\begin{definition}   \label{definition:valuation_v_r_relative_perfect}
Let $h\colon X \ra \D^d_K$ be a toric chart.
\begin{enumerate}
\item[(1)] Define $\widetilde{\mathbf{A}}_X \coloneqq W(\sq R_X)$ and $\widetilde{\mathbf{A}}_X^+ \coloneqq W(\sq R_X^+).$        
\item[(2)] For each $r > 0$, we define a norm $\lambda_r$ on $\widetilde{\mathbf{A}}^+_X$ by
$\lambda_r \left(\sum_{n=0}^{\infty} p^n [x_n]  \right) = \sup_n \left\{p^{-n} |x_n|^r\right\}$, 
which extends the norm $\lambda_r$ on $\widetilde{\mathbf{A}}_{K}^+$. 
\item[(3)] 
For each $r > 0$, define $\widetilde{\mathbf{A}}_X^r$ to be the completion of $\widetilde{\mathbf{A}}_X^+[[\varpi^\flat]^{-1}]$ under the norm $\lambda_r$. Equivalently, \[
\widetilde{\mathbf{A}}_X^r  \coloneqq  \left \{ \sum_{n =0}^{\infty} p^n [x_n] \in \widetilde{\mathbf{A}}_X \, \Bigg|   \,  p^{-n} |x_n|^r \ra 0 \textrm{ as } n \ra +\infty \right\}. \]
Let $\widetilde{\mathbf{A}}_X^{\dagger}$ be the union $\widetilde{\mathbf{A}}_X^{\dagger}  \coloneqq \varinjlim_{r>0} \widetilde{\mathbf{A}}_X^r.$  
\item[(4)] 
The period rings $\widetilde{\mathbf{B}}_X$, $\widetilde{\mathbf{B}}_X^r$, and $\widetilde{\mathbf{B}}_X^\dagger$ are obtained from the corresponding period rings of type $\mathbf{A}$ with $p$ inverted. 
\item[(5)]  
For every $0 < s \le r$, define $\widetilde{\mathbf{C}}_X^{[s, r]}$ to be the completion of the period ring $\widetilde{\mathbf{B}}_X^r$ with respect to \(\max\{\lambda_{s}, \lambda_{r}\}\) (or equivalently, completion with respect to the family of norms $\{\lambda_t\}_{t\in [s,r]}$).
\item[(6)] 
For every $r > 0$, define $\widetilde{\mathbf{C}}_X^{r}$ to be the Fr\'echet completion of  $\widetilde{\mathbf{B}}_X^r$ with respect to the family of norms $\{\lambda_{s}\}_{s \in (0, r]}$.
Finally, define $\widetilde{\mathbf{C}}_X^\dagger$ to be the union  
$\widetilde{\mathbf{C}}_X^{\dagger}  \coloneqq  \varinjlim_{r>0} \widetilde{\mathbf{C}}_X^r.$
\end{enumerate}
\end{definition}     

\begin{remark}
Strictly speaking, the perfect period rings $\widetilde{\mathbf{A}}_X$, $\widetilde{\mathbf{A}}_X^r$ etc. depend on the choice of the toric chart $h\colon X\ra \D^d_K$. In most cases, the relevant toric chart is evident from the context. For situations where this choice requires emphasis, we adopt the notation $\widetilde{\mathbf{A}}_{X,h}$, $\widetilde{\mathbf{A}}_{X,h}^r$ etc.
\end{remark}

\begin{remark}[Topology]\label{remark: topology--relative}
We endow $\widetilde{\mathbf{A}}^r_X$ with the topology induced by the $\lambda_r$-norm and endow $\widetilde{\mathbf{A}}^{\dagger}_X$ with the inductive limit topology (i.e., the finest topology making $\widetilde{\mathbf{A}}^r_X\hookrightarrow \widetilde{\mathbf{A}}^{\dagger}_X$ continuous for all $r>0$). Inverting $p$, we obtain natural topologies on $\widetilde{\mathbf{B}}^r_X$ and $\widetilde{\mathbf{B}}^{\dagger}_X$. The ring $\widetilde{\mathbf{C}}^{[s,r]}_X$ is endowed with the topology induced by $\max\{\lambda_r, \lambda_s\}$, which then induces a Fr\'echet topology on $\widetilde{\mathbf{C}}^r_X$. Finally, we equip $\widetilde{\mathbf{C}}^{\dagger}_X$ with the inductive limit topology (i.e., the LF topology).
\end{remark}

\begin{remark}\label{remark:Frobenius_relative_and_action_of_Gamma}
\begin{enumerate}
\item[(i)]
For more general adic spaces over $K$, the perfect period rings above globalize to period sheaves on $X_{\text{pro\'et}}$, which are denoted by $\widetilde{\mathbb{A}}_X, \widetilde{\mathbb{A}}^r_X, \widetilde{\mathbb{A}}^\dagger_X, \widetilde{\mathbb{B}}_X, \widetilde{\mathbb{B}}^r_X, \widetilde{\mathbb{B}}^\dagger_X,   \widetilde{\mathbb{C}}^r_X, \widetilde{\mathbb{C}}^\dagger_X$, etc. For more details, we refer the reader to \cite[Section 9.3]{Kedlaya-LiuI}. We will not use these sheaf versions of period rings in this article. 
\item[(ii)] 
There are natural surjective maps $\theta\colon \widetilde{\mathbf{A}}^+_X  \to \widehat A_{X, \infty}^+$ and $ \theta\colon \widetilde{\mathbf{A}}_X  \to \widehat A_{X,\infty}$. 
\item[(iii)] Remark~\ref{remark:Frobenius_classical} applies in the relative setting. In other words, the Frobenius maps on $\widetilde{\mathbf{A}}_X^+$ and $\widetilde{\mathbf{A}}_X$ extend to Frobenius structures on $\widetilde{\mathbf{A}}_X^r$, $\widetilde{\mathbf{A}}_X^\dagger$, $\widetilde{\mathbf{B}}_X$, $\widetilde{\mathbf{B}}_X^r$, $\widetilde{\mathbf{B}}_X^\dagger$, $\widetilde{\mathbf{C}}_X^{[s,r]}$, $\widetilde{\mathbf{C}}_X^r$, and $\widetilde{\mathbf{C}}_X^\dagger$, which we again denote as $\varphi$. 
\item[(iv)] 
The natural $\Gamma$-action on $\widetilde{R}_X$ induces an action on $\widetilde{\mathbf{A}}_X$. Notice that the $\lambda_r$-norm is preserved under the $\Gamma$-action. Consequently, for every $r > 0$, the subspace $\widetilde{\mathbf{A}}_X^r\subset \widetilde{\mathbf{A}}_X$ is stable under the $\Gamma$-action and the action is isometric with respect to $\lambda_r$. This in turn induces a continuous $\Gamma$-action on $\widetilde{\mathbf{A}}_X^{\dagger}$. Similarly, there are natural isometric (resp. continuous) $\Gamma$-actions on $\widetilde{\mathbf{B}}_X^r$ and $\widetilde{\mathbf{C}}_X^{[s, r]}$ (resp. on $\widetilde{\mathbf{B}}_X^{\dagger}$, $\widetilde{\mathbf{C}}_X^r$, and $\widetilde{\mathbf{C}}_X^{\dagger}$). 
\end{enumerate}
\end{remark}

%%% 4.3
\vspace{0.1in}
\subsection{Imperfect period rings, I: weakly decompleting data}\label{sec:weakly_decompleting_towers} 
\noindent
\vspace{0.1in}

\noindent The purpose of \S \ref{sec:weakly_decompleting_towers} is to recall the notion of \emph{weakly decompleting data}, together with the imperfect period rings associated with such data, originally introduced in \cite[\S 5]{relative}. Roughly speaking, this framework axiomatizes a set of conditions ensuring that the corresponding imperfect period rings enjoy nice properties, which in turn generalizes the classical constructions (cf. \S \ref{ss:point imperfect period ring}). One key property is that the ``completed perfections'' of these imperfect period rings should coincide with the corresponding perfect period rings. Moreover, the axioms also yield a reasonable theory of $(\varphi, \Gamma)$-modules over the associated imperfect period rings. For the setup, let $h\colon X \ra \D^d_K$ be a toric chart as in Definition~\ref{defn: toric chart}. We start with a pair $\Lambda=(R_{\Lambda}, \mathbf{A}_{\Lambda})$ where
\begin{itemize}
\item $R_{\Lambda}\subset \widetilde{R}_X$ is a $\kappa(\!(\overline{\pi}_K)\!)$-subalgebra such that $R_{\Lambda}$ is complete with respect to the norm $|\cdot|$ (note that $R_{\Lambda}$ is also complete with respect to the norm $|\cdot|^r$ for any $r>0$);
\item $\mathbf{A}_{\Lambda}\subset \widetilde{\mathbf{A}}_X$ is a $\varphi$-stable $\mathbf{A}_K$-subalgebra, which is complete with respect to the $p$-adic topology, such that the image of the composition
\[\mathbf{A}_{\Lambda}/p\rightarrow \widetilde{\mathbf{A}}_X/p\cong \widetilde{R}_X\]
lies in $R_{\Lambda}$.
\end{itemize}

\begin{definition}
Given a pair $\Lambda=(R_{\Lambda}, \mathbf{A}_{\Lambda})$ as above, we introduce the following notation.
\begin{enumerate}
\item[(1)] For any $r>0$, let $\mathbf{A}^r_{\Lambda} \coloneqq \mathbf{A}_{\Lambda}\cap \widetilde{\mathbf{A}}_X^r$ equipped with the inherited $\lambda_r$-norm from $\widetilde{\mathbf{A}}_X^r$. Let $\mathbf{B}^r_{\Lambda} \coloneqq \mathbf{A}^r_{\Lambda}[\frac{1}{p}]$.
\item[(2)] Let $\mathbf{A}^{\dagger}_{\Lambda}$ be the union $\mathbf{A}^{\dagger}_{\Lambda} \coloneqq \varinjlim_{r>0} \mathbf{A}^r_{\Lambda}$ and $\mathbf{B}^{\dagger}_{\Lambda} \coloneqq \mathbf{A}^{\dagger}_{\Lambda}[\frac{1}{p}]$. It follows that $\mathbf{A}^{\dagger}_{\Lambda}=\mathbf{A}_{\Lambda}\cap \widetilde{\mathbf{A}}^{\dagger}_X$.
\item[(3)] For $0<s\le r$, let $\mathbf{C}^{[s,r]}_{\Lambda}$ be the completion of $\mathbf{B}^r_{\Lambda}$ with respect to $\max\{\lambda_r, \lambda_s\}$ (or equivalently, completion with respect to the family of norms $\{\lambda_t\}_{t\in [s,r]}$). Let $\mathbf{C}^r_{\Lambda}$ be the Fr\'echet completion of $\mathbf{B}^r_{\Lambda}$ with respect to the family of norms $\{\lambda_s\}_{s\in (0,r]}$ and let $\mathbf{C}^{\dagger}_{\Lambda} \coloneqq  \varinjlim_{r>0}\mathbf{C}^r_{\Lambda}$.
\end{enumerate}
\end{definition}

It follows from the definition that $\mathbf{A}^r_{\Lambda}$ (resp. $\mathbf{A}^{\dagger}_{\Lambda}$, $\mathbf{B}^r_{\Lambda}$, $\mathbf{B}^{\dagger}_{\Lambda}$, $\mathbf{C}^{[s,r]}_{\Lambda}$, $\mathbf{C}^r_{\Lambda}$, $\mathbf{C}^{\dagger}_{\Lambda}$) can be naturally identified as a subring of $\widetilde{\mathbf{A}}^r_X$ (resp. $\widetilde{\mathbf{A}}^{\dagger}_X$, $\widetilde{\mathbf{B}}^r_X$, $\widetilde{\mathbf{B}}^{\dagger}_X$, $\widetilde{\mathbf{C}}^{[s,r]}_X$, $\widetilde{\mathbf{C}}^r_X$, $\widetilde{\mathbf{C}}^{\dagger}_X$). We equip them with the corresponding subspace topologies (cf. Remark~\ref{remark: topology--relative}). In light of this, we further introduce the following notation.

\begin{definition}
\begin{enumerate}
\item[(1)] For any $r>0$ and any integer $m\ge 0$, write $\varphi^m\colon \widetilde{\mathbf{A}}^r_X\isom  \widetilde{\mathbf{A}}^{r/p^m}_X$ for the $m$-th power of the Frobenius map from Remark~\ref{remark:Frobenius_relative_and_action_of_Gamma} (iii). Then define $
\mathbf{A}^r_{\Lambda, (m)} \coloneqq \varphi^{-m}\mathbf{A}^{r/p^m}_{\Lambda}\subset \widetilde{\mathbf{A}}^r_X$, 
$\breve{\mathbf{A}}^r_{\Lambda} \coloneqq \bigcup_{m=0}^{\infty}\mathbf{A}^r_{\Lambda, (m)}\subset \widetilde{\mathbf{A}}^r_X, $ and $ \breve{\mathbf{A}}^{\dagger}_{\Lambda} \coloneqq \bigcup_{r>0}\breve{\mathbf{A}}^r_{\Lambda}\subset \widetilde{\mathbf{A}}^{\dagger}_X.$
 
\item[(2)] Let $\mathbf{B}^r_{\Lambda, (m)} \coloneqq \mathbf{A}^r_{\Lambda, (m)}[\frac{1}{p}]\subset \widetilde{\mathbf{B}}^r_X$, $\: \breve{\mathbf{B}}^r_{\Lambda} \coloneqq \breve{\mathbf{A}}^r_{\Lambda}[\frac{1}{p}]\subset \widetilde{\mathbf{B}}^r_X$, and $\breve{\mathbf{B}}^{\dagger}_{\Lambda} \coloneqq \bigcup_{r>0}\breve{\mathbf{B}}^r_{\Lambda}\subset \widetilde{\mathbf{B}}^{\dagger}_X$.
\item[(3)] Similarly, let 
$
\mathbf{C}^{[s,r]}_{\Lambda, (m)} \coloneqq \varphi^{-m}\mathbf{C}^{[s/p^m, r/p^m]}_{\Lambda}\subset \widetilde{\mathbf{C}}^{[s,r]}_X $ and $ \breve{\mathbf{C}}^{[s,r]}_{\Lambda} \coloneqq \bigcup_{m=0}^{\infty}\mathbf{C}^{[s,r]}_{\Lambda, (m)}\subset \widetilde{\mathbf{C}}^{[s,r]}_X$  for $0<s\le r$.
Finally, let
$
\mathbf{C}^r_{\Lambda, (m)} \coloneqq \varprojlim_{0<s\le r} \mathbf{C}^{[s,r]}_{\Lambda, (m)}=\varphi^{-m}\mathbf{C}^{r/p^m}_{\Lambda}\subset \widetilde{\mathbf{C}}^r_X$,
$\breve{\mathbf{C}}^r_{\Lambda} \coloneqq \bigcup_{m=0}^{\infty} \mathbf{C}^r_{\Lambda, (m)}\subset \widetilde{\mathbf{C}}^r_X, $ and $  \breve{\mathbf{C}}^{\dagger}_{\Lambda} \coloneqq \bigcup_{r>0}\breve{\mathbf{C}}^r_{\Lambda}\subset \widetilde{\mathbf{C}}^{\dagger}_X.$ 
\end{enumerate}
\end{definition}

\begin{remark}[{Frobenius structures}]\label{remark: Frobenius_wdc_datum}
The Frobenius map $\varphi\colon \widetilde{\mathbf{A}}^r_X \isom  \widetilde{\mathbf{A}}^{r/p}_X$ induces a Frobenius $\varphi\colon \mathbf{A}^r_{\Lambda}\rightarrow \mathbf{A}^{r/p}_{\Lambda}$ for all $r$. Similarly, we have Frobenius maps 
$\varphi\colon \mathbf{B}^r_{\Lambda} \ra \mathbf{B}^{r/p}_{\Lambda},  \varphi\colon \mathbf{C}^r_{\Lambda} \ra \mathbf{C}^{r/p}_{\Lambda},  \varphi\colon \mathbf{C}^{[s, r]}_{\Lambda} \ra \mathbf{C}^{[s/p,\,r/p]}_{\Lambda},$ as well as on their breve variants $\breve{\mathbf{B}}^r_{\Lambda}, \breve{\mathbf{C}}^r_{\Lambda}, \breve{\mathbf{C}}^{[s,r]}_{\Lambda}$. Moreover, $\mathbf{B}^{\dagger}_{\Lambda}, \mathbf{C}^{\dagger}_{\Lambda}, \breve{\mathbf{A}}^{\dagger}_{\Lambda}, \breve{\mathbf{B}}^{\dagger}_{\Lambda}, \breve{\mathbf{C}}^{\dagger}_{\Lambda}$ are stable under Frobenius.
\end{remark}

Now we are ready to define weakly decompleting data.

\begin{definition}[Weakly decompleting datum]\label{defn: weakly decompleting datum}
A pair $\Lambda=(R_{\Lambda}, \mathbf{A}_{\Lambda})$ is called a \emph{weakly decompleting datum} for $(X,h)$ if there exists $r_*>0$ such that the following hold.
\begin{enumerate}
\item[(1)] $R_{\Lambda}$ admits a finite $p$-basis; namely, $\varphi^{-1}R_{\Lambda}$ (viewed as a subring of $\widetilde{R}_X$) is a finite free $R_{\Lambda}$-module.
\item[(2)] The colimit perfection $R^{\mathrm{perf}}_{\Lambda} \coloneqq \varinjlim \; (R_{\Lambda}\xrightarrow[]{\varphi}R_{\Lambda}\xrightarrow[]{\varphi}\cdots )$ is dense in $\widetilde{R}_X$.
\item[(3)] The subring $\mathbf{A}_{\Lambda}\subset \widetilde{\mathbf{A}}_X$ is stable under the $\Gamma$-action, and coincides with the $p$-adic completion of $\mathbf{A}^{\dagger}_{\Lambda}$.
\item[(4)] For any $r\in (0,r_*]\cap\Q$, $\mathbf{A}^r_{\Lambda}$ is complete with respect to the $\lambda_r$-norm.
\item[(5)] For any $r\in (0,r_*]\cap\Q$, the composition
\[\mathbf{A}^r_{\Lambda}/p\rightarrow \mathbf{A}_{\Lambda}/p\rightarrow R_{\Lambda}\]
is an isomorphism. (Consequently, $\mathbf{A}_{\Lambda}/p\rightarrow R_{\Lambda}$ is also an isomorphism.) 
\item[(6)] For any $r\in (0,r_*]\cap\Q$, the modulo $p$ reduction map $\mathbf{A}^r_{\Lambda}\rightarrow R_{\Lambda}$ is \emph{strict} in the sense that the quotient (semi)norm on $R_{\Lambda}$ induced from the $\lambda_r$-norm on $\mathbf{A}^r_{\Lambda}$ is equivalent to the norm $|\cdot|^r$ on $R_{\Lambda}$.
\item[(7)] $R_{\Lambda}$ is an affinoid $\kappa(\!(\overline{\pi}_K)\!)$-algebra (in the sense of Tate's rigid analytic geometry).
\item[(8)] For any $0<s\le r\le r_*$ with $s,r\in \mathbb{Q}$, the Banach $K'_0$-algebra $\mathbf{C}^{[s,r]}_{\Lambda}$ is a reduced affinoid $K'_0$-algebra (in the sense of Tate's rigid analytic geometry).
\end{enumerate}
\end{definition}

\begin{remark}[{$\Gamma$-actions}]
The condition (3) above (together with the isometric property of the $\lambda_r$-norms) implies that all of the imperfect period rings
$\mathbf{A}^r_{\Lambda}, \mathbf{A}^{\dagger}_{\Lambda}, \mathbf{B}_{\Lambda}, \mathbf{B}^r_{\Lambda}, \mathbf{B}^{\dagger}_{\Lambda}, \mathbf{C}^{[s,r]}_{\Lambda}, \mathbf{C}^r_{\Lambda}, \mathbf{C}^{\dagger}_{\Lambda}$
and their breve variants
$\breve{\mathbf{A}}^r_{\Lambda}, \breve{\mathbf{A}}^{\dagger}_{\Lambda}, \breve{\mathbf{B}}^r_{\Lambda}, \breve{\mathbf{B}}^{\dagger}_{\Lambda}, \breve{\mathbf{C}}^{[s,r]}_{\Lambda}, \breve{\mathbf{C}}^r_{\Lambda}, \breve{\mathbf{C}}^{\dagger}_{\Lambda}$ are stable under the $\Gamma$-actions.
\end{remark}

\begin{remark} Suppose $\Lambda=(R_{\Lambda}, \mathbf{A}_{\Lambda})$ is weakly decompleting.
\begin{enumerate}
\item[(i)] Since $R_{\Lambda}$ is a subring of the perfectoid ring $\widetilde{R}_{X}$, it is necessarily reduced. Notice that $|\cdot|$ is a complete power-multiplicative $\kappa(\!(\overline{\pi}_K)\!)$-algebra norm on $R_{\Lambda}$. By \cite[\S 6.2.4, Theorem 1]{BGR}, the norm $|\cdot|$ must coincide with the spectral norm on $R_{\Lambda}$ viewed as a reduced affinoid algebra.
\item[(ii)] Similarly, note that $\mathrm{max}\{\lambda_r, \lambda_s\}$ is a complete power-multiplicative $K'_0$-algebra norm on $\mathbf{C}^{[s,r]}_{\Lambda}$. Therefore, for the same reason as in (i), the norm $\mathrm{max}\{\lambda_r, \lambda_s\}$ coincides with the spectral norm on $\mathbf{C}^{[s,r]}_{\Lambda}$ viewed as a reduced affinoid $K'_0$-algebra.
\end{enumerate}
\end{remark}

Let us record the following consequences of the axioms. In particular, parts (2) and (3) of Lemma~\ref{lemma: perfection is dense} (also see Definition~\ref{defn: weakly decompleting datum} (2)) assert that the ``completed perfections'' of the imperfect period rings are precisely the perfect ones. This is essential for the decompletion theory of relative $(\varphi, \Gamma)$-modules.

\begin{lemma}[{\cite[Lemma 5.7, 5.11]{relative}}]\label{lemma: perfection is dense}
Suppose $\Lambda=(R_{\Lambda}, \mathbf{A}_{\Lambda})$ is weakly decompleting. 
\begin{enumerate}
\item[(1)] For any $r\in (0,r_*]$, the natural map $\mathbf{C}^r_{\Lambda} \rightarrow  \varprojlim_{0<s\le r} \mathbf{C}^{[s,r]}_{\Lambda}$ is necessarily an isomorphism. 
\item[(2)] For any $r\in (0,r_*]$, the subring $\breve{\mathbf{A}}^r_{\Lambda}$ is dense in $\widetilde{\mathbf{A}}^r_X$ with respect to the $\lambda_r$-norm.
\item[(3)] For any $r\in (0,r_*]$ and any $s\in (0,r]$, the subring $\breve{\mathbf{C}}^{[s,r]}_{\Lambda}$ is dense inside $\widetilde{\mathbf{C}}^{[s,r]}_X$ with respect to the $\max\{\lambda_r, \lambda_s\}$-norm.
\end{enumerate}
\end{lemma} 

The axioms in Definition~\ref{defn: weakly decompleting datum} also ensure that the intersections of the period rings behave as expected.

\begin{proposition}[Intersections of period rings]\label{prop: intersection A} Let $\Lambda=(R_{\Lambda}, \mathbf{A}_{\Lambda})$ be a weakly decompleting datum for $(X,h)$. 
Then there exists a sufficiently small $r_0 \in (0, r_*]$ such that the following statements hold.
\begin{enumerate}
\item[(1)] $\widetilde{\mathbf{A}}^s_X\cap \mathbf{C}^{[s,r]}_{\Lambda}=\mathbf{A}^r_{\Lambda}$, for any $0<s\le r<r_0$, where the intersection is taken inside $\widetilde{\mathbf{C}}^{\dagger}_X$.
\item[(2)] $\widetilde{\mathbf{A}}^{\dagger}_X\cap \mathbf{C}^{\dagger}_{\Lambda}=\mathbf{A}^{\dagger}_{\Lambda}$  where the intersection is taken inside $\widetilde{\mathbf{C}}^{\dagger}_X$.
\item[(3)] $\mathbf{C}^{[s,r]}_{\Lambda}\cap \mathbf{C}^{[s', r']}_{\Lambda}=\mathbf{C}^{[s,r']}_{\Lambda}$ for any $0<s\le s'\le r\le r'<r_0$, where the intersections are taken inside $\mathbf{C}^{[s',r]}_{\Lambda}$. 
\end{enumerate}
\end{proposition}

\begin{proof}
The analogous assertion is claimed in \cite[Lemma 5.2.10]{Kedlaya-LiuII}. However, the proof in \emph{loc.\ cit.} is incomplete. A complete argument is provided in   \cite[Proposition 5.25 \& 5.26]{relative}.
\end{proof}

Finally, we briefly discuss the behavior of weakly decompleting data under finite \'etale extensions. 

\begin{construction}\label{construction: weakly decompleting finite etale}
Let $h\colon X\rightarrow \D_K^d$ be a toric chart and let $\Lambda=(R_{\Lambda}, \mathbf{A}_{\Lambda})$ be a weakly decompleting datum for $(X,h)$. Consider a finite \'etale map $g\colon X'\rightarrow X$. Let $h'\colon X'\rightarrow \D_K^d$ be the composition $h\circ g$.  Note that the finite \'etale map $g\colon X'\rightarrow X$ induces a finite \'etale morphism $\widehat{A}_{X,\infty}\rightarrow \widehat{A}_{X',\infty}$, which corresponds to some $R'\in \left(R_{\Lambda}\right)_{\mathrm{f\acute{e}t}}$ and $\mathbf{A}'\in \left(\mathbf{A}_{\Lambda}\right)_{\mathrm{f\acute{e}t}}$ via the following natural equivalences of categories of finite \'etale algebras (cf. \cite[Proposition 5.28]{relative})
\[ 
\left(\mathbf{A}_{\Lambda}\right)_{\mathrm{f\acute{e}t}} \cong\left(R_{\Lambda}\right)_{\mathrm{f\acute{e}t}} \cong \left(\widehat{A}_{X, \infty}\right)_{\mathrm{f\acute{e}t}}.
\]

\end{construction}

\begin{proposition}[{\cite[Theorem 5.30]{relative}}]\label{prop: weakly decompleting finite etale}
The pair $(R', \mathbf{A}')$ is a weakly decompleting datum for $(X',h')$.
\end{proposition}

%%% 4.4
\vspace{0.1in}
\subsection{Relative $\varphi$-modules and $(\varphi, \Gamma)$-modules}\label{ss:relative_phi_gamma_mod}  
\noindent
\vspace{0.1in}

\noindent Throughout \S \ref{ss:relative_phi_gamma_mod}, let $h\colon X\rightarrow \D^d_K$ be a toric chart and let $\Lambda=(R_{\Lambda}, \mathbf{A}_{\Lambda})$ be a weakly decompleting datum for $(X,h)$. Let $\mathbf{A}^r_{\Lambda}, \mathbf{A}^{\dagger}_{\Lambda}, \mathbf{B}^r_{\Lambda}, \mathbf{B}^{\dagger}_{\Lambda}, \mathbf{C}^{[s,r]}_{\Lambda}, \mathbf{C}^r_{\Lambda}, \mathbf{C}^{\dagger}_{\Lambda}$  be the associated imperfect period rings. We introduce the notion of $\varphi$-modules and $(\varphi, \Gamma)$-modules over these rings.

\begin{definition}[$\varphi$-modules]  
\label{def:phi_modules}
Let $r_*$ be as in Definition~\ref{defn: weakly decompleting datum}.
\begin{enumerate} 
\item  
Let $0<r\le r_*$ and let $\mathbf{D}^r$ be one of the rings from $ \{\mathbf{A}^r_{\Lambda}, \mathbf{B}^r_{\Lambda},\mathbf{C}^r_{\Lambda},
 \widetilde{\mathbf{A}}^r_X, \widetilde{\mathbf{B}}^r_X, \widetilde{\mathbf{C}}^r_X\}.$ 
A \emph{$\varphi$-module} over $\mathbf{D}^r$ is a finite projective $\mathbf{D}^r$-module $M$ together with an isomorphism
\[
\varphi_{M}\colon M\otimes_{{\mathbf{D}}^r,\varphi} {\mathbf{D}}^{r/p} \isom M\otimes_{{\mathbf{D}}^r,\iota}{\mathbf{D}}^{r/p}
\]
of ${\mathbf{D}}^{r/p}$-modules, where $\varphi\colon \mathbf{D}^{r} \ra \mathbf{D}^{r/p}$ is the Frobenius map (cf. Remark~\ref{remark:Frobenius_relative_and_action_of_Gamma} (iii) and Remark~\ref{remark: Frobenius_wdc_datum}), and $\iota$ is the inclusion map. 
\item Similarly, let $\mathbf{D}$ be one of the rings from  
$\{ \mathbf{A}_{\Lambda}, \mathbf{A}^{\dagger}_{\Lambda}, \mathbf{B}_{\Lambda}, \mathbf{B}^{\dagger}_{\Lambda}, \mathbf{C}^{\dagger}_{\Lambda}, 
\widetilde{\mathbf{A}}_X,
\widetilde{\mathbf{A}}^{\dagger}_X, \widetilde{\mathbf{B}}_X, \widetilde{\mathbf{B}}^{\dagger}_X, \widetilde{\mathbf{C}}^{\dagger}_X\}.$ A \emph{$\varphi$-module} over $\mathbf{D} $ is a finite projective $\mathbf{D}$-module $M$ together with an isomorphism
$ \varphi_{M}\colon M\otimes_{{\mathbf{D}},\varphi} \mathbf{D} \isom M $ 
of $\mathbf{D}$-modules.
\end{enumerate}
\end{definition}

\begin{definition} %[$(\varphi, \Gamma)$-modules]
\label{def:relative_phi_gamma_modules}   
\begin{enumerate}
    \item 
Let $\mathbf{D}$ be one of the rings from Definition~\ref{def:phi_modules}.  A \emph{$\Gamma$-module} over $\mathbf{D}$ is a finite projective $\mathbf{D}$-module equipped with a continuous semilinear action of $\Gamma$.\footnote{For a finite projective $\mathbf{D}$-module $M$, we may identify $M$ as a direct summand of some finite free $\mathbf{D}$-module and equip $M$ with the subspace topology. This topology is independent of the choice of the presentation. On the other hand, $\Gamma$ is equipped with the usual profinite topology.}  A \emph{$(\varphi, \Gamma)$-module} over $\mathbf{D}$ is a $\varphi$-module $M$ over $\mathbf{D}$, equipped with a continuous semilinear action of $\Gamma$ that commutes with the action of $\varphi$. 

\item A $\varphi$-module (resp.  $(\varphi, \Gamma)$-module) over one of the following rings $\{\mathbf{B}_{\Lambda}, \mathbf{B}^{\dagger}_{\Lambda}, \mathbf{C}^{\dagger}_{\Lambda}, \widetilde{\mathbf{B}}_X, \widetilde{\mathbf{B}}_X^{\dagger}, \widetilde{\mathbf{C}}_X^{\dagger} \}$ 
is  \emph{\'etale} if it is the base change of a $\varphi$-module (resp.  $(\varphi, \Gamma)$-module) over  $ \{ \mathbf{A}_{\Lambda}, \mathbf{A}^{\dagger}_{\Lambda}, \mathbf{A}^{\dagger}_{\Lambda}, \widetilde{\mathbf{A}}_X, \widetilde{\mathbf{A}}_X^{\dagger}, \widetilde{\mathbf{A}}_X^{\dagger} \}$ respectively.  
\end{enumerate}
\end{definition}

\begin{remark}
Our notion of \'etale $(\varphi, \Gamma)$-modules requires the corresponding lattice over period rings of type $\mathbf{A}$ to be stable under both $\varphi$ and $\Gamma$. This is similar to Definition \cite[Definition I.4.2]{Cherbonnier_Colmez}. In contrast, \cite{Kedlaya-LiuII} only requires the lattice to be stable under $\varphi$. In other words, the definition we use is more restrictive, but seems more suitable for the purpose of decompletion and descent. 
\end{remark}

\begin{notation} 
Let $\mathbf{D}$ be one of the rings in Definition~\ref{def:phi_modules}.  Let  $\mathrm{Mod}_{/\mathbf{D}}^{\varphi} $ (resp. $\mathrm{Mod}_{/\mathbf{D}}^{\Gamma},$  resp. $\mathrm{Mod}_{/\mathbf{D}}^{(\varphi, \Gamma)})$ denote the category of $\varphi$-modules (resp. the category of $\Gamma$-modules, resp. the category of $(\varphi, \Gamma)$-modules) over $\mathbf{D}$.  
Similarly,  we write $\mathrm{Mod}_{/\mathbf{D}}^{\varphi, \ett}$ (resp. $\mathrm{Mod}_{/\mathbf{D}}^{(\varphi, \Gamma), \ett}$) for the category of \'etale $\varphi$-modules (resp. the category of \'etale $(\varphi, \Gamma)$-modules) over $\mathbf{D}$ whenever appropriate.    
\end{notation}

Given a weakly decompleting datum $\Lambda=(R_{\Lambda}, \mathbf{A}_{\Lambda})$ for $(X,h)$, we frequently need to understand the behavior of the base change functors for $(\varphi,\Gamma)$-modules across the following diagram of period rings (for instance, whether these functors are equivalences or fully faithful).

\begin{equation} \label{diagram:big_fucking_diagram_of_rings}  
\begin{tikzcd}[column sep=2em, row sep=0.4em, scale cd=1]
    \mathbf{A}^{\dagger}_{\Lambda} \arrow[rr, "(1)"]   \arrow[ddd] & &   \widetilde{\mathbf{A}}_X^{\dagger} \arrow[rrrdd, "(4)"] \arrow[ddd]   
    \\   \\ 
    &  &  &
    \mathbf{A}_{\Lambda} \arrow[from=llluu, crossing over, swap,  "(3)"] \arrow[rr, swap, "(2)"] & &\widetilde{\mathbf{A}}_X   \arrow[ddd] \\ 
    \mathbf{B}^{\dagger}_{\Lambda} \arrow[rr, "(5)"] \arrow[ddd, "(9)"]  & &   \widetilde{\mathbf{B}}_X^{\dagger} \arrow[ddd,  near end, "(10)"] \arrow[rrrdd, "(8)"] \\  \\ 
       &  & &  \mathbf{B}_{\Lambda} \arrow[from=uuu, crossing over] \arrow[from=llluu, swap,  "(7)", crossing over] \arrow[rr, swap, "(6)"] & & \widetilde{\mathbf{B}}_X   \\
      \mathbf{C}^{\dagger}_{\Lambda} \arrow[rr, "(11)"]& & \widetilde{\mathbf{C}}_X^{\dagger}  
\end{tikzcd}  
\end{equation}
The horizontal arrows correspond to so-called ``decompletion questions'' while the diagonal arrows correspond to ``overconvergence questions''. In the rest of \S \ref{ss:relative_phi_gamma_mod}, we discuss some known results under the assumption of weakly decompleting data. We shall constantly refer to the labels of the arrows in the diagram above.

\subsubsection{Base change along the arrows (4), (8), and (10)}
\noindent
\vspace{0.1in}

\noindent We first recall the following descent theorem for perfect period rings. 

\begin{proposition}[{Arrows (4) and (10)}]\label{prop:descending_phi_mod_along_4_and_10} 
Taking base change along $\widetilde{\mathbf{A}}_X^{\dagger} \ra \widetilde{\mathbf{A}}_X$ (resp. $\widetilde{\mathbf{B}}_X^{\dagger}  \ra \widetilde{\mathbf{C}}_X^{\dagger} $) induces an equivalence of categories  
    \[ \mathrm{Mod}_{/\widetilde{\mathbf{A}}_X^{\dagger}}^{\varphi} \isom \mathrm{Mod}_{/\widetilde{\mathbf{A}}_X }^{\varphi} \qquad \left( \textrm{resp. }  \mathrm{Mod}_{/\widetilde{\mathbf{B}}_X^{\dagger}}^{\varphi, \ett}  \isom  \mathrm{Mod}_{/\widetilde{\mathbf{C}}_X^{\dagger}}^{\varphi, \ett} \right)
    \]  between $\varphi$-modules (resp. \'etale $\varphi$-modules), and thus induces an equivalence of categories
     \[ \mathrm{Mod}_{/\widetilde{\mathbf{A}}_X^{\dagger}}^{(\varphi,\Gamma)} \isom \mathrm{Mod}_{/\widetilde{\mathbf{A}}_X }^{(\varphi,\Gamma)} \qquad \left( \textrm{resp. }  \mathrm{Mod}_{/\widetilde{\mathbf{B}}_X^{\dagger}}^{(\varphi,\Gamma), \ett}  \isom  \mathrm{Mod}_{/\widetilde{\mathbf{C}}_X^{\dagger}}^{(\varphi,\Gamma), \ett} \right)
    \]
    between $(\varphi, \Gamma)$-modules (resp. \'etale $(\varphi, \Gamma)$-modules).  
\end{proposition}

\begin{proof}
The first statement follows from \cite[Theorem 8.5.3 (d)(e)]{Kedlaya-LiuI} and \cite[Theorem 8.5.6 (e)(f)]{Kedlaya-LiuI}, while the second statement is an immediate consequence.
\end{proof}

\begin{remark}[{Arrow (8)}]
As an immediate consequence of Proposition~\ref{prop:descending_phi_mod_along_4_and_10}, we obtain equivalences of categories
\[\mathrm{Mod}_{/\widetilde{\mathbf{B}}_X^{\dagger}}^{\varphi, \ett}  \isom  \mathrm{Mod}_{/\widetilde{\mathbf{B}}_X}^{\varphi, \ett}, \qquad \mathrm{Mod}_{/\widetilde{\mathbf{B}}_X^{\dagger}}^{(\varphi,\Gamma), \ett}  \isom  \mathrm{Mod}_{/\widetilde{\mathbf{B}}_X}^{(\varphi,\Gamma), \ett}.\]
\end{remark}

\subsubsection{Base change along the arrows (1), (2), (3), (5), (6), (7), (9), and (11)}

\begin{proposition}[{Arrows (2), (3), and (9)}]\label{prop:descending_phi_mod_along_2_3_9}
Let $\Lambda=(R_{\Lambda}, \mathbf{A}_{\Lambda})$ be a weakly decompleting datum for $(X,h)$.
\begin{enumerate} 
\item[(1)] Taking base change along the inclusion $\mathbf{A}_{\Lambda} \ra \widetilde{\mathbf{A}}_X$ induces an equivalence of categories of $\varphi$-modules 
$\mathrm{Mod}_{/\mathbf{A}_{\Lambda}}^{\varphi} \isom \mathrm{Mod}_{/\widetilde{\mathbf{A}}_X }^{\varphi},$ 
and thus an equivalence of categories of $(\varphi, \Gamma)$-modules
\[ 
\mathrm{Mod}_{/\mathbf{A}_{\Lambda}}^{(\varphi, \Gamma)} \isom \mathrm{Mod}_{/\widetilde{\mathbf{A}}_X }^{(\varphi, \Gamma)}.
\]
\item[(2)] The base change functor of $\varphi$-modules (resp. of \'etale $\varphi$-modules) 
\[ \mathrm{Mod}_{/\mathbf{A}^{\dagger}_{\Lambda}}^{\varphi} \lra \mathrm{Mod}_{/\mathbf{A}_{\Lambda}}^{\varphi} \qquad \Big( \textrm{resp. } \mathrm{Mod}_{/\mathbf{B}^{\dagger}_{\Lambda}}^{\varphi, \ett} \lra \mathrm{Mod}_{/\mathbf{C}^{\dagger}_{\Lambda}}^{\varphi, \ett}\Big)
\]  
along the natural map $\mathbf{A}^{\dagger}_{\Lambda} \ra \mathbf{A}_{\Lambda}$ (resp. $\mathbf{B}^{\dagger}_{\Lambda} \ra \mathbf{C}^\dagger_{\Lambda}$) is fully faithful. Consequently, we obtain fully faithful functors
\[\mathrm{Mod}_{/\mathbf{A}^{\dagger}_{\Lambda}}^{(\varphi, \Gamma)} \lra \mathrm{Mod}_{/\mathbf{A}_{\Lambda}}^{(\varphi, \Gamma)}, \qquad  \mathrm{Mod}_{/\mathbf{B}^{\dagger}_{\Lambda}}^{(\varphi, \Gamma), \ett} \lra \mathrm{Mod}_{/\mathbf{C}^{\dagger}_{\Lambda}}^{(\varphi, \Gamma), \ett}.\]
\end{enumerate}
\end{proposition}

\begin{proof}
This is \cite[Proposition 7.14 \& 7.15]{relative}.    
\end{proof}

\begin{remark}[{Arrows (1), (5), (6), (7), and (11)}]
The following are immediate consequences of Proposition~\ref{prop:descending_phi_mod_along_2_3_9}.
\begin{enumerate}
\item[(i)] The base change functor of $\varphi$-modules (resp. $(\varphi, \Gamma)$-modules) along the arrow (1) is fully faithful.
\item[(ii)] The base change functor of \'etale $\varphi$-modules (resp. \'etale $(\varphi, \Gamma)$-modules) along the arrow (6) is an equivalence of categories.
\item[(iii)] The base change functor of \'etale $\varphi$-modules (resp. \'etale $(\varphi, \Gamma)$-modules) along the arrows (5) and (7) are fully faithful. Together with Proposition~\ref{prop:descending_phi_mod_along_4_and_10}, one sees that the base change functor of \'etale $\varphi$-modules (resp. \'etale $(\varphi, \Gamma)$-modules) along the arrow (11) is also fully faithful.
\end{enumerate}
\end{remark}

%%% 4.5
\vspace{0.1in}
\subsection{Imperfect period rings, II: decompleting data}  \label{sec:decompleting_towers}
\noindent
\vspace{0.1in}

\noindent 
For our purpose, it is crucial to be able to descend \'etale $(\varphi, \Gamma)$-modules along the arrow (11); the issue at the moment is essential surjectivity. For this, we need to impose a stronger condition called \emph{decompleting data} introduced in \cite{relative}, which we now recall. The idea is to put an additional condition on the $\Gamma$-action so that $\Gamma$-modules over perfect period rings descend to imperfect period rings (hence $(\varphi, \Gamma)$-modules also descend). The main theorem is Theorem~\ref{thm:decompleting data--descent_for_type_C}.

\begin{definition}
Let $M$ be a topological abelian group whose topology is induced by a norm $|\cdot|$ (in the sense of \cite[\S 1.1.1 Definition 1]{BGR}). Suppose $M$ is equipped with a continuous action of a compact topological group $G$. Let  $\mathcal{C}^{\bullet}(G, M)$ denote the complex of inhomogeneous continuous $G$-cochains with values in $M$. Each $\mathcal{C}^i=\mathcal{C}^i(G, M)$ is equipped with the supremum norm induced from the norm $|\cdot|$ on $M$, still denoted by $|\cdot|$; namely, for $f\in \mathcal{C}^i(G, M)$, define
\[|f| \coloneqq \sup_{(g_1, \ldots, g_i)\in G^i}|f(g_1, \ldots, g_i)|.\] 
Since $G$ is compact, the supremum norm is well-defined. Let $d^{(i)}\colon \mathcal{C}^i\rightarrow \mathcal{C}^{i+1}$ denote the differentials. We say that $\mathcal{C}^{\bullet}(G, M)$ is \emph{strictly exact at $d^{(i)}$} if $\mathrm{ker}(d^{(i)})=\mathrm{im}(d^{(i-1)})$ and $d^{(i)}$ is \emph{strict} in the sense that the quotient (semi)norm on $\mathcal{C}^{i}/\mathrm{ker}(d^{(i)})\cong \mathrm{im}(d^{(i)})$ is equivalent to the subspace norm on $\mathrm{im}(d^{(i)})\subset \mathcal{C}^{i+1}$. In explicit terms, this means there exist $\mathfrak{c}, \mathfrak{c}'>0$ such that 
\begin{enumerate}
\item[(1)] for any $g\in \mathcal{C}^i$, we have
\[\left|d^{(i)}g\right|\le \mathfrak{c}\cdot\left|g+d^{(i-1)}f\right|\]
for all $f\in \mathcal{C}^{i-1}$;
\item[(2)] for any $g\in \mathcal{C}^i$, we can find $f'\in \mathcal{C}^{i-1}$ satisfying
\[\left|g+d^{(i-1)}f'\right|\le \mathfrak{c}'\cdot\left|d^{(i)}g\right|.\]
\end{enumerate}
\end{definition}

\begin{definition}\label{definition:decompleting_datum}

A weakly decompleting datum  $\Lambda=(R_{\Lambda}, \mathbf{A}_{\Lambda})$ for $(X,h)$ is called \emph{decompleting} if the following conditions hold.
\begin{enumerate}
\item For any open subgroup $\Gamma_0 \subset \Gamma$, the cochain complex
\[
\mathcal{C}^{(\bullet)}(\Gamma_0, \varphi^{-1}(R_{\Lambda}) / R_{\Lambda})
\]
is strict exact at $d^{(0)}$ and $d^{(1)}$.

\item For any finite \'etale map   $g\colon X'\rightarrow X$, Condition (1) holds with $R_{\Lambda}$ replaced by $R'$, where $R'\in \left(R_{\Lambda}\right)_{\mathrm{f\acute{e}t}}$ is given by Construction~\ref{construction: weakly decompleting finite etale}.
\end{enumerate}
\end{definition}

\begin{remark}
\begin{enumerate}
    \item[(i)]
If $\Lambda=(R_{\Lambda}, \mathbf{A}_{\Lambda})$ is decompleting, then the cochain complex $\mathcal{C}^{\bullet}(\Gamma, \varphi^{-1}R_{\Lambda}/R_{\Lambda})$ is strict exact at $d^{(0)}$ and $d^{(1)}$ (by considering the case $g_1=\cdots=g_n=1$ and $\Gamma_0=\Gamma$).
 \item[(ii)]
In fact, condition (2) in Definition~\ref{definition:decompleting_datum} is not needed for Theorem~\ref{thm:decompleting data--descent_for_type_C} or Theorem~\ref{thm:decompleting data--descent_for_type_C--finite free} below. We include this condition solely to ensure that the notion of decompleting data is stable under finite \'etale extensions (cf. Construction~\ref{construction: weakly decompleting finite etale}, Proposition~\ref{prop: weakly decompleting finite etale}).
\end{enumerate}
\end{remark}

The following result is one of the main results of \cite{relative} and plays an important role later in the article.

\begin{theorem}[{\cite[Theorem 7.18, Corollary 7.20]{relative}}]\label{thm:decompleting data--descent_for_type_C}
Let $\Lambda=(R_{\Lambda}, \mathbf{A}_{\Lambda})$ be a decompleting datum for $(X,h)$. Then the base change functors
\[
\mathrm{Mod}_{/\mathbf{A}_{\Lambda}^{\dagger}}^{(\varphi, \Gamma)} \lra \mathrm{Mod}_{/ \widetilde{\mathbf{A}}_X^{\dagger}}^{(\varphi, \Gamma)}
\]
and
\[
\mathrm{Mod}_{/\mathbf{B}_{\Lambda}^{\dagger}}^{(\varphi, \Gamma), \ett} \lra \mathrm{Mod}_{/ \widetilde{\mathbf{B}}_X^{\dagger}}^{(\varphi, \Gamma), \ett} \quad (\textrm{resp. } ~~\mathrm{Mod}_{/\mathbf{C}_{\Lambda}^{\dagger}}^{(\varphi, \Gamma), \ett} \lra \mathrm{Mod}_{/ \widetilde{\mathbf{C}}_X^{\dagger}}^{(\varphi, \Gamma), \ett})
\]
are equivalences of categories. Moreover, the same statement holds if we replace $\Gamma$ by an open subgroup.
\end{theorem}

\begin{remark}\label{remark:neither fully faithful nor essentially surjective}
By the discussion above, we have equivalences of categories of \'etale $(\varphi, \Gamma)$-modules across all the arrows (5)--(11) in the diagram (\ref{diagram:big_fucking_diagram_of_rings}). Let us note that, however, it is unclear whether the functor $\mathrm{Mod}^{(\varphi, \Gamma)}_{/\mathbf{C}^{\dagger}_{\Lambda}}\lra \mathrm{Mod}^{(\varphi, \Gamma)}_{/\widetilde{\mathbf{C}}_X}$ (for general, not necessarily \'etale, $(\varphi, \Gamma)$-modules) is fully faithful or essentially surjective. Nevertheless, we have the following partial result, which turns out to be sufficient for the purpose of this paper. 
\end{remark}

\begin{theorem}[{\cite[Theorem 7.17]{relative}}]\label{thm:decompleting data--descent_for_type_C--finite free}
Let $\Lambda=(R_{\Lambda}, \mathbf{A}_{\Lambda})$ be a decompleting datum for $(X,h)$. Then a descent of a finite free $(\varphi, \Gamma)$-module over $\widetilde{\mathbf{C}}^{\dagger}_X$  to a $(\varphi, \Gamma)$-module over $\mathbf{C}^{\dagger}_{\Lambda}$, if exists, is unique. Moreover, the same statement holds if we replace $\Gamma$ by an open subgroup.
\end{theorem}

\begin{remark}[Relation with local systems] \label{remark:local_system_relation_with_period_ring_type_A}  Let $\mathrm{Loc}_{\Z_p}(X)$  denote the category of \'etale $\Z_p$-local systems on $X$ and let $\mathrm{Loc}_{\Z_p}(X)_{\Q}$ denote its isogeny category. By \cite[Theorem 8.5.5]{Kedlaya-LiuI}, there is an equivalence of categories
\[\widetilde{D}^{\mathrm{int}}:\: \mathrm{Loc}_{\Z_p} (X)\isom\mathrm{Mod}_{/\widetilde{\mathbf{A}}_X }^{(\varphi,\Gamma)}\]
which reduces to the Riemann--Hilbert correspondence of Katz after reducing modulo $p$ (cf. \cite[Proposition 4.1.1]{Katz}, also see \cite[Proposition 3.4]{BS_crystal}). Consequently, under the assumption of decompleting data, $\mathrm{Loc}_{\Z_p}(X)$ (resp. $\mathrm{Loc}_{\Z_p}(X)_{\Q}$) is equivalent to the categories of $(\varphi, \Gamma)$-modules across all the arrows (1)-(4) (resp. categories of \'etale $(\varphi, \Gamma)$-modules across all the arrows (5)-(11)). In particular, we obtain equivalences of categories
\begin{align*}
\widetilde{D}^{\mathrm{int}, \dagger}\colon  \: \mathrm{Loc}_{\Z_p} (X) \isom & \,
     \mathrm{Mod}^{(\varphi, \Gamma)}_{/\widetilde{\mathbf{A}}_X^{\dagger}}.\\
D^{\mathrm{int},\dagger}\colon  \: \mathrm{Loc}_{\Z_p} (X) \isom & \, \mathrm{Mod}^{(\varphi, \Gamma)}_{/\mathbf{A}^{\dagger}_{\Lambda}}\\
\widetilde{D}^{\dagger}_{\mathrm{rig}}\colon  \: \mathrm{Loc}_{\Z_p}(X)_{\Q} \isom & \,
     \mathrm{Mod}^{(\varphi, \Gamma), \ett}_{/\widetilde{\mathbf{C}}_X^{\dagger}}.\\
D_{\mathrm{rig}}^{\dagger}\colon  \: \mathrm{Loc}_{\Z_p} (X)_{\Q} \isom & \, \mathrm{Mod}^{(\varphi, \Gamma), \ett}_{/\mathbf{C}^{\dagger}_{\Lambda}}
\end{align*}
where the last equivalence can be viewed as a relative version of the base change functor $D^{\dagger}_{\mathrm{rig}, K}$ in \S \ref{sss:overconvergent classical}. 
\end{remark}

%%% 4.6
\vspace{0.1in}
\subsection{Imperfect period rings, III: strongly decompleting data} \label{sec:strongly decompleting data}
\noindent
\vspace{0.1in}

\noindent Going one step further, we now recall an even more restrictive condition called \emph{strongly decompleting data} introduced in \cite{relative}. This condition will allow us to extend Berger's $\mathbf{N}_{\mathrm{dR}}$ functor to the relative setting.

Let $h\colon X\rightarrow \D^d_K$ be a toric chart. We further assume that the image of $h$ does not intersect the locus $\{T_1\cdots T_d=0\}$. In particular, the toric tower $X_{\infty}=\varprojlim_n X_n\rightarrow X$ is an affinoid perfectoid object in the pro-\'etale site $X_{\proet}$. Recall the period sheaves $\mathbb{B}^+_{\mathrm{dR}, X}$ and $\mathcal{O}\mathbb{B}^+_{\mathrm{dR}, X}$ on $X_{\proet}$ defined in \cite{Scholze_p_adic_Hodge}. The sheaf $\mathcal{O}\mathbb{B}^+_{\mathrm{dR}, X}$ is equipped with a $\mathbb{B}^+_{\mathrm{dR}}$-linear connection
\[\nabla\colon  \mathcal{O}\mathbb{B}^+_{\mathrm{dR}, X}\rightarrow \mathcal{O}\mathbb{B}^+_{\mathrm{dR}, X}\otimes_{\mathcal{O}_X}\Omega^1_X\]
and we have $\mathbb{B}^+_{\mathrm{dR}, X}= (\mathcal{O}\mathbb{B}^+_{\mathrm{dR}, X} )^{\nabla=0}$ by \cite[Corollary 6.13]{Scholze_p_adic_Hodge}. For any $r\in (0,1]$ and any integer $n\ge -\log_pr$, consider the map
$\iota_n\colon  \widetilde{\mathbf{A}}^r_X\rightarrow \mathbb{B}^+_{\mathrm{dR}, X}(X_{\infty})$ 
sending $\sum_{m=0}^{\infty}p^m[x_m]\mapsto \sum_{m=0}^{\infty}p^m[x_m^{1/p^n}]$. This  is well-defined by the same computation as in \cite[Lemme 7.2]{Andreatta-Brinon}. Inverting $p$, we obtain $\iota_n\colon  \widetilde{\mathbf{B}}^r_X\rightarrow \mathbb{B}^+_{\mathrm{dR}, X}(X_{\infty})$. By a continuity argument (\cite[Corollary 8.2]{relative}), $\iota_n$ further extends to 
\[\iota_n\colon  \widetilde{\mathbf{C}}^{[s,r]}_X\rightarrow \mathbb{B}^+_{\mathrm{dR}, X}(X_{\infty})\]
whenever $p^{-n}\in[s,r]$.
Recall the natural $\theta$-map
$\theta_X\colon  \mathcal{O}\mathbb{B}^+_{\mathrm{dR}, X}(X_{\infty})\twoheadrightarrow \widehat{A}_{X, \infty}$ with $\mathrm{ker}(\theta_X)=(t,u_1, \ldots, u_d)$ where $u_i=T_i-[T_i^{\flat}]$. Moreover, by \cite[Proposition 6.10]{Scholze_p_adic_Hodge} (also see \cite[Proposition 2.3.15]{DLLZ2}), there is a natural identification
\[\mathcal{O}\mathbb{B}^+_{\mathrm{dR}, X}(X_{\infty})\cong \mathbb{B}^+_{\mathrm{dR}, X}(X_{\infty})[\![u_1, \ldots, u_d]\!].\]
Consider the subring 
\[A_{X,n}[\![t, u_1, \ldots, u_d]\!]\subset \mathcal{O}\mathbb{B}^+_{\mathrm{dR}, X}(X_{\infty}).\]
Note that $A_{X,n}[\![t, u_1, \ldots, u_d]\!]$ is stable under the $\Gamma$-action on $\mathcal{O}\mathbb{B}^+_{\mathrm{dR}, X}(X_{\infty})$. Taking the $\nabla=0$ part, we consider the subring
\[A_{X,n}[\![t, u_1, \ldots, u_d]\!]^{\nabla=0}\subset \mathcal{O}\mathbb{B}^+_{\mathrm{dR}, X}(X_{\infty})^{\nabla=0}=\mathbb{B}^+_{\mathrm{dR}, X}(X_{\infty}).\]
It turns out $A_{X,n}[\![t, u_1, \ldots, u_d]\!]^{\nabla=0}$ is complete with respect to the $t$-adic topology (cf. \cite[\S 8.1.2]{relative}).

Now we are ready to define strongly decompleting data.
\begin{definition}\label{definition:strongly_decompleting_datum}
Let $\Lambda=(R_{\Lambda}, \mathbf{A}_{\Lambda})$ be a decompleting datum for $(X,h)$. We say that $\Lambda$ is \emph{strongly decompleting} if the following conditions hold.
\begin{enumerate}
\item[(a)] For $r$ sufficiently small, the image of the composition
\[\mathbf{A}^r_{\Lambda}\subset \widetilde{\mathbf{A}}^r_X\xrightarrow[]{\iota_n} \mathbb{B}_{\mathrm{dR}, X}^+(X_{\infty})\]
lies in $A_{X,n}[\![t, u_1, \ldots, u_d]\!]^{\nabla=0}$, for all $n\ge -\log_pr$.
\item[(b)] The composition \[ \mathbf{A}^r_{\Lambda}\xrightarrow[]{\iota_n} A_{X,n}[\![t, u_1, \ldots, u_d]\!]^{\nabla=0} \twoheadrightarrow A_{X, n}
\] is surjective for all $n\ge -\log_pr$.

\item[(c)] Let $r_n=p^{-n}$. There are an integer $n_0 \ge 0$ and a constant $C < p^{-\frac{1}{p-1}}$ such that for any $n \ge n_0$, we have
$\lambda_{r_n}\left((\gamma-1)(x)\right)\le C\cdot \lambda_{r_n}(x)$ for all $x\in \mathbf{C}_{\Lambda}^{[r_n, r_n]}$ and $\gamma\in \Gamma_n$.
\end{enumerate}
\end{definition}

\begin{remark}
\begin{enumerate}
\item[(i)] As a consequence of condition (a), the map $\iota_n\colon \mathbf{A}^r_{\Lambda}\rightarrow A_{X,n}[\![t, u_1, \ldots, u_d]\!]^{\nabla=0}$ extends to $\iota_n\colon \mathbf{C}^{[s,r]}_{\Lambda}\rightarrow A_{X,n}[\![t, u_1, \ldots, u_d]\!]^{\nabla=0}$ by continuity (cf. \cite[Corollary 8.5]{relative}), whenever $p^{-n}\in [s,r]$.
\item[(ii)] The composition in condition (b) is simply the restriction of $\theta\circ \varphi^{-n}\colon  \widetilde{\mathbf{B}}^r_X\twoheadrightarrow \widehat{A}_{X, \infty}$ on $\mathbf{B}^r_{\Lambda}$. 
\end{enumerate}
\end{remark}

\begin{lemma}[{\cite[Lemma 8.11]{relative}}]\label{lemma: dense in B_n}
Suppose $\Lambda=(R_{\Lambda}, \mathbf{A}_{\Lambda})$ is strongly decompleting. For $r$ sufficiently small, let $[s,r]$ be an interval containing $p^{-n}$. Then $\iota_n\colon \mathbf{C}^{[s,r]}_{\Lambda}\rightarrow A_{X,n}[\![t, u_1, \ldots, u_d]\!]^{\nabla=0}$ is injective with dense image with respect to the $t$-adic topology, for all $n\ge -\log_pr$.
\end{lemma}

The notion of strongly decompleting data is also stable under finite \'etale extensions.

\begin{proposition}[{\cite[Theorem 8.13]{relative}}]\label{prop: strongly decompleting datum finite etale}
Let $g\colon  X'\rightarrow X$ be a finite \'etale map and let $h'\colon X'\rightarrow \D_K^d$ be the composition $h\circ g$. Let $\Lambda=(R_{\Lambda}, \mathbf{A}_{\Lambda})$ be a strongly decompleting datum for $(X,h)$ and let $(R', \mathbf{A}')$ be the pair induced from $g$ as in Construction~\ref{construction: weakly decompleting finite etale}. Then $(R', \mathbf{A}')$ is a strongly decompleting datum for $(X',h')$.
\end{proposition}

\begin{remark}[Passing to $p$-power Kummer cover]\label{remark: passing to p-power Kummer cover}
When we prove Theorem~\ref{thm:main_intro} later in this article, we occasionally need to ``replace $X$ by a $p$-power Kummer cover''. More precisely, assume that $K$ contains $\zeta_{p^n}$ for an integer $n\ge 1$ and we consider $p^n$-th Kummer cover $g\colon X_n\rightarrow X$. The new toric chart becomes $h'\colon X_n\rightarrow \D_K$ which is the composition of $g\colon X_n\rightarrow X$ and $h\colon X\rightarrow \D_K$ while $\Gamma$ stays the same. Note that $\Gamma_n\coloneqq \Gal(X_{\infty}/X_n)$ is an open normal subgroup of $\Gamma$ of index $p^{nd}$. A key observation is that 
\[(X_n)_{\infty}=X_n\times_{h', \D_K}\D_{\infty}\cong X_{\infty}\times (\Z/p^n\Z)^d\]
is isomorphic to $p^{nd}$ copies of $X_{\infty}$. Moreover, the $\Gamma$-action on $(X_n)_{\infty}$ is given by the action of $\Gamma_n$ on $X_{\infty}$ and the finite subgroup $\Gamma/\Gamma_n\cong (\Z/p^n\Z)^d$ permutes the $p^{nd}$ copies of $X_{\infty}$. If $\Lambda=(R_{\Lambda}, \mathbf{A}_{\Lambda})$ is a strongly decompleting datum for $(X,h)$ and $(R', \mathbf{A}')$ is the associated strongly decompleting datum for $(X_n, h')$, then $R'$ (resp. $\mathbf{A}'$) is isomorphic to the direct product of $p^{nd}$ copies of $R_{\Lambda}$ (resp. $\mathbf{A}_{\Lambda}$). The $\Gamma$-action on $R'$ (resp. $\mathbf{A}'$) is given by the $\Gamma_n$-action on $R_{\Lambda}$ (resp. $\mathbf{A}_{\Lambda}$) and the finite subgroup $\Gamma/\Gamma_n\cong (\Z/p^n\Z)^d$ permutes the $p^{nd}$ copies. The same assertions hold for all the induced perfect and imperfect relative period rings.
\end{remark}

%%% 4.7
\vspace{0.1in}
\subsection{Construction of $\mathbf{N}_{\mathrm{dR}}(\L)$}\label{ss:relative NdR}
\noindent
\vspace{0.1in}

\noindent As mentioned earlier, under the assumption of strongly decompleting data, we are able to generalize Berger's $\mathbf{N}_{\mathrm{dR}}$ functor (cf. \S \ref{sss:Berger_N_dR}) to the relative setting. 
Throughout \S \ref{ss:relative NdR}, let $h\colon X\rightarrow \D^d_K$ be a toric chart and let $\Lambda=(R_{\Lambda}, \mathbf{A}_{\Lambda})$ be a strongly decompleting datum.    

Let  $\mathbb{L}$ be a de Rham local system on $X$ and let $D_{\mathrm{rig}}^{\dagger}(\mathbb{L})$ be the relative \'etale $(\varphi, \Gamma)$-module over $\mathbf{C}^{\dagger}_{\Lambda}$ from Remark~\ref{remark:local_system_relation_with_period_ring_type_A}. 
Let us briefly explain the construction of another (not necessarily \'etale) $(\varphi, \Gamma)$-module $\mathbf{N}_{\mathrm{dR}}(\mathbb{L})$ over $\mathbf{C}^{\dagger}_{\Lambda}$ by a ``modification'' along the ``$t=0$ locus''. 
A key ingredient in the construction is \emph{Beauville--Laszlo gluing} (cf. \cite{beauvillelaszlo}) for finite projective modules. It is convenient to recall the following terminology used in \cite{relative}.

\begin{definition}\label{definition: BL data}
A \emph{Beauville--Laszlo datum} is a triple $(A, A_f, \widehat{A})$ where $A$ is a commutative ring, $f\in A$ is a non-zero-divisor, $A_f \coloneqq A[\frac{1}{f}]$, and $\widehat{A} \coloneqq \varprojlim_n A/f^nA.$
\end{definition}

For simplicity, let $M^{\dagger}$ denote $D_{\mathrm{rig}}^{\dagger}(\mathbb{L})$, which corresponds to a compatible system of finite projective $\mathbf{C}^{[s,r]}_{\Lambda}$-modules $M^{[s,r]}$, for all $r$ sufficiently small and all $0<s\le r$. We want to modify each $M^{[s,r]}$ along the ``$t=0$ locus'' to obtain a compatible system of $\mathbf{C}^{[s,r]}_{\Lambda}$-modules $N^{[s,r]}$ which glue to a $(\varphi, \Gamma)$-module $N^{\dagger}$ over $\mathbf{C}^{\dagger}_{\Lambda}$. (When $[s,r]\cap p^{\mathbb{Z}}=\emptyset$, we simply put $N^{[s,r]}=M^{[s,r]}$ as $t$ is invertible on $\mathbf{C}^{[s,r]}_{\Lambda}$.)

We first consider the situation where $\#\left([s,r]\cap p^{\mathbb{Z}}\right)=1$, say $p^{-n}\in [s,r]$. By Lemma~\ref{lemma: dense in B_n}, the injection $\iota_n\colon \mathbf{C}^{[s,r]}_{\Lambda}\hookrightarrow A_{X,n}[\![t, u_1, \ldots, u_d]\!]^{\nabla=0}$ has dense image with respect to the $t$-adic topology. In other words, $A_{X,n}[\![t, u_1, \ldots, u_d]\!]^{\nabla=0}$ is precisely the $t$-adic completion of $\mathbf{C}^{[s,r]}_{\Lambda}$. Hence \[(\mathbf{C}^{[s,r]}_{\Lambda}, \,\,\mathbf{C}^{[s,r]}_{\Lambda}[\frac{1}{t}],\,\, A_{X,n}[\![t, u_1, \ldots, u_d]\!]^{\nabla=0})\] is a Beauville--Laszlo datum. Consider the finite projective $\mathbf{C}^{[s,r]}_{\Lambda}[\frac{1}{t}]$-module $M^{[s,r]}[\frac{1}{t}]$ and the finite projective $A_{X,n}[\![t, u_1, \ldots, u_d]\!]^{\nabla=0}$-module 
\[N^{\mathrm{dR}}_n \coloneqq \left(D_{\mathrm{dR}}(\mathbb{L})\otimes_{A_X}A_{X,n}[\![t,u_1,\ldots,u_d]\!]\right)^{\nabla=0}\]
of the same rank. By \cite[Proposition 10.4]{relative}, for $r$ sufficiently small (i.e., for $n$ sufficiently large), there is a canonical isomorphism \footnote{
    The existence of the canonical isomorphism (\ref{eq:d_diff_isomorphism}) boils down to a relative version of Fontaine's $D_{\mathrm{dif}}$-functor (cf. \cite[Theorem 4]{relative}) and makes full use of Condition (c) in Definition~\ref{definition:decompleting_datum}.
    }
\begin{equation}
    \label{eq:d_diff_isomorphism}
 M^{[s,r]}\otimes_{\mathbf{C}^{[s,r]}_{\Lambda}, \iota_n} A_{X,n}[\![t, u_1, \ldots, u_d]\!]^{\nabla=0}[\frac{1}{t}]\cong N^{\mathrm{dR}}_n[\frac{1}{t}].
 \end{equation}
By Beauville--Laszlo gluing, we obtain a finite projective $\mathbf{C}^{[s,r]}_{\Lambda}$-module $N^{[s,r]}$. One checks that the construction above is compatible when we vary $r$ and $s$, and hence glue to $\mathbf{C}^{[s,r]}_{\Lambda}$-modules $N^{[s,r]}$ for arbitrary $[s,r]$ with $r$ sufficiently small (resp. $\mathbf{C}^r_{\Lambda}$-modules $N^r$ for all $r$ sufficiently small). One also checks that $N^{[s,r]}$ and $N^r$ possess natural $\Gamma$-actions and that $N^r$ possesses a natural Frobenius structure $\varphi$. Taking colimit as $r\rightarrow 0^+$, we obtain a $(\varphi, \Gamma)$-module $N^{\dagger}$ over $\mathbf{C}^{\dagger}_{\Lambda}$ (cf. \cite[\S 10.4]{relative}). Consequently, we arrive at a functor \[\mathbf{N}_{\mathrm{dR}}\colon  \mathrm{Loc}^{\mathrm{dR}}_{\Z_p} (X)\lra \mathrm{Mod}^{(\varphi, \Gamma)}_{/\mathbf{C}^{\dagger}_{\Lambda}}\]sending $\mathbb{L}\mapsto N^{\dagger}$.

\begin{remark}\label{remark: classical NdR via BL gluing}
Although we always assume $d\ge 1$ in the relative setting, if we take $d=0$, the same construction above remains valid (actually much simpler) and we obtain an alternative construction of Berger's $\mathbf{N}_{\mathrm{dR}}(V)$ in the classical setting. We refer the reader to \cite[Construction 3.38, Proposition 3.43]{relative} for more details. 
\end{remark}

\begin{remark}
By \cite[Theorem 2]{relative} (also see Theorem~\ref{intro:theorem_relative_Berger}), the modified relative $(\varphi, \Gamma)$-module $\mathbf{N}_{\mathrm{dR}}(\L)$ actually carries differential operators in both the arithmetic and geometric directions, compatible with the ``normalized'' differential operators on $\mathbf{C}^{\dagger}_{\Lambda}$. This makes $\mathbf{N}_{\mathrm{dR}}(\L)$ into a ``relative $p$-adic differential equation''. 
For simplicity of exposition, we omit the precise definitions here and refer the reader to \emph{loc.\ cit.} for more details. In the case where $X$ is a circle (\S \ref{ss: annulus--pDE}), we will give an explicit description of the differential operators on $\mathbf{N}_{\mathrm{dR}}(\mathbb{L})$. In particular, as we shall see, $\mathbf{C}^{\dagger}_{\Lambda}$ carries the ``correct'' differential operators $\frac{d}{d\pi_K}$ and $\frac{d}{d[T^{\flat}]}$.
\end{remark}

\begin{remark}\label{remark: tilde NdR--relative}
Let us point out that there is also a theory of ``modified $(\varphi, \Gamma)$-modules'' over the perfect period ring $\widetilde{\mathbf{C}}^{\dagger}_{X}$. This time, one considers the Beauville--Laszlo datum \[\left(\widetilde{\mathbf{C}}^{[s,r]}_{X},\,\, \widetilde{\mathbf{C}}^{[s,r]}_{X}[\frac{1}{t}],\,\, \mathbb{B}^+_{\mathrm{dR}, X}(X_{\infty})\right)\]when $\#\left([s,r]\cap p^{\mathbb{Z}}\right)=1$. For a de Rham local system $\mathbb{L}\in \mathrm{Loc}^{\mathrm{dR}}_{\Z_p}(X)$, consider the associated \'etale $(\varphi, \Gamma)$-module $\widetilde{M}^{\dagger} \coloneqq \widetilde{D}^{\dagger}_{\mathrm{rig}}(\mathbb{L})$ (see Remark~\ref{remark:local_system_relation_with_period_ring_type_A}), which corresponds to a compatible system of finite projective modules $\widetilde{M}^{[s,r]}$ over $\widetilde{\mathbf{C}}^{[s,r]}_X$. Then, by Beauville--Laszlo gluing, the finite projective $\widetilde{\mathbf{C}}^{[s,r]}_X[\frac{1}{t}]$-module $\widetilde{M}^{[s,r]}[\frac{1}{t}]$ and the finite projective $\mathbb{B}^+_{\mathrm{dR}, X}(X_{\infty})$-module 
\[\widetilde{N}^{\mathrm{dR}} \coloneqq  (D_{\mathrm{dR}}(\mathbb{L})\otimes_{A_X}\mathcal{O}\mathbb{B}^+_{\mathrm{dR}, X}(X_{\infty}) )^{\nabla=0}
\]
glue to a finite projective $\widetilde{\mathbf{C}}^{[s,r]}_X$-module $\widetilde{N}^{[s,r]}$. These $\widetilde{N}^{[s,r]}$ then glue to a $(\varphi, \Gamma)$-module $\widetilde{N}^r$ (resp. $\widetilde{N}^{\dagger}$) over $\widetilde{\mathbf{C}}^r_X$ (resp. $\widetilde{\mathbf{C}}^{\dagger}_X)$. This yields a functor 
\[\widetilde{\mathbf{N}}_{\mathrm{dR}}\colon  \mathrm{Loc}^{\mathrm{dR}}_{\Z_p} (X)\lra \mathrm{Mod}^{(\varphi, \Gamma)}_{/\widetilde{\mathbf{C}}^{\dagger}_X}\]
sending $\mathbb{L}\mapsto \widetilde{N}^{\dagger}$. We refer the reader to \cite[\S 10.1]{relative} for more details. It turns out $\mathbf{N}_{\mathrm{dR}}(\mathbb{L})$ is a descent of $\widetilde{\mathbf{N}}_{\mathrm{dR}}(\mathbb{L})$ along the arrow (11) in (\ref{diagram:big_fucking_diagram_of_rings}). However, it is unclear whether $\mathbf{N}_{\mathrm{dR}}(\mathbb{L})$ is the \emph{unique} such descent
(cf. Remark~\ref{remark:neither fully faithful nor essentially surjective}).

Finally, let us remark that the same construction applies to $d=0$, in which case we obtain a functor
\[\widetilde{\mathbf{N}}_{\mathrm{dR}}\colon  \mathrm{Rep}^{\mathrm{dR}}_{\Q_p} (\Gal_K)\lra \mathrm{Mod}^{(\varphi, \Gamma)}_{/\widetilde{\mathbf{C}}^{\dagger}_K}\]
from the category of de Rham representations of $\Gal_K$ to the category of $(\varphi, \Gamma_K)$-modules over $\widetilde{\mathbf{C}}^{\dagger}_K$. For any $V\in \mathrm{Rep}^{\mathrm{dR}}_{\Q_p}(\Gal_K)$, $\widetilde{\mathbf{N}}_{\mathrm{dR}}(V)$ is obtained from $\widetilde{D}^{\dagger}_{\mathrm{rig}, K}(V)$ by a modification along the $t=0$ locus, and $\mathbf{N}_{\mathrm{dR}}(V)$ is a descent of $\widetilde{\mathbf{N}}_{\mathrm{dR}}(V)$ along $\mathbf{C}^{\dagger}_K\rightarrow \widetilde{\mathbf{C}}^{\dagger}_K$. For more details, see \cite[\S 3.4.5]{relative}.
\end{remark}

\begin{remark}\label{remark: NdR r}
The $(\varphi, \Gamma)$-module $\mathbf{N}_{\mathrm{dR}}(\L)$ (resp. $\widetilde{\mathbf{N}}_{\mathrm{dR}}(\L)$) descends to a $(\varphi, \Gamma)$-module $\mathbf{N}^r_{\mathrm{dR}}(\L)$ (resp. $\widetilde{\mathbf{N}}^r_{\mathrm{dR}}(\L)$) over $\mathbf{C}^r_{\Lambda}$ (resp. $\widetilde{\mathbf{C}}^r_X$) for $r$ sufficiently small. The same statement holds for $\mathbf{N}_{\mathrm{dR}}(V)$ and $\widetilde{\mathbf{N}}_{\mathrm{dR}}(V)$ in the classical setting; the resulting descents are denoted $\mathbf{N}^r_{\mathrm{dR}}(V)$ and $\widetilde{\mathbf{N}}^r_{\mathrm{dR}}(V)$, respectively.
\end{remark}

%%% 4.8
\vspace{0.1in}
\subsection{Examples of strongly decompleting data}\label{ss: examples of strongly decompleting data} 
\noindent
\vspace{0.1in}

\noindent Let us recall from \cite[\S 11]{relative} two important classes of examples of strongly decompleting data: one for closed annuli and the other for ``closed annuli with punctures'' (more precisely, closed annuli of outer radius 1 with finitely many open unit disks removed from the boundary circle $|T|=1$). In particular, we provide explicit descriptions of the associated relative period rings, which will play essential roles in the remainder of the article.

\vspace{0.05in}
\subsubsection{Closed annuli}\label{sss: example closed annuli}
\noindent
\vspace{0.1in}

\noindent
We first consider the case of closed annuli. Let $[\sigma_1, \sigma_2]\subset (0,1]$ be a closed interval with $\sigma_1, \sigma_2\in p^{\Q}$, and consider the closed annulus $A_K([\sigma_1, \sigma_2];T)$. Let $\D_K$ be the closed unit disk and take the natural inclusion $\iota\colon  A_K([\sigma_1, \sigma_2];T)\hookrightarrow \D_K$ to be the toric chart. To simplify the notation, in the subscripts of period rings we write $[\sigma_1, \sigma_2]$ instead of $A_K([\sigma_1, \sigma_2];T)$. The constructions in \S \ref{ss:relative_period_rings:perfect} produce perfectoid rings $\widehat{A}_{[\sigma_1, \sigma_2], \infty}$ and $\widetilde{R}_{[\sigma_1, \sigma_2]}$ and yield perfect period rings
\[\widetilde{\mathbf{A}}_{[\sigma_1, \sigma_2]}, \widetilde{\mathbf{A}}^r_{[\sigma_1, \sigma_2]}, \widetilde{\mathbf{A}}^{\dagger}_{[\sigma_1, \sigma_2]}, \widetilde{\mathbf{B}}_{[\sigma_1, \sigma_2]}, \widetilde{\mathbf{B}}^r_{[\sigma_1, \sigma_2]}, \widetilde{\mathbf{B}}^{\dagger}_{[\sigma_1, \sigma_2]}, \widetilde{\mathbf{C}}^{[s,r]}_{[\sigma_1, \sigma_2]}, \widetilde{\mathbf{C}}^r_{[\sigma_1, \sigma_2]}, \widetilde{\mathbf{C}}^{\dagger}_{[\sigma_1, \sigma_2]}.\]
Explicitly, let $T^{\flat} \coloneqq (T, T^{\frac{1}{p}}, T^{\frac{1}{p^2}}, \ldots)\in \widetilde{R}_{\D_K}$ using the prescribed $p$-power roots of $T$, then $\widetilde{R}_{[\sigma_1, \sigma_2]}$ can be identified as the rational localization of \[\widetilde{R}_{\D_K}=\widehat{K}^{\flat}_{\infty}\langle (T^{\flat})^{\frac{1}{p^{\infty}}}\rangle \]
given by $\sigma_1\le |T^{\flat}|\le \sigma_2$; namely, $\widetilde{R}_{[\sigma_1, \sigma_2]}$ consists of $\sum_{j\in \mathbb{Z}[\frac{1}{p}]}a_j (T^{\flat})^j$ with $a_j\in \widehat{K}^{\flat}_{\infty}$ such that for any $\varepsilon>0$, there are only finitely many indices $j$ such that $|a_j|\cdot \sigma_1^j\ge \varepsilon$ and likewise only finitely many $j$ satisfying $|a_j|\cdot \sigma_2^j\ge \varepsilon$. All these period rings are equipped with commuting actions of $\varphi$ and $\Gamma\cong \Gamma_{\mathrm{geom}}\rtimes\Gamma_K$. The choices of $p$-power roots of $T$ and $p$-power roots of unity induce an isomorphism $\eta\colon  \Gamma_{\mathrm{geom}}\cong \Z_p$ and the cyclotomic character $\chi\colon \Gamma_K\hookrightarrow \Z_p^{\times}$. Let $[T^{\flat}]$ denote the Teichm\"uller lift of $T^{\flat}$ in $\widetilde{\mathbf{A}}_{[\sigma_1, \sigma_2]}=W(\widetilde{R}_{[\sigma_1, \sigma_2]})$.

Our goal is to construct explicit imperfect period rings for the annulus $A_K([\sigma_1, \sigma_2];T)$ so that the corresponding datum is strongly decompleting.

\begin{definition}\label{definition: description period rings--annulus}
\begin{enumerate}
\item[(1)] Let $R_{[\sigma_1, \sigma_2]}=\Gamma_{\kappa(\!(\overline{\pi}_K)\!)}([\sigma_1, \sigma_2];T^{\flat})$ (cf. \S \ref{sss:notation_cyclo}) as a subring of $\widetilde{R}_{[\sigma_1, \sigma_2]}$.  
In other words, $R_{[\sigma_1, \sigma_2]}$ is the ring of functions on the closed annulus \[A_{\kappa(\!(\overline{\pi}_K)\!)}([\sigma_1, \sigma_2];T^{\flat})\subset \mathbb{D}_{\kappa(\!(\overline{\pi}_K)\!)}.\]
\item[(2)] For any $0<r<r_K$, define $\mathbf{A}^r_{[\sigma_1, \sigma_2]}$ to be the ring of formal sums $\sum_{i\in \mathbb{Z}}b_i[T^{\flat}]^i$ with $b_i\in \mathbf{A}_K^r$ such that $\lambda_r(b_i)\cdot \sigma^{ri}\rightarrow 0$ as $i\rightarrow\pm\infty$, for all $\sigma_1\le \sigma\le \sigma_2$. Define a norm $\tilde{\lambda}_r$ on $\mathbf{A}^r_{[\sigma_1, \sigma_2]}$ by
\[\tilde{\lambda}_r\left(\sum_{i\in \mathbb{Z}}b_i[T^{\flat}]^i\right)=\sup_{\sigma\in [\sigma_1, \sigma_2]}\sup_{i\in \mathbb{Z}} \left\{\lambda_r(b_i)\cdot \sigma^{ri}\right\}.\]
\item[(3)] Let $\mathbf{A}_{[\sigma_1, \sigma_2]}^{\dagger}$ be the union $\mathbf{A}_{[\sigma_1, \sigma_2]}^{\dagger}  \coloneqq  \varinjlim_{0<r<r_K}\mathbf{A}^r_{[\sigma_1, \sigma_2]}$ and let $\mathbf{A}_{[\sigma_1, \sigma_2]}$ be the $p$-adic completion of $\mathbf{A}_{[\sigma_1, \sigma_2]}^{\dagger}$. Note that every element in $\mathbf{A}_{[\sigma_1, \sigma_2]}$ can be uniquely expressed as $\sum_{i\in \mathbb{Z}}b_i[T^{\flat}]^i$ with $b_i\in \mathbf{A}_K$.
\item[(4)] Let $\mathbf{B}_{[\sigma_1, \sigma_2]}$, $\mathbf{B}^r_{[\sigma_1, \sigma_2]}$, $\mathbf{B}^{\dagger}_{[\sigma_1, \sigma_2]}$ be the rings $\mathbf{A}_{[\sigma_1, \sigma_2]}$, $\mathbf{A}^r_{[\sigma_1, \sigma_2]}$, $\mathbf{A}^{\dagger}_{[\sigma_1, \sigma_2]}$ with $p$ inverted.
\item[(5)] For $0<s\le r<r_K$, let $\mathbf{C}^{[s,r]}_{[\sigma_1, \sigma_2]}$ be the completion of $\mathbf{B}^r_{[\sigma_1, \sigma_2]}$ with respect to the norm $\max\{\tilde{\lambda}_r, \tilde{\lambda}_s\}$ (or equivalently, with respect to the family of norms $\{\tilde{\lambda}_t\}_{t\in [s,r]}$). Let
\[\mathbf{C}^r_{[\sigma_1, \sigma_2]} \coloneqq \varprojlim_{0<s\le r}\mathbf{C}^{[s,r]}_{[\sigma_1, \sigma_2]},\qquad \textrm{and}\qquad \mathbf{C}^{\dagger}_{[\sigma_1, \sigma_2]}=\varinjlim_{0<r<r_K} \mathbf{C}^r_{[\sigma_1, \sigma_2]}.\]
\end{enumerate}
\end{definition}

\begin{proposition}[{\cite[Proposition 11.11, Lemma 11.2]{relative}}]\label{prop: annulus case is strongly decompleting}
\begin{enumerate}
\item[(1)] The datum $\left(R_{[\sigma_1, \sigma_2]}, \mathbf{A}_{[\sigma_1, \sigma_2]}\right)$ is strongly decompleting. The associated imperfect period rings are precisely 
\[\mathbf{A}^r_{[\sigma_1, \sigma_2]}, \mathbf{A}^{\dagger}_{[\sigma_1, \sigma_2]}, \mathbf{B}_{[\sigma_1, \sigma_2]}, \mathbf{B}^r_{[\sigma_1, \sigma_2]}, \mathbf{B}^{\dagger}_{[\sigma_1, \sigma_2]}, \mathbf{C}^{[s,r]}_{[\sigma_1, \sigma_2]}, \mathbf{C}^r_{[\sigma_1, \sigma_2]}, \mathbf{C}^{\dagger}_{[\sigma_1, \sigma_2]}.\]
\item[(2)] The norm $\tilde{\lambda}_r$ coincides with $\lambda_r$, for all $0<r<r_K$.
\end{enumerate}
\end{proposition}

We also have the following explicit descriptions of $\mathbf{C}^{[s,r]}_{[\sigma_1, \sigma_2]}$ and $\mathbf{C}^r_{[\sigma_1, \sigma_2]}$.

\begin{lemma}[{\cite[Corollary 11.5]{relative}}]\label{lemma: explicit description C--case of annulus}
\begin{enumerate}
\item[(1)] For $0<s\le r<r_K$, the ring $\mathbf{C}^{[s,r]}_{[\sigma_1, \sigma_2]}$ consists of $\sum_{i\in \mathbb{Z}} b_i[T^{\flat}]^i$ with $b_i\in \mathbf{C}^{[s,r]}_K$ such that $\lambda_{r'}(b_i)\cdot \sigma^{r'i}\rightarrow 0$ 
as $i\rightarrow \pm\infty$, for all $r'\in [s,r]$ and for all $\sigma\in [\sigma_1, \sigma_2]$.\footnote{In fact, one only needs to verify this condition for $r'=s,r$ and $\sigma=\sigma_1, \sigma_2$.}
\item[(2)] For $0<r<r_K$, the ring $\mathbf{C}^r_{[\sigma_1, \sigma_2]}$ consists of $\sum_{i\in \mathbb{Z}} b_i[T^{\flat}]^i$ with $b_i\in \mathbf{C}^r_K$ such that $\lambda_{r'}(b_i)\cdot \sigma^{r'i}\rightarrow 0$ 
as $i\rightarrow \pm\infty$, for all $r'\in (0,r]$ and for all $\sigma\in [\sigma_1, \sigma_2]$.
\end{enumerate}
\end{lemma}

\begin{remark}\label{remark: explicit description C--case of annulus}
By Proposition~\ref{prop:identification_C_with_annulus}, $\mathbf{C}^{[s,r]}_K$ can be identified with the ring of analytic functions $\Gamma_{K'_0}([\delta_1, \delta_2];\pi_K)$ on the closed annulus $A_{K'_0}([\delta_1, \delta_2];\pi_K)$ where $\delta_1=p^{-\frac{p}{p-1}\cdot\frac{r}{e_K}}$, $\delta_2=p^{-\frac{p}{p-1}\cdot\frac{s}{e_K}}$. Moreover, for any $r'\in [s,r]$, the $\lambda_{r'}$-norm on $\mathbf{C}^{[s,r]}_K$ coincides with the $\rho$-Gauss norm on $A_{K'_0}([\delta_1, \delta_2];\pi_K)$ with $\rho=p^{-\frac{p}{p-1}\cdot\frac{r'}{e_K}}$. Therefore, by Lemma~\ref{lemma: explicit description C--case of annulus}, the ring $\mathbf{C}^{[s,r]}_{[\sigma_1, \sigma_2]}$ consists of two-variable power series $\sum_{i,j\in \Z} c_{i,j}\pi_K^j[T^{\flat}]^i$ with $c_{i,j}\in K'_0$ such that
\[|c_{i,j}|\cdot p^{-\frac{p}{p-1}\cdot \frac{r'j}{e_K}}\cdot \sigma^{r'i}\rightarrow 0\]
as $|i|+|j|\rightarrow +\infty$, for all $r'\in [s,r]$ and all $\sigma\in [\sigma_1, \sigma_2]$ (again, one only needs to verify this condition for $r'=s,r$ and $\sigma=\sigma_1, \sigma_2$). 
\end{remark}

\subsubsection{Closed annuli with punctures} \label{sss: example closed annuli with punctures}
\noindent
\vspace{0.05in}

\noindent
We now consider a closed annulus of outer radius $1$, with finitely many open disks removed.

\begin{notation} \label{notation:closed unit disk removing open unit disks} 
Let $\delta\in (0,1]\cap p^{\Q}$ and let $f(T)\in \mathcal{O}_{K_0}[T]$ be a monic polynomial that factors as $f(T)=\prod_{i=1}^m (T-[\beta_i])$ for some $\beta_i\in \overline{k}_K^{\times}$. Let $X_{\delta,f}$ be the rational subset\footnote{Roughly speaking, $X_{\delta, f}$ is $\mathbb{D}_K$ minus finitely many open unit disks on the boundary circle $\{|T|=1\}$ (these disks are centered at $[\beta_i]$'s), then further remove an open disk of radius $\delta$ centered at $0$.} of $\D_K$ defined by
\[X_{\delta, f}\coloneqq \{|T|\ge \delta, \: |f(T)|=1\}\subset \D_K.\]
When $f$ is understood from the context, we often write $X_1$ for $X_{1, f}$ and $X_{\delta}$ for $X_{\delta, f}$. 
\end{notation}

\begin{remark} \label{rem: open-nbhds-type-II-V}
Let $\xi$ denote the Gaussian point of $\D_K$ and let $\xi$ denote the type V point corresponding to the open unit circle centered at 0. Then the collection $\Sigma=\{X_{1,f}\}$ forms a basis of open neighborhoods of $\xi_{\mathrm{I\!I}}$ inside $\mathbb{D}_K$; while the collection $\Sigma'=\{X_{\delta, f}\}$ forms a basis of open neighborhoods of $\xi$ inside $\mathbb{D}_K$.
\end{remark}

Let $\delta \in (0,1]\cap p^{\Q}$ and let $f \in \mathcal{O}_{K_0}[T]$ be a polynomial which takes the form $f(T) = \prod_i (T-[\beta_i]) \in \overline{K}[T]$ with $\beta_i \in \overline{k}_K^{\times}$. Take the natural inclusion $\iota: X_{\delta, f} \rightarrow \D_K$ to be the toric chart. To simplify the notation, in the subscripts of period rings we write $\delta, f$ instead of $X_{\delta, f}$. From the construction in \S \ref{ss:relative_period_rings:perfect}, we have perfectoid rings $\widehat{A}_{\delta, f, \infty}$ and $\widetilde{R}_{\delta, f}$ and associated perfect period rings
\[
\widetilde{\mathbf{A}}_{\delta, f}, \widetilde{\mathbf{A}}^r_{\delta, f}, \widetilde{\mathbf{A}}^{\dagger}_{\delta, f}, \widetilde{\mathbf{B}}_{\delta, f}, \widetilde{\mathbf{B}}^r_{\delta, f}, \widetilde{\mathbf{B}}^{\dagger}_{\delta, f}, \widetilde{\mathbf{C}}^{[s,r]}_{\delta, f}, \widetilde{\mathbf{C}}^r_{\delta, f}, \widetilde{\mathbf{C}}^{\dagger}_{\delta, f}
\]
equipped with commuting actions of $\varphi$ and $\Gamma$. Let $\overline{f}$ denote the reduction of $f$ modulo $p$. Then $\widetilde{R}_{\delta, f}$ is the rational localization of $\widetilde{R}_{[\delta, 1]}$ given by $|\overline{f}(T^{\flat})| = 1$.

We have the following explicit imperfect relative period rings for $X_{\delta, f}$ where the corresponding datum is strongly decompleting. Let $h$ denote the degree of the polynomial $f(T)$.

\begin{definition}
\begin{enumerate}
\item Let $R_{\delta, f}$ be the subring of $\widetilde{R}_{\delta, f}$ defined by
\[
R_{\delta, f}\coloneqq \left\{\sum_{i \in \Z} a_i(T^\flat)^i+ \sum_{j \ge 1}\sum_{k=0}^{h-1} \frac{b_{jk}(T^{\flat})^k}{\overline{f}(T^{\flat})^j} ~\Big|~ a_i, b_{jk} \in \kappa(\!(\overline{\pi}_K)\!), ~~|a_i|\sigma^i \rightarrow 0 \text{ as } i \rightarrow \pm \infty \text{ for all } \delta \le \sigma \le 1, ~~|b_{jk}| \rightarrow 0 \text{ as } j \rightarrow \infty \right\}.
\]

\item For each $0 < r < r_K$, let $\fA_{\delta, f}^r$ be a subring of $\widetilde{\fA}_{\delta, f}$ given by
\[
\fA_{\delta, f}^r\coloneqq \left\{\sum_{i\in \Z} a_i[T^{\flat}]^i+\sum_{j \ge 1}\sum_{k=0}^{h-1} \frac{b_{jk} [T^{\flat}]^k}{f([T^{\flat}])^j} ~\Big|~ a_i, b_{jk} \in \fA_K^r, \lambda_r(a_i)\sigma^{ri} \rightarrow 0 \text{ as } i \rightarrow \pm \infty \text{ for all } \delta \le \sigma \le 1, \lambda_r(b_{jk}) \rightarrow 0 \text{ as } j \rightarrow \infty \right\}.
\]
Consider the norm $\tilde{\lambda}_r$ on $\fA_{\delta, f}^r$ given by
\[
\tilde{\lambda}_r\left(\sum_{i\in \Z} a_i[T^{\flat}]^i+\sum_{j \ge 1}\sum_{k=0}^{h-1}b_{jk}\frac{[T^{\flat}]^k}{f([T^{\flat}])^j}\right)\coloneqq \max\left\{ \sup_{\delta \le \sigma \le 1} \sup_{i\in \Z} \{\lambda_r(a_i)\sigma^{ri}\}, \sup_{j \ge 1, 0 \le k \le h-1}\{\lambda_r(b_{jk})\} \right\}. 
\]

\item Let 
\[
\fA_{\delta, f}^{\dagger} \coloneqq \varinjlim_{0 < r < r_K}  \fA_{\delta, f}^r,
\]
and let $\fA_{\delta, f}$ be the $p$-adic completion of $\fA_{\delta, f}^{\dagger}$. Note that $\fA_{\delta, f}$ is a subring of $\widetilde{\fA}_{\delta, f}$.

\item Let $\mathbf{B}_{\delta, f}$, $\mathbf{B}^r_{\delta, f}$, $\mathbf{B}^{\dagger}_{\delta, f}$ be the rings $\mathbf{A}_{\delta, f}$, $\mathbf{A}^r_{\delta, f}$, $\mathbf{A}^{\dagger}_{\delta, f}$ with $p$ inverted.

\item For $0<s\le r<r_K$, let $\mathbf{C}^{[s,r]}_{\delta, f}$ be the completion of $\mathbf{B}^r_{\delta, f}$ with respect to the norm $\max\{\tilde{\lambda}_r, \tilde{\lambda}_s\}$ (or equivalently, with respect to the family of norms $\{\tilde{\lambda}_t\}_{t\in [s,r]}$). Let
\[
\mathbf{C}^r_{\delta, f} \coloneqq \varprojlim_{0<s\le r}\mathbf{C}^{[s,r]}_{\delta, f},\qquad \textrm{and}\qquad \mathbf{C}^{\dagger}_{\delta, f}=\varinjlim_{0<r<r_K} \mathbf{C}^r_{\delta, f}.
\]
\end{enumerate}    
\end{definition}

\begin{proposition}[{\cite[Proposition~11.19, Lemma~11.14]{relative}}]
\begin{enumerate}
\item The datum $\left(R_{\delta, f}, \mathbf{A}_{\delta, f}\right)$ is strongly decompleting. The associated imperfect period rings are precisely 
\[
\mathbf{A}^r_{\delta, f}, \mathbf{A}^{\dagger}_{\delta, f}, \mathbf{B}_{\delta, f}, \mathbf{B}^r_{\delta, f}, \mathbf{B}^{\dagger}_{\delta, f}, \mathbf{C}^{[s,r]}_{\delta, f}, \mathbf{C}^r_{\delta, f}, \mathbf{C}^{\dagger}_{\delta, f}.
\]

\item The norm $\tilde{\lambda}_r$ coincides with $\lambda_r$, for all $0<r<r_K$.
\end{enumerate}
\end{proposition}

We also have the explicit descriptions of $\mathbf{C}^{[s,r]}_{\delta, f}$ and $\mathbf{C}^r_{\delta, f}$.

\begin{lemma}[{\cite[Corollary~11.15]{relative}}] \label{lem: puncture-C-description}
\begin{enumerate}
\item For $0 < s \le r < r_K$, the ring $\fC_{\delta, f}^{[s, r]}$ consists of elements
\[
\sum_{i\in \Z} a_i[T^{\flat}]^i+\sum_{j \ge 1}\sum_{k=0}^{h-1}b_{jk}\frac{[T^{\flat}]^k}{f([T^{\flat}])^j}, \quad a_i, b_{jk} \in \fC_K^{[s, r]}
\]
such that $\lambda_{r'}(a_i)\sigma^{r'i} \rightarrow 0$ as $i \rightarrow \pm \infty$ for all $r' \in [s, r]$ and $\delta\le \sigma \le 1$, and $\lambda_{r'}(b_{jk}) \rightarrow 0$ as $j \rightarrow \infty$ for all $r' \in [s, r]$.

\item For $0 < r < r_K$, the ring $\fC_{\delta, f}^r$ consists of
\[
\sum_{i\in \Z} a_i[T^{\flat}]^i+\sum_{j \ge 1}\sum_{k=0}^{h-1}b_{jk}\frac{[T^{\flat}]^k}{f([T^{\flat}])^j}, \quad a_i, b_{jk} \in \fC_K^r
\]
such that $\lambda_{r'}(a_i)\sigma^{r'i} \rightarrow 0$ as $i \rightarrow \pm \infty$ for all $r' \in (0, r]$ and $\delta\le \sigma \le 1$, and $\lambda_{r'}(b_{jk}) \rightarrow 0$ as $j \rightarrow \infty$ for all $r' \in (0, r]$.
\end{enumerate}    
\end{lemma}

\vspace{0.1in}
%%% 4.9
\subsection{From de Rham local systems to $p$-adic differential equations}\label{ss: annulus--pDE}
\noindent
\vspace{0.1in}

\noindent Keep the notations in \S \ref{ss: examples of strongly decompleting data}. Let $X$ be either $X_{[\sigma_1, \sigma_2]}$ or $X_{\delta, f}$, and let $\L$ be a de Rham local system on $X$. We recall the construction of the associated $p$-adic differential equation $\fN_{\mathrm{dR}}(\L)$ from \cite{relative}. For simplicity, we mainly explain the case $X = X_{[\sigma_1, \sigma_2]}$; analogous statements hold for $X = X_{\delta, f}$.

Let $X = X_{[\sigma_1, \sigma_2]} = A_K([\sigma_1, \sigma_2]; T)$, and recall the perfect and imperfect period rings from \S \ref{ss: examples of strongly decompleting data}. The machinery in \S \ref{section:relative_period_rings} yields finite projective $(\varphi, \Gamma)$-modules $\widetilde{D}^{\dagger}_{\mathrm{rig}}(\L)$ and $D^{\dagger}_{\mathrm{rig}}(\L)$ over $\widetilde{\fC}_{[\sigma_1, \sigma_2]}^{\dagger}$ and $\fC_{[\sigma_1, \sigma_2]}^{\dagger}$ respectively. Applying the machinery in \S \ref{ss:relative NdR} and the fact that $(R_{[\sigma_1, \sigma_2]}, \mathbf{A}_{[\sigma_1, \sigma_2]})$ is strongly decompleting, we obtain $(\varphi, \Gamma)$-modules $\widetilde{\mathbf{N}}_{\mathrm{dR}}(\L)$ and $\mathbf{N}_{\mathrm{dR}}(\L)$ over $\widetilde{\mathbf{C}}^{\dagger}_{[\sigma_1, \sigma_2]}$ and $\mathbf{C}^{\dagger}_{[\sigma_1, \sigma_2]}$. We first equip $\mathbf{N}_{\mathrm{dR}}(\L)$ with natural differential operators following \cite{relative}, making it into a ``relative $p$-adic differential equation''. 

Recall that $\Gamma\cong \Gamma_{\mathrm{geom}}\rtimes \Gamma_K$ acts continuously on the $\mathbf{C}$-type imperfect period rings $\mathbf{C}^{[s,r]}_{[\sigma_1, \sigma_2]}$, $\mathbf{C}^r_{[\sigma_1, \sigma_2]}$, and $\mathbf{C}^{\dagger}_{[\sigma_1, \sigma_2]}$. For every $\gamma\in \Gamma_{\mathrm{geom}}$ and any integer $n\ge 1$, there is a unique $\eta_n(\gamma)\in \mathbb{Z}/p^n\mathbb{Z}$ such that $\gamma(T^{1/p^n})=\zeta_{p^n}^{\eta_n(\gamma)}\cdot T^{1/p^n}$; these give rise to an element $\eta(\gamma)\in \mathbb{Z}_p$ such that $\eta(\gamma)\!\mod p^n=\eta_n(\gamma)$ for all $n$. We obtain a continuous group homomorphism $\eta\colon  \Gamma_{\mathrm{geom}}\rightarrow \mathbb{Z}_p$ which is an isomorphism.%\footnote{The composition of $\eta$ with the injection $\Gal (K_{a, \infty}/K_{\infty}) \ra \Gamma_{\mathrm{geom}}$ (cf. Lemma~\ref{lemma: Gamma_K,a vs Gamma} and its proof) is simply the map $c$ from (\ref{eq:cocycle_becomes_homomorphism}.)} 
On the other hand, we have the usual cyclotomic character $\chi\colon \Gamma_K\hookrightarrow \Z_p^{\times}$. The action of $\mathrm{Lie}\,\Gamma$ on these imperfect relative period rings induces differential operators
\[d_{\mathrm{geom}} \coloneqq \frac{\log(\gamma_{\mathrm{geom}})}{\eta(\gamma_{\mathrm{geom}})}, \qquad \textrm{and}\qquad d_{\mathrm{arith}} \coloneqq \frac{\log(\gamma_{\mathrm{arith}})}{\log(\chi(\gamma_{\mathrm{arith}}))}\]
where $\gamma_{\mathrm{geom}}$ (resp. $\gamma_{\mathrm{arith}}$) is any element in $\Gamma_{\mathrm{geom}}$ (resp. in $\Gamma_K$) sufficiently close to 0 (resp. to 1). A direct computation shows that
$d_{\mathrm{geom}}=t\cdot[T^{\flat}]\cdot \frac{\partial}{\partial [T^{\flat}]}$ 
and
$d_{\mathrm{arith}}=t\cdot (1+\pi)\cdot\frac{\partial}{\partial \pi}=t\cdot (1+\pi)\cdot u_0\cdot\frac{\partial}{\partial \pi_K}$ 
for some unit $u_0\in \mathbf{C}^{\dagger}_K$ (cf. Remark~\ref{remark:partial vs d/d pi_K}). By shrinking $r$, we may assume that $u_0$ is a unit in $\mathbf{C}^r_K$. We refer the reader to \cite[\S 11.3]{relative} for more details. 

Given any finite projective module $M$ over $\mathbf{C}^{\dagger}_{[\sigma_1, \sigma_2]}$ (or $\mathbf{C}^{[s,r]}_{[\sigma_1, \sigma_2]}$, $\mathbf{C}^r_{[\sigma_1, \sigma_2]}$) equipped with a continuous $\Gamma$-action, the action of $\mathrm{Lie}\,\Gamma$ also induces differential operators $\nabla_{\mathrm{geom}}$ and $\nabla_{\mathrm{arith}}$ on $M$ defined by
\[\nabla_{\mathrm{geom}}=\frac{\log(\gamma_{\mathrm{geom}})}{\eta(\gamma_{\mathrm{geom}})}, \qquad \textrm{and}\qquad \nabla_{\mathrm{arith}}=\frac{\log(\gamma_{\mathrm{arith}})}{\log(\chi(\gamma_{\mathrm{arith}}))}.\]
For our purpose, we normalize these differential operators by putting
\begin{align*}
    D_{[T^{\flat}]} & \coloneqq t^{-1}\cdot [T^{\flat}]^{-1}\cdot \nabla_{\mathrm{geom}} \\
    D_{\pi_K} & \coloneqq t^{-1}\cdot ((1+\pi)\cdot u_0)^{-1}\cdot \nabla_{\mathrm{arith}}
\end{align*}
We caution the reader that, \emph{a priori}, $(M,D_{\pi_K})$ (resp. $(M, D_{[T^{\flat}]})$) could have poles along the vanishing locus of $t=\log(1+\pi)$ when viewed as a differential module over the differential ring $\left(\mathbf{C}^{\dagger}_{[\sigma_1, \sigma_2]}, \frac{d}{d\pi_K}\right)$ (resp. $\left(\mathbf{C}^{\dagger}_{[\sigma_1, \sigma_2]},\frac{d}{d[T^{\flat}]}\right)$).

Recall that $\mathbf{N}_{\mathrm{dR}}(\mathbb{L})$ is obtained from $D^{\dagger}_{\mathrm{rig}}(\mathbb{L})$ by modifying along the ``$t=0$ locus''. In fact, these operators $D_{[T^{\flat}]}$ and $D_{\pi_K}$ are precisely the geometric and arithmetic differential operators alluded to in the statement of 
\ref{intro:theorem_relative_Berger} in the introduction (cf. \cite[Theorem 2, Proposition 10.9]{relative}), which asserts that $\left(\mathbf{N}_{\mathrm{dR}}(\mathbb{L}), D_{[T^{\flat}]}, D_{\pi_K}\right)$ does not have poles, hence it is a genuine differential module (i.e., a relative $p$-adic differential equation) over the differential ring $\left(\mathbf{C}^{\dagger}_{[\sigma_1, \sigma_2]},\frac{d}{d[T^{\flat}]},\frac{d}{d\pi_K}\right)$. In other words, the modification along the ``$t=0$ locus'' removes the singularities from $D^{\dagger}_{\mathrm{rig}}(\mathbb{L})$.

When $X = X_{\delta, f}$, we similarly obtain $\fN_{\mathrm{dR}}(\L)$ from $D_{\mathrm{rig}}^{\dagger}(\L)$, and $\fN_{\mathrm{dR}}(\L)$ is a differential module over the differential ring $\left(\mathbf{C}^{\dagger}_{\delta, f},\frac{d}{d[T^{\flat}]},\frac{d}{d\pi_K}\right)$ (see \cite[\S 11.3]{relative}).

\vspace{0.1in}
\subsection{A criterion for semistability}\label{ss:unipotence_and_st}
\noindent
\vspace{0.1in}

\noindent For later use, we establish a criterion for semistability of the $p$-adic local system $\L$. Throughout \S \ref{ss:unipotence_and_st}, we assume $(X,h)$ is one of the following:
    \begin{enumerate}
        \item Let $X$ be a closed annulus $X_{[\sigma_1, \sigma_2]} = \{\sigma_1 \le |T| \le \sigma_2\}$, with $\sigma_1, \sigma_2 \in p^\Q$ and $0<\sigma_1\le \sigma_2\le 1$, as in \S \ref{sss: example closed annuli}. Let $h\colon X_{[\sigma_1, \sigma_2]}\hookrightarrow \D_K$ be the natural embedding.
        \item Let $X$ be a closed annulus with punctures $X_{\delta, f} = \{ \delta \le |T| \le 1,\, |f(T)| = 1 \}$ as in \S \ref{sss: example closed annuli with punctures}. Let $h\colon X_{\delta, f}\hookrightarrow \D_K$ be the natural embedding.
        \item Let $X$ be a finite \'etale cover of either $X_{[\sigma_1, \sigma_2]}$ or $X_{\delta, f}$ and let $h$ be the composition of this finite \'etale map and the natural embedding into $\D_K$.
     \end{enumerate}     
In \cite[\S 11]{relative}, we have explicitly constructed a strongly decompleting datum for $(X,h)$ (also see \S \ref{ss: examples of strongly decompleting data} and Proposition~\ref{prop: strongly decompleting datum finite etale}). 
Let $\mathbf{A}_{X}^\dagger$, $\mathbf{B}_{X}^\dagger$, $\mathbf{C}_{X}^\dagger$, etc. denote the corresponding imprefect relative period rings. Let $\L$ be a de Rham $\Z_p$-local system of rank $d$ on $X$ and let $D^{\dagger}_{\mathrm{rig}}(\L)$ be the associated \'etale $(\varphi, \Gamma)$-module over $\mathbf{C}^{\dagger}_X$.
  
\begin{proposition} \label{prop:unipotence_and_st}
Suppose that the $\mathbf{C}_{X}^\dagger [\log [T^\flat], \log \pi_K, 1/t]$-module\footnote{Note that the $(\varphi, \Gamma)$-actions on $\mathbf{C}^{\dagger}_X$ naturally extends to $\mathbf{C}_{X}^\dagger [\log [T^\flat], \log \pi_K, 1/t]$.}
\[D^{\dagger}_{\mathrm{rig}}(\L)\otimes_{\mathbf{C}^{\dagger}_X} \mathbf{C}_{X}^\dagger [\log [T^\flat], \log \pi_K, 1/t]\]
is free and admits a basis $w_1, \ldots, w_d$ such that $\gamma(w_i)=w_i$ for all $\gamma\in \Gamma$ and all $i=1,\ldots, d$. Then $\L$ is pointwise semistable on $X$ (i.e., semistable at every classical point on $X$).  
\end{proposition}

We point out that there is no direct ``specialization map'' relating $D^{\dagger}_{\mathrm{rig}}(\L)$ and $D^{\dagger}_{\mathrm{rig}}(\L_x)$ for a classical point $x$. Instead, we directly relate $D^{\dagger}_{\mathrm{rig}}(\L)$ with $D_{\mathrm{st}}(\L_x)$ and use the $\Gamma$-invariant basis to produce enough semistable periods at each classical point.

\begin{proof} 
Let $x\in X$ be a $L$-rational classical point for some finite extension $L/K$. Then $x$ is represented by a ring map $f_x\colon A_X\rightarrow L$ (actually $x$ is determined by the $\Gal(L/K)$-orbit of such a ring map). Consider the ring map $h^*\colon K\langle T\rangle\rightarrow A_X$ induced by the toric chart $h$. Let $a\coloneqq f_x\circ h^*(T)\in L^{\times}$ where $T$ is the coordinate of $\D_K$. Pick a compatible $p$-power roots $\{a^{1/p^m}\}_{m \ge 0}$ of $a$. This induces a map $\eta\colon \widehat{K}_{\infty}\langle T^{1/p^{\infty}}\rangle\rightarrow \C_K$ sending $T^{1/p^m} \mapsto a^{1/p^m}$. Since 
\[\widehat{A}_{X, \infty}=A_X\widehat{\otimes}_{K\langle T\rangle}\widehat{K}_{\infty}\langle T^{1/p^{\infty}}\rangle,\]
the maps $f_x$ and $\eta$ induces a ring map $e_x\colon \widehat A_{X, \infty} \ra \C_K$.

In fact, the map $e_x$ can be made $\Gal_L$-equivariant in a suitable sense. First, the action of $\Gal_L$ on $a^{1/p^m}$, given by $g (a^{1/p^m}) = \zeta_{p^m}^{c_m (g)} a^{1/p^m}$ for a unique $c_m(g) \in \Z/p^m$, determines a continuous cocycle $c\colon  \Gal_L \ra \Z_p$ relative to the cyclotomic character; namely, $c(g_1g_2) = c(g_1) + \chi(g_1) c(g_2)$. Next, recall that the choice of $\{T^{1/p^m}\}_{m \ge 0}$ determines a splitting $\Gamma\cong \Gamma_{\mathrm{geom}}\rtimes \Gamma_K$ (cf. \S \ref{ss:toric_tower}). Then consider a group homomorphism $\rho_x\colon  \Gal_L \ra \Gamma$ sending $g \mapsto (c(g), g|_{K_\infty})$. Now, consider the action of $\Gal_L$ on $\widehat{A}_{X, \infty}$ given by the natural $\Gamma$-action and $\rho_x\colon \Gal_L\rightarrow \Gamma$. Then $e_x\colon \widehat{A}_{X, \infty}\rightarrow \C_K$ is $\Gal_L$-equivariant. Indeed, it suffices to check $\Gal_L$-equivariance for $\eta$. On $\widehat{K}_{\infty}$, the map $\eta$ is just the identity and $\Gal_L$ acts naturally on both the source and target; on the other hand, we have 
\[\eta (\rho_x(g) \cdot (T^{1/p^m})) = \eta (\zeta_{p^m}^{c(g)} T^{1/p^m}) =  \zeta_{p^m}^{c(g)} a^{1/p^m} = g \cdot (a^{1/p^m}) = g \cdot \eta (T^{1/p^m}),\]
as desired.

The map $e_x$ induces a $\Gal_L$-equivariant map $e^\flat_x\colon  \sq R_X = \widehat A_{X, \infty}^\flat \ra \C_K^\flat$ which sends $T^\flat \mapsto a^\flat \coloneqq (a, a^{1/p}, \ldots )$ and thus a $\Gal_L$-equivariant map $\tilde{e}_x\colon \widetilde{\mathbf{A}}_X\rightarrow \widetilde{\mathbf{A}}$. It further induces $\Gal_L$-equivariant maps $\tilde{e}_x\colon \sq{\mathbf{A}}_{X}^\dagger \ra  \sq{\mathbf{A}}^\dagger$ and $\tilde{e}_x\colon \sq{\mathbf{C}}_{X}^\dagger \ra  \sq{\mathbf{C}}^\dagger$. Composing with the natural inclusion $\mathbf{C}^{\dagger}_X\subset \widetilde{\mathbf{C}}^{\dagger}_X$, we arrive at the ``specialization map'' $s_x\colon \mathbf{C}^{\dagger}_X\rightarrow \widetilde{\mathbf{C}}$. This in turn induces a $\Gal_L$-equivariant specialization map
\[ s_x\colon \mathbf{C}_{X}^\dagger [\log [T^\flat], \log \pi_K, 1/t] \ra  \sq{\mathbf{C}}_{\log}^\dagger [1/t]\]
sending $\log [T^\flat] \mapsto \log [a^\flat]$ and $\log \pi_K \mapsto \log \pi_K$.  

Observe that if $\mathbb M$ is a finite locally free $\Z/p^m$-local system on $\spec \sq R_X$ and $D_m(\mathbb M)$ is the associated $\varphi$-module over $W_m(\sq R_X)$ by \cite[Proposition 3.2.7]{Kedlaya-LiuI}, then we have a canonical $\varphi$- and $\Gal_L$-equivariant isomorphism
\[
D_m(\mathbb M) \otimes_{W_m(\sq R_X), e^\flat_x} W_m (\C_K^\flat) \cong D_m ((e^\flat_x)^* \mathbb M)
\]
where $(e^\flat_x)^* \mathbb M$ denotes the local system on $\spec \C_K^\flat$ obtained by base change along $e^\flat_x\colon  \sq R_{X} \ra \C_K^\flat$. Consequently, by unwinding the construction of $\widetilde{D}^{\mathrm{int}} (\L)$ (cf. Remark~\ref{remark:local_system_relation_with_period_ring_type_A} and the references therein), we obtain a canonical $\varphi$- and $\Gal_L$-equivariant isomorphism
$
    \widetilde{D}^{\mathrm{int}}(\L) \otimes_{\sq{\mathbf{A}}_{X}, \tilde{e}_x}  \sq{\mathbf{A}}  \cong  \L_{x} \otimes_{\Z_p} \sq{\mathbf{A}}, 
$  
which further induces an isomorphism 
$ 
\widetilde{D}^{\mathrm{int}, \dagger}(\L) \otimes_{\sq{\mathbf{A}}_{X}^\dagger, \tilde{e}_x}  \sq{\mathbf{A}}^\dagger 
\cong \L_{x} \otimes_{\Z_p} \sq{\mathbf{A}}^\dagger
$ by \cite[Theorem 8.5.3 (d)(e)]{Kedlaya-LiuI}. This further induces a $\varphi$- and $\Gal_L$-equivariant isomorphism 
\[
D^{\dagger}_{\mathrm{rig}}(\L) \otimes_{\mathbf{C}_{X}^\dagger, s_x}  \sq{\mathbf{C}}^\dagger=\widetilde{D}^{\dagger}_{\mathrm{rig}}(\L) \otimes_{\sq{\mathbf{C}}_{X}^\dagger, \tilde{e}_x}  \sq{\mathbf{C}}^\dagger
\cong \L_{x} \otimes_{\Z_p} \sq{\mathbf{C}}^{\dagger}.
\]
Finally, joining $\log[T^{\flat}]$, $\log\pi_K$, $1/t$ and using the $\Gamma$-invariant basis in the statement of the proposition, we arrive at $\Gal_L$-equivariant isomorphisms
\[
\L_{x} \otimes_{\Z_p} \sq{\mathbf{C}}_{\log}^\dagger [1/t] \cong \left( \bigoplus_{i = 1}^d  \mathbf{C}_{X}^\dagger [\log [T^\flat], \log \pi_K, 1/t] \cdot w_i  \right) \otimes_{\mathbf{C}_{X}^\dagger [\log [T^\flat], \log \pi_K, 1/t], s_x}   \sq{\mathbf{C}}_{\log}^\dagger [1/t]  \cong   \bigoplus_{i = 1}^d  \sq{\mathbf{C}}_{\log}^\dagger [1/t] \cdot w_i.
\]
Here, $\Gal_L$ acts on $w_i$ via $\rho_x\colon \Gal_L \ra \Gamma$, thus $\Gal_L$ fixes each $w_i$. By \cite[Proposition 3.4]{Berger-differential} (combined with a Tate twist argument), we have
\[
D_{\mathrm{st}, L} (\L_x) \cong \left( \L_{x} \otimes_{\Z_p} \sq{\mathbf{C}}_{\log}^\dagger [1/t]  \right)^{\Gal_L}.
\]
Therefore, $D_{\mathrm{st}, L} (\L_x)$ has $L_0$-dimension $d$ (where $L_0$ denotes the maximal unramified subfield in $L$) and thus $\L_x$ is semistable, as desired. 
\end{proof}

\vspace*{0.3in}
%%% Section 5
\section{Around points of type I and I\!I}\label{section: type 1 and 2}

As mentioned in the introduction, to prove Theorem~\ref{thm:main_intro}, we first prove the following local version of the $p$-adic monodromy theorem.

\begin{theorem}[The local $p$-adic monodromy theorem for curves]\label{thm: local main thm}
Let $X$ be a smooth rigid analytic curve over $K$, viewed as an adic space over $\spa(K, \mathcal{O}_K)$. Let $\L$ be a de Rham local system on $X$. Let $\xi\in X_{\C_K}$ be a point (of any type). Then there exists a finite extension $L/K$, an open neighborhood $U$ of the image of $\xi$ in $X_L$, and a finite \'etale cover $f\colon V\rightarrow U$ such that the pullback $f^*(\mathbb{L}|_U)$ is pointwise semistable (i.e., semistable at every classical point in $V$).
\end{theorem}

We prove Theorem~\ref{thm: local main thm} by a case-by-case analysis based on the type of the point $\xi$; in fact, points of type I, I\!I, I\!I\!I, I\!V, V are treated in \S \ref{ss: type 1}, \S \ref{ss: type 2}, \S \ref{section: type 3}, \S \ref{section: type 4}, \S \ref{section: type 5}--\ref{section: type V general}, respectively.

\vspace{0.1in}
%% 5.1
\subsection{Around type I points}\label{ss: type 1}
\noindent
\vspace{0.1in}

\noindent The case of classical type I points is handled in \cite{Shimizu_2}. In fact, we have the following generalization of \cite[Theorem 9.2]{Shimizu_2}, which also deal with non-classical type I points.

\begin{proposition}\label{prop:Shimizu_generalization}
Let $X$ be a smooth rigid analytic space over $K$ (viewed as an adic space) and let $\L$ be a de Rham $\Z_p$-local system on $X$. Let $x\in X_{\C_K}$ be a type I point. Then there exist a finite extension $L/K$ and an open neighborhood $U\subset X_{L}$ of the image $x'$ of $x$ in $X_L$ such that $\L|_U$ is semistable at every classical point in $U$. 
\end{proposition}

\begin{proof}
With some modifications, the proof of \cite[Theorem 9.2]{Shimizu_2} remains valid. First, by the same argument as in \cite[Proposition 9.6]{Shimizu_2}, $x$ admits an open neighborhood $U_x$ that is isomorphic to a closed unit polydisk over $\C_K$. The $p$-adic Riemann-Hilbert functor of Liu--Zhu \cite{Liu_Zhu} yields a vector bundle with connection $(D_{\mathrm{dR}}(\L), \nabla)$ over $X$. By base change, we obtain a vector bundle with connection $(D_{\mathrm{dR}}(\L)_{\C_K}, \nabla)$ over $X_{\C_K}$. By \cite[Theorem 9.7]{Shimizu_2}, up to shrinking the polydisk, we may assume that the connection is trivial on $U_x$.

Notice that $U_x$ contains a $K'$-rational point for some finite extension $K'/K$. In particular, the image $U'_x$ of $U_x$ in $X_{K'}$ is isomorphic to the closed unit polydisk over $K'$. By construction, the vector bundle with connection $(D_{\mathrm{dR}}(\L)_{K'}, \nabla)$ is trivial over $U'_x$, hence $\L|_{U'_x}$ is horizontally de Rham. By shrinking, we may replace $U'_x$ by some $U''_x$ isomorphic to 
\[\mathbb{T}^n_{K'}\cong \spa\left(K'\langle T_1^{\pm}, \ldots, T_n^{\pm}\rangle, \mO_{K'}\langle T_1^{\pm}, \ldots, T_n^{\pm}\rangle\right).\]
Finally, by the same argument as in \cite[\S 9.4]{Shimizu_2}, there exists a finite extension $L/K'$ such that $\L|_{U''_{x, L}}$ is horizontally de Rham and semistable at every classical point.
\end{proof}

\vspace{0.1in}
%% 5.2
\subsection{Around type I\!I points}\label{ss: type 2}
\noindent
\vspace{0.1in}

\noindent Recall that every type I\!I point is $K'$-rational for some finite extension $K'/K$. Therefore, when $\xi$ is a type I\!I point, Theorem~\ref{thm: local main thm} follows from the following proposition.

\begin{proposition}\label{prop: local main thm type II}
Let $X$ and $\L$ be as in Theorem~\ref{thm: local main thm}. Let $\xi$ be a $K$-rational type I\!I point on $X$.   Then, after replacing $K$ by a finite extension if necessary,\footnote{Since $\xi$ is $K$-rational, for any finite extension $K'/K$, the preimage of $\xi$ in $X_{K'}$ is a single point of type I\!I (cf. Remark~\ref{remark:K_rational_point_under_field_extension}). Hence, there is no ambiguity when enlarging $K$.} there exists an open neighborhood $U$ of $\xi$ and a finite \'etale cover $U'\rightarrow U$ such that the pullback of $\mathbb{L}$ to $U'$ is pointwise semistable (i.e., semistable at every classical point on $U'$).
\end{proposition}

\begin{proof}
We first claim that, by enlarging $K$ if necessary, $\xi$ admits a connected affinoid open neighborhood $U$ which admits a smooth formal model over $\mathcal{O}_K$. Indeed, by \cite[Proposition 2.2.1]{Berkovich_differential} (which corrects the original statements in \cite[Proposition 3.6.1]{Berkovich_etale}), by enlarging $K$ if necessary, there exist
\begin{itemize}
\item a smooth connected projective curve $\mathcal{Y}$ over $\mathcal{O}_K$ with generic fiber $Y$ (viewed as an adic space over $\spa(K, \mathcal{O}
_K)$);
\item finitely manly closed points $z_1, \ldots, z_n$ on the special fiber $Y_s$; and  
\item an open neighborhood $U$ of $\xi$,
\end{itemize} 
such that $U$ is isomorphic to $\mathrm{sp}^{-1}\left(Y_s-\{z_1, \ldots, z_n\}\right)$, where $\mathrm{sp}\colon  Y\rightarrow Y_s$ is the specialization map, and that $\xi$ is the preimage of the generic point of $Y_s$. It follows from the construction that $U$ is affinoid connected and it admits a smooth affine formal integral model.

Given this, the rest of the proof follows from the arguments in \cite[\S 6 \& 7]{GY}. For the reader's convenience, we recall the proof here (using terminology from \emph{loc.\ cit.}). Notice that $\xi$ is the unique Shilov point of $U$. Let $L=\mathscr{H}(\xi)$ be the completed residue field of $\xi$, which is a complete discrete valuation field. Write $U=\mathcal{U}_{\eta}$ for some smooth connected formal scheme $\mathcal{U}=\spf R$ over $\mathcal{O}_K$. The restriction $\mathbb{L}|_{\xi}$ corresponds to a de Rham representation $V$ of $ \Gal_L$. By either \cite[Main Theorem]{Ohkubo-p-adic-monodromy} or \cite[Corollary 1.2]{Morita-p-adic-monodromy}, there exists a finite Galois extension $L'/L$ such that $V|_{\Gal_{L'}}$ is a semistable representation. According to \cite[Theorem 6.11, Proposition 6.16, and Theorem 6.17]{GY}, by shrinking $\mathcal{U}$ and replacing $L'$ by a further finite extension if necessary, there exists a finite Galois extension $K'/K$ and a finite cover $\mathcal{U}'\rightarrow \mathcal{U}$ such that
\begin{itemize}
\item $\mathcal{U}'$ is affine smooth connected over $\mathcal{O}_{K'}$;
\item $\mathcal{U}'\rightarrow \mathcal{U}$ factors as $\mathcal{U}'\rightarrow \mathcal{U}_{\mathcal{O}_{K'}}\rightarrow \mathcal{U}$ where $\mathcal{U}_{\mathcal{O}_{K'}}\rightarrow \mathcal{U}$ is the extension of scalars and $\mathcal{U}'\rightarrow \mathcal{U}$ induces a finite \'etale map $\mathcal{U}'_{\eta}\rightarrow \mathcal{U}_{\eta, K'}$ on the generic fibers;
\item the inverse image of $\xi$ in $\mathcal{U}'_{\eta}$ is the unique Shilov point $\xi'$ of $\mathcal{U}'_{\eta}$;
\item the extension $\mathscr{H}(\xi)\rightarrow \mathscr{H}(\xi')$ coincides with $L\rightarrow L'$.
\end{itemize}
In particular, the restriction $\mathbb{L}|_{\xi'}$ is a semistable representation of $\Gal_{L'}$. Let $U'=\mathcal{U}'_{\eta}$. By the purity theorem of Du--Liu--Moon--Shimizu (Theorem~\ref{thm: purity DLMS}), we conclude that $\mathbb{L}|_{U'}$ is a semistable local system; in particular, the pullback of $\mathbb{L}$ to $U'$ is pointwise semistable.
\end{proof}

\vspace{0.3in}
%%% Section 6
\section{Around points of type I\!I\!I}\label{section: type 3}

The goal of this section is to prove the local $p$-adic monodromy theorem (Theorem~\ref{thm: local main thm}) around type I\!I\!I points. Since each type I\!I\!I point in $X_{\C_K}$ is $K'$-rational for some finite extension $K'/K$, it is equivalent to prove the following theorem.

\begin{theorem}\label{thm: local main thm type III}
Let $X$ and $\L$ be as in Theorem~\ref{thm: local main thm}. Let $\xi$ be a $K$-rational type I\!I\!I point on $X$. Then, after replacing $K$ by a finite extension if necessary, there exists an open neighborhood $U$ of $\xi$ and a finite \'etale cover $U'\rightarrow U$ such that the pullback of $\mathbb{L}$ to $U'$ is pointwise semistable (i.e., semistable at every classical point on $U'$).
\end{theorem}

To prove Theorem~\ref{thm: local main thm type III}, we will carefully study relative $p$-adic differential equations constructed in \cite{relative} (recalled in \S \ref{section:relative_period_rings}) over suitable open neighborhoods of the type I\!I\!I points $\xi$. Our strategy is to relate a completed colimit of the imperfect relative period rings over neighborhoods of $\xi$ to the \emph{fake Robba ring} studied by Kedlaya in \cite{kedlaya_fake} and \cite{kedlaya_st3}. Using a form of $p$-adic local monodromy theorem proved in \emph{loc.\ cit.} for these fake Robba rings, we are able to achieve unipotence of the differential operators ``at $\xi$'' after pulling back to a finite \'etale cover, and then we ``propagate'' this unipotence to an open neighborhood of $\xi$. 
In a similar spirit as in the proof of the classical $p$-adic monodromy theorem (cf. \S \ref{ss: proof of classical p-adic monodromy theorem}), we shall see that unipotence of differential operators allows us to produce enough semistable periods in a neighborhood of $\xi$, by passing to a further finite \'etale cover if necessary (see \S \ref{ss:contagion_of_goem_unipotence}). Since this propagation process is quite technical, we defer its discussion to \S \ref{ss:geometric_uni_special}--\ref{ss:geometric_uni_general}: the special case is treated in \S \ref{ss:geometric_uni_special}, and the general case in \S \ref{ss:geometric_uni_general}.

\vspace{0.1in}
%%% 6.1
\subsection{The fake Robba ring (\`a la Kedlaya)} 
\label{ss:fake_Robba_ring} 
\noindent 
\vspace{0.1in}

\noindent
By \cite[Proposition~2.2.1]{Berkovich_differential}, up to replacing $K$ by a finite extension, we may assume that $\xi$ admits an open neighborhood of the form $X_{[\sigma_1, \sigma_2]} = A_{K} ([\sigma_1, \sigma_2])$ (i.e. the closed annulus $\{\sigma_1 \le |T| \le \sigma_2 \}$) with $\sigma_1, \sigma_2 \in p^\Q$ and $\xi$ is the type I\!I\!I point corresponding to the circle of radius $\sigma_{\xi}$ where $\sigma_1<\sigma_{\xi}< \sigma_2$ and $\sigma_{\xi}\not\in p^{\mathbb{Q}}$. Write $\sigma_\xi = p^{-\frac{p}{p-1} \frac{1}{e_K} \cdot \alpha}$. Then
\[ \alpha = - \frac{(p-1)e_K}{p} \cdot\log_p (\sigma_\xi)  = - \frac{1}{c_K} \cdot \log_p (\sigma_\xi)\] is an irrational number, where $c_K = \frac{p}{(p-1)e_K}$. Similarly, consider the rational numbers $\alpha_i = - \frac{1}{c_K} \log_p (\sigma_i)$ for $i = 1, 2$. If $[\sigma_1, \sigma_2]$ runs through closed intervals containing $\sigma_{\xi}$ with $\sigma_1, \sigma_2\in p^{\Q}$, then the closed annuli $X_{[\sigma_1, \sigma_2]}$'s form a basis of open neighborhoods of $\xi$. Recall the imperfect relative period rings
\begin{equation}\label{eq: imperfect relative period rings sigma_1 sigma_2}
\mathbf{A}^r_{[\sigma_1, \sigma_2]}, \mathbf{A}^{\dagger}_{[\sigma_1, \sigma_2]}, \mathbf{B}_{[\sigma_1, \sigma_2]}, \mathbf{B}^r_{[\sigma_1, \sigma_2]}, \mathbf{B}^{\dagger}_{[\sigma_1, \sigma_2]}, \mathbf{C}^{[s,r]}_{[\sigma_1, \sigma_2]}, \mathbf{C}^r_{[\sigma_1, \sigma_2]}, \mathbf{C}^{\dagger}_{[\sigma_1, \sigma_2]}
\end{equation}
from \S \ref{ss: examples of strongly decompleting data}. Let us introduce some auxiliary period rings $\mathbf{A}^r_{\xi}$, $\mathbf{A}^{\dagger}_{\xi}$, $\mathbf{B}^r_{\xi}$, $\mathbf{B}^{\dagger}_{\xi}$, $\mathbf{C}^{[s,r]}_{\xi}$, $\mathbf{C}^r_{\xi}$, and $\mathbf{C}^{\dagger}_{\xi}$, which are certain ``completed colimits'' of the period rings in \eqref{eq: imperfect relative period rings sigma_1 sigma_2}.

\begin{definition}
    \begin{enumerate}
        \item On each $\mathbf{A}^r_{[\sigma_1, \sigma_2]}$, we define a norm $\lambda_{r,\xi}$ by setting
\[\lambda_{r,\xi}\left(\sum_{i\in \mathbb{Z}}c_i[T^{\flat}]^i\right)=\sup_{i\in \mathbb{Z}}\left\{\lambda_r(c_i)\cdot \sigma_{\xi}^{ri}\right\}\]
for each $\sum_{i\in \mathbb{Z}}c_i[T^{\flat}]^i\in \mathbf{A}^r_{[\sigma_1, \sigma_2]}$ with $c_i\in \mathbf{A}^r_K$. Then we define $\mathbf{A}^r_{\xi}$ to be the completion of the colimit 
\[\mathbf{A}^r_{\xi, \mathrm{colim}}=\varinjlim_{\sigma_1<\sigma_{\xi}<\sigma_2} \mathbf{A}^r_{[\sigma_1, \sigma_2]}\]
with respect to $\lambda_{r,\xi}$. It follows that $\mathbf{A}^r_{\xi}$ consists of $\sum_{i\in \mathbb{Z}}c_i[T^{\flat}]^i$ with $c_i\in \mathbf{A}_K^r$ such that $\lambda_r(c_i)\cdot \sigma_{\xi}^{ri}\rightarrow 0$ as $i\rightarrow\pm\infty$. 
        \item Define $\mathbf{B}^r_{\xi}=\mathbf{A}^r_{\xi}[\frac{1}{p}]$. 
        \item Define $\mathbf{A}_{\xi}^{\dagger}$ (resp. $\mathbf{B}_{\xi}^{\dagger}$) as the colimit $\varinjlim_{r \ra 0^+} \mathbf{A}_{\xi}^{r}$ (resp. $\varinjlim_{r \ra 0^+}\mathbf{B}_{\xi}^{r}$).
    \end{enumerate}
\end{definition}
  The norm $\lambda_{r,\xi}$ naturally extends to $\mathbf{A}^r_{\xi}$ and $\mathbf{B}^r_{\xi}$. If we express an element in $\mathbf{B}^r_{\xi}$ as $x = \sum_{i, j} c_{i, j} \pi_K^j [T^\flat]^i$ with $c_{i,j}\in K'_0$, then
\[\lambda_{r, \xi} (x) = \sup_{i, j \in \Z} \left\{|c_{i, j}| \cdot p^{-\frac{p}{p-1} \frac{r j}{e_K}} \cdot \sigma_{\xi}^{ri}\right\}
= \sup_{i, j \in \Z} \left\{ |c_{i, j}| \cdot p^{-\frac{p}{p-1} \frac{1}{e_K} \cdot r \cdot (j +  \alpha i)}\right\}.  
\]

\begin{definition}
    \begin{enumerate} 
        \item Define  $\mathbf{C}^{[s, r]}_{\xi}$ as the completion of $\mathbf{B}_{\xi}^r$ under the norm $\max\{\lambda_{s, \xi}, \lambda_{r, \xi}\}$ and let $\mathbf{C}^r_{\xi}$ be the Fr\'echet completion of $\mathbf{B}_{\xi}^r$ with respect to the family of norms $\{\lambda_{s, \xi}\}_{s\in (0, r]}$.
        \item Define $\mathbf{C}_{\xi}^{\dagger}=\varinjlim_{r \ra 0^+}\mathbf{C}_{\xi}^{r}$.
    \end{enumerate}
\end{definition}

Note that $\mathbf{C}^{\dagger}_{\xi}$ and $\mathbf{C}^r_{\xi}$ are equipped with Frobenius structures, as well as geometric and arithmetic differential operators, inherited from those of $\mathbf{C}^{\dagger}_{[\sigma_1, \sigma_2]}$ and $\mathbf{C}^r_{[\sigma_1, \sigma_2]}$. 

\begin{remark} \label{remark:fake_robba_ring}
Equivalently, we can define $\mathbf{C}^{[s, r]}_{\xi}$  as the completion of $K_0' [\pi_K^{\pm 1}, [T^\flat]^{\pm 1}]$ under the norm $\max\{\lambda_{s, \xi}, \lambda_{r, \xi}\}$, that is, it consists of $\sum_{i\in \mathbb{Z}}c_i[T^{\flat}]^i$ with $c_i\in \mathbf{C}_K^{[s, r]}$ such that $\lambda_t(c_i)\cdot \sigma_{\xi}^{ti}\rightarrow 0$ as $i\rightarrow\pm\infty$ for both $t = s$ and $t = r$. This is precisely the \emph{fake annulus} $\mathscr R^{\lambda}_{[p^{-c_K r}, p^{- c_K s}]}$ in the language of \cite[Definition~5.1.1]{kedlaya_st3}\footnote{Fake annuli and fake Robba rings are first introduced in \cite{kedlaya_fake} using a different notation.} where $\lambda$ stands for the linear functional $\lambda\colon \mathbb{R}^2\rightarrow \mathbb{R}$ sending $(x,y)\mapsto x+\alpha y$. Similarly, $\mathbf{C}^{r}_{\xi}$ (resp. $\mathbf{C}^{\dagger}_{\xi}$) can be identified with the fake annulus $\mathscr R^{\lambda}_{[p^{-c_K r}, 1)}$ (resp. the \emph{fake Robba ring} $\mathscr R^{\lambda}$) in \emph{loc.\ cit.} Moreover, under this identification, the fake Robba ring $\mathscr{R}^{\lambda}$ is equipped with a Frobenius structure in the sense of \cite[Definition~5.1.6]{kedlaya_st3}.
\end{remark}

\begin{remark}  \label{remark:description_R_alpha}
The mod $p$ reduction $R_{\xi}\coloneqq \mathbf{A}_{\xi}^{\dagger}/p$ can be described as the completion of $\kappa [(\cl \pi_K)^{\pm 1}, (T^\flat)^{\pm 1}]$ under the valuation 
\begin{equation}
   \label{eq:v_alpha_on_mod_p} 
v_{\alpha} \left(\sum_{i, j} b_{i, j} (\cl \pi_K)^j (T^\flat)^i \right) =  \frac{p}{p-1} \cdot \frac{1}{e_K} \cdot \inf_{b_{i, j} \ne 0} (j+ \alpha i)
\end{equation}
(or equivalently, completion with respect to the norm $p^{-v_{\alpha}}$). Note that $R_{\xi}$ is a field. In fact, if we consider the type I\!I\!I point $\cl \xi$ on the annulus $A_{\kappa (\!(\cl \pi_K)\!)} ([\sigma_1, \sigma_2]; T^\flat)$ (in characteristic $p$) corresponding to the circle of radius $\sigma_\xi$, then $R_{\xi}$ agrees with the completed residue field $\sH(\overline{\xi})$. 
\end{remark}

\vspace{0.1in}
%%% 6.2
\subsection{Geometric realization of extensions of fake Robba rings} \label{ss:geometric_realization_fake} 
\noindent 
\vspace{0.1in}

\noindent
For any interval $[\sigma_1, \sigma_2]$ as above, the machinery in \S \ref{ss: annulus--pDE} produces a finite projective module $\mathbf{N}_{\mathrm{dR},[\sigma_1, \sigma_2]}=\mathbf{N}_{\mathrm{dR}}(\L|_{X_{[\sigma_1, \sigma_2]}})$ over $\mathbf{C}^{\dagger}_{[\sigma_1, \sigma_2]}$, equipped with both an arithmetic operator $D_{\pi_K}$ and a geometric differential operator $D_{[T^{\flat}]}$, as well as a Frobenius structure. We also consider the
integrable and Frobenius compatible connection (``the full connection'')
\[
\nabla = D_{\pi_K} d \pi_K + D_{[T^\flat]} d [T^\flat], 
\]
on $\mathbf{N}_{\mathrm{dR},[\sigma_1, \sigma_2]}$. Moreover, for $r$ sufficiently small, we have a finite projective $(\varphi, \nabla)$-module $\mathbf{N}_{\mathrm{dR}, [\sigma_1, \sigma_2]}^r$ over $\mathbf{C}^r_{[\sigma_1, \sigma_2]}$, also equipped with differential operators $D_{\pi_K}$ and $D_{[T^{\flat}]}$ (cf. Remark~\ref{remark: NdR r}). The formations of $\mathbf{N}_{\mathrm{dR},[\sigma_1, \sigma_2]}$ and $\mathbf{N}^r_{\mathrm{dR},[\sigma_1, \sigma_2]}$ are compatible with restriction to a smaller $[\sigma_1, \sigma_2]$ according to \cite[Proposition~10.13]{relative}. Taking colimits as $\sigma_1, \sigma_2 \ra \sigma_\xi$, we obtain finite projective $(\varphi, \nabla)$-modules
\[\mathbf{N}_{\mathrm{dR}, \mathrm{colim}}=\varinjlim_{\sigma_1 < \sigma_\xi < \sigma_2}\mathbf{N}_{\mathrm{dR}, [\sigma_1, \sigma_2]},\qquad \mathbf{N}^r_{\mathrm{dR}, \mathrm{colim}}=\varinjlim_{\sigma_1 < \sigma_\xi < \sigma_2}\mathbf{N}^r_{\mathrm{dR}, [\sigma_1, \sigma_2]}\]
over $\mathbf{C}_{\xi, \mathrm{colim}}^{\dagger} = \varinjlim_{\sigma_1 < \sigma_\xi < \sigma_2} \mathbf{C}_{[\sigma_1, \sigma_2]}^{\dagger}$ and $\mathbf{C}_{\xi, \mathrm{colim}}^{r} = \varinjlim_{\sigma_1 < \sigma_\xi < \sigma_2} \mathbf{C}_{[\sigma_1, \sigma_2]}^{r}$, respectively. Note that $\mathbf{C}_{\xi}^{r}$ is the Fr\'echet completion of $\mathbf{C}_{\xi, \mathrm{colim}}^{r}$ under the family of norms $\{\lambda_{s, \xi}\}_{s \in (0, r]}$. 
Write 
\[
\mathbf{N}_{\mathrm{dR},\xi}^r \coloneqq \mathbf{N}_{\mathrm{dR}, \mathrm{colim}}^r \otimes_{\mathbf{C}_{\xi, \mathrm{colim}}^{r}} \mathbf{C}_{\xi}^{r} =\mathbf{N}_{\mathrm{dR}, [\sigma_1, \sigma_2]}^r \otimes_{\mathbf{C}_{\xi, [\sigma_1, \sigma_2]}^{r}} \mathbf{C}_{\xi}^{r}  
\]
and 
\[
\mathbf{N}_{\mathrm{dR},\xi} \coloneqq 
\varinjlim_{r \ra 0^+} \mathbf{N}_{\mathrm{dR},\xi}^r= \mathbf{N}_{\mathrm{dR}, \mathrm{colim}} \otimes_{\mathbf{C}_{\mathrm{colim}}^{\dagger}} \mathbf{C}_{\xi}^{\dagger}=\mathbf{N}_{\mathrm{dR}, [\sigma_1, \sigma_2]} \otimes_{\mathbf{C}_{[\sigma_1, \sigma_2]}^{\dagger}} \mathbf{C}_{\xi}^{\dagger}. 
\]
Recall that $\mathbf{C}_{\xi}^{\dagger}$ is isomorphic to the fake Robba ring $\mathscr R^{\lambda}$. We view $(\mathbf{N}_{\mathrm{dR},\xi}, \nabla)$ as a finite projective $(\phi, \nabla)$-module over $\mathscr R^{\lambda}$, also equipped with differential operators $D_{\pi_K}$ and $D_{[T^\flat]}$. 

The following is Kedlaya's $p$-adic local monodromy theorem for fake Robba rings.

\begin{theorem}[{\cite[Theorem~5.2.4]{kedlaya_fake}}] \label{thm:monodromy_fake_annuli}
There exists a finite separable extension $R'/R_{\xi}$ with associated fake Robba ring $\mathscr{R}'$ such that $\mathbf{N}_{\mathrm{dR},\xi}\otimes_{\mathscr{R}^{\lambda}}\mathscr{R}'$ becomes unipotent for the full connection $\nabla$; that is, it is a successive extension of trivial connections over $\mathscr{R}'$ of rank 1. 
\end{theorem}

For the reader's convenience, let us explain the construction of the fake Robba ring $\mathscr R'$ associated to $R'$. Firstly, the valuation $v_\alpha$ in \eqref{eq:v_alpha_on_mod_p} uniquely extends to $R'$, which we still denote as $v_a$ by a slight abuse of notation. The value group of $v_{\alpha}$ contains $\Z + \alpha \Z$ as a subgroup of finite index. In particular, $R'$ is a \emph{nearly monomial field} in the sense of \cite{kedlaya_fake} (see \cite[Definition~2.1.1, Remark~2.1.4]{kedlaya_fake}). Let $(\mathscr R^{\mathrm{int}\prime}, R')$ be the henselian pair that is finite \'etale over the henselian pair $(\mathbf{A}_{\xi}^{\dagger}, R_{\xi})$. 
Then $\mathscr R'$ is defined as
\[
\mathscr{R}' \coloneqq  \mathbf{C}_{\xi}^{\dagger} \otimes_{\mathbf{A}_{\xi}^{\dagger}}  \mathscr R^{\mathrm{int}\prime}.
\]
Moreover, for sufficiently small $r$, the finite \'etale map $\mathbf{A}_{\xi}^{\dagger}\rightarrow \mathscr{R}^{\mathrm{int}\prime}$ descends to a finite \'etale map $\mathbf{A}^r_{\xi}\rightarrow {\mathscr{R}^{\mathrm{int}\prime}}^r$. Then we put 
\[
{\mathscr R'}^r \coloneqq \mathbf{C}_{\xi}^r \otimes_{\mathbf{A}_{\xi}^r}  {\mathscr{R}^{\mathrm{int}\prime}}^r,
\]
and we have ${\mathscr R'}=\varinjlim_{r\rightarrow 0^+}{\mathscr R'}^r$. By \cite[\S 2]{kedlaya_fake} (in particular, see Proposition~2.3.4 and Remark~2.3.5), $\mathscr{R}^{\mathrm{int}\prime}$ (resp. $\mathscr{R}'$) is indeed an integral fake Robba ring (resp. a fake Robba ring). Since $\mathbf{C}^{\dagger}_{\xi}\rightarrow \mathscr{R}'$ (resp. $\mathbf{C}^r_{\xi}\rightarrow {\mathscr{R}'}^r$) is finite \'etale, the differential operators $\frac{d}{d[T^{\flat}]}$ and $\frac{d}{d\pi_K}$ uniquely extend to $\mathscr{R}'$ (resp. ${\mathscr{R}'}^r$). As a corollary of Theorem~\ref{thm:monodromy_fake_annuli}, there exists a sufficiently small $r$, such that the full connection $\nabla$ on $\mathbf{N}_{\mathrm{dR},\xi}^r \otimes_{\mathbf{C}_{\xi}^{r}} {\mathscr R}'^r$ is unipotent.

In the remainder of \S \ref{ss:geometric_realization_fake}, let us realize this cover of fake Robba rings $\mathbf{C}^{\dagger}_{\xi} \ra \mathscr R'$ in geometric terms. First, write $R_{[\sigma_1, \sigma_2]} \coloneqq  \mathbf{A}^\dagger_{[\sigma_1, \sigma_2]} /p = \mathbf{A}^r_{[\sigma_1, \sigma_2]}/p$ for $r$ sufficiently small, as in \S \ref{ss: examples of strongly decompleting data} (also see \cite[\S 11]{relative}). Recall that $R_{[\sigma_1, \sigma_2]}$ is the ring of analytic functions on the closed annulus over $\kappa (\!(\cl \pi_K)\!)$ with radius in $[\sigma_1, \sigma_2]$ and variable $T^\flat$; in particular, it is an affinoid $\kappa (\!(\cl \pi_K)\!)$-algebra.  Write $\cl \lambda_{[\sigma_1, \sigma_2]}$ for the corresponding spectral norm. Equivalently, define $v_{\alpha_i}$ as in \eqref{eq:v_alpha_on_mod_p} but use $\alpha_i$ instead of $\alpha$. Then $R_{[\sigma_1, \sigma_2]}$ is the completion of $\kappa [(\cl \pi_K)^{\pm 1}, (T^\flat)^{\pm 1}]$ with respect to the valuation $\min\{v_{\alpha_1}, v_{\alpha_2}\}$. Moreover, write $R_{\xi, \mathrm{colim}}$ for the colimit $\varinjlim_{\sigma_\xi \in [\sigma_1, \sigma_2]} R_{[\sigma_1, \sigma_2]}$. Then $R_{\xi}$ agrees with the completion of $R_{\xi, \mathrm{colim}}$ under the valuation $v_{\alpha}$ from \eqref{eq:v_alpha_on_mod_p}.

\begin{lemma} \label{lemma:geometric_realization_fake_Robba}
     There exist $\sigma_1, \sigma_2$ sufficiently close to $\sigma_{\xi}$, and a finite \'etale map $R_{[\sigma_1, \sigma_2]} \ra R'_{[\sigma_1, \sigma_2]}$ such that 
     \[R'=R_{\xi}\otimes_{R_{[\sigma_1, \sigma_2]}} R'_{[\sigma_1, \sigma_2]}.\]
\end{lemma}

\begin{proof}
    The assertion is an immediate consequence of Corollary~\ref{cor: equiv fet R_alpha} below. 
\end{proof}

\begin{lemma}\label{lemma: colimit of R[sigma1, sigma2]}
\begin{enumerate}
        \item For $\sigma_1 \le \sigma_1' < \sigma_2' \le \sigma_2$, the restriction map $R_{[\sigma_1, \sigma_2]} \ra R_{[\sigma_1', \sigma_2']}$ is submetric. 
        \item  The colimit $R_{\xi, \mathrm{colim}}\coloneqq \varinjlim_{\sigma_1<\sigma_{\xi}<\sigma_2} R_{[\sigma_1, \sigma_2]}$ is a field. 
        \item Consider the seminorm $\cl \lambda_{\alpha}$ on $\varinjlim_{\sigma_1<\sigma_{\xi}<\sigma_2} R_{[\sigma_1, \sigma_2]}$ defined as the infimum of the quotient seminorms induced by  $\cl \lambda_{[\sigma_1, \sigma_2]}$'s. Then $\cl \lambda_{\alpha}$ agrees with $p^{-v_{\alpha}}$. In particular, it is multiplicative and the completion of $\varinjlim_{\sigma_1<\sigma_{\xi}<\sigma_2} R_{[\sigma_1, \sigma_2]}$ with respect to $\cl \lambda_{\alpha}$ is precisely $R_{\xi}$. 
    \end{enumerate}
\end{lemma}

\begin{proof}
   This is standard. In particular, $R_{\xi, \mathrm{colim}}$ is simply the (uncompleted) residue field $\mathbf{k}(\overline{\xi})$ of the type I\!I\!I point $\overline{\xi}$ in Remark~\ref{remark:description_R_alpha}.
\end{proof}

\begin{corollary}\label{cor: equiv fet R_alpha}
    The base functor induces an equivalence of categories of finite \'etale algebras.
    \[
    \left(\varinjlim_{\sigma_1<\sigma_{\xi}<\sigma_2}  R_{[\sigma_1, \sigma_2]} \right)_{\textup{f\'et}} \isom \left(R_{\xi}\right)_{\textup{f\'et}}.
    \]
\end{corollary}

\begin{proof}
    This follows from Lemma~\ref{lemma: colimit of R[sigma1, sigma2]} and \cite[Lemma~2.2.4]{Kedlaya-LiuI}. 
\end{proof}

\begin{construction}[Geometric realization of $\mathbf{C}^{\dagger}_{\xi} \ra \mathscr R'$] 
\label{construction:geometric_realization}
Write $A_{[\sigma_1, \sigma_2]} = \mO({X_{[\sigma_1, \sigma_2]}})$ for the ring of analytic functions on the annulus $X_{[\sigma_1, \sigma_2]}$ and recall the perfectoid algebra $\widehat{A}_{[\sigma_1, \sigma_2], \infty}$ from \S \ref{ss: examples of strongly decompleting data}. By the diagram of equivalences of categories in \cite[Proposition~5.28]{relative}, the finite \'etale map $R_{[\sigma_1, \sigma_2]} \ra R'_{[\sigma_1, \sigma_2]}$ from Lemma~\ref{lemma:geometric_realization_fake_Robba} corresponds to a finite \'etale map $\widehat{A}_{[\sigma_1, \sigma_2], \infty} \ra B$. By \cite[\S 7]{perfectoid} (also see \cite[Proposition~2.6.8]{Kedlaya-LiuI}), up to replacing $K$ by $K (\zeta_{p^n})$, we may assume that this finite \'etale map descends to a finite \'etale map $A_{[\sigma_1, \sigma_2],n} \ra B_0$ where $A_{[\sigma_1, \sigma_2],n}$ is the ring of analytic functions on the $p^n$-th Kummer cover $X_{[\sigma_1, \sigma_2],n}$ of $X_{[\sigma_1, \sigma_2]}$. Note that, for the purpose of proving Theorem~\ref{thm: local main thm type III}, we are always allowed to extend the base field and replace $X_{[\sigma_1, \sigma_2]}$ by a finite \'etale cover. By Remark~\ref{remark: passing to p-power Kummer cover}, we may further assume that $n=0$; namely, there exists a finite \'etale cover $g\colon  X'_{[\sigma_1, \sigma_2]}=\spa(B_0, B_0^{\circ}) \ra X_{[\sigma_1, \sigma_2]}$ so that the corresponding strongly decompleting data given by \cite[Proposition~5.28]{relative} (also see Proposition~\ref{prop: strongly decompleting datum finite etale}) gives rise to relative imperfect period rings ${\mathbf{A}'}_{[\sigma_1, \sigma_2]}^{r}$, ${\mathbf{B}'}_{[\sigma_1, \sigma_2]}^{r}$, ${\mathbf{C}'}_{[\sigma_1, \sigma_2]}^{r}$, etc., satisfying ${\mathbf{A}'}_{[\sigma_1, \sigma_2]}^{r}/p = R'_{[\sigma_1, \sigma_2]}$ for $r$ sufficiently small. Let $\xi'$ denote the pre-image of $\xi$ in $X_{[\sigma_1, \sigma_2]}'$. This is again a point of type I\!I\!I. By Lemma~\ref{lemma:geometric_realization_fake_Robba} and the henselianity of $\mathscr{R}^{\mathrm{int}\prime}$, we have
\[\mathscr{R}^{\mathrm{int}\prime} = \mathbf{A}_{\xi}^\dagger \otimes_{\mathbf{A}_{[\sigma_1, \sigma_2]}^r} {\mathbf{A}'}_{[\sigma_1, \sigma_2]}^{r} \] for $r$ sufficiently small, and thus  ${\mathscr R'}^r = \mathbf{C}_{\xi}^{r}
 \otimes_{\mathbf{C}^{r}_{[\sigma_1, \sigma_2]}}  {\mathbf{C}'}_{[\sigma_1, \sigma_2]}^{r}$. Let \[\mathbf{C}_{\xi', \mathrm{colim}}^{r} \coloneqq \varinjlim_{\sigma_1<\sigma_{\xi}<\sigma_2} {\mathbf{C}'}_{[\sigma_1, \sigma_2]}^{r}.\] Then ${\mathscr R'}^r$ can be identified as the completion $\mathbf{C}^r_{\xi'}$ of $\mathbf{C}_{\xi', \mathrm{colim}}^{r}$ with respect to the unique extension $\lambda_{r, \xi}'$ of the norm $\lambda_{r, \xi}$ on $\mathbf{C}_{\xi, \mathrm{colim}}^{r}$, and thus $\mathscr{R}'=\varinjlim_{r\rightarrow 0^+} \mathbf{C}^r_{\xi'}$.
 
 Moreover, write $ {\mathbf{N}'}_{\mathrm{dR}, [\sigma_1, \sigma_2]}^r =  \mathbf{N}_{\mathrm{dR}}^r ( \L|_{X'_{[\sigma_1, \sigma_2]}})$ and $ \mathbf{N}'_{\mathrm{dR}, [\sigma_1, \sigma_2]} =  \mathbf{N}_{\mathrm{dR}}( \L|_{X'_{[\sigma_1, \sigma_2]}})$. Then for $r$ sufficiently small, we have 
\begin{align*}
    \mathbf{N}^r_{\mathrm{dR},\xi} \otimes_{\mathbf{C}_{\xi}^{r}} {\mathscr R'}^r  & =   \left(\mathbf{N}_{\mathrm{dR}, \mathrm{colim}}^r \otimes_{\mathbf{C}_{\xi, \mathrm{colim}}^{r}} \mathbf{C}_{\xi}^{r}   \right)  \otimes_{\mathbf{C}_{\xi}^{r}} {\mathscr R'}^r \\ 
    & =  \left( \varinjlim_{\sigma_1<\sigma_{\xi}<\sigma_2}   \mathbf{N}_{\mathrm{dR}, [\sigma_1, \sigma_2]}^r \otimes_{\mathbf{C}_{[\sigma_1, \sigma_2]}^r}  {\mathbf{C}'}_{[\sigma_1, \sigma_2]}^{r} \right) \otimes_{\mathbf{C}_{\xi', \mathrm{colim}}^r} {\mathscr R'}^r 
    \\
    & = \left( \varinjlim_{\sigma_1<\sigma_{\xi}<\sigma_2} {\mathbf{N}'}_{\mathrm{dR}, [\sigma_1, \sigma_2]}^r \right) \otimes_{\mathbf{C}_{\xi', \mathrm{colim}}^r} {\mathscr R'}^r. 
\end{align*}
Note that, in the last identification, we make use of the compatibility of the $\mathbf{N}_{\mathrm{dR}}$ functor under finite \'etale extensions (cf. \cite[Proposition~10.13]{relative}). 
\end{construction}

\vspace{0.1in}
%%% 6.3
\subsection{Subrings of constants}\label{ss: subring of constants}
\noindent
\vspace{0.1in}

\noindent For later use, we need to compute the subrings of ${\mathbf{C}'}^{\dagger}_{[\sigma_1, \sigma_2]}$ and $\mathscr{R}'$ consisting of horizontal elements with respect to $\frac{d}{d[T^{\flat}]}$. Let $\kappa_1$ denote the maximal algebraic extension of $\kappa (\!(\cl \pi_K)\!)$ contained in $R'$. In particular, $\kappa_1$ is a finite separable extension of $\kappa (\!(\cl \pi_K)\!)$. 

\begin{lemma}  \label{lemma:relatively_algebraically_closed}
    There exists $\sigma_1, \sigma_2$ sufficiently close to $\sigma_\xi$, such that the inclusion $\kappa_1 \subset R'$ factors as 
    \[
    \kappa_1 \hooklongrightarrow R'_{[\sigma_1, \sigma_2]} \hooklongrightarrow R'.
    \]
\end{lemma}

\begin{proof} 
Recall that the colimit $R_{\xi, \mathrm{colim}}$ is precisely the uncompleted stalk $\mathbf{k} (\cl \xi)$ at the type I\!I\!I point $\cl \xi$ on the annulus $A_{\kappa (\!(\cl \pi_K)\!)} ([\sigma_1, \sigma_2]; T^\flat)$ corresponding to the circle of radius $\sigma_\xi$ (cf. Remark~\ref{remark:description_R_alpha}, Lemma~\ref{lemma: colimit of R[sigma1, sigma2]}).\footnote{Here, the stalk can be taken in either the Berkovich or the adic topology. See \S \ref{subsection: comparison_between_stalks}.} For simplicity, let $Y $ denote $A_{\kappa (\!(\cl \pi_K)\!)} ([\sigma_1, \sigma_2]; T^\flat)$ so that $\mO (Y) = R_{[\sigma_1, \sigma_2]}$. Let $h\colon  Y' \ra Y$ be the finite \'etale cover given by the finite \'etale map $R_{[\sigma_1, \sigma_2]} \ra R'_{[\sigma_1, \sigma_2]}$ and let $\cl \xi' \in Y'$ denote the preimage of $\cl \xi$. Since $h$ is finite \'etale and $h^{-1} (\cl \xi)$ consists of a singleton, the preimage of a system of open neighborhoods of $\cl \xi$ in $Y$ is cofinal among the neighborhoods of $\cl \xi'$ in $Y'$. In particular, the colimit $\varinjlim_{\sigma_1<\sigma_{\xi} <\sigma_2} R'_{[\sigma_1, \sigma_2]}$ is the uncompleted stalk $\mathbf{k}(\cl \xi')$ of the point $\cl \xi'$, and $R' = \sH (\cl \xi')$ is the completed stalk. 
By \cite[Theorem~2.3.3, Proposition~2.4.1]{Berkovich_etale}, we know that $\mathbf{k}(\cl \xi')$ is quasicomplete and is a henselian subfield of $\sH(\cl \xi')$. By \cite[Theorem~1.9.5]{Temkin_notes}, we conclude that $\mathbf{k}(\cl \xi')$ is relatively algebraically closed in $\sH(\cl \xi')$. This implies that $\kappa_1 \subset \mathbf{k}(\cl \xi')$ and thus proves the lemma. 
\end{proof}

\begin{remark} \label{remark:tensor_still_a_field_exercise}
    By construction, $\kappa_1$ is relatively algebraically closed in $R'$. Therefore, by elementary field theory, if $L/\kappa_1$ is any finite separable extension, then $R' \otimes_{\kappa_1} L$ is a field. 
\end{remark}

By the classical theory of Robba ring, the finite separable extension $\kappa_1/\kappa (\!(\cl \pi_K)\!)$ corresponds to a finite \'etale extension of Robba rings which we denote by $\mathbf{C}^{\dagger}_K\rightarrow \mathscr{R}(\kappa_1)$. In terms of the notation in \S \ref{sss:notation_cyclo}, $\mathscr{R}(\kappa_1)$ is abstractly isomorphic to $\mathscr{R}_{K''_0}$ for some finite unramified extension $K''_0$ of $K'_0$. Note that the differential operator $\frac{d}{d\pi_K}$ on $\mathbf{C}^{\dagger}_K$ uniquely extends to $\mathscr{R}(\kappa_1)$. Moreover, $\mathscr{R}(\kappa_1)$ can be naturally identified as a subring of ${\mathbf{C}'}^{\dagger}_{[\sigma_1, \sigma_2]}$ for $\sigma_1, \sigma_2$ sufficiently close to $\sigma_{\xi}$.

\begin{proposition}\label{prop:geometric_horizontal_section}
Let ``$\log[T^{\flat}]$'' be a free variable with $\frac{d}{d[T^\flat]} (\log[T^{\flat}]) = \frac{1}{[T^\flat]}$ and $\frac{d}{d\pi_K}(\log[T^{\flat}])=0$. Then, after shrinking $[\sigma_1, \sigma_2]$ if necessary, we have 
\[
 \left( {\mathbf{C}'}_{[\sigma_1, \sigma_2]}^{\dagger}\right)^{\frac{d}{d[T^\flat]} = 0 }=\left( {\mathbf{C}'}_{[\sigma_1, \sigma_2]}^{\dagger} [\log[T^{\flat}]]\right)^{\frac{d}{d[T^\flat]} = 0 } = \mathscr{R}(\kappa_1).\]
 identifying the differential operators $\frac{d}{d\pi_K}$ on both sides. As an immediate corollary, we have \[\left(\mathscr{R}'\right)^{\frac{d}{d[T^{\flat}]}=0}=\left(\mathscr{R}'[\log[T^{\flat}]]\right)^{\frac{d}{d[T^{\flat}]}=0}=\mathscr{R}(\kappa_1).\]
\end{proposition}

\begin{proof} Unwinding the constructions and applying Lemma~\ref{lemma:relatively_algebraically_closed}, we know that \[\mathscr{R}(\kappa_1) \subset \left( {\mathbf{C}'}_{[\sigma_1, \sigma_2]}^{\dagger} \right)^{\frac{d}{d[T^\flat]} = 0}\subset \left( {\mathbf{C}'}_{[\sigma_1, \sigma_2]}^{\dagger} [\log[T^{\flat}]]\right)^{\frac{d}{d[T^\flat]} = 0 }\] for $\sigma_1$, $\sigma_2$ sufficiently close to $\sigma_{\xi}$. It remains to show the reverse inclusions. It suffices to show that $\big( {\mathbf{C}'}_{[\sigma_1, \sigma_2]}^{\dagger} \big)^{\frac{d}{d[T^\flat]} = 0 } \subset \mathscr{R}(\kappa_1)$. Indeed, using the explicit description of $ {\mathbf{C}}_{[\sigma_1, \sigma_2]}^{[s, r]}$ in Lemma~\ref{lemma: explicit description C--case of annulus}, we have
 \[
 \left( {\mathbf{C}}_{[\sigma_1, \sigma_2]}^{[s, r]}\right)^{\frac{d}{d[T^\flat]} = 0 } = \mathbf{C}_K^{[s, r]}\]
for $r$ sufficiently small and for all $0<s\le r$. This implies
$ \big( {\mathbf{C}}_{[\sigma_1, \sigma_2]}^{\dagger} \big)^{\frac{d}{d[T^\flat]} = 0 } = \mathbf{C}_K^{\dagger}$, and hence 
\[
 \left( {\mathbf{C}}_{[\sigma_1, \sigma_2]}^{\dagger}  \otimes_{\mathbf{C}^{\dagger}_K} \mathscr{R}(\kappa_1) \right)^{\frac{d}{d[T^\flat]} = 0 } = \mathscr{R}(\kappa_1).\]
 
Consider the finite \'etale map 
 \[g\colon  {\mathbf{C}}_{[\sigma_1, \sigma_2]}^{\dagger} \otimes_{\mathbf{C}^{\dagger}_K} \mathscr{R}(\kappa_1) \ra {\mathbf{C}'}_{[\sigma_1, \sigma_2]}^{\dagger}. 
 \] 
The map $g$ realizes the target as a finite projective module over the source, say of rank $m$. Take an element $\alpha \in {\mathbf{C}'}_{[\sigma_1, \sigma_2]}^{\dagger}$ such that $\frac{d}{d[T^\flat]} (\alpha ) = 0$. Let $T_\alpha\colon  {\mathbf{C}'}_{[\sigma_1, \sigma_2]}^{\dagger} \ra {\mathbf{C}'}_{[\sigma_1, \sigma_2]}^{\dagger}$ denote the linear map given by multiplication by $\alpha$ and let $P_\alpha (x)  \in  ({\mathbf{C}}_{[\sigma_1, \sigma_2]}^{\dagger} \otimes_{\mathbf{C}^{\dagger}_K} \mathscr{R}(\kappa_1))[x]$ denote the characteristic polynomial. Since $\frac{d}{d[T^\flat]} (\mathrm{tr} (T_\alpha^i )) = \mathrm{tr} ( \frac{d}{d[T^\flat]} (T_\alpha^i )) =0$ for all $i\ge 1$, we know that all coefficients of $P_\alpha (x)$ are $\frac{d}{d[T^{\flat}]}$-horizontal sections, and thus $P_{\alpha} (x) \in \mathscr{R}(\kappa_1)[x]$.

Next, we construct a finite separable extension $L/\kappa_1$ (with associated extension of Robba rings $\mathscr{R}(L)/\mathscr{R}(\kappa_1)$) together with a ring map $\eta \colon \mathbf{C}^{\dagger}_{[\sigma_1, \sigma_2]}\otimes_{\mathbf{C}^{\dagger}_K} \mathscr{R}(\kappa_1)\rightarrow \mathscr{R}(L)$ such that the base change of $g$ along $\eta$ splits into a direct product $\prod_{i=1}^m \mathscr{R}(L)$. For this, we first choose any finite separable extension $L/\kappa_1$ such that $\kappa_1\rightarrow L$ factors through natural inclusion $\kappa_1\hookrightarrow R'_{[\sigma_1, \sigma_2]}$. This induces a finite \'etale map $\delta\colon \mathscr{R}(\kappa_1)\rightarrow \mathscr{R}(L)$ which factors as $\mathscr{R}(\kappa_1)\hookrightarrow {\mathbf{C}'}^{\dagger}_{[\sigma_1, \sigma_2]}\xrightarrow[]{\mu}\mathscr{R}(L)$. Now, the composition $\mathbf{C}^{\dagger}_{[\sigma_1, \sigma_2]}\rightarrow {\mathbf{C}'}^{\dagger}_{[\sigma_1, \sigma_2]}\xrightarrow[]{\mu}\mathscr{R}(L)$ and the map $\delta$ induce a ring map
$\eta \colon \mathbf{C}^{\dagger}_{[\sigma_1, \sigma_2]}\otimes_{\mathbf{C}^{\dagger}_K} \mathscr{R}(\kappa_1)\rightarrow \mathscr{R}(L)$. The base change of $g$ along $\eta$ is a finite \'etale cover of $\mathscr{R}(L)$. Replacing $L$ by a further separable extension, we can make sure that this finite \'etale cover splits into direct products $\prod_{i=1}^m\mathscr{R}(L)$. To sum up, we obtain a commutative diagram

\[ 
\begin{tikzcd}
         {\mathbf{C}'}^{\dagger}_{[\sigma_1, \sigma_2]}\arrow[r] & {\mathbf{C}'}^{\dagger}_{[\sigma_1, \sigma_2]}\otimes_{\mathscr{R}(\kappa_1)}\mathscr{R}(L) \arrow[r] & \prod_{i=1}^m\mathscr{R}(L) \\
         \mathbf{C}^{\dagger}_{[\sigma_1, \sigma_2]} \otimes_{\mathbf{C}^{\dagger}_K}\mathscr{R}(\kappa_1)\arrow[u, "g"] \arrow[r]\arrow[rr, "\eta"' {xshift=0.3in}, bend right=13, crossing over] & \mathbf{C}^{\dagger}_{[\sigma_1, \sigma_2]} \otimes_{\mathbf{C}^{\dagger}_K}\mathscr{R}(L) \arrow[u, "g'"]\arrow[r] & \mathscr{R}(L)\arrow[u]\\
         \mathscr{R}(\kappa_1)\arrow[u]\arrow[r, "\delta"] & \mathscr{R}(L)\arrow[u]
\end{tikzcd} 
\]
where each square is a push-out square.

Let $(\alpha_1,\ldots, \alpha_m)\in \prod_{i=1}^m\mathscr{R}(L)$ be the image of $\alpha$ under the composition of the top row of the diagram. Also consider term ${\mathbf{C}'}^{\dagger}_{[\sigma_1, \sigma_2]}\otimes_{\mathscr{R}(\kappa_1)}\mathscr{R}(L)$ in the middle column. Since taking base change commutes with taking characteristic polynomials, the multiplication-by-$\alpha$ map on ${\mathbf{C}'}\otimes_{\mathscr{R}(\kappa_1)}\mathscr{R}(L)$ and the multiplication-by-$(\alpha_1,\ldots, \alpha_m)$ map on $\prod_{i=1}^m\mathscr{R}(L)$ both have characteristic polynomial $P_{\alpha}(x)\in \mathscr{R}(\kappa_1)[x]$. In particular, we have
\[
P_{\alpha}(x) = (x-\alpha_1)\cdots (x-\alpha_m) \in \mathscr{R}(L)[x]    
\]
By Cayley--Hamilton, we have an identity $(\alpha-\alpha_1)\cdots (\alpha-\alpha_m) = 0$ on ${\fC'}_{[\sigma_1, \sigma_2]}^{\dagger}\otimes_{\mathscr{R}(\kappa_1)} \mathscr{R}(L)$.

Now, consider $R'_L\coloneqq R'\otimes_{\kappa_1} L$, which is a field by Remark~\ref{remark:tensor_still_a_field_exercise}. Let $\mathscr{R}'\rightarrow\mathscr{R}'_L$ be the extension of fake Robba rings associated with the finite separable extension $R'_L/R'$. Then there is a natural embedding ${\fC'}_{[\sigma_1, \sigma_2]}^{\dagger}\otimes_{\mathscr{R}(\kappa_1)} \mathscr{R}(L) \rightarrow \mathscr{R}'_L$, and thus ${\fC'}_{[\sigma_1, \sigma_2]}^{\dagger}\otimes_{\mathscr{R}(\kappa_1)} \mathscr{R}(L)$ is an integral domain. Consequently, $\alpha = \alpha_i$ as an element in ${\fC'}_{[\sigma_1, \sigma_2]}^{\dagger}\otimes_{\mathscr{R}(\kappa_1)} \mathscr{R}(L)$, for some $i$. By the injectivity of $\mathscr{R}(\kappa_1) \rightarrow {\fC'}_{[\sigma_1, \sigma_2]}^{\dagger}$ and the faithful flatness of $\mathscr{R}(\kappa_1) \rightarrow \mathscr{R}(L)$, we conclude that
\[\alpha\otimes 1\in \big(\{1\}\otimes \mathscr{R}(L)\big)\cap \left({\mathbf{C}'}^{\dagger}_{[\sigma_1, \sigma_2]}\otimes\{1\}\right)=\mathscr{R}(\kappa_1);\]
namely, 
$\alpha \in \mathscr{R}(\kappa_1)$. This finishes the proof of the proposition.
\end{proof}

From now on, we shrink $[\sigma_1, \sigma_2]$ so that Lemma~\ref{lemma:relatively_algebraically_closed} and Proposition~\ref{prop:geometric_horizontal_section} apply.

\vspace{0.1in}
%%% 6.4
\subsection{Unipotence implies semistability} 
\label{ss:contagion_of_goem_unipotence} 
\noindent 
\vspace{0.1in}
 
\noindent Recall from Construction~\ref{construction:geometric_realization} that we have realized the finite extension $\mathbf{C}^{\dagger}_{\xi}\rightarrow \mathscr R'$ of fake Robba rings via a finite \'etale cover $f\colon  X_{[\sigma_1, \sigma_2]}' \ra X_{[\sigma_1, \sigma_2]}$. The finite projective module $\mathbf{N}'_{\mathrm{dR}, [\sigma_1, \sigma_2]} = \mathbf{N}_{\mathrm{dR}} (\L|_{X'_{[\sigma_1, \sigma_2]}})$ over ${\mathbf{C}'}_{[\sigma_1, \sigma_2]}^{\dagger}$ is equipped with an arithmetic differential operator $D_{\pi_K}$ and a geometric differential operator $D_{[T^{\flat}]}$, as well as a Frobenius structure. We also consider the full connection $\nabla = D_{\pi_K}  d \pi_K + D_{[T^\flat]} d [T^\flat]$ so that $(\mathbf{N}'_{\mathrm{dR}, [\sigma_1, \sigma_2]}, \nabla)$ is a finite projective $(\varphi, \nabla)$-module over ${\mathbf{C}'}_{[\sigma_1, \sigma_2]}^{\dagger}$. Similarly, for $r$ sufficiently small, we have a finite projective $(\varphi, \nabla)$-module $ {\mathbf{N}'}_{\mathrm{dR}, [\sigma_1, \sigma_2]}^r = \mathbf{N}_{\mathrm{dR}}^r (\L|_{X'_{[\sigma_1, \sigma_2]}})$ over ${\mathbf{C}'}_{[\sigma_1, \sigma_2]}^{r}$, equipped with differential operators $D_{\pi_K}$ and $D_{[T^{\flat}]}$. Moreover, we know that
\[
\mathbf{N}'_{\mathrm{dR}, \xi}= \mathbf{N}'_{\mathrm{dR}, [\sigma_1, \sigma_2]}\otimes_{{\mathbf{C}'}^{\dagger}_{[\sigma_1, \sigma_2]}}\mathscr{R}'
\]
is unipotent with respect to $\nabla$. A key ingredient for the proof of Theorem~\ref{thm: local main thm type III} is the following proposition, which states that the unipotence of $(\mathbf{N}'_{\mathrm{dR}, \xi}, \nabla)$ can be spread out to $(\mathbf{N}'_{\mathrm{dR}, [\sigma_1, \sigma_2]}, \nabla)$ for some sufficiently small interval $[\sigma_1, \sigma_2]$.

\begin{proposition} \label{prop:geometric_unipotence_over_a_region}
After shrinking $[\sigma_1, \sigma_2]$ if necessary, the differential module $\mathbf{N}'_{\mathrm{dR}, [\sigma_1, \sigma_2]}$ is strongly unipotent for $\nabla$; that is, $\mathbf{N}'_{\mathrm{dR}, [\sigma_1, \sigma_2]}$ is free over ${\mathbf{C}'}_{[\sigma_1, \sigma_2]}^{\dagger}$ and admits a basis $w_1, \ldots, w_d$ such that the matrices of $[T^{\flat}]D_{[T^{\flat}]}$ and $\pi_KD_{\pi_K}$ relative to this basis are both upper-triangular nilpotent with entries in $K''_0$. In particular, $\mathbf{N}'_{\mathrm{dR}, [\sigma_1, \sigma_2]}$ is unipotent for the full connection $\nabla$.
\end{proposition}

The proof of Proposition~\ref{prop:geometric_unipotence_over_a_region} will occupy \S \ref{ss:geometric_uni_special}--\ref{ss:geometric_uni_general}. Granting Proposition~\ref{prop:geometric_unipotence_over_a_region}, let us finish the proof of Theorem~\ref{thm: local main thm type III}. 

\begin{proof}[Proof of Theorem~\ref{thm: local main thm type III} assuming Proposition~\ref{prop:geometric_unipotence_over_a_region}]
\noindent 

\noindent 
The first step is to compute the geometric horizontal sections after adjoining ``$\log [T^\flat]$''. Let $\sigma_1, \sigma_2$ be sufficiently close to $\sigma_{\xi}$ as in Proposition~\ref{prop:geometric_unipotence_over_a_region}. Consider \[ 
N' \coloneqq \left(\mathbf{N}'_{\mathrm{dR}, [\sigma_1, \sigma_2]} [\log [T^\flat]] \right)^{D_{[T^\flat]} = 0}.
\]
By Proposition~\ref{prop:geometric_horizontal_section} and Proposition~\ref{prop:geometric_unipotence_over_a_region}, and noting that $D_{[T^\flat]}$ commutes with $D_{\pi_K}$ and is compatible with Frobenius, we deduce that $N'$ is a finite free $\mathscr{R}(\kappa_1)$-module of rank $d$ (where $d$ is the rank of $\L$) and $(N', D_{\pi_K})$ is a differential module over $(\mathscr{R}(\kappa_1), \frac{d}{d\pi_K})$ equipped with Frobenius structure. By the classical $p$-adic local monodromy theorem of Andr\'e--Kedlaya--Mebkhout (see \cite{Andre_p-adic_monodormy, Kedlaya_p-adic_monodormy, Mebkhout_p-adic_monodormy}), there exists a finite separable extension $\kappa_2/\kappa_1$ (with the associated Robba ring extension denoted by $\mathscr{R}(\kappa_2)/\mathscr{R}(\kappa_1)$), such that $N'\otimes_{\mathscr{R}(\kappa_1)}\mathscr{R}(\kappa_2)$ becomes unipotent with respect to $D_{\pi_K}$. By replacing $K$ by finite extension $K'$ so that $\mathbf{C}^{\dagger}_{K'}$ contains $\mathscr{R}(\kappa_2)$, we may assume that $N'$ is already unipotent with respect to $D_{\pi_K}$. Hence, the $\mathscr{R}(\kappa_1)$-module $N'[\log \pi_K]$ admits a horizontal basis $w_1, \ldots, w_d$ with respect to $D_{\pi_K}$. Note that $w_1, \ldots, w_d$ are also geometrically horizontal as $\frac{d}{d[T^\flat]} (\log \pi_K) = 0$. We arrive at a decomposition 
\[
\mathbf{N}'_{\mathrm{dR}, [\sigma_1, \sigma_2]}[\log [T^\flat], \log \pi_K]  \cong \bigoplus_{i = 1}^{d}  {\mathbf{C}'}_{[\sigma_1, \sigma_2]}^{\dagger} [\log [T^\flat], \log \pi_K] \cdot w_i
\]
with each $w_i$ satisfying $D_{[T^\flat]} (w_i) = D_{\pi_K} (w_i) = 0$. Note that all these $w_i$'s are fixed by an open subgroup of $\Gamma$. By enlarging $K$ and passing to a $p$-power Kummer cover of $X'_{[\sigma_1, \sigma_2]}$ (cf. Remark~\ref{remark: passing to p-power Kummer cover}), we may further assume that all $w_i$'s are fixed by $\Gamma$. 
This allows us to apply Proposition~\ref{prop:unipotence_and_st} (note that the functors $\mathbf{N}_{\mathrm{dR}}$ and $D^{\dagger}_{\mathrm{rig}}$ coincide after inverting $t$) and conclude the proof. 
\end{proof}

\vspace{0.1in}
%%% 6.5
\subsection{Unipotence over annuli: special case} \label{ss:geometric_uni_special}
\noindent 
\vspace{0.1in}
 
\noindent In the remainder of this section, we prove Proposition~\ref{prop:geometric_unipotence_over_a_region}. Thanks to Theorem~\ref{thm:monodromy_fake_annuli}, we already achieve unipotence of $\nabla$ ``at $\xi$'' by passing to a finite \'etale cover. We want to propagate this unipotence to an open neighborhood of $\xi$. This is akin to the process of extending the region of convergence for solutions of $p$-adic differential equations over a (fake) poly-annulus.

We first prove Proposition~\ref{prop:geometric_unipotence_over_a_region} in the special case where $\mathscr{R}'=\mathbf{C}^{\dagger}_{\xi}$. In this case, $\mathbf{N}^r_{\mathrm{dR}, \xi}$ is already unipotent with respect to $\nabla$, for $r$ sufficiently small. Our goal is to show that
\[\mathbf{N}_{\mathrm{dR}, \mathrm{colim}} = \varinjlim_{r \ra 0^+} \mathbf{N}_{\mathrm{dR}, \mathrm{colim}}^r =\varinjlim_{r\rightarrow 0^+}\varinjlim_{\sigma_1<\sigma_{\xi<\sigma_2}}\mathbf{N}^r_{\mathrm{dR}, [\sigma_1, \sigma_2]}\] 
is strongly unipotent. Since 
\[\mathbf{N}_{\mathrm{dR}, \xi} = \mathbf{N}_{\mathrm{dR}, \mathrm{colim}} \otimes_{\mathbf{C}^\dagger_{\xi, \mathrm{colim}}} \mathbf{C}^\dagger_{\xi},\]
this can be viewed as a ``decompletion question'' from $\mathbf{C}^{\dagger}_{\xi}$ to $\mathbf{C}^{\dagger}_{\xi, \mathrm{colim}}$.

\vspace{0.1in}
%%% 6.5.1
\subsubsection{Intrinsic generic radii over log RP subsets} 
\noindent
\vspace{0.1in}

\noindent
We first recall the notion of generic radius of convergence, originally developed by Dwork, Robba, Kedlaya, and others. In the classical setting, suppose $(V, D)$ is a finite nonzero differential module over $(F_\rho, \frac{d}{dt})$ where $\rho > 0$ and $F_\rho$ denotes the completion of $K(t)$ under the $\rho$-Gauss norm (cf. \cite[Definition~9.4.1]{kedlaya_book}). Consider the spectral radius of the operator $D$ defined by $|D|_{\mathrm{sp}, V} = \lim_{s \ra +\infty} |D^s|_V^{1/s}$, where $|D^s|_V$ denotes the operator norm of $D^s$. Then the \emph{generic radius of convergence} of $(V, D)$ is defined to be $R(V) = \omega \cdot |D|_{\mathrm{sp}, V}^{-1}$, where $\omega = p^{-1/(p-1)}$. The \emph{intrinsic generic radius of convergence} is defined to be \[I\!R(V) = \rho^{-1} R(V)=\frac{|\frac{d}{dt}|_{\mathrm{sp}, F_{\rho}}}{|D|_{\mathrm{sp}, V}},\] which is a real number in $(0, 1]$ (cf. \cite[Definition~9.4.4, 9.4.7, 6.1.2, 6.1.3]{kedlaya_book}). 

To treat the type I\!I\!I point of interest, we will need generalizations of these notions from \cite{kedlaya_st3}. For the setup, let us describe $\mathbf{C}_{[\sigma_1, \sigma_2]}^{[s, r]}$ (resp.  $\mathbf{C}_{[\sigma_1, \sigma_2]}^{r}$) in terms of \emph{log-RP subset} (resp.  \emph{ind-log-RP subset}) of $(\R_{>0})^2$ in the sense of \cite[Definition~4.1.2]{kedlaya_st3}. For an interval $I = [\sigma_1, \sigma_2]$, consider the subset
\begin{align*} 
S_{[\sigma_1, \sigma_2]}^{[s, r]}  \coloneqq \: & \left\{(p^{- c_K \cdot u}, \sigma^{u}) \: ~\vline~  \: u \in [s, r],\, \sigma \in [\sigma_1, \sigma_2] \right \}   \\
 = \: & \left\{(x, x^\beta) \: ~\vline~  \: x \in [p^{-c_K \cdot r}, p^{- c_K \cdot s} ] , \:\beta \in [- \frac{1}{c_K}\log_p\sigma_2, - \frac{1}{c_K}\log_p \sigma_1]  \right\} \subset (\R_{>0})^2. 
\end{align*}
Here $c_K = \frac{p}{(p-1)} \frac{1}{e_K}$ as in Notation~\ref{sss:notation_cyclo}. The subset $S_{[\sigma_1, \sigma_2]}^{(0, r]}$ is defined in a similar way.

\begin{remark} 
In the language of \cite{kedlaya_st3}, $S_{[\sigma_1, \sigma_2]}^{[s, r]}$ (resp. $S_{[\sigma_1, \sigma_2]}^{(0, r]}$) is referred to as a log-RP subset (resp. ind-log-RP subset) because its image under the logarithmic map is a rational polygon (resp. limit of rational polygons); namely,
 \[
 \log_p (S_{[\sigma_1, \sigma_2]}^{[s, r]}) = \left\{(x, \beta x) \: \vline  \: x \in [ {-c_K \cdot r},  -c_K \cdot s] , \beta \in [- \frac{1}{c_K}\log_p\sigma_2, - \frac{1}{c_K}\log_p \sigma_1]  \right\}\subset \R^2. 
 \]
See Figure \ref{fig:regions} for an illustration of the log-RP subset. (In the figure, $\beta_- = - \frac{1}{c_K}\log_p \sigma_2$ and $\beta_+ = - \frac{1}{c_K}\log_p \sigma_1$). 
\begin{figure} 
    \centering
  \scalebox{0.7}{
\begin{tikzpicture}[x=5cm,y=5cm,font=\large,
  every node/.style={text=sectoraxis,inner sep=2pt}]
  \draw[sectoraxis,line width=.65pt,-{Stealth[length=5pt,width=5pt]}] (-0.2,0) -- (1.10,0)
    node[below] { \qquad \quad $x = |\pi_K|$};
  \draw[sectoraxis,line width=.65pt,-{Stealth[length=5pt,width=5pt]}] (0,-0.15) -- (0,1.10)
    node[left] {$y = |[T^\flat]]|$};
 \path[pattern=north east lines,pattern color=sectorink]    plot[variable=\t,domain=.5:.833333333,samples=160]     (\t,{pow(\t,.12)})--plot[variable=\t,domain=.833333333:.5,samples=160]      (\t,{pow(\t,.45)}) -- cycle;

% Greenboundary 
  \foreach \beta in {.12,.45} {
    % Parameterizing the left extension by y keeps its steep end smooth.
    \draw[sectorboundary,line width=.75pt,dash pattern=on 3pt off 3pt]      plot[variable=\t,domain=0:{pow(.5,\beta)},samples=200]        ({pow(\t,1/\beta)},\t);
    \draw[sectorboundary,line width=1.25pt]      plot[variable=\t,domain=.5:.833333333,samples=160]        (\t,{pow(\t,\beta)});
    \draw[sectorboundary,line width=.75pt,dash pattern=on 3pt off 3pt]      plot[variable=\t,domain=.833333333:1,samples=100]       (\t,{pow(\t,\beta)}); }

  \foreach \xx/\period in {.5/r,.833333333/s} {
    \draw[sectorink,line width=.85pt]
      (\xx,{pow(\xx,.45)}) -- (\xx,{pow(\xx,.12)});
    \foreach \beta in {.12,.45} {
      \fill[sectorink] (\xx,{pow(\xx,\beta)}) circle[radius=1.4pt];
    }
    \draw[sectorguide,line width=.45pt,dash pattern=on 1pt off 3pt] (\xx,0) -- (\xx,{pow(\xx,.45)});
    \draw[sectoraxis,line width=.6pt] (\xx,-.006) -- (\xx,.006);
    \node[below=3pt] at (\xx,0) {$p^{-c_K\period}$};
  }

  \draw[sectorred,line width=.95pt,line cap=round,dash pattern=on .1pt off 3.6pt]
    plot[variable=\t,domain=0:{pow(.5,.316227766)},samples=200]
      ({pow(\t,1/.316227766)},\t);
  \draw[sectorred,line width=1.2pt]
    plot[variable=\t,domain=.5:.833333333,samples=160]
      (\t,{pow(\t,.316227766)});
  \draw[sectorred,line width=.95pt,line cap=round,dash pattern=on .1pt off 3.6pt]
    plot[variable=\t,domain=.833333333:1,samples=100]
      (\t,{pow(\t,.316227766)});

  \draw[sectorink,fill=white,line width=.8pt]
    (1,1) circle[radius=1.8pt];
  \draw[sectoraxis,line width=.6pt] (1,-.006) -- (1,.006);
  \draw[sectoraxis,line width=.6pt] (-.006,1) -- (.006,1);
  \node[below=5pt] at (-0.07,0) {$0$};
  \node[below=5pt] at (1,0) {$1$};
  \node[left=3pt] at (0,1) {$1$};
%  \node[anchor=south east] at (.985,1.035) {$(1,1)$};
  \node[sectorboundary,anchor=west] at (.58,1.02) {$y=x^{\beta_-}$};
  \node[sectorboundary,anchor=west] at (.58,.75) {$y=x^{\beta_+}$};
  \node[sectorred,anchor=west] at (.33,.752) {$\sigma_\xi$};
\end{tikzpicture} \qquad \qquad 
\raisebox{0mm}{
\begin{tikzpicture}[x=4cm,y=6.5cm,font=\large,
  every node/.style={text=sectoraxis,inner sep=2pt}]
  \draw[sectoraxis,line width=.65pt,-{Stealth[length=5pt,width=5pt]}] (-1.55,0) -- (.25,0)
    node[right=3pt] {$\log_p x$};
  \draw[sectoraxis,line width=.65pt,-{Stealth[length=5pt,width=5pt]}] (0,-.8) -- (0,.17)
    node[right=4pt] {$\log_p y$};

  \path[pattern=north east lines,pattern color=sectorink]
    (-1,{.45*-1})
    -- (-1,{.12*-1})
    -- (-.263034406,{.12*-.263034406})
    -- (-.263034406,{.45*-.263034406}) -- cycle;

  \foreach \beta in {.12,.45} {
    \draw[sectorboundary,line width=.75pt,dash pattern=on 3pt off 3pt]
      (-1.55,{\beta*-1.55}) -- (-1,{\beta*-1});
    \draw[sectorboundary,line width=1.25pt]
      (-1,{\beta*-1}) -- (-.263034406,{\beta*-.263034406});
    \draw[sectorboundary,line width=.75pt,dash pattern=on 3pt off 3pt] (-.263034406,{\beta*-.263034406}) -- (0,0);
  }

  \foreach \uu/\period in {-1/r,-.263034406/s} {
    \draw[sectorink,line width=.85pt]
      (\uu,{.45*\uu}) -- (\uu,{.12*\uu});
    \foreach \beta in {.12,.45} {
      \fill[sectorink] (\uu,{\beta*\uu}) circle[radius=1.4pt];
    }
    \draw[sectorguide,line width=.45pt,dash pattern=on 1pt off 3pt] (\uu,{.12*\uu}) -- (\uu,0);
    \draw[sectoraxis,line width=.6pt] (\uu,-.008) -- (\uu,.008);
    \node[above=3pt] at (\uu,0) {$-c_K\period$};
  }

  \draw[sectorred,line width=.95pt,line cap=round,dash pattern=on .1pt off 3.6pt]
    (-1.55,{.316227766*-1.55}) -- (-1,{.316227766*-1});
  \draw[sectorred,line width=1.2pt]
    (-1,{.316227766*-1}) -- (-.263034406,{.316227766*-.263034406});
  \draw[sectorred,line width=.95pt,line cap=round,dash pattern=on .1pt off 3.6pt] (-.263034406,{.316227766*-.263034406}) -- (0,0);

  \draw[sectorink,fill=white,line width=.8pt]
    (0,0) circle[radius=1.8pt];
  \node[above right=2pt] at (0,0) {$0$};
  \node[sectorboundary,anchor=west] at (-1.31,-.06) {$\log_p y =\beta_- \log_p x $};
  \node[sectorboundary,anchor=west] at (-.7,-.35) {$\log_p y =\beta_+ \log_p x$};
  \node[sectorred,anchor=west] at (-1.38,-.35) {$\sigma_\xi$};
\end{tikzpicture}
}% raisebox
}%scalebox
    \caption{A cartoon illustration of $S_{[\sigma_1, \sigma_2]}^{[s, r]}$ (left) and $\log_p S_{[\sigma_1, \sigma_2]}^{[s, r]}$ (right) }
    \label{fig:regions}
\end{figure}
\end{remark}

We adopt the notation from \cite[Notation~4.1.3]{kedlaya_st3}. Let $\mathbb{A}^{2,\mathrm{an}}_{K'_0}$ denote the analytic affine plane with coordinates $\pi_K$ and $[T^{\flat}]$. For any (ind-)log-RP subset $S\subset \mathbb{R}_{>0}^2$, let $A_{K'_0}(S)$ denote the open subset of $\mathbb{A}^{2,\mathrm{an}}_{K'_0}$ defined by the condition $(|\pi_K|, |[T^{\flat}]|)\in S$. Then we have
\[\mathbf{C}_{[\sigma_1, \sigma_2]}^{[s, r]} = \Gamma(A_{K_0'} (S_{[\sigma_1, \sigma_2]}^{[s, r]}), \mO),\qquad \mathbf{C}_{[\sigma_1, \sigma_2]}^r = \Gamma(A_{K_0'} (S_{[\sigma_1, \sigma_2]}^{(0, r]}), \mO)\]
(note that $\lambda_u (\pi_K) = p^{- c_K u}$). For any $\sigma \in [\sigma_1, \sigma_2]$ and $s \in (0, r]$, consider the point
\[
R(\sigma, s) = (p^{-c_K \cdot s}, \sigma^{s}) = (p^{-c_K \cdot s}, (p^{-c_K \cdot s})^{\beta})\in S_{[\sigma_1, \sigma_2]}^{(0,r]}
\]
where $\beta = -\frac{1}{c_K}\cdot\log_p \sigma$. This point gives rise to a norm $\lambda_{s, \sigma}$ on $\mathbf{C}_{[\sigma_1, \sigma_2]}^{(0, r]}$ given by\footnote{This norm is denoted as $|\cdot|_{R(\sigma, s)}$ in \cite[Definition~4.1.6]{kedlaya_st3}.}
\[\lambda_{s,\sigma}\left(\sum_{i\in \mathbb{Z}}c_i[T^{\flat}]^i\right)=\sup_{i\in \mathbb{Z}}\left\{\lambda_{s}(c_i)\cdot \sigma^{si}\right\}.\]

\begin{definition}[{\cite[Definition~4.2.1]{kedlaya_st3}}]
    Let $(M, D_{\pi_K}, D_{[T^\flat]})$ be a finite projective  differential module over $(\mathbf{C}_{[\sigma_1, \sigma_2]}^{(0, r]}, \frac{d}{d\pi_K}, \frac{d}{d[T^\flat]})$ and let $R(\sigma, s) \in S_{[\sigma_1, \sigma_2]}^{(0, r]}$ be a point as above. Let $F_{\sigma, s}$ denote the completion of $\mathrm{Frac}(\mathbf{C}_{[\sigma_1, \sigma_2]}^{(0, r]})$ under $\lambda_{\sigma, s}$, and consider the base change
    $V_{\sigma, s} \coloneqq M \otimes_{\mathbf{C}_{[\sigma_1, \sigma_2]}^{(0, r]}} F_{\sigma, s}$. Define the \emph{intrinsic generic radius} of $M$ at $R(\sigma, s)$ to be 
    \[
    I\!R(M, R(\sigma, s)) \coloneqq \min \left\{1, \,\,\frac{|\frac{d}{d[T^\flat]}|_{\mathrm{sp}, F_{\sigma, s}}}{|D_{[T^\flat]}|_{\mathrm{sp}, V_{\sigma, s}}}, \,\,\frac{|\frac{d}{d\pi_K}|_{\mathrm{sp}, F_{\sigma, s}}}{|D_{\pi_K}|_{\mathrm{sp}, V_{\sigma, s}}} \right\}.
    \]
\end{definition}

By assumption, we know that $\mathbf{N}_{\mathrm{dR}, \xi}^r$ is unipotent for some sufficiently small $r$; that is, it admits a filtration whose graded pieces are trivial differential modules. In particular, for $\sigma_1, \sigma_2$ sufficiently closed to $\sigma_{\xi}$, there exists $r_0>0$ such that
\[
I\!R(\mathbf{N}_{\mathrm{dR}, [\sigma_1, \sigma_2]}^{r_0}, R(\sigma_\xi, s)) = 1
\]
for all $s \in (0, r_0]$ (cf. \cite[Remark~4.2.4]{kedlaya_st3} and \cite[Lemma~1.2.9]{Kedlaya-Xiao}). Note that $\alpha = - \frac{1}{c_K}\log_p \sigma_\xi$ is irrational by the assumption. Therefore, according to \cite[Proposition~6.2.1]{kedlaya_st3}, by decreasing $r_0$ and shrinking the interval $[\sigma_1, \sigma_2]$ if necessary, we may assume that 
\begin{equation}\label{eq: Robba condition NdR}
I\!R(\mathbf{N}_{\mathrm{dR}, [\sigma_1, \sigma_2]}^{r_0}, R(\sigma, s)) = 1
\end{equation}
for all $\sigma\in [\sigma_1, \sigma_2]$ and $s\in (0,r_0]$.

\vspace{0.1in}
%%% 6.5.2
\subsubsection{A monomial change of coordinates} 
\noindent
\vspace{0.1in}

\noindent
From the above discussion, we may assume that $\mathbf{N}_{\mathrm{dR}, [\sigma_1, \sigma_2]}^{r_0}$ has intrinsic generic radius equal to $1$ at all points $R(\sigma, s) \in S_{[\sigma_1, \sigma_2]}^{(0, r_0]}$. To deduce unipotence from this condition, we would like to put ourselves in a setting where we can apply the multivariable $p$-adic Fuchs theorem\footnote{The single variable $p$-adic Fuchs theorem is due to Christol and Mebkhout (see \cite[Theorem 13.6.1]{kedlaya_book}).} of Wang \cite[Corollary~3.20]{Wang} (which slightly generalizes the work of Gachet \cite{Gachet}). However, Wang's theorem only applies to honest poly-annuli rather than the fake ones. The key is to use a change of coordinates to transfer the ind-log-RP subsets into a ``standard rectangle''.

Let $A = \left( \begin{smallmatrix} a & b \\ c & d \end{smallmatrix} \right)$ be a matrix with $a, b, c, d \in \Z$ and $\det A = \pm 1$. Consider the map $\mu_A\colon  \R_{>0}^2 \ra \R_{>0}^2$ sending $(x, y) \mapsto (x^a y^b, x^c y^d)$. We refer to $\mu_A$ as a \emph{monomial change of coordinates}.\footnote{This is called a \emph{toric coordinate change} in \cite[Definition~4.1.4]{kedlaya_st3}.} In fact, if we take logarithms on both the source and target, the resulting map
\[
\log \mu_A\colon  \log_p(\R^2_{>0})=\R^2 \ra \log_p(\R^2_{>0})=\R^2
\]
is simply the affine transformation given by the multiplication of the matrix $A$. Now, for any $\delta>0$, let $S^{\delta}$ denote the ind-log-RP subset $[p^{-\delta}, 1) \times [p^{-\delta}, 1) \subset (\R_{>0})^2$. Note that $\log_pS^{\delta}=[-\delta, 0)\times[-\delta, 0)$. For a suitable choice of $A$ and shrinking $[\sigma_1, \sigma_2]$ if necessary, we can ensure that the image of $S^{\delta}$ under $\mu_A$ lies in $S^{(0,r]}_{[\sigma_1, \sigma_2]}$ for sufficiently small $\delta$. Indeed, we first choose $A$ such that the image of the third quadrant (i.e., $(-\infty,0)\times(-\infty,0)$) under $\log\mu_A$ lies entirely inside the cone
\[\big\{(x,\beta x)\,\big|\, x<0,\,\, \beta\in [-\frac{\log_p\sigma_2}{c_K}, -\frac{\log_p\sigma_1}{c_K}]\big\}\]
and contains the ray 
\[\big\{(x,\beta x)\,\big|\, x<0,\,\, \beta=-\frac{\log_p\sigma_{\xi}}{c_K}\big\}.\]
Shrinking the interval $[\sigma_1, \sigma_2]$ and multiplying $A$ by $\left( \begin{smallmatrix} 0 & 1 \\ 1 & 0 \end{smallmatrix} \right)$ if necessary, we may assume that the images of the rays $(-\infty,0)\times\{0\}$ and $\{0\}\times (-\infty, 0)$ under the linear transformation $\log\mu_A$ are precisely the rays $\{(x,\beta x)\,|\,x<0,\,\beta=- \frac{1}{c_K} \log_p (\sigma_1)\}$ and $\{(x,\beta x)\,|\,x<0,\,\beta=- \frac{1}{c_K} \log_p (\sigma_2)\}$, respectively. For any given $r\in (0, r_0]$, by making $\delta$ sufficiently small, we can ensure that the image of $\log_p S^{\delta}$ under $\log\mu_A$ is entirely contained in $\log_p (S_{[\sigma_1, \sigma_2]}^{(0, r]})$. Consequently, we obtain an injective map
\[\mu_A \colon S^{\delta}\rightarrow S^{(0,r]}_{[\sigma_1, \sigma_2]}.\]
%
%
% Schematic layout, not a specified integral matrix.
% Both panels are 5.95cm square; x and y units agree in every scope.
%
\begin{figure}  
    \centering
  \scalebox{0.65}{
\begin{tikzpicture}[font=\large,every node/.style={text=sectoraxis,inner sep=2pt}]
\fill[white] (0,0) rectangle (17.2,7.3);
\begin{scope}[x=16.913727cm,y=16.913727cm,
  xshift=-10.855394cm,yshift=-11.305394cm]
  \draw[sectoraxis,-{Stealth}] (.706845610,.706845610)--(1.058630878,.706845610) node[below right] {$x$};
  \draw[sectoraxis,-{Stealth}] (.706845610,.706845610)--(.706845610,1.058630878) node[above] {$y$};
  \path[pattern=north east lines,pattern color=orange]
    (.824149364,.824149364) rectangle (1,1);
  \draw[excluded,line width=.85pt]
    (.824149364,1)--(.824149364,.824149364)--(1,.824149364);
  \draw[sectorink,line width=.85pt,dashed]
    (.824149364,1)--(1,1)--(1,.824149364);
  \draw[sectorred,line width=1.2pt]
    plot[domain=.876475194:1,samples=120] (\x,{pow(\x,1.466879637)});
  \draw[sectorred,line width=.95pt,line cap=round,dash pattern=on .1pt off 3.6pt]
    plot[domain=.789371904:.876475194,samples=100] (\x,{pow(\x,1.466879637)});
  \node[sectorred,anchor=west,xshift=6pt] at (.831784162,.763247247) {$\sigma_\xi$};
  \foreach \xx/\yy in {.824149364/1,1/1,1/.824149364}
    \draw[sectoraxis,fill=white,line width=.8pt] (\xx,\yy) circle[radius=1.8pt];
  \fill[sectoraxis] (.824149364,.824149364) circle[radius=1.8pt];
  \foreach \t/\lab in {.824149364/{p^{-\delta}},1/1} {
    \draw[sectoraxis] (\t,.706845610) ++(0,-.002)--++(0,.004);
    \node[below=3pt] at (\t,.706845610) {$\lab$};
    \draw[sectoraxis] (.706845610,\t) ++(-.002,0)--++(.004,0);
    \node[left=3pt] at (.706845610,\t) {$\lab$};
  }
\end{scope}
\draw[sectoraxis,-{Stealth},line width=.85pt]
  (7.55,3.56)--node[above=5pt] {$\mu_A$}(9.85,3.56);
\begin{scope}[x=5.409091cm,y=5.409091cm,xshift=10.5cm,yshift=.65cm]
  \draw[sectoraxis,-{Stealth}] (0,0)--(1.1,0) node[below right] {$x$};
  \draw[sectoraxis,-{Stealth}] (0,0)--(0,1.1) node[left] {$y$};
  \path[fill=sectorboundary!14]
    plot[domain=.5:1,samples=160] (\x,{pow(\x,.12)})
    -- plot[domain=1:.5,samples=160] (\x,{pow(\x,.45)}) -- cycle;
  \foreach \beta in {.12,.45} {
    \draw[sectorboundary,line width=.75pt,dashed]
      plot[variable=\t,domain=0:{pow(.5,\beta)},samples=160] ({pow(\t,1/\beta)},\t);
    \draw[sectorboundary,line width=.85pt]
      plot[domain=.5:1,samples=160] (\x,{pow(\x,\beta)});
  }
  \draw[sectorboundary,line width=.85pt] (.5,{pow(.5,.45)})--(.5,{pow(.5,.12)});
  \draw[sectorguide,densely dotted] (.5,0)--(.5,{pow(.5,.45)});
  \node[below=3pt] at (.5,0) {$p^{-c_Kr}$};
  \coordinate (O) at (1,1);
  \coordinate (A) at (.714142843,.960404554);
  \coordinate (B) at (.714142843,.859415713);
  \coordinate (C) at (.51,.825386764);
  \path[pattern=north east lines,pattern color=orange]
    plot[domain=1:.714142843,samples=100] (\x,{pow(\x,.12)})
    -- plot[domain=.714142843:.51,samples=100] (\x,{pow(.714142843,-.33)*pow(\x,.45)})
    -- plot[domain=.51:.714142843,samples=100] (\x,{pow(.714142843,.33)*pow(\x,.12)})
    -- plot[domain=.714142843:1,samples=100] (\x,{pow(\x,.45)}) -- cycle;
  \draw[sectoraxis,line width=.85pt]
    plot[domain=.714142843:.51,samples=100] (\x,{pow(.714142843,-.33)*pow(\x,.45)})
    -- plot[domain=.51:.714142843,samples=100] (\x,{pow(.714142843,.33)*pow(\x,.12)});
  \foreach \beta in {.12,.45}
    \draw[sectorink,line width=.85pt,dashed]
      plot[domain=.714142843:1,samples=100] (\x,{pow(\x,\beta)});
  \draw[sectorred,line width=.95pt,line cap=round,dash pattern=on .1pt off 3.6pt]
    plot[variable=\t,domain=0:{pow(.567685206,.316227766)},samples=160]
      ({pow(\t,1/.316227766)},\t);
  \draw[sectorred,line width=1.2pt]
    plot[domain=.567685206:1,samples=120] (\x,{pow(\x,.316227766)});
  \foreach \pt in {O,A,B}
    \draw[sectorink,fill=white,line width=.8pt] (\pt) circle[radius=1.8pt];
  \fill[sectoraxis] (C) circle[radius=1.8pt];
%  \node[above left=2pt] at (A) {$A'$};
%  \node[below right=4pt] at (B) {$B'$};
%  \node[left=10pt] at (C) {$C'$};
%  \node[above left=4pt] at (O) {$(1,1)$};
  \node[sectorred,anchor=west] at (.33,.752) {$\sigma_\xi$};
  \node[sectorboundary,anchor=west] at (.36,1) {$y=x^{\beta_-}$};
  \node[sectorboundary,anchor=west] at (.62,.69) {$y=x^{\beta_+}$};
  \foreach \t in {.5,1} \draw[sectoraxis] (\t,-.006)--(\t,.006);
  \draw[sectoraxis] (-.006,1)--(.006,1);
  \node[below=3pt] at (0,0) {$0$};
  \node[below=3pt] at (1,0) {$1$};
  \node[left=3pt] at (0,1) {$1$};
\end{scope}
\end{tikzpicture}
}%scalebox
    \caption{A cartoon illustration of $\mu_A: S^{\delta}\rightarrow S^{(0,r]}_{[\sigma_1, \sigma_2]}$}
    \label{fig:region2}
\end{figure}

Passing to the ring of functions, we obtain an injective ring homomorphism
\[
\mu_{A}^*\colon \mathbf{C}^r_{[\sigma_1, \sigma_2]}=\Gamma(A_{K_0'} (S_{[\sigma_1, \sigma_2]}^{(0, r]}), \mO) \ra \Gamma (A_{K_0'} (S^{\delta}), \mO) 
\]
sending $\pi_K \mapsto u^a v^b, [T^\flat] \mapsto u^c v^d$. Here, we have identified $A_{K_0'} (S^{\delta})$ with the half-open poly-annulus $A_{K_0'} ([p^{-\delta}, 1)) \times A_{K_0'} ([p^{-\delta}, 1))$ with coordinate functions $u,v$. The merit here is that $A_{K_0'} (S^{\delta})$ is an honest poly-annulus. In what follows, we shall study the $p$-adic differential equation $\mathbf{N}_{\mathrm{dR}, [\sigma_1, \sigma_2]}^{r}$ over $\mathbf{C}^r_{[\sigma_1, \sigma_2]}$ using the coordinate $u,v$ by taking base change to $A_{K_0'} ([p^{-\delta}, 1)) \times A_{K_0'} ([p^{-\delta}, 1))$ along $\mu_A^*$, so that the multivariable $p$-adic Fuchs theorem applies.

\begin{remark} \label{remark:restrict_to_fake_robba_1}
The map $\mu_A$ sends the rectangle $S^{\delta}$ to a ``tilted rectangular region'' $\mu_A(S^{\delta})$ inside $S^{(0,r]}_{[\sigma_1, \sigma_2]}$. See Figure~\ref{fig:region2} for an illustration of $\mu_A$, where both the rectangular and the tilted rectangular regions are colored in orange. The curve $y=x^{\beta_{-}}$ (resp. $y=x^{\beta_+}$) on the right-hand side of the illustration corresponds to the line $y=1$ (resp. $x=1$) on the left-hand side. Moreover, the curvy ray (colored in red)
\[\big\{(x, x^{\beta})\,\big|\,p^{-c_Kr}\le x<1,\,\,\beta=-\frac{1}{c_K}\log_p\sigma_{\xi}\big\}\]
on the right-hand side corresponds to the fake Robba ring $\mathscr{R}^{\lambda}$ (as $r\rightarrow 0^+$). Its preimage is a line segment (colored in red) with irrational slope on the left-hand side, corresponding to another fake Robba ring $\mathscr{R}^{A^{-1}\lambda}$. The map $\mu_A^*$ induces an isomorphism $\mathscr R^{\lambda} \cong \mathbf{C}_{\xi}^{\dagger}\isom \mathscr R^{A^{-1} \lambda}$ that fits into a commutative diagram
\[ 
\begin{tikzcd}
         \Gamma(A_{K_0'} (S_{[\sigma_1, \sigma_2]}^{(0, r]}), \mO) \arrow[r, "\mu_A^*"] \arrow[d] & \Gamma (A_{K_0'} (S^{\delta}), \mO) \arrow[d]  \\
         \mathscr R^{\lambda} \arrow[r, "\sim"]  & \mathscr R^{A^{-1} \lambda}. 
     \end{tikzcd} 
     \]
The Frobenius structure on $\mathscr{R}^{\lambda}$ transfers to a Frobenius structure on $\mathscr{R}^{A^{-1} \lambda}$.
\end{remark}

\begin{remark}\label{remark: two dimensional Robba ring}
Let $A^{-1} = \left( \begin{smallmatrix} a' & b' \\ c' & d' \end{smallmatrix} \right)$ be the inverse of $A$. Then $u\mapsto \pi_K^{a'}[T^{\flat}]^{b'}$, $v\mapsto \pi_K^{c'}[T^{\flat}]^{d'}$ induces an injective ring map
\[\nu_A^*\colon \Gamma (A_{K_0'} (S^{\delta}), \mO)\rightarrow \Gamma(A_{K_0'} (S_{[\sigma'_1, \sigma'_2]}^{(0, {r_1}])}, \mO)\]
for some $0<r_1<r$ and for any $\sigma'_1, \sigma'_2\in p^{\Q}$ satisfying $\sigma_1<\sigma'_1<\sigma_{\xi}<\sigma'_2<\sigma_2$, such that the composite $\nu_A^*\circ \mu_A^*$ is simply the natural inclusion $\mathbf{C}^r_{[\sigma_1, \sigma_2]}\subset\mathbf{C}^{r_1}_{[\sigma'_1, \sigma'_2]}$. Taking colimit as $r$ and $s$ go to 0, we obtain injective ring maps
\[\mathbf{C}^{\dagger}_{[\sigma_1, \sigma_2]}\xrightarrow[]{\mu_A^*} \varinjlim_{\delta\rightarrow 0^+}\Gamma (A_{K_0'} (S^{\delta}), \mO)\xrightarrow[]{\nu_A^*}\mathbf{C}^{\dagger}_{[\sigma'_1, \sigma'_2]}\]
such that the composite is the natural inclusion. Note that the term in the middle is just the (honest) two-dimensional Robba ring \[\mathscr{R}_{K'_0}\widehat{\otimes}_{K'_0}\mathscr{R}_{K'_0}=\varinjlim_{r\rightarrow 0^+}\varprojlim_{0<s\le r}\mathscr{R}^{[s,r]}_{K'_0}\widehat{\otimes}_{K'_0}\mathscr{R}^{[s,r]}_{K'_0}\] with coordinates $u,v$. Moreover, since both $\mu_A^*$ and $\nu_A^*$ are given by monomial change of coordinates, it is straightforward to check that the Frobenius structures on $\mathbf{C}^r_{[\sigma_1, \sigma_2]}$ and $\mathbf{C}^{\dagger}_{[\sigma_1, \sigma_2]}$ induce Frobenius structures on the poly-annulus $\Gamma (A_{K_0'} (S^{\delta}), \mO)$ and the two-dimensional Robba ring $\mathscr{R}_{K'_0}\widehat{\otimes}_{K'_0}\mathscr{R}_{K'_0}$.
\end{remark}

\vspace{0.1in}
%%% 6.5.3
\subsubsection{Strong unipotence of $\nabla$} \label{sss:unipotence_special}
\noindent
\vspace{0.1in}

\noindent
For $r\in (0,r_0]$ and $\delta$ as above, consider 
\[
N_{\delta} \coloneqq \mathbf{N}_{\mathrm{dR}, [\sigma_1, \sigma_2]}^{r} \otimes_{\mathbf{C}_{[\sigma_1, \sigma_2]}^r, \mu_A^*} \Gamma(A_{K_0'} (S^{\delta}), \mO),
\]
equipped with the full connection
\[
\nabla = D_u du + D_v dv = \Delta_u \cdot  d\log u + \Delta_v \cdot  d\log v 
\]
where $\Delta_u = u D_u$ and $\Delta_v = v D_v$. By Remark~\ref{remark:restrict_to_fake_robba_1}, the map $\mu_A^*$ induces an isomorphism
\begin{equation}\label{eq:N_delta}
\mathbf{N}_{\mathrm{dR}, [\sigma_1, \sigma_2]}^r \otimes_{\mathbf{C}_{[\sigma_1, \sigma_2]}^r}  \mathbf{C}_{\xi}^\dagger=\mathbf{N}_{\mathrm{dR}, \xi} \isom N_{\delta}\otimes_{\Gamma(A_{K_0'} (S^{\delta}), \mO)} \mathscr R^{A^{-1} \lambda}
\end{equation}
of $(\varphi, \nabla)$-modules, under the identification $\mathscr{R}^{\lambda}\cong \mathscr{R}^{A^{-1}\lambda}$.

By \eqref{eq: Robba condition NdR}, Remark~\ref{remark: two dimensional Robba ring}, and \cite[Proposition~4.2.5]{kedlaya_st3}, we know that $N_{\delta}$ satisfies the Robba condition. Moreover, since $N_{\delta}$ is equipped with a Frobenius structure (cf. Remark~\ref{remark: two dimensional Robba ring}), we know that $N_{\delta}$ admits a set of $p$-adic exponents $\{a_1, \ldots, a_d\}$ with $a_i \in (\Q\cap \Z_p)^2$ such that they have $p$-adically non-Liouville differences. This is precisely the situation where we can apply the multivariable $p$-adic Fuchs theorem to $N_{\delta}$. By \cite[Corollary~3.20]{Wang},\footnote{\cite[Corollary~3.20]{Wang} is an enhancement of \cite[Th\'eor\`eme pg.216]{Gachet}, removing the freeness condition from the latter.} we conclude that $N_{\delta}$ is free and admits a basis such that the matrices of $\Delta_u$ and $\Delta_v$ both have entries in $K_0'$ and their eigenvalues represent (the equivalence class of) the exponents $a_i$'s. With respect to this basis, the full connection can be written as  
\[
\nabla = d + T_1 \cdot d \log u + T_2 \cdot d \log v,
\]
with $T_i \in \mathrm{Mat}_{d\times d} (K_0')$ satisfying $T_1 T_2 = T_2 T_1$.

We claim that there exist a basis as above so that $T_1$ and $T_2$ are upper-triangular nilpotent. Via a change of basis by a matrix in $\GL_d(K'_0)$, we may assume that $N_{\delta}$ decomposes into simultaneous generalized eigenspaces of $T_1$ and $T_2$. Let $M$ be any one of these direct summands, with corresponding basis $e_1, \ldots, e_h$. With respect to this basis, we may assume that $T_i=\lambda_iI_d+N_i$ with $N_i$ upper-triangular nilpotent, for $i=1,2$; namely, 
\[\Delta_{u} = u \frac{d}{du} + \lambda_1 I + N_1,\qquad \Delta_{v} = v \frac{d}{dv} + \lambda_2 I + N_2.\]
Note that $\lambda_i \in \Q \cap \Z_p$. We claim that $\lambda_1$ and $\lambda_2$ are actually integers. Granting this claim, we simply modify the basis on $M$ by multiplying $u^{-\lambda_1} v ^{- \lambda_2}$ (and repeat this on each direct summand), and we are done.

It remains to prove the claim. Suppose that $(\lambda_1, \lambda_2) \notin \Z^2$. Consider
\[M' \coloneqq M \otimes_{\Gamma(A_{K_0'} (S^{\delta}), \mO)} \mathscr R^{A^{-1} \lambda}\]
(cf. Remark~\ref{remark:restrict_to_fake_robba_1}), which is a direct summand of $N_{\delta} \otimes_{\Gamma(A_{K_0'} (S^{\delta}), \mO)} \mathscr R^{A^{-1} \lambda}\cong \mathbf{N}_{\mathrm{dR}, \xi}$ (cf. \eqref{eq:N_delta}). By assumption, $\mathbf{N}_{\mathrm{dR}, \xi}$ is unipotent, and hence so is $M'$. In particular, $M'$ admits a nonzero horizontal vector for the full connection, say $s = \sum_{i = 1}^{h} s_i e_i \in M'$. We write each $s_i = \sum_{m, n \in \Z} s^{mn}_i u^m v^n$ in terms of power series expansions in the fake Robba ring. Let $\mathbf{s} = (s_1, \cdots, s_h)$ as a column vector and write $\mathbf{s} = \sum_{m, n \in \Z} \mathbf{s}^{mn} u^m v^n$ accordingly. Since $\Delta_u(\mathbf{s})=\Delta_v(\mathbf{s})=0$, we obtain
\[
\Big( (m+\lambda_1) I_d + N_1 \Big) \mathbf{s}^{mn} = 0, \qquad \Big( (n+\lambda_2) I_d + N_2 \Big) \mathbf{s}^{mn} = 0. 
\]
If $\lambda_1 \notin \Z$, then $(m+\lambda_1) I + N_1 $ is invertible for all $m$, so $s = 0$. Similarly, if $\lambda_2 \notin \Z$, then $s = 0$. This is a contradiction and justifies the claim.

So far, we have achieved strong unipotence for $N^{\delta}$. Using $\nu_A^*$, we obtain strong unipotence for $\mathbf{N}^r_{\mathrm{dR}, [\sigma_1, \sigma_2]}$ for sufficiently small $r$ and for $\sigma_1, \sigma_2$ sufficiently close to $\sigma_{\xi}$ (note that $K''_0$ is just $K'_0$ here). This settles the special case of Proposition~\ref{prop:geometric_unipotence_over_a_region}.

\vspace{0.1in}
%%% new 6.6
\subsection{Unipotence over annuli: general case}\label{ss:geometric_uni_general}
\noindent
\vspace{0.1in}

\noindent Finally, we treat the general case of Proposition~\ref{prop:geometric_unipotence_over_a_region}. The idea is to reduce to the special case handled in \S \ref{ss:geometric_uni_special}.

\vspace{0.05in}
%%% 6.6.1
\subsubsection{Some preparations}\label{sss: some prep}
\noindent
\vspace{0.1in}

\noindent Recall from \S \ref{ss:geometric_realization_fake}--\ref{ss:contagion_of_goem_unipotence} the finite \'etale extension $\mathbf{C}^{\dagger}_{\xi}\rightarrow \mathscr{R}'$ of fake Robba rings induced from the finite separable extension $R'/R_{\xi}$. We have also constructed a finite \'etale cover $f\colon X_{[\sigma_1, \sigma_2]}' \ra X_{[\sigma_1, \sigma_2]}$ that realizes the extension $\mathbf{C}^{\dagger}_{\xi}\rightarrow \mathscr{R}'$. For $r>0$ sufficiently small, the map $f$ induces a finite \'etale map $\mathbf{C}_{[\sigma_1, \sigma_2]}^{r} \ra {\mathbf{C}'}_{[\sigma_1, \sigma_2]}^{r}$ and a differential module ${\mathbf{N}'}_{\mathrm{dR}, [\sigma_1, \sigma_2]}^r$ over ${\mathbf{C}'}_{[\sigma_1, \sigma_2]}^{r}$, equipped with a Frobenius structure. 

To proceed, we first introduce some more notation. Recall that $R_{\xi}$ (resp. $R'$) coincides with the completed stalk $\mathscr{H}(\overline{\xi})$ (resp. $\mathscr{H}(\overline{\xi}')$) (see Remark~\ref{remark:description_R_alpha} and the proof of Lemma~\ref{lemma:relatively_algebraically_closed}). In particular, the residue field of $R_{\xi}$ (resp. $R'$) is $\kappa$ (resp. the algebraic closure $\kappa'$ of $\kappa$ in $R'$). Note that $K''_0$ in the statement of Proposition~\ref{prop:geometric_unipotence_over_a_region} coincides with $W(\kappa')[1/p]$. The following is an analogue of $\mathbf{A}^r_{[\sigma_1, \sigma_2]}$, but with variables $x,y$ instead of $\pi_K$ and $[T^{\flat}]$.

\begin{definition}
Let $\tau_0>0$ such that $\tau_0\not\in p^{\Q}$.
\begin{enumerate}
    \item For any $0<\tau_1 < \tau_0 < \tau_2$ with $\tau_1, \tau_2 \in p^\Q$ and any sufficiently small $s>0$, define 
\[
\mathbf{A}^{s}_{[\tau_1, \tau_2]} (x, y) \coloneqq  \left\{ \sum_{i, j \in \Z} c_{ij} x^i y^j \: ~\vline~ \:  c_{ij} \in W(\kappa'), |c_{ij}|\cdot p^{-c_K \cdot s i} \tau^{sj} \ra  0 \textrm{ for all } \tau\in [\tau_1, \tau_2], \textrm{ as } |i|+|j| \ra +\infty\right\}. 
\]
For any $s'\in (0,s]$, we can define a norm $\lambda_{s'}$ on $\mathbf{A}^{s}_{[\tau_1, \tau_2]} (x, y)$ by
\[\lambda_{s'} \left(\sum_{i,j\in \Z} c_{ij} x^i y^j\right) = \sup_{i,j\in \Z}\sup_{\tau \in [\tau_1, \tau_2]} \{|c_{ij}|\cdot p^{-c_K \cdot s'i} \tau^{s'j}\}.\] 
\item Let
\[
\mathbf{A}^s_{\tau_0, \mathrm{colim}} (x, y) \coloneqq  \varinjlim_{\tau_1<\tau_0< \tau_2}\mathbf{A}^{s}_{[\tau_1, \tau_2]} (x, y)
\]
and
\[
\mathbf{A}^{\dagger}_{\tau_0, \mathrm{colim}} (x, y) \coloneqq  \varinjlim_{s \ra 0^+}\mathbf{A}^s_{\tau_0, \mathrm{colim}} (x, y). 
\]
Let $\mathbf{A}^s_{\tau_0}(x,y)$ be the completion of $\mathbf{A}^s_{\tau_0, \mathrm{colim}} (x, y)$ with respect to the norm $\lambda_{s, \tau_0}$ defined by
\[\lambda_{s, \tau_0} \left(\sum_{i,j\in \Z} c_{ij} x^i y^j\right) = \sup_{i,j\in \Z}\{|c_{ij}|\cdot p^{-c_K \cdot si} \tau_0^{sj}\}\]
and let $\mathbf{A}^{\dagger}_{\tau_0}(x,y)=\varinjlim_{s\rightarrow 0^+}\mathbf{A}^s_{\tau_0}(x,y)$.
\item Let $\mathbf{B}^{s}_{[\tau_1, \tau_2]} (x, y)=\mathbf{A}^{s}_{[\tau_1, \tau_2]} (x, y)[1/p]$ and $\mathbf{B}^s_{\tau_0}(x,y)=\mathbf{A}^s_{\tau_0}(x,y)[1/p]$.
\item Let $\mathbf{C}^{[s',s]}_{[\tau_1, \tau_2]} (x, y)$ (resp. $\mathbf{C}^{s}_{[\tau_1, \tau_2]} (x, y)$) be the completion of $\mathbf{B}^{s}_{[\tau_1, \tau_2]} (x, y)$ with respect to the family of norms $\{\lambda_t\}_{t\in [s', s]}$ (resp. family of norms $\{\lambda_t\}_{t\in (0, s]}$). Let $\mathbf{C}^{\dagger}_{[\tau_1, \tau_2]} (x, y)=\varinjlim_{s>0}\mathbf{C}^{s}_{[\tau_1, \tau_2]} (x, y)$.
\item Let
\[
\mathbf{C}^s_{\tau_0, \mathrm{colim}} (x, y) \coloneqq  \varinjlim_{\tau_1<\tau_0< \tau_2}\mathbf{C}^{s}_{[\tau_1, \tau_2]} (x, y).
\]
and
\[
\mathbf{C}^{\dagger}_{\tau_0, \mathrm{colim}} (x, y) \coloneqq  \varinjlim_{s \ra 0^+}\mathbf{C}^s_{\tau_0, \mathrm{colim}} (x, y). 
\]
Let $\mathbf{C}^s_{\tau_0} (x, y)$ be the Fr\'echet completion of $\mathbf{C}^s_{\tau_0, \mathrm{colim}} (x, y)$ with respect to $\{\lambda_{t,\tau_0}\}_{t \in (0, s]}$ and let \[\mathbf{C}^{\dagger}_{\tau_0}(x,y)=\varinjlim_{s\rightarrow 0^+}\mathbf{C}^s_{\tau_0}(x,y).\]
\end{enumerate}
\end{definition}

The main idea to prove Proposition~\ref{prop:geometric_unipotence_over_a_region} is to construct continuous maps 
\[{\mathbf{C}'}^r_{[\sigma_1, \sigma_2]}\xrightarrow[]{\psi} \mathbf{C}^s_{[\tau_1, \tau_2]}(x,y), \qquad \mathbf{C}^s_{[\tau_1, \tau_2]}(x,y)\xrightarrow[]{\phi} {\mathbf{C}'}^{r'}_{[\sigma'_1, \sigma'_2]}\]
for suitable $r,s,r',\sigma_1, \sigma_2, \tau_1, \tau_2, \sigma'_1, \sigma'_2$ and coordinates $x,y$, such that $\phi, \psi$ are compatible with Frobenius and differential structures, and such that $\phi\circ\psi$ and $\psi\circ \phi$ coincide with natural inclusions. Given such (systems of) maps, we can translate the question of $p$-adic differential equations over ${\mathbf{C}'}^r_{[\sigma_1, \sigma_2]}$ to a question over $\mathbf{C}^s_{[\tau_1, \tau_2]}(x,y)$; but the latter is already analyzed in \S \ref{ss:geometric_uni_special}. As we shall see, to construct these maps, we first construct them on the integral level; that is, we construct continuous maps
\[{\mathbf{A}'}^r_{[\sigma_1, \sigma_2]}\xrightarrow[]{\psi} \mathbf{A}^s_{[\tau_1, \tau_2]}(x,y), \qquad \mathbf{A}^s_{[\tau_1, \tau_2]}(x,y)\xrightarrow[]{\phi} {\mathbf{A}'}^{r'}_{[\sigma'_1, \sigma'_2]}\]
with similar properties. The construction of these maps will occupy \S \ref{sss: algebraization of analytic covers}--\ref{sss: construction of phi and psi}.

\vspace{0.05in}
%%% 6.6.2
\subsubsection{Algebraization of analytic covers}\label{sss: algebraization of analytic covers}
\noindent
\vspace{0.1in}

\noindent
Let ${\mathbf{A}'}_{\xi, \mathrm{colim}}^\dagger  \coloneqq   \varinjlim_{r \ra 0^+}\varinjlim_{\sigma_1<\sigma_{\xi} < \sigma_2} {\mathbf{A}'}_{[\sigma_1, \sigma_2]}^{r}$, which is finite \'etale over $\mathbf{A}_{\xi, \mathrm{colim}}^\dagger$. Recall that $\mathbf{A}_{\xi, \mathrm{colim}}^\dagger/p=\mathbf{k} (\cl \xi)$ and ${\mathbf{A}'}_{\xi, \mathrm{colim}}^\dagger/p = \mathbf{k} (\cl \xi')$. Let $|\cdot|_{\overline{\xi}}$ (resp. $|\cdot|_{\overline{\xi}'}$) denote the corresponding norm on the uncompleted stalk $\mathbf{k}(\overline{\xi})$ (resp. $\mathbf{k}(\overline{\xi}')$). Also note that the finite separable extension $ \mathbf{k} (\cl \xi')/  \mathbf{k} (\cl \xi)$ induces $ \mathbf{k} (\cl \xi') \otimes_{ \mathbf{k} (\cl \xi)} R_{\xi} \cong R'$.  

\begin{lemma} \label{lemma:colim_A_hensel}
    \begin{enumerate}
\item Both $(\mathbf{A}_{\xi, \mathrm{colim}}^\dagger,p)$ and $({\mathbf{A}'}_{\xi, \mathrm{colim}}^\dagger, p)$ are henselian pairs. 
\item Let $f \in {\mathbf{A}'}_{\xi, \mathrm{colim}}^\dagger$ be an element whose mod $p$ reduction $\cl f$ is nonzero. Then $f$ is invertible. If moreover $f$ satisfies $|\cl f|_{\cl \xi'} < 1$, then there exist a sufficiently small $r$ and $\sigma_1, \sigma_2$ sufficiently close to $\sigma_\xi$ such that $f$ is topologically nilpotent in ${\mathbf{A}'}_{[\sigma_1, \sigma_2]}^r$. 
    \end{enumerate}
\end{lemma}

\begin{proof}
\begin{enumerate}
\item Each $(\mathbf{A}_{[\sigma_1, \sigma_2]}^\dagger , p)$ is henselian by the proof of \cite[Proposition~5.28]{relative}, so the filtered colimit $(\mathbf{A}_{\xi, \mathrm{colim}}^\dagger , p)$ is again henselian. The henselianity of $({\mathbf{A}'}_{\xi, \mathrm{colim}}^\dagger, p)$ follows from finite \'etaleness. 

\item The first claim on the invertibility of $f$ is a consequence of henselianity. 

For the second claim, we first show that the claim holds for $\mathbf{A}_{\xi, \mathrm{colim}}^\dagger$ in place of ${\mathbf{A}'}_{\xi, \mathrm{colim}}^\dagger$. Suppose $f_0 \in \mathbf{A}_{\xi, \mathrm{colim}}^\dagger$ such that $|\cl f_0|_{\cl \xi} < 1$. Then, up to shrinking $[\sigma_1, \sigma_2]$, we may assume that $\overline{f}_0\in R_{[\sigma_1, \sigma_2]}$ and $q \coloneqq \overline{\lambda}_{[\sigma_1, \sigma_2]}(\overline{f}_0) < 1$ where $\overline{\lambda}_{[\sigma_1, \sigma_2]}$ is the spectral norm on $R_{[\sigma_1, \sigma_2]}$. By \cite[Proposition~5.18(1)]{relative}, there exists a sufficiently small $r_0$ and a lift $\widetilde{f}_0\in \mathbf{A}^{r_0}_{[\sigma_1, \sigma_2]}$ of $\cl f_0$ such that $\lambda_r (\widetilde{f}_0) = q^r$ for all $r \in (0, r_0]$. If we write $f_0 = \widetilde{f}_0 + p g$ for $g \in \mathbf{A}_{[\sigma_1, \sigma_2]}^{r_0}$, then $\lambda_r (g) \le \max \{1, \lambda_{r_0} (g)^{r/r_0}\}$. Hence, there exists $r_1 < r_0$ such that  $\lambda_r (p g) < 1$ for all $r \in (0, r_1]$. This implies $\lambda_r (f_0)  < 1 $ for all $r \in (0, r_1]$. 

Now we return to the statement for ${\mathbf{A}'}_{\xi, \mathrm{colim}}^\dagger$. View ${\mathbf{A}'}_{\xi, \mathrm{colim}}^\dagger$ as a finite projective module of rank $m$ over $\mathbf{A}_{\xi, \mathrm{colim}}^\dagger$, and write $P_f(T) = T^m + a_{m-1} T^{m-1} + \cdots + a_0 \in \mathbf{A}_{\xi, \mathrm{colim}}^\dagger [T]$ for the characteristic polynomial of the multiplication-by-$f$ map. The roots of the mod $p$ reduction $\cl P_f (T)$ are Galois conjugates $\tau (\cl f)$ of $\cl f \in \mathbf{k}(\cl \xi') =  {\mathbf{A}'}_{\xi, \mathrm{colim}}^\dagger/p$ in some fixed algebraic closure $ \mathbf{k}(\cl \xi)^{\tu{alg}}$, and each of them satisfies $|\tau (\cl f)|_{\cl \xi'} < 1$. This implies $|\cl a_i|_{\cl \xi} < 1$ for all $i$. By the discussion above, we may assume that $a_i\in \mathbf{A}^{r_1}_{[\sigma_1, \sigma_2]}$ and $f\in {\mathbf{A}'}^{r_1}_{[\sigma_2, \sigma_2]}$ for some $r_1, \sigma_1, \sigma_2$, and that $\lambda_r (a_i) < 1$ for all $i$ and all $r\in (0,r_1]$. Finally, we show that $\lambda_r(f)<1$ for all $r\in (0,r_1]$. By Cayley--Hamilton, we have $f^m + a_{m-1} f^{m-1} \cdots + a_0 = 0$. If $\lambda_r (f) \ge 1$, then 
\[
\lambda_r(f)^m = \lambda_r (f^m) \le \max_{i < m} \{\lambda_r(a_i) \lambda_r(f)^i\} < \lambda_r(f)^m
\]
which is impossible. Here, we have used the fact that $\lambda_r$ is power-multiplicative on ${\mathbf{A}'}^r_{[\sigma_1, \sigma_2]}$.
\end{enumerate}
\end{proof}

By Krasner's lemma, there is a finite separable extension $F'/\kappa (\cl \pi_K, T^\flat)$ such that $R' = F' \otimes_{\kappa (\cl \pi_K, T^\flat)} R_{\xi}$. We view $F'$ as a subfield of $R'$. Let \[\cl G (X) = X^{d_0} + \sum_{i=0}^{d_0-1} a_i X^i \in \kappa' (\cl \pi_K, T^\flat) [X]\] be the minimal polynomial for a primitive element of the finite separable extension $F'/\kappa' (\cl \pi_K, T^\flat)$. Let $\cl g \in \kappa' [\cl \pi_K, T^\flat]$ be a nonzero polynomial that clears the denominators of the coefficients of $\cl G$, that is, $\cl g \cdot a_i \in  \kappa' [\cl \pi_K, T^\flat]$ for all $i=0,\ldots, d_0-1$. Pick an (arbitrary) $g \in W(\kappa') [\pi_K, [T^\flat]]$ lifting $\cl g$ and pick a lift
\[G [X] \in W(\kappa') [\pi_K, [T^\flat], 1/g] [X] \]
of $\cl G$. Let $\Delta(G) \in W(\kappa') [\pi_K, [T^\flat], 1/g] $ denote the discriminant of $G$ and consider
\[
C_0 \coloneqq W(\kappa')  [\pi_K^{\pm}, [T^\flat]^{\pm}, 1/g, 1/\Delta(G)], \qquad D_0 \coloneqq  C_0[X]/G(X).  
\]
By construction, $C_0 \ra D_0$ is finite \'etale. The fraction fields of $C_0/p$ and $D_0/p$ are precisely $\kappa' (\cl \pi_K, T^\flat)$ and $F'$, respectively. Moreover, we have $\cl \Delta(G) \ne 0$, and is thus a unit in the fraction field $ \kappa' (\cl \pi_K, T^\flat)$. By Lemma~\ref{lemma:colim_A_hensel}, there is a natural map $C_0 \ra \mathbf{A}^{\dagger}_{\xi, \mathrm{colim}}$ as both $g$ and $\Delta (G)$ map to units. In particular, we obtain a map $C_0 \ra  \mathbf{A}^{r}_{[\sigma_1, \sigma_2]}$ after decreasing $r$ and shrinking $[\sigma_1, \sigma_2]$ if necessary.

\begin{lemma} \label{lemma:D_0_maps_to_colimA'}
For sufficiently small $r$ and for $\sigma_1$, $\sigma_2$ sufficiently close to $\sigma_{\xi}$, we have 
    \[
      D_0 \otimes_{C_0}  \mathbf{A}^{r}_{[\sigma_1, \sigma_2]}\cong {\mathbf{A}'}^{r}_{[\sigma_1, \sigma_2]}. 
    \]
\end{lemma}

\begin{proof}
By construction, there is an isomorphism 
\[
  D_0/p \otimes_{C_0/p} R_{\xi} \cong  R'.
    \]
By quasi-completeness of $\mathbf{k}(\cl \xi)$ (cf. \cite[Theorem~2.3.3, Proposition~2.4.1]{Berkovich_etale}), we have an equivalence of categories $\big( \mathbf{k} (\cl \xi) \big)_{\mathrm{f\'et}} \isom \big( R_{\xi} \big)_{\mathrm{f\'et}}$. This produces an isomorphism 
\[
 D_0/p \otimes_{C_0/p}  \mathbf{A}^{\dagger}_{\xi, \mathrm{colim}}/p \cong {\mathbf{A}'}_{\xi, \mathrm{colim}}^\dagger/p. 
\]
Finally, by the henselianity of $(\mathbf{A}^{\dagger}_{\xi, \mathrm{colim}}, p)$ proved in Lemma~\ref{lemma:colim_A_hensel}, we obtain an isomorphism 
\[
 D_0\otimes_{C_0}  \mathbf{A}^{\dagger}_{\xi, \mathrm{colim}} \cong {\mathbf{A}'}_{\xi, \mathrm{colim}}^\dagger, 
\]
from which the lemma follows. 
\end{proof}

\vspace{0.05in}
%%% 6.6.3
\subsubsection{The coordinates $x, y$}
\noindent
\vspace{0.1in}

\noindent
To construct $\psi$ and $\phi$, we need a further approximation step to relate $D_0$ with a suitable polynomial ring $W(\kappa')[x,y]$. Recall from \S \ref{ss:geometric_realization_fake} that the valuation $v_\alpha$ on $R_{\xi}$ uniquely extends to a valuation (again denoted as $v_\alpha$) on $R'$. Since $\kappa' (\cl \pi_K, T^\flat) \subset F'\subset R'$, the value group $\Gamma'$ of $F'$ contains $\Z + \alpha \Z$ as a subgroup of finite index. In particular, the valuation $v_\alpha$ on $F'$ is \emph{monomial} in the sense of \cite[Definition~6.1.1]{kedlaya_st3}. Pick a $\Z$-basis $\{\beta_1, \beta_2\}$ of $\Gamma'$ with $\beta_1, \beta_2>0$ (note that $\beta_2/\beta_1 \notin \Q$), and then pick elements $t_1, t_2 \in (F')^\times$ such that $v_{\alpha} (t_i) = \beta_i$.   

\begin{proposition}
There exists a smooth algebraic surface $S$ over $\spec \kappa'$ with function field $F(S) = F'$, a closed point $s_0 \in S$ at which the valuation $v_{\alpha}$ is centered, and regular parameters $\cl x, \cl y \in \mO_{S, s_0}$ such that $t_i = \cl u_i \cdot \cl x^{a_i} \cl y^{b_i}$ for $i=1,2$, with $\cl u_i \in \mO_{S, s_0}^\times$  and $a_i, b_i \in \Z_{\ge 0}$. Here, we say that the valuation is \emph{centered} at $s_0$ if $\fm_{S, s_0}=\mO_{S, s_0}\cap \{a\in F'\;|\; v_{\alpha}(a)>0\}$ where the intersection is taken inside $F'=\mathrm{Frac}\,\mO_{S, s_0}$ (cf. \cite[Definition~2.4.1]{kedlaya_st2}). 
\end{proposition}

\begin{proof}
This follows from the local uniformization theorem of monomial valuations (cf. \cite[Theorem~1.1]{Knaf_Kuhlmann}, also see \cite[Proposition~6.1.3]{kedlaya_st3}). 
\end{proof}  
 
Write $\gamma_1 = v_{\alpha}(\cl x)$ and $\gamma_2 = v_{\alpha}(\cl y)$ for the valuations of $\cl x$ and $\cl y$. Then $\beta_1 = a_1 \gamma_1 + b_1 \gamma_2$ and $\beta_2 = a_2 \gamma_1 + b_2 \gamma_2$. In particular, $\gamma_1$ and $\gamma_2$ are rationally independent. 
Since $S$ is smooth at $s_0$ with regular parameters $\cl x, \cl y$, we have an isomorphism $(\mO_{S, s_0})_{\fm_{S,s_0}}^{\wedge} \cong \kappa' [\![\cl x, \cl y]\!]$. This implies that $\kappa' [\cl x, \cl y]$ is dense in $\mO_{S, s_0}$ for the topology induced by the valuation $v_{\alpha}$. Therefore, $\kappa' (\cl x, \cl y)$ is dense in $F'$, and thus is dense in $R'$. %
Also note that, since $\cl x, \cl y$ are regular parameters at $s_0$, we have an isomorphism 
\[
F' \otimes_{\kappa' (\cl x, \cl y)} \Omega^1_{\kappa' (\cl x, \cl y)/\kappa'} \isom \Omega^1_{F'/\kappa'}  
\]
as $ \Omega^1_{F'/\kappa'}  \cong F' \otimes_{\mO_{S, s_0}} \Omega^1_{\mO_{S, s_0}/\kappa'} \cong F' d \cl x \oplus F' d \cl y$. This implies $\Omega^1_{F'/\kappa' (\cl x, \cl y)} = 0$, so the extension $F'/\kappa' (\cl x, \cl y)$ is finite separable (it is finite because it is finitely generated and has transcendental degree $0$). 

\begin{remark}
    We point out that, the fact that $F'/\kappa' (\cl x, \cl y)$ is a finite separable extension and that $\kappa' (\cl x, \cl y)$ is dense in $F'$ implies that they have the same completion, but does not imply that they are equal. 
\end{remark}

Since $F'=\mathrm{Frac} (D_0/p)$, we can write $\cl x = \cl h_{1}/\cl g_{1}$ and $\cl y = \cl h_2/\cl g_2$ with $h_1, h_2, g_1, g_2\in D_0/p$. Pick liftings $h_1, h_2, g_1, g_2\in D_0$ and consider a ring map $W(\kappa') [x, y] \ra D_0 [1/(g_1 g_2)]$ sending $x \mapsto h_1/g_1$, $y \mapsto h_2/g_2$. This is a lift of $\kappa' [x, y] \lra \mO_{S, s_0} \subset  F' = \mathrm{Frac} (D_0/p)$. By further inverting the Jacobian $J \in D_0 [1/(g_1 g_2)]$ of $x, y$ with respect to $\pi_K, [T^\flat]$, we arrive at an \'etale ring map (``the framing'')
\[
 \square\colon W(\kappa') [x, y] \lra D \coloneqq  D_0[1/(g_1 g_2), 1/J]. 
\]

\vspace{0.05in}
%%% 6.6.4
\subsubsection{Constructions of $\psi$ and $\phi$}\label{sss: construction of phi and psi}
\noindent
\vspace{0.1in}

\noindent \textit{Construction of $\psi$}. 
\vspace{0.05in}

Let $\tau_0 \coloneqq  p^{- c_K \cdot\gamma_2/\gamma_1}$ (note that $\tau_0\not\in p^{\Q}$) and let $[\tau_1, \tau_2] $ be an interval containing $\tau_0$, with $\tau_1, \tau_2 \in p^\Q$. Consider $\mathbf{A}^{\dagger}_{\tau_0, \mathrm{colim}} (x, y)$ and $\mathbf{A}^{\dagger}_{\tau_0} (x, y)$ defined in \S \ref{sss: some prep}, and note that $(\mathbf{A}^{\dagger}_{\tau_0, \mathrm{colim}} (x, y), p)$ is a henselian pair. Let $R_{\tau_0,\mathrm{colim}}\coloneqq \mathbf{A}^{\dagger}_{\tau_0, \mathrm{colim}} (x, y)/p$ and $R_{\tau_0}\coloneqq \mathbf{A}^{\dagger}_{\tau_0} (x, y)/p$. Similarly to Remark~\ref{remark:description_R_alpha} and Lemma~\ref{lemma: colimit of R[sigma1, sigma2]}, $R_{\tau_0,\mathrm{colim}}$ (resp. $R_{\tau_0}$) coincides with the uncompleted residue field (resp. completed residue field) of a type I\!I\!I point corresponding to a circle of radius $\tau_0$ in characteristic $p$, and we have an isomorphism $R_{\tau_0}\cong R'$ by sending $x\mapsto \cl x$, $y\mapsto \cl y$. 

By construction, we have inclusions $F'\subset R'$ and $R_{\tau_0, \mathrm{colim}}\subset R_{\tau_0} \cong R'$. We claim that the first inclusion factors through the second via an injection $F' \hookrightarrow R_{\tau_0,\mathrm{colim}}$. To see this, consider $E \coloneqq F' \otimes_{\kappa' (\cl x, \cl y)}  R_{\tau_0,\mathrm{colim}}$. Note that $E$ is a finite \'etale $R_{\tau_0,\mathrm{colim}}$-algebra since $F'/\kappa' (\cl x, \cl y)$ is finite separable. Consider the composition $m \colon E \otimes_{ R_{\tau_0,\mathrm{colim}}} R'\cong F' \otimes_{\kappa' (\cl x, \cl y)} R' \rightarrow R'$ where the second map is the multiplication map. By the fully faithfulness in \cite[Lemma 2.2.4(a)]{Kedlaya-LiuI}, $m$ descends to a unique $ R_{\tau_0,\mathrm{colim}}$-algebra map $m_0 \colon E \ra  R_{\tau_0,\mathrm{colim}}$. Finally, composing with $F' \ra E$ sending $a \mapsto a \otimes 1$, we get the desired inclusion $F' \hookrightarrow  R_{\tau_0,\mathrm{colim}}$. 

Given this, we deduce that the natural map $W(\kappa') [x, y] \ra \mathbf{A}^{\dagger}_{\tau_0, \mathrm{colim}} (x, y)$ factors as the composition
\[W(\kappa') [x, y]\xrightarrow[]{\square} D \xrightarrow[]{\eta} \mathbf{A}^{\dagger}_{\tau_0, \mathrm{colim}} (x, y)\]
for a unique map 
$\eta\colon D \lra \mathbf{A}^{\dagger}_{\tau_0, \mathrm{colim}} (x, y)$. Indeed, 
since $(\mathbf{A}^{\dagger}_{\tau_0, \mathrm{colim}} (x, y), p)$ is henselian, it suffices to observe that the modulo $p$ reduction $\kappa'[x,y]\rightarrow R_{\tau_0, \mathrm{colim}}$ factors as $\kappa'[x,y]\rightarrow D/p\rightarrow F'\rightarrow R_{\tau_0, \mathrm{colim}}$. Consequently, we obtain a map
\[
\eta \colon D \lra \mathbf{A}^{s}_{[\tau_1, \tau_2]} (x, y)
\]
for sufficiently small $s>0$ and for $\tau_1, \tau_2$ sufficiently close to $\tau_0$.
 
To construct the desired $\psi \colon {\mathbf{A}'}^r_{[\sigma_1, \sigma_2]}\rightarrow \mathbf{A}^s_{[\tau_1, \tau_2]}(x,y)$, it suffices to produce a map $\psi_1 \colon D_0 \ra 
\mathbf{A}^{s}_{[\tau_1, \tau_2]} (x, y) $ and a map $\psi_2 \colon {\mathbf{A}}_{[\sigma_1, \sigma_2]}^{r} \ra \mathbf{A}^{s}_{[\tau_1, \tau_2]} (x, y)$ that agrees on $C_0$ (cf. Lemma~\ref{lemma:D_0_maps_to_colimA'}). In fact, we take $\psi_1$ to be the composition
\[
D_0\rightarrow D \xrightarrow[]{\eta} \mathbf{A}^{s}_{[\tau_1, \tau_2]} (x,y).
\] 

To construct $\psi_2$, we need the following lemma. For the setup, let us write the rational numbers $c_K r$ and $\alpha_i = - \frac{1}{c_K}\log_p \sigma_i$ (for $i =1, 2$) in their reduced rational forms
\[
c_K r = \frac{h_1}{k_1}, \qquad  - \frac{1}{c_K}\log_p \sigma_2 = \frac{e_1}{m_1}, \qquad  - \frac{1}{c_K}\log_p \sigma_1 = \frac{f_1}{n_1}.
\]
Note that $e_1/m_1 < \alpha <  f_1/n_1$ by assumption. 

\begin{lemma} \label{lemma:map_from_A_ring}
    Let $B$ be a Banach $W(\kappa)$-algebra where $W(\kappa)$ is equipped with the standard $p$-adic norm. Suppose that $x, y \in B^\times$ are topologically nilpotent units and suppose that $p^{h_1} x^{-k_1}, y^{m_1}x^{-e_1}, x^{f_1}y^{-n_1}$ are power-bounded. Then there exists a unique continuous ring map $\mathbf{A}^{r}_{[\sigma_1, \sigma_2]} \ra B$ sending $\pi_K \mapsto x$, $[T^\flat] \mapsto y$.
\end{lemma}

\begin{proof}
From the description of $\mathbf{A}^{r}_{[\sigma_1, \sigma_2]}$ in Definition~\ref{definition: description period rings--annulus}(2) (cf. Proposition~\ref{prop:identificaiton_for_A^r}), we know that 
    \[
        \mathbf{A}^{r}_{[\sigma_1, \sigma_2]} = \left( W(\kappa) [\![\pi_K]\!] \Big\langle[T^\flat],  \frac{p^{h_1}}{\pi_K^{k_1}}, \frac{ [T^\flat]^{m_1}}{\pi_K^{e_1}},  \frac{\pi_K^{f_1}}{ [T^\flat]^{n_1}} \Big\rangle \right) [\frac{1}{\pi_K}]
    \]
where completion is taken with respect to the $(p, \pi_K)$-adic topology (equivalently, the $p$-adic topology). In other words, we have 
\[ \mathbf{A}^{r}_{[\sigma_1, \sigma_2]}  \cong  \left( \frac{W(\kappa) [\![\pi_K]\!] \gr{[T^\flat], X_0, X_1, X_2} }{ \big( p^{h_1} - \pi_K^{k_1} X_0, [T^\flat]^{m_1} - \pi_K^{e_1} X_1, \pi_K^{f_1} - [T^\flat]^{n_1} X_2 \big)} \right) [\frac{1}{\pi_K}].\] 
The map to $B$ is then given by sending $\pi_K\mapsto x, [T^\flat] \mapsto y, X_0 \mapsto  p^{h_1}/x^{k_1}, X_1 \mapsto y^{m_1}/x^{e_1}, X_2 \mapsto x^{f_1} / y^{n_1}$. It is straightforward to check that the map is indeed well-defined and continuous. 
\end{proof}

Consider elements $x_1 \coloneqq  \eta(\pi_K)$ and $y_1 \coloneqq \eta([T^\flat])$ in $\mathbf{A}^s_{[\tau_1, \tau_2]}(x,y)$. Note that the images of $x_1$, $y_1$, $p^{h_1} x_1^{-k_1}$, $y_1^{m_1}x_1^{-e_1}$, $x_1^{f_1}y_1^{-n_1}$ in $\mathbf{A}^{\dagger}_{\tau_0, \mathrm{colim}}(x,y)$ all have mod $p$ reductions with norm $<1$. Hence, by decreasing $s$ and shrinking $[\tau_1, \tau_2]$ if necessary, we may assume that 
$x_1$, $y_1$, $p^{h_1} x_1^{-k_1}$, $y_1^{m_1}x_1^{-e_1}$, $x_1^{f_1}y_1^{-n_1}$ are all topologically nilpotent, by applying Lemma~\ref{lemma:colim_A_hensel}. Applying Lemma~\ref{lemma:map_from_A_ring} to $B=\mathbf{A}^s_{[\tau_1, \tau_2]}(x,y)$, we obtain the desired map $\psi_2 \colon {\mathbf{A}}_{[\sigma_1, \sigma_2]}^{r} \ra \mathbf{A}^{s}_{[\tau_1, \tau_2]} (x, y)$. It is routine to check that these two maps agree on $C_0$. \\

\noindent \textit{Construction of $\phi$}.
\vspace{0.05in}

Write $c_K s = h_2/k_2$, $- \frac{1}{c_K}\log_p \tau_2 = e_2/m_2$, $- \frac{1}{c_K}\log_p \tau_1 = f_2/n_2$ in their reduced rational forms. Since the images of $1/(g_1 g_2)$ and $1/J$ in $\mathrm{Frac}(D_0/p) = F'$ are nonzero, the map $D_0 \ra  {\mathbf{A}'}^{r}_{[\sigma_1, \sigma_2]}$ given by Lemma~\ref{lemma:D_0_maps_to_colimA'} extends to a map 
\[
\theta \colon D \lra  {\mathbf{A}'}^{r}_{[\sigma_1, \sigma_2]}
\]
by decreasing $r$ and shrinking $[\sigma_1, \sigma_2]$ if necessary, using the henselianity of $({\mathbf{A}'}_{\xi, \mathrm{colim}}^\dagger, p)$ (see Lemma~\ref{lemma:colim_A_hensel}). Write $x_2 \coloneqq  \theta(x)$ and $y_2 \coloneqq  \theta(y)$. Since $x_2$ and $y_2$ are invertible, and $x_2$, $y_2$, $p^{h_2} x_2^{-k_2}$, $y_2^{m_2}x_2^{-e_2}$, $x_2^{f_2}y_2^{-n_2}$ all have reduction of absolute value $< 1$, (after further decreasing $r'$ and shrinking $[\sigma_1', \sigma_2']$ if necessary), we obtain the desired map $\phi \colon \mathbf{A}^{s}_{[\tau_1, \tau_2]}(x,y)\rightarrow {\mathbf{A}'}^{r'}_{[\sigma'_1, \sigma'_2]}$ by applying Lemma~\ref{lemma:map_from_A_ring}. \\

\noindent \textit{Compositions of $\psi$ and $\phi$}.
\vspace{0.05in}

For the reader's convenience, let us summarize our constructions. Given any sufficiently small $r$ and sufficiently small interval $[\sigma_1, \sigma_2]$ containing $\sigma_{\xi}$, we have a continuous ring map $\psi \colon {\mathbf{A}'}^r_{[\sigma_1, \sigma_2]}\rightarrow \mathbf{A}^s_{[\tau_1, \tau_2]}(x,y)$ for some $s$ and some sufficiently small interval $[\tau_1, \tau_2]$ containing $\tau_0$. Note that if we decrease $r$, then $\psi$ remains valid as long as we decrease $s$ accordingly, without changing $[\sigma_1, \sigma_2]$ or $[\tau_1, \tau_2]$.\footnote{This follows from the construction. Note that in Lemma~\ref{lemma:colim_A_hensel}(2), the choices of $\sigma_1$ and $\sigma_2$ are independent of $r$.} Therefore, we actually obtain a system of maps $\psi$ as $r$ varies. Likewise, we obtain a system of continuous ring maps $\phi \colon \mathbf{A}^{s}_{[\tau_1, \tau_2]}(x,y)\rightarrow {\mathbf{A}'}^{r'}_{[\sigma'_1, \sigma'_2]}$ as $s$ varies.

Unwinding the construction, it is straightforward to check that the compositions 
\[{\mathbf{A}'}^r_{[\sigma_1, \sigma_2]}\xrightarrow[]{\psi} \mathbf{A}^s_{[\tau_1, \tau_2]}(x,y)\xrightarrow[]{\phi}{\mathbf{A}'}^{r'}_{[\sigma'_1, \sigma'_2]}\]
and
\[\mathbf{A}^s_{[\tau_1, \tau_2]}(x,y)\xrightarrow[]{\phi}{\mathbf{A}'}^{r'}_{[\sigma'_1, \sigma'_2]}\xrightarrow[]{\psi} \mathbf{A}^{s'}_{[\tau'_1, \tau'_2]}(x,y)\]
are both natural restriction maps (for suitable $r'<r$, $s'<s$, $[\sigma'_1, \sigma'_2]\subset [\sigma_1, \sigma_2]$, and $[\tau'_1, \tau'_2]\subset [\tau_1, \tau_2]$). 

Lastly, by inverting $p$ and passing to Fr\'echet completions, we obtain systems of continuous ring maps
\[{\mathbf{C}'}^r_{[\sigma_1, \sigma_2]}\xrightarrow[]{\psi} \mathbf{C}^s_{[\tau_1, \tau_2]}(x,y), \qquad \mathbf{C}^s_{[\tau_1, \tau_2]}(x,y)\xrightarrow[]{\phi} {\mathbf{C}'}^{r'}_{[\sigma'_1, \sigma'_2]}\]
such that both $\phi \circ \psi$ and $\psi \circ \phi$ are natural restriction maps.

\vspace{0.05in}
%%% 6.6.5
\subsubsection{Proof of Proposition~\ref{prop:geometric_unipotence_over_a_region}}
\noindent
\vspace{0.1in}

\noindent Passing to colimits, we obtain maps
\[{\mathbf{C}'}^{\dagger}_{[\sigma_1, \sigma_2]}\xrightarrow[]{\psi}\mathbf{C}^{\dagger}_{[\tau_1, \tau_2]}(x,y)\xrightarrow[]{\phi}{\mathbf{C}'}^{\dagger}_{[\sigma'_1, \sigma'_2]},\qquad \mathbf{C}^{\dagger}_{[\tau_1, \tau_2]}(x,y)\xrightarrow[]{\phi}{\mathbf{C}'}^{\dagger}_{[\sigma'_1, \sigma'_2]}\xrightarrow[]{\psi}\mathbf{C}^{\dagger}_{[\tau'_1, \tau'_2]}(x,y)\]
whose compositions are both natural restriction maps. For $\tau_1, \tau_2$ sufficiently close to $\tau_0$, consider
\[N_{[\tau_1, \tau_2]}\coloneqq \mathbf{N}'_{\mathrm{dR}, [\sigma_1, \sigma_2]}\otimes_{{\mathbf{C}'}^{\dagger}_{[\sigma_1, \sigma_2]}, \psi} \mathbf{C}^{\dagger}_{[\tau_1, \tau_2]}(x,y).\]
Using $\psi$ and $\phi$, the Frobenius structure on ${\mathbf{C}'}^{\dagger}_{[\sigma_1, \sigma_2]}$ induces a Frobenius structure on $\mathbf{C}^{\dagger}_{[\tau_1, \tau_2]}(x,y)$. This is essentially due to the construction; indeed, since $\phi$ and $\psi$ are constructed from the integral level, the induced Frobenius map is indeed a lift of the absolute Frobenius map in characteristic $p$. Moreover, passing the completed colimits as $\sigma_1, \sigma_2\rightarrow \sigma_{\xi}$ and $\tau_1, \tau_2\rightarrow \tau_0$, we obtain an isomorphism $\mathscr{R}'\cong \mathbf{C}^{\dagger}_{\tau_0}(x,y)$ of fake Robba rings, compatible with Frobenius structures. Consequently, $N_{[\tau_1, \tau_2]}$ is a finite projective $(\phi, \nabla)$-module over $\mathbf{C}^{\dagger}_{[\tau_1, \tau_2]}(x,y)$ equipped with differential operators $D_x$ and $D_y$ such that
\[N_{[\tau_1, \tau_2]}\otimes_{\mathbf{C}^{\dagger}_{[\tau_1, \tau_2]}(x,y)} \mathbf{C}^{\dagger}_{\tau_0}(x,y)\]
is already unipotent. Finally, applying the arguments in the special case (see \S \ref{ss:geometric_uni_special}), we deduce that $N_{[\tau_1, \tau_2]}$ is strongly unipotent, by shrinking the interval $[\tau_1, \tau_2]$ if necessary. Using $\psi$ and $\phi$, we conclude that $\mathbf{N}'_{\mathrm{dR}, [\sigma_1, \sigma_2]}$ is also strongly unipotent, for $\sigma_1, \sigma_2$ sufficiently close to $\sigma_{\xi}$. This finishes the proof of Proposition~\ref{prop:geometric_unipotence_over_a_region}.

\vspace*{0.3in}

%%% Section 7
\section{Around points of type I\!V}\label{section: type 4}

In this section, we prove the following local version of the $p$-adic monodromy theorem (Theorem~\ref{thm: local main thm}) for a type I\!V point. 

\begin{theorem} \label{thm:horizontal_st_type4}
Let $X$ be a smooth rigid analytic curve over $K$ (viewed as an adic space over $\spa(K, \mathcal{O}_K)$) and let $\L$ be a de Rham $\Z_p$-local system on $X$. Let $\xi \in X_{\C_K}$ be a point of type I\!V. Then there exist a finite extension $L/K$, an open neighborhood $U$ of the image of $\xi$ in $X_L$, and a finite \'etale cover $f\colon  V \ra U$, such that $f^*(\L|_U)$ is horizontally semistable in the sense of \cite[Definition 4.13]{Shimizu_2}; in particular, $f^*(\L|_U)$ is pointwise semistable (i.e. semistable at every classical point in $V$). 
\end{theorem}

To prove Theorem~\ref{thm:horizontal_st_type4}, we first show that every $p$-adic differential equation is \'etale locally trivial around any type I\!V point. 

\begin{theorem} \label{thm:pDE_near_type_4_trivial}
    Let $\xi \in X_{\C_K}$ be a point of type I\!V and let $(M, \nabla)$ be a vector bundle with connection on a neighborhood of $\xi$. Then there exists an open neighborhood $U$ of $\xi$ and a finite \'etale cover $f\colon  V \ra U$ such that $f^* M$ is a trivial connection.   
\end{theorem}

Applying Theorem~\ref{thm:pDE_near_type_4_trivial} to $D_{\mathrm{dR}}(\L)$ (together with a descent argument), we reduce to the case where $\L$ is horizontally de Rham. Then Theorem~\ref{thm:horizontal_st_type4} follows from the work of Shimizu \cite{Shimizu_2} on potential semistability of horizontally de Rham local systems. 

Throughout the section, let $\omega=p^{-\frac{1}{p-1}}$.

\begin{comment}
\begin{corollary} \label{cor:pDE_near_type_4_trivial}
    Let $\xi \in X$ be a $K$-rational point of type \textup{I\!V} and let $(M, \nabla)$ be a vector bundle with connection on a neighborhood of $\xi$. Then, after replacing $K$ by a finite extension, there exists an open neighborhood $U$ of $\xi$ and a finite \'etale cover $f: V \ra U$ such that $f^* M$ is a trivial connection.   
\end{corollary}
\end{comment}

\vspace{0.1in}
%% 7.1
\subsection{Some preparations} \label{ss: reduction_of_type_4}
\noindent 
\vspace{0.1in}

\noindent 
We first make some elementary observations before proceeding to the proof of Theorem~\ref{thm:pDE_near_type_4_trivial}. By \cite[Proposition~2.2.1]{Berkovich_differential}, a type I\!V point on $X_{\C_K}$ admits a basis of open neighborhoods consisting of closed disks. Therefore, without loss of generality we may assume that $\xi$ is a type I\!V point on a closed unit disk with coordinate $T$. Also recall that $\mathbf{k}(\xi)$ denotes the (uncompleted) residue field at $\xi$, which coincides with the stalk $\mO_{X_{\C_K}, \xi}$, and $\mathscr{H}(\xi)$ stands for the completed residue field. Moreover, $\mathbf{k}(\xi)$ is a differential field with derivation $\partial = d/{dT}$. Note that the constant field $\mathbf{k}(\xi)^{\partial=0}$ is $\C_K$, and $\partial\colon  \mathbf{k}(\xi)\rightarrow \mathbf{k}(\xi)$ is surjective.\footnote{To see the surjectivity of $\partial$, suppose $a\in \mathbf{k}(\xi)$ is defined on some closed disk $U$ containing $\xi$ with coordinate $T$. Integrating $a$ with respect to $T$, we obtain an analytic function $b$ that converges on any strictly smaller closed disk containing $\xi$ and such that $\partial b=a$.} For later use, we fix an algebraic closure $\mathscr{K}(\xi)$ of $\mathbf{k}(\xi)$. When the type I\!V point $\xi$ is understood, we simply write $\mathscr K$ for $\mathscr K(\xi)$.

\begin{lemma} \label{lemma:partial_is_surjective}
    Let $\mathbf k/\mathbf{k}(\xi)$ be a finite extension. Then there exist an open neighborhood $U$ of $\xi$ and a finite \'etale cover $f\colon  V \ra U$ with a unique preimage $\xi'$ lying above $\xi$ such that $\mathbf k = \mO_{V, \xi'}$; in particular, $\xi'$ has type I\!V. Moreover, the derivation $\partial = d/dT$ uniquely extends to a derivation $\partial$ on $\mathbf k$ with constant field $\C_K$, and $\partial\colon  \mathbf k \ra \mathbf k$ is surjective. 
\end{lemma}

\begin{proof}
    The first statement is straightforward. The derivation uniquely extends to $\mathbf{k}$ because $\mathbf{k}/\mathbf{k}(\xi)$ is a finite separable extension. Pick a closed disk $V'\subset V$ containing $\xi'$ with coordinate $S$. Then we have $\partial=g \frac{d}{dS}$ where $g=\frac{dS}{dT}\in \mathbf{k}^{\times}$. The desired surjectivity of $\partial\colon\mathbf{k}\rightarrow \mathbf{k}$ follows from the surjectivity of $\frac{d}{dS}$ on $\mathbf{k}$.
\end{proof}

As an immediate consequence, $\partial=d/dT$ uniquely extends to a derivation $\partial\colon \mathscr{K}\rightarrow \mathscr{K}$. The following corollary is a special case of Theorem~\ref{thm:pDE_near_type_4_trivial} when the vector bundle has rank $1$. 

\begin{corollary} \label{cor:1_dim_differential_equation_all_solvable}
    Every rank $1$ differential module over $(\mathscr{K}, \partial)$ is trivial. 
\end{corollary}

\begin{proof}
    Let $e$ be a basis with $\nabla(e) = a e$ for some $a \in \sK$. It suffices to find $c \in \sK^\times$ such that $\partial c = ac$. Given such $c$, we have $\nabla(c^{-1} e) = - c^{-2} (ac) e + c^{-1} \cdot ae = 0$ as desired. 
    
    Let $\mathbf k/\mathbf k(\xi)$ be a finite extension that contains $a$. Applying Lemma~\ref{lemma:partial_is_surjective}, we obtain a finite \'etale cover $V \ra U$ and a type I\!V point $\xi'\in V$ such that $\xi'$ is the unique preimage of $\xi$ and $\mathbf{k} = \mO_{V, \xi'}$. By shrinking, we may assume that $V$ is a closed disk. Moreover, there exists $b \in \mathbf k$ such that $\partial b = a$, which we may assume to be an analytic function on $V$. Write $A = \mO(V)$, equipped with its spectral norm $|\cdot|_V$. Let $m$ be a sufficiently large integer such that $|p^m b|_V < \omega$. Then the exponential $u = \exp (p^m b)$ converges on $V$ and $u \in A^\times$. Note that the cover of $V$ given by $A[z]/(z^{p^m} - u)$ is finite \'etale by the Jacobian criterion. Therefore,  a root $z$ of $z^{p^m} = u$ in $\sK$ satisfies $\partial z = a z$ and is precisely the solution we are looking for.   
\end{proof}

\vspace{0.1in}
%% 7.2
\subsection{Differential modules over $(\mathbf{k}(\xi), \partial)$} \label{ss:type_4_prime_to_p_rank}
\noindent 
\vspace{0.1in}

\noindent 
Theorem~\ref{thm:pDE_near_type_4_trivial} boils down to potential triviality of differential modules over $(\mathbf{k}(\xi), \partial)$. To prove this, we need the following key technical result.

\begin{proposition} \label{prop:type_4_key}
Let $\xi\in X_{\C_K}$ be a point of type I\!V and let $V$ be an irreducible differential module over $(\mathbf{k}(\xi), \partial)$ of rank $n$ with trivial determinant (i.e., the determinant $\wedge^nV$ is a trivial differential module). Then either $p |n$ or $n =1$. 
\end{proposition}

\subsubsection{Radii of optimal convergence}\label{subsubsection: radii of optimal convergence}
\noindent
\vspace{0.05in}

\noindent
To prove this key proposition, let us recall some basic notions and results from the theory of $p$-adic differenial equations (for example, see \cite{kedlaya_book}). Just for the purpose of proving Proposition~\ref{prop:type_4_key}, we temporarily treat all disks and annuli below as Berkovich spaces instead of as adic spaces. This is mainly because the theory of $p$-adic differential equations are most developed over Berkovich spaces in standard literature. Note that for a type I\!V point $\xi$ on a rigid analytic curve $X$, its Berkovich stalk $\mO_{X^{\mathrm{Ber}},\xi}$ coincides with its adic stalk $\mO_{X^{\mathrm{ad}},\xi}$ (cf. \S \ref{subsection: comparison_between_stalks}). In particular, there is no ambiguity in writing $\mathbf{k}(\xi)$ and $\mathscr{H}(\xi)$. 

For the setup, let $(M, \nabla)$ be a differential module of rank $n$ over the closed disk $\mathbb{D}_{\C_K} (0, r_0)$ of radius $r_0\in p^{\Q}$ with coordinate $T$. For any Berkovich point $z \in \mathbb D_{\C_K} (0, r_0)$ (that is, a point of type I, I\!I, I\!I\!I, or I\!V), define 
\[M_{z,  0} \coloneqq M \otimes_{\mO(\D_{\C_K} (0, r_0))} \sH (z) [\![T - t_z]\!]\] 
where $t_z$ denotes the image of $T$ in $\sH(z)$. This is a differential module over $(\sH (z) [\![T - t_z]\!], \frac{d}{d(T- t_z)}= \frac{d}{dT})$ with full (formal) horizontal sections; namely, 
\[(M_{z,  0} )^{\nabla = 0 } \otimes_{\sH(z)} \sH (z) [\![T - t_z]\!] \isom M_{z,  0}.\] For any $\rho \in (0, r_0]$, let 
\[M_{z,  \rho} \coloneqq 
 M \otimes_{\mO(\D_{\C_K} (0, r_0))} \sH (z) \gr{\frac{T - t_z}{\rho}}\] 
 and let  
$\mathrm{Sol}_z (\rho) \coloneqq (M_{z,  0} )^{\nabla = 0 } \cap M_{z, \rho}$ denote the $\sH(z)$-vector space of horizontal sections of $\nabla$ that converge on the closed disk of radius $\rho$. 

\begin{definition}
For any Berkovich point $z\in \D_{\C_K}(0,r_0)$, the \emph{radii of optimal convergence}\footnote{The normalization we use here is slightly different from \cite[Definition 11.9.1]{kedlaya_book}. Also note that there is a typo in \emph{loc.\ cit.}: $n - i$ should be $n + 1 -i$.}  $0 < s_1 (M, z) \le \cdots \le s_n (M, z) \le r_0 $ are defined by 
\[
s_i (M, z) \coloneqq \mathrm{sup} \{\rho < r_0 \: : \: \dim_{\sH(z)}   \mathrm{Sol}_z (\rho)  \ge n + 1 -i \}.
\]
In other words, $s_i (M, z)$ is the radius of the maximal open disk around $t_z$ on which there exist $n-i+1$ linearly independent horizontal sections of $M$. 
\end{definition}

\begin{definition}
For any $\mu\in \C_K$ with $|\mu|\le r_0$ and any $0<\rho\le r_0$, let $x_{\mu, \rho}$ denote the point of type I\!I or I\!I\!I corresponding to the disk $\D_{\C_K}(\mu, \rho)$ of radius $\rho$ centered at $\mu$ (see \S \ref{subsection: points on closed unit disks}). We put 
\[
R_i (M, x_{\mu, \rho}) = \min \{\rho, s_i(M, x_{\mu, \rho})\}, \qquad 
I\!R_i (M, x_{\mu, \rho}) = \min \{1, s_i(M, x_{\mu,\rho})/\rho\}
\]
for $i=1, \ldots, n$. 
\end{definition}

\begin{remark}
The quantities $R_i (M, x_{\mu, \rho})$ and $I\!R_i (M, x_{\mu, \rho})$ are closely related to the notions of \emph{generic subsidiary radii of convergence} and \emph{intrinsic generic subsidiary radii of convergence} (cf. \cite[Definition~9.4.4, 9.4.7, \& 9.8.1]{kedlaya_book}). Indeed, let $\mathscr{H}_{\mu, \rho}$ denote the completed residue field at $x_{\mu, \rho}$ and let $M_{\mu, \rho}$ be the base change of $(M,\nabla)$ to $\mathscr{H}_{\mu, \rho}$. Then $R_i (M, x_{\mu, \rho})$ (resp. $I\!R_i (M, x_{\mu, \rho})$) is precisely the generic subsidiary radii of convergence $R_i(M_{\mu, \rho})$ (resp. the intrinsic generic subsidiary radii of convergence $I\!R_i(M_{\mu, \rho})$) defined in \emph{loc.\ cit.}
\end{remark}

The following local constancy result will be a key ingredient of the proof of Proposition~\ref{prop:type_4_key}.

\begin{theorem} \label{thm:local_constancy_radii_type4}
    Let $\xi \in \D_{\C_K} (0, r_0)$ be a type I\!V point and let $(M, \nabla)$ be a differential module of rank $n$ over $\D_{\C_K} (0, r_0)$. Then there exists a closed disk $U\subset \D_{\C_K} (0, r_0)$ containing $\xi$ such that $s_i (M, z)$ is a constant function on $U$ for all $i=1,\ldots, n$, where the variable $z$ ranges over all Berkovich points on $U$. 
\end{theorem}

\begin{proof}
    This is \cite[Lemma~4.5.14]{kedlaya_local_global} (also see \cite{Kedlaya_Shiho}). 
    %Also see \cite{Pulita}. 
\end{proof}

\vspace{0.05in}
\subsubsection{Some reduction steps}\label{subsubsection: some reduction steps}
\noindent
\vspace{0.05in}

\noindent
To prove Proposition~\ref{prop:type_4_key}, let us start with some reduction steps.

First, we may assume that the differential module $V$ is represented by a finite free irreducible differential module $(M, \nabla)$ over a closed disk $\mathbb{D}_{\C_K} (0, r_0)$ centered at $0$ with radius $r_0\in p^{\Q}$, over which the determinant of $M$ is trivial. By Theorem~\ref{thm:local_constancy_radii_type4}, we may further assume that all $s_i(M, z)$'s are constant on $\mathbb{D}_{\C_K} (0, r_0)$, which we simply denote as $s_1, \ldots, s_n$.   

Let $r \coloneqq r(\xi)$ denote the radius of the type I\!V point $\xi$. Note that $r < r_0$. We may further assume that $s_1 \le r$. Otherwise, there exists a closed disk $\mathbb{D}_{\C_K}(a, d)$ centered at some point $a$ with radius $d$ satisfying $r < |T - a|_{\xi} < d < s_1$, where $|\cdot|_{\xi}$ is the continuous valuation associated with $\xi$. This disk contains the type I\!V point $\xi$ and $M$ admits a full set of horizontal solutions over $\mathbb{D}_{\C_K}(a, d)$. The irreducibility then forces $n = 1$ and there is nothing to prove. From now on, we assume that $s_1 \le r$.  

\begin{lemma} \label{lemma:annulus_splitting}
    We have $s_1 = s_2 = \cdots =s_n$.
\end{lemma}

\begin{proof}
   Suppose $s_1 = \cdots = s_j < s_{j+1}$ for some $j < n$. By assumption, we have $R_1 (M, x_{0,\rho}) = \cdots = R_j (M, x_{0,\rho}) = s_1$, and $R_{j+1} (M, x_{0,\rho}) = \min \{\rho, s_{j+1}\} > s_1$, for all $\rho\in (r,r_0)$. By \cite[Theorem~12.5.2]{kedlaya_book}, there exists a decomposition $M = M_1 \oplus M_2$ of differential modules over the open annulus $A_{\C_K} (r, r_0)$, where $M_1$ has rank $j$ while $M_2$ has rank $n - j$. This contradicts the irreducibility assumption. 
\end{proof}

Let $s$ denote the value $s_1=\cdots =s_n$. The following is an immediate consequence of Lemma~\ref{lemma:annulus_splitting}.

\begin{corollary}
For any $\mu\in \C_K$ with $|\mu|<r_0$ and any $\rho\in [s,r_0)$, we have $R_i(M, x_{\mu, \rho})=s$ for all $i=1, \ldots, n$.
\end{corollary}

Now, we pick $d \in p^\Q\cap (r, r_0)$ such that $d \notin \{ s \cdot \omega^{- 1/p^j} ~|~ j \in \Z_{\ge 0}\}$, and pick $a, c \in \C_K$ with $|c| = d$ and $|T-a|_{\xi} < d$.  Let $b = a + c$. Note that $|T-b|_{\xi} = d$. In other words, $\xi$ lies on the circle centered at $b$ with radius $d$. In fact, $\xi$ lies inside the open disk $\D_{\C_K}(a, < d)$. By a change of coordinates, we may assume that $b = 0$. Under this new coordinate, the open disk $\mathbb D_{\C_K} (a, < d)$ (and thus $\xi$) is contained in the circle $\{|T|=d\}$, and hence contained in the closed annulus $\{ \alpha \le |T| \le \beta  \}$ for any $\alpha, \beta\in p^{\Q}$ satisfying 
\[r < \alpha < d < \beta < r_0.\] In other words, we have 
\begin{equation} \label{eq:nesting_of_disks_and_annuli_for_type4_point}
\xi \in \mathbb{D}_{\C_K}(a, < d)\subset A_{\C_K}([d, d]) \subset A_{\C_K} ([\alpha,\beta]) \subset \D_{\C_K} (0, r_0)
\end{equation}
as illustrated in Figure~\ref{figure:nbhd_type4}. By scaling the coordinate $T$, we may assume that $d=1$ from now on. Note that the previous assumption on $d$ implies that $s\neq \omega^{1/p^j}$ for all integers $j\ge 0$. 

\begin{figure} \label{figure:nbhd_type4}
    \centering
  \scalebox{0.7}{
\begin{tikzpicture}[x=1cm,y=1cm]
 \path[use as bounding box] (-6.2,-3.55) rectangle (7.4,3.2);
 \coordinate (b) at (-1,0);
 \coordinate (residuecenter) at (-2.64,1.1);
 \coordinate (a) at (-2.64,1.10);
 \coordinate (x) at (-2.62,1.36);
 \fill[ann!10,even odd rule] (b) circle (2.92) (b) circle (1.03);
 \fill[level!12,even odd rule] (b) circle (2.52) (b) circle (1.43);
 \draw[draw=excluded,line width=1pt,dashed] (b) circle (3.5);
 \draw[draw=ann,line width=1.3pt,dashed] (b) circle (2.92);
 \draw[draw=ann,line width=1.3pt,dashed] (b) circle (1.03);
 \draw[draw=level,line width=1.25pt] (b) circle (2.52);
 \draw[draw=level,line width=1.25pt] (b) circle (1.43);
 \fill[fill=disk!15] (residuecenter) circle (.46);
 \draw[draw=disk,line width=1.3pt,dashed] (residuecenter) circle (.46);
 \node[dot] at (b) {};
 \node[anchor=north] at (-1,-.15) {\(0\)};
 \node[dot,fill=disk] at (a) {};
 \node[anchor=north east,inner sep=2pt] at (a) {\(a\)};
 \node[xdot] at (x) {};
 \node[anchor=west,text=point,font=\normalsize,inner sep=4pt] at (x) {\(\xi\)};
 \node[text=disk,anchor=west,align=left] at (-7.2,2.2)
      {\small The open disk \(\mathbb D(a, < d)\)\\ \small   containing $\xi$\\ };
 \draw[-{Stealth[length=2mm]},draw=disk,line width=.8pt]
      (-5.07,2.26)--(-3.04,1.36);
 \node[anchor=west,text=excluded] at (2.92,3.1) {\small The closed disk  \(\mathbb{D}_{\C_K}(0,r_0)\)};
 \draw[draw=excluded,line width=.7pt] (2.79,3.04)--(1.47,2.47);
 \node[anchor=west,text=ann,align=left] at (2.92,1.85)
      {\small The outer boundary \(|T|=\beta\)};
 \draw[draw=ann,line width=.7pt] (2.78,1.88)--(1.78,.89);
 \node[anchor=west,text=level,align=left] at (2.92,.2)
      {\small The closed circle (amplified in picture) \\ \small  \( A_{\C_K}([d, d])\coloneqq \{|T|=d\}\)};
 % The leader ends inside the band representing C_b(d).
 \draw[draw=level,line width=.7pt] (2.78,.2)--(.96,.2);
 \node[anchor=west,text=ann,align=left] at (2.92,-1.5)
      {\small The inner boundary \(|T|=\alpha\)};
 \draw[draw=ann,line width=.7pt] (2.78,-1.5)--(-.21,-.66);
 \node[anchor=west,text=ann,align=left] at (2.92,-2.98)
      {\small The closed annulus \( A_{\C_K} ([\alpha,\beta])\)};
 \draw[draw=ann,line width=.7pt] (2.78,-2.98)--(.64,-2.18);
\end{tikzpicture}
}%scalebox
    \caption{A Cartoon of annuli and disks containing $\xi$}
    \label{fig:placeholder}
\end{figure}

Let us make one final reduction before proving Proposition~\ref{prop:type_4_key}.  

\begin{lemma} \label{lemma:type_4_key_prop_reduction3}
It suffices to prove Proposition~\ref{prop:type_4_key} under the additional assumption that $s<\omega$.
\end{lemma}

\begin{proof}
    Let $h \ge 0$ be the least integer such that $s^{p^h} < \omega$. If $h = 0$, there is nothing to prove, so we assume $h > 0$. Note that $s^{p^h} < \omega$ and $s^{p^{h - 1}} > \omega$ (here we use $s\neq \omega^{1/p^j}$ for all $j\ge 0$). In particular, $s > \omega$. Pick an interval $[\alpha, \beta]$ as above so that 
    \begin{equation} \label{eq:choice_of_integer_h}
            (s/\rho)^{p^h} < \omega \textup{ and } (s/\rho)^{p^{h - 1}} > \omega \quad  \textup{for all } \rho \in [\alpha, \beta].
    \end{equation}
The idea, which dates back to the works of Christol and Dwork, is to apply Frobenius antecedents to make the generic radii visible. More precisely, consider the ``Froebnius map''
\[\phi\colon  \{\alpha \le |T| \le \beta\} \ra \{\alpha^p \le |T_1| \le \beta^p \}\]
given by $T_1 \mapsto T^p$. By \cite[Theorem~10.4.4, 10.4.2]{kedlaya_book}, there exists a Frobenius antecedent $M_1$ of $M$ along $\phi$; that is, a differential module $M_1$ over the annulus $\{ \alpha^p \le  |T_1| \le \beta^p \}$ such that $M$ is the pullback of $M_1$ along the Frobenius map $\phi$. For this new annulus, the open disk that contains the type I\!V point $\phi(\xi)$ becomes $\mathbb D_{\C_K}(\phi(a), < 1)$, while the radii of optimal convergence all become $s^p$. Hence, for any $\rho_1 \in [\alpha^p, \beta^p]$, we have $I\!R_i(M_1, x_{0,\rho_1})=s^p/\rho_1$ for all $i=1, \ldots, n$, and we have $(s^p/\rho_1)^{p^{h-1}} < \omega$ and $(s^p/\rho_1)^{p^{h-2}} > \omega$ by \eqref{eq:choice_of_integer_h}. Moreover, $M_1$ is irreducible and has trivial determinant. The latter follows from the uniqueness of Frobenius antecedent.  
Applying this construction $h$-times, we obtain $M_h$ (which is irreducible with trivial determinant) over the annulus $\{\alpha^{p^h} \le |T_h| \le \beta^{p^h} \}$ with all radii of optimal convergence equal to $s^{p^h}<\omega$. This finishes the proof of the lemma.
\end{proof}

From now on, we assume $s<\omega$. Pick $r<\alpha<d=1<\beta<r_0$ as above such that $s/\alpha <\omega$. In particular, we have $s/\rho < \omega$ for all $\rho \in [\alpha, \beta]$. 

\vspace{0.05in}
\subsubsection{Proof of Proposition~\ref{prop:type_4_key}}\label{subsubsection: proof of 7.5}
\noindent
\vspace{0.05in}

\noindent
We are ready to prove Proposition~\ref{prop:type_4_key}. Thanks to the reduction steps, we are in the following setup: 
\begin{itemize}
\item $(M, \nabla)$ is a finite free differential module over $\mathbb{D}_{\C_K}(0, r_0)$ of rank $n$;
\item For each $i=1, \ldots, n$, the radii of optimal convergence $s_i(M,z)$ are constant on $\mathbb{D}_{\C_K}(0, r_0)$, and are all equal to some $s<\omega$;
\item For any $\mu\in \C_K$ with $|\mu|<r_0$ and any $\rho\in [s, r_0)$, we have $R_i(M, x_{\mu, \rho})=s$ for all $i=1, \ldots, n$;
\item We have $r<\alpha<1<\beta<r_0$ with $\alpha, \beta\in p^{\Q}$ such that $s/\rho<\omega$ for all $\rho \in [\alpha, \beta]$, and such that
\[\xi \in \mathbb{D}_{\C_K}(a, < 1)\subset A_{\C_K}([1, 1]) \subset A_{\C_K} ([\alpha,\beta]) \subset \D_{\C_K} (0, r_0)\]
as illustrated in Figure (\ref{figure:nbhd_type4}).
\end{itemize}
Let us assume that $n>1$ and $p \nmid n$. Our goal is to find a contradiction. The proof is organized in several steps.

\vspace{0.1in}
\noindent \textbf{Step 1}.

\vspace{0.05in}

\noindent  
Consider the base change $(M_F, \nabla)$ of $(M,\nabla)$ to the fraction field $F$ of the ring of analytic functions on the annulus $\{ \alpha \le |T| \le \beta \}$. By \cite[Theorem~5.4.2]{kedlaya_book}, the differential module $(M_F, \nabla)$ admits a cyclic vector $v$. Write $\nabla^n (v) = \sum_{i = 0}^{n-1} c_i \cdot \nabla^i (v)$, where each $c_i = c_i (T)$ is a meromorphic function on the annulus, and let $P (X) \coloneqq   X^n - \sum_{i=0}^{n-1} c_i \cdot X^i$ be the characteristic polynomial of $v$. For any $\mu\in \C_K$ with $|\mu|<r_0$ and any $\rho\in [\alpha, \beta]$, let $\sH_{\mu, \rho}=\sH(x_{\mu, \rho})$ denote the completed residue field of $x_{\mu, \rho}$ (equipped with the continuous valuation $|\cdot|_{\mu,\rho}$) and let $M_{\mu, \rho}$ denote the base change of $M$ to $\sH_{\mu, \rho}$. By construction, $v$ is also a cyclic vector of $M_{\mu, \rho}$, again with characteristic polynomial $P(X)$. Let $\beta_{i, \mu, \rho}$ denote the roots of $P(X)$ in an algebraic closure of $\sH_{\mu,\rho}$, for $i=1,\ldots, n$. By the classical result of Christol and Dwork (cf. \cite[Theorem~6.5.3, Corollary~6.5.4]{kedlaya_book}), we have $|\beta_{i, \mu,\rho}|_{\mu,\rho} = \omega/s$ for all $i$ (in particular, the $|\cdot|_{\mu, \rho}$-norm of $\beta_{i, \mu, \rho}$ are independent of $\rho$). Note that the entire spectrum of $M_{\mu,\rho}$ is visible in the sense of \cite[Definition~6.5.1]{kedlaya_book} because $\omega/s > 1/\rho$.

Note that $\omega/s \in p^{\Q}$. Pick $\gamma \in \C_K^\times$ with $|\gamma|= \omega/s$ and write $b_i \coloneqq \gamma^{-(n-i)} \cdot c_i\in F$. Consider a new polynomial \[\sq P(X) \coloneqq \gamma^{-n}P(\gamma X)= X^n - \sum_{i=0}^{n-1} b_i X^i\in F[X].\] Over $\sH_{\mu,\rho}$, it has roots $\gamma^{-1} \beta_{i, \mu, \rho}$, which all have norm $1$. In particular, we have $|b_i|_{\mu,\rho} \le 1$ for all $i=1, \ldots, n-1$, and $|b_0|_{\mu,\rho} = 1$. 

Now we consider the type I\!I point $x_{0,1}$ and the associated completed residue field $\sH_{0,1}=\sH(x_{0,1})$. It is a complete non-archimedean field with respect to the $1$-Gauss norm $|\cdot|_{0,1}$, with an imperfect residue field $\overline{k}(T)$ (recall that $k=k_K$ is the residue field of $K$). Then $b_i = b_i (T)$ lies in $\mO_{\sH_{0,1}}$ for all $0 \le i \le n-1$. For simplicity, let us write $\beta_i = \beta_{i, 0,1}$ for the roots of $P(X)$ over an algebraic closure of $\sH_{0,1}$, which we fix in the rest of the proof.  \\

\noindent \textbf{Sublemma A}. The image of $\sq P(X)$ under the reduction map $\mO_{\sH_{0,1}} [X] \ra \overline{k}(T)[X]$ lies in $\overline{k}[X]$. 

\begin{proof}[Proof of Sublemma A] Let $\cl \mu \in \overline{k}^{\times}$ be an element and let $\mu \in \mO_{\C_K}^\times$ be a lift of $\cl \mu$ (for example one may take $\mu = [\cl \mu]$). The residue disk %inside $\D_{\C_K} (0, 1)$ %(recall that we have normalized the radius so that $d = 1$) 
corresponding to $\cl \mu$ is $\D_{\C_K} (\mu, < 1)$. Inside this residue disk, we have type I\!I and I\!I\!I points $x_{\mu, \rho}$ with $\rho < 1$, corresponding to the closed disk $\D_{\C_K}(\mu, \rho) \subset \D_{\C_K}(\mu, < 1)$. Let $f$ be a meromorphic function on the annulus $\{\alpha \le |T| \le \beta\}$. Recall the number $s_{\cl \mu} (f)$ defined in \cite[Definition~11.1.1]{kedlaya_book} 
\[
s_{\cl \mu} (f) \coloneqq \frac{d}{dr} \vline_{ \raisebox{-0.24cm}{\tiny \;$r\ra0^+$}}\! \left(-\log |f|_{x_{\mu, p^{-r}}} \right).
\]
It is independent of the lift $\mu$ and measures the change in the norm of $f$ when ``entering'' into the residue disk $\D_{\C_K} (\mu, < 1)$. Moreover, if $|f|_{0,1} = 1$, we know that $s_{\cl \mu} (f) = \mathrm{ord}_{\cl \mu} (\cl f)$ is simply the order of the zeros of $\cl f \in \overline{k} (T)$ at $\cl \mu$ (cf. the proof of \cite[Proposition~11.1.2]{kedlaya_book}). There are similar definitions and statements for $\cl \mu = 0$ and $\infty$; see \emph{loc.\ cit.}

We apply the discussion above to the meromorphic functions $b_i$. Recall that $|b_i|_{0,1} \le 1$. If $|b_i|_{0,1} < 1$ then $\cl b_i = 0$ and there is nothing to prove. So we assume that $|b_i|_{0,1} = 1$. If $\cl \mu \in \overline{k}^{\times}$ is a pole of $\cl b_i$, then $s_{\cl \mu} (b_i)  = -m$ for some $m \ge 1$ and hence $|b_i|_{\mu, \rho} = \rho^{-m}>1$ for $\rho$ sufficiently close to $1^-$, which is a contradiction. Therefore $\cl b_i$ has no poles at any $\cl \mu \in \overline{k}^{\times}$. A similar argument implies that $\cl b_i$ has no poles at $0$ and $\infty$ (by considering $x_{0, \rho}$ as $\rho \ra 1^-$ and $x_{0, \rho'}$ as $\rho' \ra 1^+$, respectively). Therefore, each $\cl b_i$ must be a constant, as desired.
\end{proof}

\noindent \textbf{Step 2}.

\vspace{0.05in}

\noindent  
By Sublemma A, any root of the reduction of $\sq P(X)$ lies in $\cl{k}$. Pick a root $\cl \eta \in \cl{k}$ of the reduction of $\sq P(X)$, and let $\eta = [\cl \eta]$ be its Teichmuller lift in $\mO_{\C_K}^\times$. Let $\lambda = \eta \gamma$, which has norm $|\lambda| = |\gamma| = \omega/s$. By the construction, we know that $|\beta_i - \lambda|_{0,1} \le \omega/s$ for all $i=1, \ldots, n$. Moreover, we have $|\beta_{i_0} - \lambda|_{0,1} <  \omega/s$ for at least one of the indices $i_0$. \\

\noindent \textbf{Sublemma B}. In fact, we have $|\beta_i - \lambda|_{0,1} <  \omega/s$ for all $i = 1, \ldots  , n$.

\begin{proof}[Proof of Sublemma B] Suppose $|\beta_i - \lambda|_{0,1} = \omega/s$ for some indices $i$. Let $l$ be the number of such indices; in particular, $1 \le l \le n-1$. Consider the free rank-1 differential module $L_{-\lambda}$ over the annulus $A_{\C_K}([\alpha, \beta])$ with basis $e_{-\lambda}$ defined by $D (e_{-\lambda}) =  - \lambda e_{-\lambda}$. Then consider the differential module \[M_{-\lambda} \coloneqq M \otimes_{\Gamma_{\C_K}([\alpha, \beta])} L_{-\lambda}\] where $\Gamma_{\C_K}([\alpha, \beta])$ denotes analytic functions on the annulus $A_{\C_K}([\alpha, \beta])$ (cf. \S \ref{sss:notation_cyclo}). If we write $\nabla_{-\lambda}$ for the resulting differential on $M_{-\lambda}$, then $\nabla_{-\lambda} = \nabla - \lambda$. That is, for all $x \in M$, we have 
 \[\nabla_{-\lambda} (x \otimes e_{-\lambda}) = (\nabla (x) - \lambda x) \otimes e_{-\lambda}.\] From now on, we shall identify an element $x \in M$ with its image $x \otimes e_{-\lambda} \in M_{-\lambda}$. In particular, for the cyclic vector $v$ chosen previously, %(over the fraction field of the annulus), 
we have $P(\nabla_{-\lambda} + \lambda) v = 0$. Moreover, $\{v, \nabla_{-\lambda}(v), \cdots, \nabla^{n-1}_{-\lambda}(v)\}$ form a cyclic basis of $M_{-\lambda}$ with characteristic polynomial $P(X+ \lambda)$, whose roots in the algebraic closure of $\sH_{0,1}$ are $\beta_i - \lambda$. 

Since $|\lambda| = \omega/s$, $\exp (\lambda T)$ has convergence radius $s$ and thus $R_1(L_{-\lambda}, x_{0,\rho})  = \min \{\rho, \frac{\omega}{|\lambda|}\} = s$ for all $\rho \in [\alpha, \beta]$. This in turn implies that $R_1 (M_{-\lambda}, x_{0,\rho}) \ge s$ for all $\rho \in [\alpha, \beta]$ by \cite[Lemma~9.4.6]{kedlaya_book}. We claim that $R_n (M_{-\lambda}, x_{0,1}) > s$. Otherwise, we must have 
\[R_1(M_{-\lambda}, x_{0,1}) = \cdots = R_n (M_{-\lambda}, x_{0,1}) = s.\]
Since the entire spectrum of $M_{-\lambda}$ is visible at $x_{0,1}$, by the result of Christol and Dwork again (cf. \cite[Theorem~6.5.3, Corollary~6.5.4]{kedlaya_book}), the roots $\beta_i - \lambda$ of $P(X+\lambda)$ all have norms $|\beta_i - \lambda|_{0,1} = \omega/s$. This contradicts the fact that $|\beta_{i_0} - \lambda|_{0,1} < \omega/s$ for at least one of the indices $i_0$. 

Consequently, up to shrinking the interval $[\alpha, \beta]$ if necessary, we have $R_n (M_{-\lambda}, x_{0,\rho}) > s$ for all $\rho \in [\alpha, \beta]$. But this would imply that $M_{-\lambda}$ (and thus $M$) is not irreducible, by the same splitting argument in the proof of Lemma~\ref{lemma:annulus_splitting} (using  \cite[Theorem~12.5.2]{kedlaya_book}). This finishes the proof of the lemma. 
\end{proof}

\noindent \textbf{Step 3}.

\vspace{0.05in}

\noindent  
Now we finish the proof of Proposition~\ref{prop:type_4_key}. By Sublemma B, we have
\[
\big|c_{n-1}-n\lambda\big|_{0,1}=\big|\sum_{i = 1}^n \beta_i - n \lambda \big|_{0,1} < \omega/s.
\]
Since $|n\lambda| = |\lambda| = \omega/s$ (this is where we use $p \nmid n$), we deduce that $|c_{n-1}|_{0,1} = \omega/s > 1$. Now, consider the element $w^{\det} \coloneqq v \wedge \nabla(v) \wedge \cdots \wedge \nabla^{n-1} (v) $ in the determinant of $M_{0,1}=M\otimes_{\Gamma_{\C_K}([\alpha, \beta])} \sH_{0,1}$. Differentiating $w^{\det}$ yields $\nabla (w^{\det}) = c_{n-1} w^{\det}$. By assumption, $M_{0,1}$ has trivial determinant, so there exists some element $\delta \in \sH_{0,1}^\times$ such that $\delta w^{\det}$ is horizontal; namely, $\nabla(\delta w^{\det}) = (\frac{d}{dT}\delta) w^{\det} + \delta\cdot c_{n-1} w^{\det} = 0$. This implies $c_{n-1} = - (\frac{d}{dT} \delta) \cdot \delta^{-1}$. However, $|\frac{d}{dT} \delta|_{0,1} \le |\delta|_{0,1}$ for any $\delta\in \sH_{0,1}^{\times}$. This forces $|c_{n-1}|_{0,1} \le 1$, and we find the desired contradiction.

\vspace{0.1in}
%% 7.3
\subsection{Proof of Theorem~\ref{thm:pDE_near_type_4_trivial}} \label{ss:finish_the_proof_of_trivial_connections}
\noindent 
\vspace{0.1in}

\noindent 
We will reduce Theorem~\ref{thm:pDE_near_type_4_trivial} to Proposition~\ref{prop:type_4_key} using a similar strategy as in \cite[\S 21.5]{kedlaya_book}. Let us start with the following observation on differential modules over $\sK$. Recall that $\sK=\sK(\xi)$ is a fixed algebraic closure of $\mathbf{k}(\xi)$.

\begin{lemma} \label{lemma:trivial_subquotient}
Let $V$ be a differential module over $(\sK, \partial)$ that admits a subquotient of dimension $1$. Then $V$ admits a nontrivial horizontal section. 
\end{lemma}

\begin{proof}
  By assumption, there exist differential submodules $V'' \subset V' \subset V$ such that $V'/V''$ is one dimensional. By Corollary~\ref{cor:1_dim_differential_equation_all_solvable}, $V'/V''$ is necessarily a trivial differential module. As in \S \ref{ss: reduction_of_type_4}, there exist an open neighborhood $U$ of $\xi$ and a finite \'etale cover $f\colon  U' \ra U$ with a unique preimage $\xi'$ lying above $\xi$ such that $V'' \subset V' \subset V$ descend to differential modules $M''\subset M'\subset M$ on $U'$. We may assume that $U'=\mathbb{D}_{\C_K} (0, r_0)$ is a closed disk, with coordinate $T$. Let $r = r(\xi')$ be the radius of the type I\!V point $\xi'$. Replacing $U'$ by a smaller disk if necessary, we may further assume that all radii of optimal convergence are constant on $U'$. 
  
  Pick $a \in \C_K$ and $d \in p^\Q$ such that $r < |T - a|_{\xi'} < d < r_0$. Then we have $\xi' \in \mathbb D_{\C_K} (a, < d) \subset U' = \mathbb D_{\C_K} (0, r_0)$. At the type I\!I point $x_{a, d}$ (namely, the Gauss point of $\D_{\C_K} (a, d)$), the existence of the trivial subquotient $M'/M''$ implies $I\!R_n (M, x_{a, d}) = 1$ by \cite[Lemma~1.2.9]{Kedlaya-Xiao}. This implies $s_n (M, x_{a, d}) \ge  d$.
  By the constancy of radii of optimal convergence, we must have $s_n(M, z) \ge d$ for all Berkovich points $z \in U'$. In particular, we have $s_n (M, a) \ge d$. Now, we pick $d'\in p^{\Q}$ such that $|T- a|_{\xi'} < d' <d$ so that $\D_{\C_K} (a, d')$ is strictly contained inside $\D_{\C_K} (a, < d)$ and $\xi'\in \D_{\C_K} (a, d')$. Then $s_n (M, a) \ge d$ implies that at least one horizontal section of $M$ converges on the entire disk $\D_{\C_K} (a, d')$. This provides the desired horizontal section.   
\end{proof}

Given this, we can prove the following upgraded version of Proposition~\ref{prop:type_4_key}. 

\begin{proposition} \label{prop:irreducible_implies_trivial!!!}
    Every irreducible differential module over $(\sK, \partial)$ is trivial (and in particular has dimension 1). 
\end{proposition}

\begin{proof}
Let $(V, \nabla)$ be an irreducible differential module over $(\sK, \partial)$ of rank $n$. First, suppose that $p \nmid n$. Note that the determinant of $V$ is trivial by 
Corollary~\ref{cor:1_dim_differential_equation_all_solvable}. Then we can descend $V$ to a differential module over $\mathbf{k} = \mathbf{k}(\xi')$ with trivial determinant, where $\xi'$ is a type I\!V point lying above $\xi$ as before. Applying Proposition~\ref{prop:type_4_key}, we conclude that $n=1$ and $V$ is trivial.

For the rest of the proof, we assume $p | n$. Our goal is to find a contradiction.  Consider $W\coloneqq \End_{\sK}(V)$ equipped with the differential $\nabla_{W} (f) = \nabla \circ f - f \circ \nabla$. By choosing a basis of $V$, we can write $\nabla = \partial + H$ for some $H\in \mathrm{Mat}_{n\times n}(\sK)=\End_{\sK} (V)$. For any matrix $B \in \mathrm{Mat}_{n\times n}(\sK)=\End_{\sK} (V)$, we have $\nabla_W (B) = \partial B + HB - BH$. Taking trace, we obtain $\mathrm{tr} (\nabla_W (B)) = \partial (\mathrm{tr} B)$. In particular, the trace-zero subspace $W_0 \coloneqq \ker (\mathrm{tr}\colon  \End_\sK(V) \ra \sK)$ is a differential submodule of $W$ of dimension $n^2-1$. Since $p$ is coprime to $n^2 - 1$, there exists an irreducible subquotient of $W_0$ with rank coprime to $p$. By the previous paragraph, this subquotient is trivial of dimension $1$. By Lemma~\ref{lemma:trivial_subquotient}, $W_0$ has a nontrivial horizontal section under $\nabla_W$. In other words, we have found a (nonzero) traceless endomorphism $B\colon  V \ra V$ that commutes with $\nabla$. In the matrix form (under the basis chosen above), we have $\partial B = B H - H B$. By induction, we deduce $\partial (B^i) = B^i H - H B^i$ for all $i \ge 1$, and hence $\partial (\mathrm{tr} (B^i)) = \mathrm{tr} (\partial (B^i)) = 0$. Therefore, the characteristic polynomial of $B$ has coefficients in $(\sK)^{\partial = 0}=\C_K$ by Lemma~\ref{lemma:partial_is_surjective}. In particular, all eigenvalues of $B$ lie in $\C_K$. Let $\lambda \in \C_K$ be one of the eigenvalues of $B$ and consider $V_{\lambda} \coloneqq \ker (B - \lambda I)$. It is clear that $B-\lambda I \ne 0$ as $B$ is traceless. Since $B$ commutes with $\nabla$, we know that $V_{\lambda} \subset V$ is a proper differential submodule. This contradicts the irreducibility of $V$ and finishes the proof of the proposition. 
\end{proof}

Now we are ready to prove Theorem~\ref{thm:pDE_near_type_4_trivial}. 

\begin{proof}[Proof of Theorem~\ref{thm:pDE_near_type_4_trivial}]
Let $M_{\mathbf{k}(\xi)}$ denote the restriction of $(M, \nabla)$ to $\mathbf{k}(\xi)$ and let $M_{\sK}$ denote its base change to $\sK=\sK(\xi)$. It suffices to show that $(M_{\sK}, \nabla)$ is a trivial differential module over $(\sK, \partial)$ (then the desired statement follows from Lemma~\ref{lemma:partial_is_surjective}). By Proposition~\ref{prop:irreducible_implies_trivial!!!}, $M_{\sK}$ admits a composition series $0 = M_0 \subset M_1 \subset \cdots \subset M_n = M_{\sK}$ where each quotient $M_i/M_{i-1}$ is a trivial differential module of dimension $1$. We need to show that this series splits. 

We prove this by induction. Suppose that $M_{i-1}$ has horizontal basis $e_1, \ldots  , e_{i-1}$. Let $\cl e_i \in M_i/M_{i-1}$ be a horizontal section under $\nabla$, and let $\sq e_i \in M_i$ be any lift of $\cl e_i$. Then we have 
$\nabla(\sq e_i) = \sum_{j = 1}^{i -1} a_j e_j$ with $a_j \in \sK$. By Lemma~\ref{lemma:partial_is_surjective}, for each $j =1, \ldots  , i - 1$, there exists $c_j \in \sK$ with $\partial c_j = a_j $. Then $e_i \coloneqq \sq e_i - \sum_{j =1}^{i-1} c_j e_j$ satisfies $\nabla (e_i) = 0$, as desired. This finishes the proof.
\end{proof}

\vspace{0.1in}
%%% 7.4
\subsection{Proof of Theorem~\ref{thm:horizontal_st_type4}} \label{ss:horizontal_semistable!}
\noindent 
\vspace{0.1in}

\noindent
To prove Theorem~\ref{thm:horizontal_st_type4}, we need a version of Theorem~\ref{thm:pDE_near_type_4_trivial} over some finite extension of $K$ rather than $\C_K$.

\begin{lemma} \label{lemma:Galois_descent}
Let $X=\mathrm{Sp}(A)$ be a smooth geometrically connected affinoid rigid analytic curve over $K$ and let $(M, \nabla)$ be a vector bundle with connection on $X$ whose global section is a finite free $A$-module. If $M_{\C_K}$ is a trivial connection, then so is $M$. 
\end{lemma}

\begin{proof}
    Let $A_{\C_K} = \mO(X_{\C_K}) = A \widehat \otimes_{K} \C_K$. By picking a basis of $M$ over $A$, we can write $\nabla = d + \Omega$ for some $\Omega\in \mathrm{Mat}_{n\times n}(A)$. Finding a horizontal basis of $M_{\C_K}$ over $A_{\C_K}$ is equivalent to finding a matrix $F \in \GL_n (A_{\C_K})$ such that $d F + \Omega F = 0$. Note that any two such matrices differ by a matrix in $\GL_n(\C_K)$.

    Choose any classical point $x_0 \in X$ defined over a finite Galois extension $L/K$. Write $A_L = \mO(X_L) = A \otimes_K L$. The evaluation map $A_{\C_K} \ra \C_K$ at the point $x_0$ yields a matrix $F(x_0) \in \GL_n (\C_K)$. Let $F_0 \coloneq F \cdot F (x_0)^{-1}$. Then we have $d F_0 + \Omega F_0 = 0$ and $F_0 (x_0) = I_n$ is the identity matrix. These two conditions uniquely characterize $F_0$. For any $g \in \Gal_L$, since $g$ fixes $\Omega$ and commutes with the differential $d$, we again have $d g(F_0) + \Omega g(F_0) = 0$ and $g(F_0)(x_0)=I_n$. The uniqueness then implies $g (F_0) = F_0$ for all $g \in \Gal_L$. By Corollary~\ref{cor: Galois orbit}, we have $(A_{\C_K})^{\Gal_L} = A_L$ and hence $F_0 \in \GL_n (A_L)$. In other words, the base change $M_L$ is a trivial connection. Finally, the horizontal sections $(M_L)^{\nabla = 0}$ of $M_L$ descend to those of $M$ by Galois descent. This finishes the proof of the lemma. 
\end{proof}

\begin{remark}
The freeness assumption on $M$ in Lemma~\ref{lemma:Galois_descent} is redundant. This simplified version is already sufficient for our purpose, and we leave the more general situation to the interested reader. 
\end{remark}

\begin{corollary}\label{cor:pDE_near_type_4_trivial--over K}
Let $X$ be a smooth rigid analytic curve over $K$ and let $\xi\in X_{\C_K}$ be a point of type I\!V. Let $(M, \nabla)$ be a vector bundle with connection defined on a neighborhood of the image of $\xi$ in $X$. Then there exist a finite extension $L/K$, an open neighborhood $U$ of the image of $\xi$ in $X_L$, and a finite \'etale cover $f\colon  V \ra U$ such that $f^* M$ is a trivial connection. 
\end{corollary}

\begin{proof}
This follows immediately from Theorem~\ref{thm:pDE_near_type_4_trivial} and Lemma~\ref{lemma:Galois_descent}.
\end{proof}

Finally, we finish the proof of Theorem~\ref{thm:horizontal_st_type4}.
 
\begin{proof}[Proof of Theorem~\ref{thm:horizontal_st_type4}] 
    Applying Corollary~\ref{cor:pDE_near_type_4_trivial--over K} to $D_{\mathrm{dR}} (\L)$, we can find a finite extension $L/K$, an open neighborhood $U$ of the image of $\xi$ in $X_L$, and a finite \'etale cover $f\colon  V \ra U$ such that $D_{\mathrm{dR}} (f^*\L) = f^*(D_{\mathrm{dR}}(\L))$ is a trivial connection. In particular, $f^* \L$ is horizontally de Rham in the sense of \cite[Definition~9.1]{Shimizu_2}. After shrinking $U$ and $V$ and enlarging $L$ if necessary, we may assume that $V$ is isomorphic to a circle over $L$. By \cite[Theorem~1.9]{Shimizu_2}, after a further finite extension of $L$, $f^*\L$ becomes horizontally semistable, thus semistable at every classical point in $V$.
\end{proof}

\vspace{0.3in}
%%%  Section 8
\section{Around points of type V: special case}\label{section: type 5}

The goal of this section is to prove Theorem~\ref{thm: local main thm} for a special class of type V points --- those type V points lying on closed unit disks. Later in \S \ref{section: type V general}, we will deduce the general case of type V points from this.

For the setup, let $\xi_{\mathrm{I\!I}}$ be the (type I\!I) Gaussian point of the closed unit disk $\mathbb{D}_K=\D_K(0,1)$ over $K$ and let $\xi$ denote the type V point corresponding to the open unit disk $\D_K(0,<1)\subset \mathbb{D}_K$ (cf. \S \ref{subsection: points on closed unit disks}). Both $\xi_{\mathrm{I\!I}}$ and $\xi$ are $K$-rational points, and $\xi$ lives in the closure of $\xi_{\mathrm{I\!I}}$. 

Let $\delta \in p^\Q \cap (0, 1)$ be a real number and $f(T) = \prod_{i=1}^n (T- [\beta_i]) \in \mO_{K_0} [T]$ be a polynomial with $\beta_i \in \cl{k}_K^\times$. Consider the closed annulus $A_{K} ([\delta, 1]) = \{\delta \le |T| \le 1\}$ and let $X_{\delta, f} \subset A_{K} ([\delta, 1])$ be the rational subset obtained by removing finitely many residue disks centered at $[\beta_i]$ on the boundary circle $\{|T| = 1\}$ (cf. Notation~\ref{notation:closed unit disk removing open unit disks}). We equip $X_{\delta, f}$ with the toric chart $X_{\delta, f}\hookrightarrow A_{K} ([\delta, 1])\hookrightarrow \D_K$.
Moreover, let $X_{1,f}\coloneqq \{|T|=1\}\cap X_{\delta, f}$ denote the ``outer boundary'', which admits a smooth integral model over $W(k_K)$. Note that $X_{\delta, f}$'s form a basis of open neighborhoods of $\xi$ while $X_{1,f}$'s form a basis of open neighborhoods of $\xi_{\mathrm{I\!I}}$. When $f$ is clear from the context, we simply write $X_{\delta}$ for $X_{\delta, f}$ and write $X_1$ for $X_{1, f}$. 

The main result of this section is the following. 

\begin{theorem} \label{thm:main_type5_disk}
Let $\L$ be a de Rham $\Z_p$-local system on $X_{\delta}=X_{\delta,f}$ for some $\delta$ and $f$. Suppose that $\L$ is semistable on the ``outer boundary'' $X_{1,f}$. Then, up to replacing $K$ by a finite extension and shrinking $X_{\delta}$ (i.e. increasing $\delta$ and removing finitely many more residue disks on the outer boundary), there exists a finite \'etale cover $g\colon X' \ra X_{\delta}$ such that $g^* \L$ is pointwise semistable. 
\end{theorem}

Since $\L$ is already semistable on the outer boundary, we are mainly interested in the half-open annulus $A_K([\delta, 1); T)$. To prove Theorem~\ref{thm:main_type5_disk}, we consider the $p$-adic differential equation $\mathbf{N}_{\mathrm{dR}}(\L)$ over $\mathbf{C}^{\dagger}_{X_{\delta}}$, equipped with both arithmetic and geometric differential operators (cf. \S \ref{ss: annulus--pDE}). As we shall see, the semistability of $\L$ on the boundary $X_{1,f}$ will force $\mathbf{N}_{\mathrm{dR}}(\L)$ to be ``strongly unipotent'' with respect to the arithmetic differential operator $D_{\pi_K}$. For a more precise statement, consider the rational subset
\[V_{\delta, f} \coloneqq \{ \delta \le |[T^\flat]| \le 1, ~~|f([T^\flat])| = 1 \}\subset A_{K_0'}([\delta, 1]; [T^{\flat}])\]
obtained by removing finitely many residue disks on the circle $\{|[T^{\flat}]|=1\}$, and let $A_{\delta, f}$ denote the ring of analytic functions on $V_{\delta, f}$. From the explicit description of $\fC_{X_{\delta}}^r$ (Lemma~\ref{lem: puncture-C-description}), we have $\left(\mathbf{C}^r_{X_{\delta, f}}\right)^{\frac{d}{d\pi_K}=0}=A_{\delta^r, f}$. The key technical result is the following.

\begin{theorem} \label{thm:spread_arithmetic_unipotence}
Let $\L$ be a de Rham $\Z_p$-local system of rank $d$ on $X_{\delta}=X_{\delta,f}$ for some $\delta$ and $f$. Suppose that $\L$ is semistable on $X_{1,f}$. Then, up to enlarging $K$ and shrinking $X_{\delta}$, there exists $r>0$ such that the $\mathbf{C}_{X_{\delta}}^r$-module ${\mathbf{N}}^r_{\mathrm{dR}}(\L)$ is free and admits a basis $e_1, \ldots, e_d$ such that the matrix of $\pi_KD_{\pi_K}$ relative to this basis is upper-triangular nilpotent with entries in $A_{\delta^r, f}$. In other words, ${\mathbf{N}}^r_{\mathrm{dR}}(\L)$ is strongly unipotent with respect to the arithmetic differential operator $D_{\pi_K}$.
\end{theorem}

Given this, we may assume that the restriction of $\mathbf{N}_{\mathrm{dR}}(\L)$ on the half-open annulus $A_K([\delta, 1); T)$ is strongly unipotent with respect to $D_{\pi_K}$. Then, similarly to our treatment in \S \ref{section: type 3}, we can produce enough semistable periods on $A_K([\delta, 1); T)$ after pulling back to some finite \'etale cover. The proof of Theorem~\ref{thm:spread_arithmetic_unipotence} will occupy \S \ref{ss:type5_special_1} -- \S \ref{ss:proof_of_arithmetic_unipotence}, from which we deduce Theorem~\ref{thm:main_type5_disk} in \S \ref{ss:arithmetic_hs} -- \S \ref{ss:Katz_special_algebraization}.

\vspace{0.1in}
%%% 8.1
\subsection{Completed perfection at the type I\!I point} \label{ss:type5_special_1}
\noindent 
\vspace{0.1in}

\noindent 
For the rest of the section, let $\L$ be a de Rham $\Z_p$-local system of rank $d$ on $X_{\delta}=X_{\delta,f}$ as in the statement of Theorem~\ref{thm:main_type5_disk} and assume that $\L$ is semistable on $X_1=X_{1,f}$.

The first step is to consider the completed stalk $\sH(\xi_{\mathrm{I\!I}})$ at the type I\!I point, which is given by the completion of $K(T)$ under the 1-Gauss norm. Write $\sH(\xi_{\mathrm{I\!I}})^{(\mathrm{pf})} \coloneqq \cup_{m \ge 0} \sH(\xi_{\mathrm{I\!I}}) (T^{1/p^m})$ and let $\sH(\xi_{\mathrm{I\!I}})^{\mathrm{pf}}$ denote its completion under the $p$-adic topology. Then $\sH(\xi_{\mathrm{I\!I}})^{\mathrm{pf}}$ a complete discrete valuation field, with residue field being the perfection $k(T)^{\mathrm{perf}}$ of $k(T)$. Write
\[
\sH(\xi_{\mathrm{I\!I}})^{\mathrm{pf}}_{\infty} \coloneqq \bigcup_{n \ge 1} \sH(\xi_{\mathrm{I\!I}})^{\mathrm{pf}}(\zeta_{p^n}),
\]
and let $\sH'$ denote the maximal unramified extension of $W(k(T)^{\mathrm{perf}})[p^{-1}]$ inside $\sH(\xi_{\mathrm{I\!I}})^{\mathrm{pf}}_{\infty}$. Let $\fC_{\sH(\xi_{\mathrm{I\!I}})^{\mathrm{pf}}}^{\dagger}$ denote the imperfect period ring associated with the $p$-adic field $\sH(\xi_{\mathrm{I\!I}})^{\mathrm{pf}}$ in the classical theory. Using the notation in \S \ref{sss:notation_cyclo}, $\fC_{\sH(\xi_{\mathrm{I\!I}})^{\mathrm{pf}}}^{\dagger}$ is simply the Robba ring $\mathscr{R}_{\sH'}$ over $\sH'$, with coordinate $\pi_K$. Consider the restriction of $\L$ to the type I\!I point $\xi_{\mathrm{I\!I}}$, and then to $\sq{\xi}_{\mathrm{I\!I}}=\spa (\sH(\xi_{\mathrm{I\!I}})^{\mathrm{pf}}, \sH(\xi_{\mathrm{I\!I}})^{\mathrm{pf}, \circ})$. This gives rise to a Galois representation over $\Gal_{\sH(\xi_{\mathrm{I\!I}})^{\mathrm{pf}}}$, which we shall denote by $\L|_{\sq{\xi}_{\mathrm{I\!I}}}$.

\begin{lemma}\label{lemma: base change to tilde_xi_II}
 There is a natural isomorphism 
  \[\mathbf{N}_{\mathrm{dR}} (\L |_{X_1}) \otimes_{\mathbf{C}_{ X_1 }^{\dagger}} 
\mathbf{C}_{\sH(\xi_{\mathrm{I\!I}})^{\mathrm{pf}}}^{\dagger} \isom  \mathbf{N}_{\mathrm{dR}} (\L |_{\sq{\xi}_{\mathrm{I\!I}}}) \]
compatible with arithmetic differential operators $D_{\pi_K}$, where $\mathbf{N}_{\mathrm{dR}} (\L |_{X_1}) $ is the relative $p$-adic differential equation constructed in \cite{relative} (cf. \S \ref{ss: annulus--pDE}), while $ \mathbf{N}_{\mathrm{dR}} (\L |_{\sq{\xi}_{\mathrm{I\!I}}}) $ is the classical $p$-adic differential equation of Berger \cite{Berger-differential} (cf. \S \ref{sss:Berger_N_dR}). 
\end{lemma}

\begin{proof}
Let $A_{X_1}$ denote the ring of analytic functions on $X_1$. The map $A_{X_1} \rightarrow \sH(\xi_{\mathrm{I\!I}})^{\mathrm{pf}}$ together with a choice of a common geometric point (i.e. a choice of a map $\overline{A}_{X_1} \rightarrow \overline{\sH(\xi_{\mathrm{I\!I}})^{\mathrm{pf}}}$, cf. \cite[0BND]{stacks-project}) yields a continuous map of Galois groups $\Gal_{\sH(\xi_{\mathrm{I\!I}})^{\mathrm{pf}}} \rightarrow \pi_1^{\text{\'et}}(X_1)$, which induces a map $\Gamma_{\sH(\xi_{\mathrm{I\!I}})^{\mathrm{pf}}} \rightarrow \Gamma$. Furthermore, the map $A_{X_1} \rightarrow \sH(\xi_{\mathrm{I\!I}})^{\mathrm{pf}}$ induces natural maps $\widetilde{\mathbf{A}}_{X_1} \rightarrow \widetilde{\mathbf{A}}_{\sH(\xi_{\mathrm{I\!I}})^{\mathrm{pf}}}$ and $\mathbf{C}_{X_1}^{\dagger} \rightarrow \mathbf{C}_{\sH(\xi_{\mathrm{I\!I}})^{\mathrm{pf}}}^{\dagger}$ compatible with $\varphi$- and $\Gamma_{\sH(\xi_{\mathrm{I\!I}})^{\mathrm{pf}}}$-actions. Now, consider the natural isomorphism
\[
(\mathbb{L}|_{X_1}\otimes_{\Z_p} W(\overline{A}_{X_1}^{\wedge, \flat})) \otimes_{W(\overline{A}_{X_1}^{\wedge, \flat}))} W(\overline{\sH(\xi_{\mathrm{I\!I}})^{\mathrm{pf}}}^{\wedge, \flat}) \cong \mathbb{L}|_{\widetilde{\xi}_{\mathrm{I\!I}}}\otimes_{\Z_p} W(\overline{\sH(\xi_{\mathrm{I\!I}})^{\mathrm{pf}}}^{\wedge, \flat})
\]
compatible with $\varphi$- and $\Gal_{\sH(\xi_{\mathrm{I\!I}})^{\mathrm{pf}}}$-actions. Unwinding the construction of $D_{\mathrm{rig}}^{\dagger}(\cdot)$, we deduce a natural isomorphism
\[
D_{\mathrm{rig}}^{\dagger}(\mathbb{L}|_{X_1})\otimes_{\mathbf{C}_{X_1}^{\dagger}} \mathbf{C}_{\sH(\xi_{\mathrm{I\!I}})^{\mathrm{pf}}}^{\dagger} \cong D_{\mathrm{rig}}^{\dagger}(\mathbb{L}|_{\widetilde{\xi}_{\mathrm{I\!I}}})
\]
compatible with $\varphi$- and $\Gamma_{\sH(\xi_{\mathrm{I\!I}})^{\mathrm{pf}}}$-actions. The desired statement then follows from \cite[Th\'eor\`eme 5.10]{Berger-differential}.
\end{proof}

Since $\L|_{\sq{\xi}_{\mathrm{I\!I}}}$ is semistable, as an immediate consequence of Lemma~\ref{lemma: base change to tilde_xi_II}, we know that 
\[
\mathbf{N}_{\mathrm{dR}, \sq{\xi}_{\mathrm{I\!I}}} \coloneqq  \mathbf{N}_{\mathrm{dR}} (\L |_{X_1}) \otimes_{\mathbf{C}_{ X_1 }^{\dagger}} 
\fC_{\sH(\xi_{\mathrm{I\!I}})^{\mathrm{pf}}}^{\dagger}
\] 
is strongly unipotent with respect to $D_{\pi_K}$;\footnote{According to \cite[Proposition~5.2.6]{Kedlaya_coefficients}, there is no difference between unipotence and strongly unipotence.} that is, there exists a basis of $\mathbf{N}_{\mathrm{dR}, \sq{\xi}_{\mathrm{I\!I}}}$ over $\fC_{\sH(\xi_{\mathrm{I\!I}})^{\mathrm{pf}}}^{\dagger}$, with respect to which the differential $\Delta_{\pi_K} \coloneqq \pi_K D_{\pi_K}$ acts through an upper-triangular nilpotent matrix with coefficients in $\sH'=\left(\fC_{\sH(\xi_{\mathrm{I\!I}})^{\mathrm{pf}}}^{\dagger}\right)^{\frac{d}{d\pi_K}=0}$.

\vspace{0.1in}
%%% 8.2
\subsection{Arithmetic unipotence at the type I\!I point} \label{ss:type5_special_localize_at_type2}
\noindent 
\vspace{0.1in}

\noindent 
Let $L$ denote the completion of $K_0' ([T^\flat])$ with respect to the 1-Gauss norm and let $|\cdot|_L$ denote the resulting norm on $L$. The Robba ring $\mathscr{R}_L$ with coordinate $\pi_K$ (see \S \ref{sss:notation_cyclo}) can be naturally viewed as a subring of $\fC_{\sH(\xi_{\mathrm{I\!I}})^{\mathrm{pf}}}^{\dagger}$, and the map
$\mathbf{C}_{X_1}^{\dagger} \ra  \fC_{\sH(\xi_{\mathrm{I\!I}})^{\mathrm{pf}}}^{\dagger}$
factors as $\mathbf{C}_{X_1}^{\dagger} \ra \mathscr{R}_L \ra 
\fC_{\sH(\xi_{\mathrm{I\!I}})^{\mathrm{pf}}}^{\dagger}$, compatibly with differential operators $\frac{d}{d\pi_K}$.

\begin{lemma}\label{lemma:Galois_descent_of_unipotence}
Let $N_0$ be a finite projective differential module over $(\mathscr{R}_L, \frac{d}{d\pi_K})$. Let $\C_{L}$ be a completed algebraic closure of $L$ containing $\sH'$ and write $N_{\C_L} \coloneqq N_0 \otimes_{\mathscr{R}_L} \mathscr{R}_{\C_L}$.
If $N_{\C_L}$ is strongly unipotent with respect to $D_{\pi_K}$, then so is $N_0$.  
\end{lemma}

\begin{proof}
By assumption, $N_{\C_L}$ admits a basis $e_1, \ldots, e_d$ over $\mathscr{R}_{\C_L}$, with respect to which $\Delta_{\pi_K}=\pi_K D_{\pi_K}$ takes the form 
\[
\Delta_{\pi_K} = \pi_K \cdot \frac{d}{d \pi_K} + Z
\]
where $Z \in \mathrm{Mat}_{d\times d} (\C_L)$ is upper-triangular nilpotent. Consider the $\C_L$-vector space $\ker (\Delta_{\pi_K}^d) \subset N_{\C_L}$. We claim that $\ker (\Delta_{\pi_K}^d)$ is precisely the $\C_L$-span of $e_1, \ldots, e_d$. It suffices to show that, if $x = \sum_{i =1, \ldots, d} a_i e_i $ (with $a_i \in \mathscr R_{\C_L}$) lies in $\ker (\Delta_{\pi_K}^d)$, then $a_i \in \C_L$ for all $1 \le i \le d$.  We first show $a_d \in \C_L$. Write $a_d = \sum_{j \in \Z} b_j \pi_K^{j}$ with $b_j \in \C_L$, and write $\Delta_{\pi_K}^d (x) = \sum_{i} \sq a_i e_i$. Note that $\widetilde{a}_i = 0$ for all $i$ since $x \in \ker (\Delta_{\pi_K}^d)$. To proceed, we compute $\sq a_d$ in terms of the $a_i$'s. Since $Z(e_i)$ lies in the span of $\{e_1, \ldots, e_{i-1} \}$, we have  
\[
\sq a_d = (\pi_K \frac{d}{d \pi_K})^d (a_d) = \sum_{j \ne 0} j^d b_j \pi_K^j = 0.
\]
Thus, $b_j = 0$ whenever $j \neq 0$, and $a_d \in \C_L$. Next, since $a_d \in \C_L$, we have $x - a_d e_d \in \ker(\Delta_{\pi_K}^d)$. Applying the same argument as above inductively, we conclude that $a_i \in \C_L$ for all $i$, as desired. %In particular, $\ker(\Delta_{\pi_K}^d)$ has dimension $d$ and $\Delta_{\pi_K}$ acts on it via $T$ under the basis chosen above. 

The rest of the proof is essentially Galois descent. First note that $(\C_L)^{\Gal_L} = L$ \`a la Ax--Tate--Sen. Consider the $\Gal_L$-action on the Robba ring $\mathscr{R}_{\C_L}$ given by the usual $\Gal_L$-action on the coefficients in $\C_L$ and the trivial action on $\pi_K$. We equip $\mathscr{R}_L$ and $N_0$ with the trivial $\Gal_L$-actions and equip $N_{\C_L}=N_0\otimes_{\mathscr{R}_L}\mathscr{R}_{\C_L}$ with the induced action. Note that $\Gal_L$ commutes with $\Delta_{\pi_K} = \pi_K D_{\pi_K}$. In particular, $\ker(\Delta_{\pi_K}^d)$ is $\Gal_L$-stable. 

Choose a basis $f_1, \ldots, f_d$ of $N_0$. For each $1 \le i \le d$ and $j \in \Z$, consider the $\C_L$-linear $\Gal_L$-equivariant projections 
\[
\lambda_{i, j}\colon N_{\C_L} \ra \C_L
\]
given by sending $\sum_{i=1}^d (\sum_{j\in \Z} b_{i, j} \pi_{K}^j) \cdot f_i$ (with $b_{i, j} \in \C_L$) to $b_{i, j}$. Since $\ker(\Delta_{\pi_K}^d)$ is $d$-dimensional $\Gal_L$-stable $\C_L$-sub-vector space in $N_{\C_L}$, there exist $d$ such linear maps among $\lambda_{i, j}$'s such that their restrictions to $\ker(\Delta_{\pi_K}^d)$ give rise to a $\Gal_L$-equivariant isomorphism 
\[
\lambda\coloneqq \ker(\Delta_{\pi_K}^d) \isom \C_L^{d}.
\]
For each $1 \le i \le d$, let $\sq f_i$ be the pre-image of the $i$-th standard basis of $\C_L^{d}$ under $\lambda$. Then $\Gal_{L}$ fixes $\widetilde{f}_i$, so $\widetilde{f}_i \in N_0$ for all $1 \le i \le d$. Note that $\widetilde{f}_1, \ldots, \widetilde{f}_d$ form another basis of $N_0$. Write $W_0 \subset \ker(\Delta_{\pi_K}^d)$ for the $L$-span of $\sq f_1, \ldots, \sq f_d$, which is a $d$-dimensional $L$-subspace. Note that $\Delta_{\pi_K}$ preserves $W_0$ and acts nilpotently on $W_0$. So we can pick an $L$-basis $g_1, \ldots, g_d$ with respect to which $\Delta_{\pi_K}$ acts through an upper-triangular nilpotent matrix $Z_0$. This then yields an $\mathscr R_{L}$-basis of $N_0$ such that $\Delta_{\pi_K} = \pi_K \frac{d}{d\pi_K} + Z_0$, as desired.
\end{proof}

As an immediate consequence of Lemma~\ref{lemma:Galois_descent_of_unipotence}, we know that 
\[  
\mathbf{N}_{\mathrm{dR}, \xi_{\mathrm{I\!I}}} \coloneqq \mathbf{N}_{\mathrm{dR}}(\L|_{X_1}) \otimes_{\mathbf{C}_{X_1}^{\dagger}} \mathscr R_{L} = \mathbf{N}_{\mathrm{dR}}(\L|_{X_\delta}) \otimes_{\mathbf{C}_{X_\delta}^{\dagger}} \mathscr R_{L}
\] 
is strongly unipotent with respect to $D_{\pi_K}$. Therefore, for $r$ sufficiently small, we have
\[  
\mathbf{N}^r_{\mathrm{dR}, \xi_{\mathrm{I\!I}}} \coloneqq \mathbf{N}^r_{\mathrm{dR}}(\L|_{X_1}) \otimes_{\mathbf{C}_{X_1}^{r}} \mathscr R_{L}^r = \mathbf{N}^r_{\mathrm{dR}}(\L|_{X_\delta}) \otimes_{\mathbf{C}_{X_\delta}^{r}} \mathscr R_{L}^r
\] 
is strongly unipotent with respect to $D_{\pi_K}$.

\vspace{0.1in}
%% 8.3
\subsection{Propagation of arithmetic unipotence} \label{ss:type5_special_half_dagger}
\noindent 
\vspace{0.1in}

\noindent 
To prove Theorem~\ref{thm:spread_arithmetic_unipotence}, it remains to propagate the unipotence along the map $\mathbf{C}_{X_\delta}^{r} \ra \mathscr R_{L}^r$, for some sufficiently small $r$. To this end, we first prove a propagation result on the level of relative Robba rings. 

Let us introduce some notation. For any $\delta\in (0,1)\cap p^{\Q}$ and any monic polynomial $g=g(T) \in \mO_{K'_0} [T]$ of the form $g(T) = \prod_i (T-[\gamma_i])$ with $\gamma_i \in \overline{k}_K^{\times}$ (so that $|g| = 1$ for the 1-Gauss norm), consider the rational subset 
\[
V_{\delta, g} \coloneqq \{ \delta \le |[T^\flat]| \le 1, ~~|g([T^\flat])| = 1 \}\subset A_{K_0'}([\delta, 1]; [T^{\flat}]),
\]
by removing finitely many residue disks from the closed annulus $A_{K_0'}([\delta, 1]; [T^{\flat}])$. Let $A_{\delta, g}$ be the ring of analytic functions on $V_{\delta, g}$, equipped with the spectral norm $|\cdot|_{\delta, g}$. 
Note that $K_0'([T^\flat]) \subset \mathrm{Frac} (A_{\delta, g}) \subset L$, and $L$ coincides with the completion of $\mathrm{Frac} (A_{\delta, g})$ under the 1-Gauss norm.

\begin{remark} 
Our strategy is inspired by \cite[Proposition~5.3.3]{Kedlaya_coefficients}, where a similar result is proved for certain over-convergent algebras (which are called \emph{dagger algebras} in \textit{loc. cit.}). If we ignore the removed residue disks on $X_{1}$ for the moment, then a tempting choice of such a dagger algebra (in order to apply \cite[Proposition~5.3.3]{Kedlaya_coefficients}) could be 
\[
A^{\mathrm{naive}, \dagger} \coloneqq K_0' \gr{[T^\flat]^{\pm 1}}^{\dagger}, 
\]
where the over-convergence condition is both in the inward direction (i.e. towards radius $\delta$) and the outward  direction (i.e. towards radius $1+\epsilon$ for some $\epsilon > 0$). However, our $p$-adic differential equation is only defined over $\mathbf{C}_{X_\delta}^r$, whose coefficients in $[T^\flat]$ are not over-convergent in the outward direction. Nevertheless, we will show that the argument of \cite[Proposition~5.3.3]{Kedlaya_coefficients} can be adapted to our setting by considering a ``one-sided'' dagger algebra $\varinjlim_{\delta \ra 1^-} A_{\delta, g}$ instead. For the precise statement, see Proposition~\ref{prop:analogue_Kedlaya_5.3.3} below.
\end{remark}

\begin{notation} Let $\delta\in (0,1)\cap p^{\Q}$ and $g\in \mathcal{O}_{K'_0}[T]$ as above.
\begin{enumerate}
    \item For each $0 < s \le r$, consider the ring $\mathscr{R}_{A_{\delta, g}}^{[s, r]}$ with coordinate $\pi_K$ as defined in \S \ref{sss:notation_cyclo}. In particular, $\mathscr{R}_{A_{\delta, g}}^{[s, r]}$ is the affinoid $K'_0$-algebra of analytic functions on the affinoid space $V_{\delta, g} \times A_{K_0'} ([p^{-c_K \cdot r}, p^{-c_K \cdot s}])$. 
    
    \item Consider the following norms on $\mathscr{R}_{A_{\delta, g}}^{[s, r]}$:
    \begin{enumerate}
    \item[(i)] The spectral norm $|\!| \cdot |\!|_{\delta, g, [s, r]}$ on the product $V_{\delta, g} \times A_{K_0'} ([p^{-c_K \cdot r}, p^{-c_K \cdot s}])$ is given by
    \[
    |\!| \sum_{i \in \Z} a_i \pi_K^i |\!|_{\delta, g, [s, r]} \coloneqq \sup_{i \in \Z}\sup_{t \in [s, r]}\{ |a_i|_{\delta, g}\cdot p^{-c_K \cdot ti}\}.
    \]
    with $a_i\in A_{\delta, g}$.
    \item[(ii)] The ``relative 1-Gauss norm'' $|\!| \cdot |\!|_{L,  [s, r]}$ is defined to be
    \[
    |\!| \sum_{i \in \Z} a_i \pi_K^i |\!|_{L,  [s, r]} \coloneqq \sup_{i \in \Z}\sup_{t \in [s, r]} \{|a_i|_L \cdot p^{-c_K \cdot ti}\}.
    \]
    where we view each $a_i$ as an element of $L$ via the inclusion $A_{\delta,g}\subset L$ and $|\cdot|_L$ is the 1-Gauss norm on $L$. 
\end{enumerate}
\end{enumerate}  
\end{notation}

\begin{proposition} \label{prop:analogue_Kedlaya_5.3.3}
Let $\delta\in (0,1)\cap p^{\Q}$ and $g\in \mathcal{O}_{K'_0}[T]$ as above. Let $(M, \Delta)$ be a finite projective differential module over $(\mathscr R_{A_{\delta, g}}^{[s_0, r_0]}, \pi_K \frac{d}{d\pi_K})$ of rank $d$ for some $s_0, r_0$ satisfying $s_0 < r_0/p$. Suppose that the base change 
\[
M_L \coloneqq M \otimes_{\mathscr{R}_{A_{\delta, g}}^{[s_0, r_0]}} \mathscr{R}_L^{[s_0, r_0]}
\]
admits a strongly unipotent basis (i.e. with respect to this basis, $\Delta$ takes the form $\pi_K \frac{d}{d\pi_K} + Z_L$ for some upper-triangular nilpotent matrix $Z_L \in \mathrm{Mat}_{d\times d}(L)$). Then $M$ becomes strongly unipotent after shrinking $V_{\delta, g}$. More precisely, given any interval $[s, r] \subset (s_0, r_0)$ with $s<r/p$, up to increasing $\delta$, there exists a monic polynomial $g' \in \mathcal{O}_{K'_0}[T]$ of the form $g'(T)=\prod_{i=1}^m(T-[\gamma'_i])$ with $\gamma'_i\in \overline{k}^{\times}_K$, such that the base change 
\[
M_{\delta, g'g} \coloneqq  M \otimes_{\mathscr{R}_{A_{\delta, g}}^{[s_0, r_0]}} \mathscr{R}_{A_{\delta, g'g}}^{[s, r]}
\]
admits a strongly unipotent basis with respect to which $\Delta=\pi_KD_{\pi_K}$ acts through an upper-triangular nilpotent matrix with entries in $A_{\delta, g'g}$.
\end{proposition}

In the remainder of this subsection, we explain how the arguments in \cite{Kedlaya_coefficients} carry over to our setting, pointing out differences along the way. We start with the following analogue of \cite[Proposition~2.4.2]{Kedlaya_coefficients}, which is an essential ingredient in the proof of Proposition~\ref{prop:analogue_Kedlaya_5.3.3}.  

\begin{lemma} \label{lemma:analogue_kedlaya_2.4.2}
For any $x \in A_{\delta, g}$ and any $\epsilon \in (0, 1]$, we have 
\[
|x|_{\delta^{\epsilon}, g} \le |x|^{1-\epsilon}_{L} \cdot |x|_{\delta, g}^{\epsilon}
\]
\end{lemma}

\begin{proof}
We first claim that it suffices to prove the inequality for elements of the form $x = \frac{h}{g^N}$ for some $h \in K_0' [[T^\flat]^{\pm 1}]$ and $N \in \Z_{\ge 0}$. Indeed, these rational functions are dense in $A_{\delta, g}$ under the spectral norm. Suppose that $x \in A_{\delta, g}$ and that $x_n$ is a sequence of such rational functions with $|x_n - x|_{\delta, g} \ra 0$. Then, since $V_{\delta^\epsilon, g} \subset V_{\delta, g}$, we have $|x_n - x|_{\delta^{\epsilon}, g} \le |x_n - x|_{\delta, g} \ra 0$  and $|x_n - x|_{L} \le |x_n - x|_{\delta, g} \ra 0$, and the claim follows. 

For the rational function $x = \frac{h}{g^N}$, we have $|x|_{\delta, g} = \max \{|h|_{\delta}, |h|_{1}\}$ and $|x|_L = |h|_1$. Here, $|\cdot|_{\rho}$ stands for the $\rho$-Gauss norm; that is, for $h = \sum_{j \in \Z} a_j [T^\flat]^j$ with $a_j \in K_0'$, we put $|h|_{\rho} = \sup_{j} \{|a_j| \cdot \rho^j\}$. For each monomial $a_j [T^\flat]^j$, we have $|a_j| \delta^{\epsilon \cdot i} = |a_j|^{1-\epsilon} (|a_j|\delta^{i})^{\epsilon}$. Taking the supremum over $j$ yields the desired inequality  
$|x|_{\delta^{\epsilon}, g} \le |x|^{1-\epsilon}_{L} \cdot |x|_{\delta, g}^{\epsilon}$.  
\end{proof}

\begin{corollary}\label{cor:analogue_kedlaya_2.4.2}
For every element $x \in \mathscr{R}_{A_{\delta, g}}^{[s, r]}$ and every $\epsilon \in (0, 1]$, we have 
\[
|\!|x|\!|_{\delta^{\epsilon}, g, [s, r]} \le |\!|x|\!|^{1-\epsilon}_{L, [s, r]} \cdot |\!|x|\!|_{\delta, g, [s, r]}^{\epsilon}.
\]
\end{corollary}

\begin{proof}
This follows by applying Lemma~\ref{lemma:analogue_kedlaya_2.4.2} coefficient-wise to $x= \sum_{i} x_i \cdot \pi_K^i$. 
\end{proof}

\begin{lemma} \label{lemma:type_5_reduce_to_free}
Keep the assumption of Proposition~\ref{prop:analogue_Kedlaya_5.3.3}. Then, after shrinking $V_{\delta,g}$ (i.e., increasing $\delta$ and replacing $g$ by $g'g$ as in the statement of Proposition~\ref{prop:analogue_Kedlaya_5.3.3}), $M$ becomes finite free over $\mathscr{R}_{A_{\delta, g}}^{[s_0, r_0]}$. 
\end{lemma}

\begin{proof}
Note that $M_L$ is free by assumption. Choosing a basis $e_1, \ldots, e_d$ of $M_L$ and consider the isomorphism $S_L\colon (\mathscr{R}_L^{[s_0, r_0]})^d \xrightarrow[]{\sim} M_L$ sending $(a_1, \ldots, a_d)\mapsto \sum_{i=1}^d a_ie_i$. Pick a finite set of generators $\{m_j\}_{j\in J}$ of the finite projective $\mathscr{R}^{[s_0, r_0]}_{A_{\delta, g}}$-module $M$ and write $e_i=\sum_{j\in J}c_jm_j$ for each $i$, with $c_j\in \mathscr{R}^{[s_0, r_0]}_L$. Such an expression is not unique in general.

Observe that every element in $\mathscr{R}^{[s_0, r_0]}_L$ can be approximated by elements in $\mathscr{R}^{[s_0, r_0]}_{A_{\delta, g}}$ (up to shrinking $V_{\delta,g}$). To see this, it suffices to notice that $K_0'([T^\flat]) [\pi_K^{\pm 1}] \subset \mathscr R_{L}^{[s_0, r_0]}$ is dense under the relative 1-Gauss norm $|\!|\cdot|\!|_{L, [s_0, r_0]}$, and then note that each element in $K_0'([T^\flat])$ is defined on some $V_{\delta, g'g}$ (by increasing $\delta$ if necessary). 

Now, we approximate each $c_j$ by $c'_j\in \mathscr{R}^{[s_0, r_0]}_{A_{\delta, g}}$ as above, and let $e'_i=\sum_{j\in J}c'_j m_j\in M$. We obtain a linear map $S: \left(\mathscr{R}^{[s_0, r_0]}_{A_{\delta, g}}\right)^d\rightarrow M$ sending $(a_1, \ldots, a_d)\mapsto \sum_{i=1}^d a_ie'_i$, such that $S$ approximates $S_L$. Similarly, let us consider the dual basis $e_1^{\vee}, \ldots, e_d^{\vee}$ of $M_L^{\vee}$ and consider $S'_L\colon (\mathscr{R}_L^{[s_0, r_0]})^d \xrightarrow[]{\sim} M_L^{\vee}$ sending $(a_1, \ldots, a_d)\mapsto \sum_{i=1}^d a_ie_i^{\vee}$. By choosing a finite set of generators of $M^{\vee}$ and repeating the approximation process above, we arrive at a linear map $S':\left(\mathscr{R}^{[s_0, r_0]}_{A_{\delta, g}}\right)^d\rightarrow M^{\vee}$ that approximates $S'_L$. Taking duals, we obtain ${S'}^{\vee}: M\rightarrow \left(\mathscr{R}^{[s_0, r_0]}_{A_{\delta, g}}\right)^d$.

Choosing sufficiently accurate approximates, we can ensure that $P\coloneqq {S'}^{\vee}\circ S-\mathrm{id}$ satisfies $|\!| P |\!|_{L, [s_0, r_0]} < 1$. Applying Corollary~\ref{cor:analogue_kedlaya_2.4.2} entry-wise to $P$, we have 
\[ 
|\!| P |\!|_{\delta^{\epsilon}, g,  [s_0, r_0]} \le  |\!| P|\!|^{1-\epsilon}_{L, [s_0, r_0]} \cdot |\!|P|\!|_{\delta, g, [s_0, r_0]}^{\epsilon}. 
\]
The right hand side is $\le 1$ for $\epsilon$ sufficiently close to $0^+$. Hence, by choosing such an $\epsilon$ and replacing $\delta$ by $\delta^{\epsilon}$, we obtain $|\!| P |\!|_{\delta, g, [s_0, r_0]} < 1$ and so ${S'}^{\vee}\circ S$ is invertible. Therefore, $S$ is a split injection and $M \cong (\mathscr{R}_{A_{\delta, g}}^{[s_0, r_0]})^d  \oplus M'$ for a finite projective module $M'$ over $\mathscr{R}_{A_{\delta, g}}^{[s_0, r_0]}$. Since $M$ has rank $d$, we must have $M' = 0$. This finishes the proof.      
\end{proof}

Given the preparations above, Proposition~\ref{prop:analogue_Kedlaya_5.3.3} follows from the same argument as in the proof of \cite[Proposition~5.3.3]{Kedlaya_coefficients}, using Lemma~\ref{lemma:analogue_kedlaya_2.4.2} and Corollary~\ref{cor:analogue_kedlaya_2.4.2} instead of Proposition~2.4.2 in \textit{loc. cit}. Let us briefly sketch out the argument below. 

\begin{proof}[Proof of Proposition~\ref{prop:analogue_Kedlaya_5.3.3}]
Fix $[s, r] \subset (s_0, r_0)$ with $s < r/p$. By assumption, there exists a strongly unipotent basis of $M_L$ over $\mathscr R_{L}^{[s_0, r_0]}$ with respect to which $\Delta$ acts as $\Delta = \pi_K \frac{d}{d\pi_K} + Z_L$,  where $Z_L$ is an upper-triangular nilpotent matrix with coefficients in $L$. Let $e$ be the nilpotency index of $Z_L$. By Lemma~\ref{lemma:type_5_reduce_to_free}, we may assume that $M$ is free over $\mathscr{R}_{A_{\delta, g}}^{[s_0, r_0]}$. Consider the operators $h_l\colon M\rightarrow M$ defined by
\[ 
h_0 = \Delta^{e-1}, \qquad h_l = (1 - \frac{\Delta^2}{l^2})^e h_{l -1} \: \: \:  \textrm{ for } l \ge 1 
\]
as in \cite[Lemma~5.3.1]{Kedlaya_coefficients}. Then by the same argument as in the proof of \cite[Lemma~5.3.1]{Kedlaya_coefficients}, we have that for any $w \in M_L$, the sequence $\{h_l (w)\}_{\ell = 0}^{\infty}$ converges to an element $h_\infty (w)$ in $M_L$ under $|\!|\cdot|\!|_{L, [s_0, r_0]}$, satisfying $\Delta(h_\infty(w)) = 0$. The key is then to prove the following analogue of \cite[Lemma~5.3.2]{Kedlaya_coefficients}. %

\vspace{0.05in}
    
\noindent \textbf{Claim.} Let $w \in M$ and $[s', r'] \subset (s_0, r_0)$. There exists $\delta' \ge \delta$ such that $h_\infty (w)$ lies in 
\[
M \otimes_{ \mathscr{R}_{A_{\delta, g}}^{[s_0, r_0]}} \mathscr{R}_{A_{\delta', g}}^{[s', r']}.
\] 

\noindent \textit{Proof of Claim}. We follow a similar idea as in the proof of \cite[Lemma~5.3.2]{Kedlaya_coefficients}. Write $w_l \coloneqq  h_{l+1}(w) - h_l (w)$. Then by the proof of \textit{loc. cit.}, we have  
$|\!| w_l |\!|_{L, [s', r']} \le c \cdot \alpha^l$ and 
$|\!| w_l |\!|_{\delta, g, [s', r']} \le c' \cdot (\alpha')^l$ 
for some $c, c', \alpha, \alpha'>0$ with $\alpha <1$ and $\alpha' \ge 1$ (here $[s', r'] \subset (s_0, r_0)$ is used to ensure that we can find such $\alpha < 1$). Now by Corollary~\ref{cor:analogue_kedlaya_2.4.2} (this replaces the role of \cite[Proposition~2.4.2]{Kedlaya_coefficients} in the proof of  \cite[Lemma~5.3.2]{Kedlaya_coefficients}), we have 
\[
|\!|w_l|\!|_{\delta^\epsilon, g, [s', r']} \le c^{1-\epsilon} (c')^\epsilon (\alpha^{1-\epsilon} (\alpha')^\epsilon)^l 
\]
for all $l$. Choosing $\epsilon$ such that $\alpha^{1-\epsilon} (\alpha')^\epsilon < 1$ and taking $\delta' = \delta^\epsilon$, we obtain the desired statement in the claim.

\vspace{0.05in}

Thanks to the claim, for a given $w \in M$, we may assume that $h_\infty(w) \in M$ (after increasing $\delta$ if necessary), which again satisfies $\Delta (h_\infty (w)) = 0$. This puts us in a situation where we can apply the induction argument as in the proof of \cite[Proposition~5.3.3]{Kedlaya_coefficients}. Now the argument in the proof of \cite[Proposition~5.2.6]{Kedlaya_coefficients} carries over \emph{verbatim} and produces the desired strongly unipotent basis. 
\end{proof}

\vspace{0.1in}
%%% 8.4
\subsection{Proof of Theorem~\ref{thm:spread_arithmetic_unipotence}} \label{ss:proof_of_arithmetic_unipotence}
\noindent 
\vspace{0.1in}

\noindent 
Using Proposition~\ref{prop:analogue_Kedlaya_5.3.3}, we now prove Theorem~\ref{thm:spread_arithmetic_unipotence}. The first step is to achieve unipotence over $\mathbf{C}_{X_\delta}^{[s, r]}$ for some suitable choice of $[s, r]$. For any $[s,r]$ with $r$ sufficiently small, let
\[
\mathbf{N}^{[s,r]}_{\mathrm{dR}}(\L|_{X_{\delta}}):=\mathbf{N}^r_{\mathrm{dR}}(\L|_{X_{\delta}})\otimes_{\mathbf{C}^r_{X_{\delta}}}\mathbf{C}^{[s,r]}_{X_{\delta}}.
\]
Recall that we already know that $\mathbf{N}^{r_0}_{\mathrm{dR}, \xi_{\mathrm{I\!I}}}=\mathbf{N}^{r_0}_{\mathrm{dR}}(\L|_{X_{\delta}})\otimes_{\mathbf{C}^{r_0}_{X_{\delta}}} \mathscr{R}_L^{r_0}$ is strongly unipotent with respect to $D_{\pi_K}$, for some sufficiently small $r_0$. We further assume that $\lambda_r(u_{\pi_K}-1) < p^{-\frac{1}{p-1}}$ for all $0<r\le r_0$, where $u_{\pi_K} = \varphi(\pi_K)/\pi_K^p$. In particular, we have a well-defined element $\log (u_{\pi_K})$ in $\mathbf{C}_K^{r_0}$.

\begin{lemma}  \label{lemma:before_Frob_spreading}
Up to replacing $K$ by an unramified extension and shrinking $X_{\delta}=X_{\delta,f}$, there exists $0<s< r<r_0$ with $s<r/p$, such that $\mathbf{N}^{[s,r]}_{\mathrm{dR}}(\L|_{X_{\delta}})$ is free over $\mathbf{C}_{X_\delta}^{[s, r]}$ and admits a basis with respect to which $\Delta_{\pi_K}=\pi_KD_{\pi_K}$ acts through an upper-triangular nilpotent matrix with entries in $A_{\delta^r, f}$.
\end{lemma}

\begin{proof}
Pick $s_0\in (0,r_0/p]$. Then $\mathbf{N}^{[s_0,r_0]}_{\mathrm{dR}}(\L|_{X_{\delta}})\otimes_{\mathbf{C}^{[s_0,r_0]}_{X_{\delta}}} \mathscr{R}_L^{[s_0,r_0]}$ is strongly unipotent with respect to $D_{\pi_K}$.

By Lemma~\ref{lem: puncture-C-description}, $\mathbf{C}_{X_{\delta}}^{[s_0, r_0]}$ is precisely the ring of analytic functions on the rational subdomain  
\[
W_{\delta, f}^{ [s_0, r_0]} = \{|\pi_K| = p^{-c_K \cdot t},  ~~\delta^{t} \le  |[T^{\flat}]| \le 1, ~~|f([T^{\flat}])| = 1, ~~t \in [s_0, r_0]\}
\]
inside the two-dimensional poly-annulus
\[A_{K'_0}([\delta^{r_0},1];[T^{\flat}])\times A_{K'_0}([p^{-c_Kr_0}, p^{-c_Ks_0}];\pi_K).\]
Therefore, we have containment of rational subsets
\[ 
V_{\delta^{s_0}, f} \times A_{K_0'} ([p^{-c_K \cdot r_0}, p^{-c_K \cdot s_0}]) \: \subset \: W_{\delta, f}^{ [s_0, r_0]} \:  \subset  \: V_{\delta^{r_0}, f} \times A_{K_0'} ([p^{-c_K \cdot r_0}, p^{-c_K \cdot s_0}]).
\]
In particular, for any choice of $\delta'\in [\delta^{s_0},1)$, we have 
\[
V_{\delta', f} \times A_{K_0'} ([p^{-c_K \cdot r_0}, p^{-c_K \cdot s_0}]) \: \subset \: W_{\delta, f}^{ [s_0, r_0]},
\]
which induces an injection
\[
\iota_1\colon \mathbf{C}_{X_{\delta}}^{[s_0, r_0]} \ra \mathscr R_{A_{\delta', f}}^{[s_0, r_0]}
\]
compatible with the arithmetic differential operator $D_{\pi_K}$. The base change of $\mathbf{N}_{\mathrm{dR}}^{[s_0, r_0]} (\L |_{X_{\delta}})$ along $\iota_1$ then produces a differential module satisfying the assumption of Proposition~\ref{prop:analogue_Kedlaya_5.3.3}. Consequently, by choosing $[s,r]\subset(s_0,r_0)$ with $s<r/p$ and increasing $\delta'$, we know that $\mathbf{N}_{\mathrm{dR}}^{[s_0, r_0]} (\L |_{X_{\delta}})$ becomes strongly unipotent with respect to $D_{\pi_K}$ after further base change to $\mathscr{R}^{[s,r]}_{A_{\delta', g'f}}$ for some $g'\in \mathcal{O}_{K'_0}[T]$. Up to replacing $K$ by an unramified extension\footnote{Note that replacing $K$ by an unramified extension does not change $\pi_K$.} (so that the coefficients of $g'$ lie in $K_0$) and shrinking $X_{\delta}$, we may assume that 
\[
\mathbf{N}_{\mathrm{dR}}^{[s_0, r_0]} (\L |_{X_{\delta}}) \otimes_{ \mathbf{C}_{X_{\delta}}^{[s_0, r_0]} } \mathscr R_{A_{\delta', f}}^{[s, r]}
\]
admits a basis with respect to which $\Delta_{\pi_K}$ acts through an upper-triangular matrix with entries in $A_{\delta', f}$.
    
Finally, we pick $\delta'' \in ((\delta')^{1/r}, 1)$ so that
\[
W_{\delta'', f}^{[s, r]} \:  \subset  \: V_{\delta', f} \times A_{K_0'} ([p^{-c_K \cdot r}, p^{-c_K \cdot s}]),
\]
which induces an injection
\[
\iota_2\colon \mathscr R_{A_{\delta', f}}^{[s, r]} \ra  \mathbf{C}_{X_{\delta''}}^{[s, r]} 
\]
compatible with $D_{\pi_K}$. Since the composition $\iota_2 \circ \iota_1$ is the natural inclusion and we have 
\[\left(\mathbf{C}^{[s,r]}_{X_{\delta'', f}}\right)^{\frac{d}{d\pi_K}=0}=A_{{\delta''}^r, f}\]
from the explicit description of $\fC_{X_{\delta}}^{[s,r]}$ (cf. Lemma~\ref{lem: puncture-C-description}), we obtained the desired strong unipotence.
\end{proof}

We are now ready to prove Theorem~\ref{thm:spread_arithmetic_unipotence}. The idea is to use Frobenius to translate the arithmetic unipotence from $\mathbf{C}_{X_\delta}^{[s, r]}$ to $\mathbf{C}_{X_\delta}^{r}$. 

\begin{proof}[Proof of Theorem~\ref{thm:spread_arithmetic_unipotence}] To simplify notations, write $N^{[s, r]} \coloneqq \mathbf{N}_{\mathrm{dR}}^{[s, r]} (\L |_{X_{\delta}})$. By Lemma~\ref{lemma:before_Frob_spreading}, we may assume that $0<s<r/p$ and $N^{[s, r]}$ is free over $\mathbf{C}_{X_\delta}^{[s, r]}$ and admits a basis with respect to which $\Delta_{\pi_K}$ acts as  
\begin{equation} \label{eq:strongly_unipotent_for_s_r}
\Delta_{\pi_K} = \pi_K \frac{d}{d\pi_K} + Z,
\end{equation}
for some upper-triangular nilpotent matrix $Z$ with coefficients in $A_{\delta^r, f}$.

Write $s_n \coloneqq  s/p^n$ and $r_n \coloneqq  r/p^n$. By \cite[Construction~10.7]{relative}, there are compatible isomorphisms 
\[
\phi_n\colon N^{[s_n, r_n]} \otimes_{\mathbf{C}_{X_{\delta}}^{[s_n, r_n]}, \varphi} \mathbf{C}_{X_{\delta}}^{[s_{n+1}, r_{n+1}]} \isom  N^{[s_{n+1}, r_{n+1}]}. 
\] 
Our goal is to construct a basis $\mathscr{B}_n$ of $N^{[s_n, r_n]}$ for each $n \ge 1$, such that $\mathscr{B}_n$'s agree on the overlap and $\Delta_{\pi_K} = \pi_K \frac{d}{d\pi_K} + Z$ relative to these bases for the same $Z$ as in \eqref{eq:strongly_unipotent_for_s_r}. We will construct $\mathscr{B}_n$ inductively. 

Recall that we have a well-defined element $\log(u_{\pi_K})\in \mathbf{C}^r_K$. Assuming $\mathscr{B}_n$ is given, we construct $\mathscr{B}_{n+1}$ as follows. Write $\mathscr{B}_n = \{e_1, \ldots, e_d\}$, and consider $\mathscr{B}_{n+1}'= \{\phi_n(e_1), \ldots, \phi_n(e_d) \}$ which is a basis of $N^{[s_{n+1}, r_{n+1}]}$. Since $\Delta_{\pi_K}$ is compatible with Frobenius (cf. \cite[Construction~10.7]{relative}), we have \[D_{\pi_K} (\phi(x)) = \frac{d\varphi(\pi_K)}{d\pi_K}\cdot \phi(D_{\pi_K} (x))\] for all $x \in N^{[s_n, r_n]}$. Hence
\[\Delta_{\pi_K} (\phi(x)) = \frac{\pi_K}{\varphi(\pi_K)}\cdot\frac{d\varphi(\pi_K)}{d\pi_K}\cdot \phi(\Delta_{\pi_K} (x)).\]
This implies 
\[
\Delta_{\pi_K} (\mathscr{B}_{n+1}') =  \mathscr{B}_{n+1}' \cdot \varphi(Z) \cdot \big(p+ \Delta_{\pi_K} (\log u_{\pi_K}) \big).
\]
We modify $\mathscr B_{n+1}'$ into a new basis 
\[ 
\mathscr{B}_{n+1}'' \coloneqq  \mathscr{B}_{n+1}' \cdot \exp \big(- \log u_{\pi_K} \cdot \varphi(Z) \big),
\]
which is well-defined as $Z$ is nilpotent. It is straightforward to check that $\Delta_{\pi_K} (\mathscr{B}_{n+1}'') = p \cdot  \mathscr{B}_{n+1}'' \cdot \varphi(Z)$. In particular, $\mathscr{B}_{n+1}''$ is a strongly unipotent basis for $N^{[s_{n+1}, r_{n+1}]}$. We will further modify it so that it agrees with $\mathscr{B}_n$ on the overlap $\mathbf{C}_{X_{\delta}}^{[s_n, r_{n+1}]}$. To this end, let $C_n \in \mathrm{GL}_d (\mathbf{C}_{X_{\delta}}^{[s_n, r_{n+1}]})$ be the change of basis matrix so that  
\begin{equation} \label{eq:change_basis_C_n}
\mathscr B_{n+1}''  = \mathscr B_{n} \cdot C_n
\end{equation} 
on the overlap. Differentiating both sides of \eqref{eq:change_basis_C_n} gives 
$\Delta_{\pi_K}(C_n) = p \cdot C_n \cdot \varphi(Z) - Z \cdot C_n$. Since $\Delta_{\pi_K} (\varphi(Z)) = \Delta_{\pi_K} (Z) = 0$ and $Z$ is nilpotent with $Z^d=0$, we deduce that $\Delta_{\pi_K}^{2d-1}(C_n) = 0$. On the other hand, using the explicit description of $\fC_{X_{\delta}}^{[s, r]}$ (cf. Lemma~\ref{lem: puncture-C-description}), one computes that   
\[
\ker (\Delta_{\pi_K}^{m}\colon \mathbf{C}_{X_{\delta}}^{[s, r]} \ra \mathbf{C}_{X_{\delta}}^{[s, r]}) = A_{\delta^r, f}
\]
for all $m \ge 1$. Hence, the coefficients of $C_n$ all lie in $A_{\delta^{r_{n+1}}, f} \subset \fC_{X_{\delta}}^{[s_{n+1}, r_{n+1}]}$. Applying the same argument to $C_n^{-1}$, we deduce that the coefficients of $C_n^{-1}$ also lie in $A_{\delta^{r_{n+1}}, f} \subset \fC_{X_{\delta}}^{[s_{n+1}, r_{n+1}]}$. Therefore, 
\[
\mathscr{B}_{n+1}  = \mathscr{B}_{n+1}'' \cdot C_n^{-1}, 
\]
is a well-defined basis of $N^{[s_{n+1}, r_{n+1}]}$ that coincides with $\mathscr{B}_n$ on $\mathbf{C}_{X_{\delta}}^{[s_n, r_{n+1}]}$, as desired. This finishes the proof of the theorem.    
\end{proof}

\vspace{0.1in}
%%% 8.5
\subsection{Arithmetic horizontal sections} \label{ss:arithmetic_hs}
\noindent 
\vspace{0.1in}

\noindent 
In the remainder of the section, we prove Theorem~\ref{thm:main_type5_disk}. We want to show that, up to enlarging $K$ and shrinking $X_{\delta}$, there exists a finite \'etale cover over which $\L$ becomes pointwise semistable. By  Theorem~\ref{thm:spread_arithmetic_unipotence}, we may assume that $\mathbf{N}^r_{\mathrm{dR}}(\L)$ is strongly unipotent with respect to the arithmetic differential operator $D_{\pi_K}$ over $\mathbf{C}^r_{X_{\delta}}$ for some sufficiently small $r$. For any $\sigma\in [\delta, 1)\cap p^{\Q}$, consider the embedding $A_{K'_0}([\delta, \sigma];T)\subset X_{\delta}$ which induces a natural map $\mathbf{C}^r_{X_{\delta}}\rightarrow \mathbf{C}^r_{[\delta, \sigma]}$. Then the base change
\[
N_{[\delta, \sigma]}^{r} \coloneqq \mathbf{N}^r_{\mathrm{dR}} (\L) \otimes_{\mathbf{C}_{X_{\delta}}^{r}} \mathbf{C}_{[\delta, \sigma]}^{r}\cong \mathbf{N}^r_{\mathrm{dR}}(\L|_{A_{K'_0}([\delta, \sigma])})
\]
is a finite free $\mathbf{C}_{[\delta, \sigma]}^{r}$-module equipped with commuting differential operators $D_{[T^\flat]}$ and $D_{\pi_K}$, as well as a Frobenius structure, where the last isomorphism follows from \cite[Proposition 10.13]{relative}. Moreover, $N_{[\delta, \sigma]}^{r}$ is strongly unipotent with respect to $D_{\pi_K}$ for all $\sigma\in [\delta, 1)\cap p^{\Q}$.

To prove Theorem~\ref{thm:main_type5_disk}, we will first study the $D_{\pi_K}$-horizontal sections after adjoining $\log \pi_K$. Then we study the geometric differential operator on these arithmetic horizontal sections using the classical $p$-adic local monodromy theorem of Andr\'e--Kedlaya--Mebkhout. As a preparation, we first compute the arithmetic horizontal sections in relevant relative period rings.

\begin{lemma} \label{lemma:computing_D_pi=0} For any $0 < \sigma_1 \le \sigma_2 \le 1$ and $0 < s \le r$, we have 
\[
\left( \mathbf{C}_{[\sigma_1, \sigma_2]}^{[s, r]}\right)^{
\frac{d}{d\pi_K} = 0}=\left( \mathbf{C}_{[\sigma_1, \sigma_2]}^{[s, r]} [\log \pi_K] \right)^{
\frac{d}{d\pi_K} = 0} = \Gamma_{K_0'}([\sigma_1^r, \sigma_2^s]; [T^\flat]) 
\]
where the right-hand side is the ring of analytic functions on the annulus $A_{K_0'}([\sigma_1^r, \sigma_2^s]; [T^\flat])$ (cf. \S \ref{sss:notation_cyclo}). Taking limit as $s\rightarrow 0^+$, we obtain
\[\left( \mathbf{C}_{[\sigma_1, \sigma_2]}^{r}\right)^{\frac{d}{d\pi_K} = 0}=\left( \mathbf{C}_{[\sigma_1, \sigma_2]}^{r}  [\log \pi_K] \right)^{\frac{d}{d\pi_K} = 0} = \Gamma_{K_0'}([\sigma_1^r, 1); [T^\flat]).\]
\end{lemma}

\begin{proof}
Since $\log \pi_K$ is transcendental over $\mathbf{C}_{[\sigma_1, \sigma_2]}^{[s, r]}$, it suffices to show \[\left( \mathbf{C}_{[\sigma_1, \sigma_2]}^{[s, r]} \right)^{
\frac{d}{d\pi_K} = 0} = \Gamma_{K_0'}([\sigma_1^r, \sigma_2^s]; [T^\flat]).\]
But this follows from the explicit description of the period ring $\mathbf{C}_{[\sigma_1, \sigma_2]}^{[s, r]}$ in Lemma~\ref{lemma: explicit description C--case of annulus}.
\end{proof}

\begin{remark}
Let us emphasize that the arithmetic horizontal sections of $\mathbf{C}_{[\sigma_1, \sigma_2]}^{r}$ and $\mathbf{C}_{[\sigma_1, \sigma_2]}^{r} [\log \pi_K]$ are actually independent of the outer radius $\sigma_2$. 
\end{remark}

In what follows, we will mostly use ``$\mathscr{R}^{?}_{?}$'' (cf. \S \ref{sss:notation_cyclo}) to denote the corresponding rings in coordinate ``$[T^{\flat}]$''. More precisely, by a slight abuse of notation, we write
\[
\mathscr{R}_{K_0'}^{[s, r]}\coloneqq \Gamma_{K_0'}([p^{-r\cdot c_K}, p^{-s\cdot c_K}]; [T^{\flat}]) ~~\text{ and }~~ \mathscr{R}_{K_0'}^r\coloneqq \Gamma_{K_0'}( [p^{-r\cdot c_K}, 1); [T^{\flat}] ).
\]
We also write $w_{\delta}\coloneqq -\frac{\log_p \delta}{c_K}$ so that
\[
\mathscr{R}_{K_0'}^{w_{\delta} r} = \Gamma_{K_0'}( [\delta^r, 1); [T^{\flat}] ).
\]

Now, we return to the situation at the beginning of \S \ref{ss:arithmetic_hs}. In particular, $N^r_{[\delta, \sigma]}$ is strongly unipotent with respect to $D_{\pi_K}$.

\begin{corollary} \label{cor:N^r_arith_computed}
After joining $\log \pi_K$, the arithmetic horizontal sections 
\[
N^{r, \nabla_{\mathrm{arith}}}\coloneqq \left( N_{[\delta, \sigma]}^{r} [\log \pi_K] \right)^{D_{\pi_K} = 0} 
\] 
(equipped with the geometric differential operator $D_{[T^{\flat}]}$) is a finite free differential module over $(\mathscr{R}_{K_0'}^{w_{\delta} r},\frac{d}{d[T^{\flat}]})$, equipped with a Frobenius structure. Moreover, it is independent of $\sigma$.
\end{corollary}

\begin{proof}
Note that $D_{[T^\flat]}$ commutes with $D_{\pi_K}$. The statements (except those regarding Frobenius structures) follow from the strong unipotence assumption and the Lemma~\ref{lemma:computing_D_pi=0}. 

For the Frobenius structure, first note that Frobenius preserves the arithmetic horizontal sections. Indeed, we have an identity
\begin{equation}\label{eq: D vs phi--log}
D_{\pi_K} \circ \phi = \frac{d \varphi(\pi_K)}{d\pi_K} \cdot \phi \circ D_{\pi_K}
\end{equation}
on $N^r_{[\delta, \sigma]}$, by the compatibility of the differential operators and Frobenius (cf. \cite[Construction~10.7]{relative}). Since $\frac{d}{d\pi_K}\circ \varphi (\log \pi_K) = \frac{d \varphi(\pi_K)}{d\pi_K}\cdot \varphi(\pi_K)^{-1}$, the identity \eqref{eq: D vs phi--log} remains valid on $N^r_{[\delta, \sigma]}[\log\pi_K]$. Therefore, for $x \in N_{[\delta, \sigma]}^{r}[\log \pi_K]$, if $D_{\pi_K} (x) = 0$, then $D_{\pi_K} (\phi(x)) = 0$. 

Finally, the original Frobenius structure $\phi$ on $N^r_{[\delta, \sigma]}$ induces a Frobenius isomorphism
\[
\phi\colon N^{r, \nabla_{\mathrm{arith}}} \otimes_{\mathscr{R}_{K_0'}^{w_{\delta} r} \varphi} \mathscr{R}_{K_0'}^{w_{\delta} , \frac{r}{p}} \isom N^{\frac{r}{p}, \nabla_{\mathrm{arith}}} 
\]
over $\mathscr{R}_{K_0'}^{w_{\delta} \frac{r}{p}}$, as desired.
\end{proof}

\vspace{0.1in}
%%% 8.6
\subsection{Achieving semistability} \label{ss:Katz_special_algebraization} 
\noindent 
\vspace{0.1in}

\noindent
We are ready to prove Theorem~\ref{thm:main_type5_disk}. By Corollary~\ref{cor:N^r_arith_computed}, 
\[
N^{\dagger, \nabla_{\mathrm{arith}}} \coloneqq \varinjlim_{r \ra 0^+} N^{r, \nabla_{\mathrm{arith}}} 
\] 
is a finite free differential module (of rank $d$) over $(\mathscr{R}_{K_0'}^{\dagger}, \frac{d}{d[T^{\flat}]})$ equipped with a Frobenius structure, where $\mathscr{R}_{K_0'}^{\dagger} \coloneqq \varinjlim_{r \rightarrow 0^+} \mathscr{R}_{K_0'}^{w_{\delta} r}$ is the \emph{geometric} Robba ring in coordinate $[T^{\flat}]$. Note that the mod $p$ reduction of the integral Robba ring inside $\mathscr{R}_{K_0'}^{\dagger}$ is precisely $\kappa (\!(T^\flat)\!)$. By the classical $p$-adic monodromy theorem of Andr\'e--Kedlaya--Mebkhout, there exists a finite separable extension $E'/\kappa (\!(T^\flat)\!)$ with associated Robba ring extension $\mathscr{R}_{K_0'}^{\dagger}\rightarrow \mathscr{R}'$, such that $N^{\dagger, \nabla_{\mathrm{arith}}}$ becomes (strongly) unipotent for the geometric differential operator $D_{[T^\flat]}$ after tensoring with $\mathscr{R}'$. Our goal is to find a finite \'etale cover of $X_{\delta}$ (after enlarging $K$ and shrinking $X_{\delta}$ if necessary) that realizes this extension of Robba rings. 

By \cite[Theorem~1.4.1]{Katz_local_global}, there exists a finite \'etale $\kappa[(T^\flat)^{\pm 1}]$-algebra $Q$\footnote{This corresponds to what Katz calls a \emph{special} finite \'etale cover of $\mathbb{G}_{m,\kappa}$} together with an isomorphism 
\[
\alpha\colon Q \otimes_{\kappa[(T^\flat)^{\pm 1} ]}  \kappa (\!(T^\flat)\!) \isom E'. 
\]
Recall from \S \ref{sss: example closed annuli with punctures} the perfectoid algebra $\widehat{A}_{X_{\delta}, \infty}$ with its tilt $\widetilde{R}_{X_{\delta}}$. Consider the finite \'etale map
\[
\widetilde{R}_{X_{\delta}} \rightarrow Q \otimes_{\kappa [(T^\flat)^{\pm 1}]} \widetilde{R}_{X_{\delta}}.
\]
By \cite[Proposition~5.28]{relative}, this corresponds to some finite \'etale maps $R_{X_{\delta}} \rightarrow R'$ and $\widehat{A}_{X_{\delta}, \infty} \rightarrow B$. By the same argument as in Construction~\ref{construction:geometric_realization}, after replacing $K$ by $K(\zeta_{p^n})$ and $X_{\delta}$ by its $p^n$-th Kummer cover $X_{\delta, n}$ for some $n \ge 0$ if necessary, we may assume that the finite \'etale map $\widehat{A}_{X_{\delta}, \infty} \rightarrow B$ descends to a finite \'etale map $A_{X_{\delta}} \rightarrow B_0$. Let $\mathrm{pr}:X'_{\delta} \ra X_{\delta}$ be the corresponding finite \'etale cover, and let ${\mathbf{A}}_{X_{\delta}'}^{r}$, ${\mathbf{B}}_{X_{\delta}'}^{r}$, ${\mathbf{C}}_{X_{\delta}'}^{r}$, etc. denote the corresponding imperfect relative period rings (cf. Construction~\ref{construction: weakly decompleting finite etale}). In particular, we have ${\mathbf{A}}_{X_{\delta}'}^{r}/p = R'$ for $r$ sufficiently small. 

Theorem~\ref{thm:main_type5_disk} is an immediate consequence of the following proposition. 

\begin{proposition} 
After enlarging $K$ if necessary, there exists a finite \'etale cover $X_{\delta}'' \ra X_{\delta}'$ such that $\L$ becomes pointwise semistable after pulling back along the composition $X_{\delta}'' \ra X'_{\delta}\xrightarrow[]{\mathrm{pr}} X_{\delta}$. 
\end{proposition}

\begin{proof} 
By assumption, $\L$ is already pointwise semistable on the outer boundary $X_1\subset X_{\delta}$. So, it suffices to construct a finite \'etale cover $X_{\delta}'' \rightarrow X_{\delta}'$ (after enlarging $K$ if necessary) such that for any closed annulus $A_K([\delta, \sigma]) \subset X_{\delta}$ with $\sigma \in [\delta, 1)\cap p^{\Q}$, the pullback of $\L|_{A_K([\delta, \sigma])}$ over $X''_{[\delta, \sigma]}\coloneqq A_K([\delta, \sigma])\times_{X_{\delta}} X''_{\delta}$ is pointwise semistable. 

Write $X'_{[\delta, \sigma]}\coloneqq A_K([\delta, \sigma])\times_{X_{\delta}} X'_{\delta}$. By \cite[Proposition~5.28, Corollary 5.31]{relative}, the isomorphism
\[
\widetilde{R}_{X'_{[\delta, \sigma]}} = \widetilde{R}_{X_{\delta}}\otimes_{\kappa [(T^\flat)^{\pm 1}]} Q \cong  \widetilde{R}_{X_{\delta}}\otimes_{\kappa(\!(T^{\flat})\!))} E'
\]
induces an isomorphism
\[
\fC_{X'_{[\delta, \sigma]}}^{\dagger} \cong \fC_{X_{[\delta, \sigma]}}^{\dagger} \otimes_{\mathscr{R}_{K_0'}} \mathscr{R}',
\]
which is compatible with both differential operators $D_{\pi_K}$ and $D_{[T^{\flat}]}$. Thus, by Lemma~\ref{lemma:computing_D_pi=0}, we have
\[
\left( \fC_{X'_{[\delta, \sigma]}}^{\dagger}[\log \pi_K] \right)^{\frac{d}{d\pi_K} = 0} = \mathscr{R}'.
\]

Now, since $\fN_{\mathrm{dR}}(\cdot)$ is compatible with finite \'etale covers by \cite[Proposition~10.13]{relative}, we have
\begin{align*}
\left(  \mathbf{N}_{\mathrm{dR}} (\mathrm{pr}^* \L|_{X'_{[\delta, \sigma]}})  [\log \pi_K] \right)^{D_{\pi_K} = 0} &  
= \left(  \mathbf{N}_{\mathrm{dR}} ( \L|_{X_{[\delta, \sigma]}})  [\log \pi_K]  \otimes_{ \mathbf{C}^{\dagger}_{X_{[\delta, \sigma]}} } \mathbf{C}^{\dagger}_{X'_{[\delta, \sigma]}} \right)^{D_{\pi_K} = 0} \\
& 
= \left(    N_{[\delta, \sigma]}^{\dagger} [\log \pi_K] \otimes_{\mathscr{R}_{K_0'}} \mathscr{R}'   \right)^{D_{\pi_K} = 0} = N^{\dagger, \nabla_{\mathrm{arith}}}  \otimes_{\mathscr{R}_{K_0'}} \mathscr{R}'. 
\end{align*} 
which is (strongly) unipotent with respect to the geometric differential $D_{[T^\flat]}$.
Therefore, we obtain a decomposition
\[
\mathbf{N}_{\mathrm{dR}} (\mathrm{pr}^* \L|_{X'_{[\delta, \sigma]}})  [\log \pi_K, \log [T^{\flat}]] \cong \bigoplus_{i=1}^d \fC_{X'_{[\delta, \sigma]}}^{\dagger} [\log \pi_K, \log [T^{\flat}]]\cdot w_i
\]
compatible with $D_{\pi_K}$ and $D_{[T^{\flat}]}$ on both sides, where each $w_i$ satisfies $D_{[T^\flat]} (w_i) = D_{\pi_K} (w_i) = 0$. %
Finally, note that all these $w_i$'s are fixed by an open subgroup of $\Gamma$. By enlarging $K$ and passing to a $p$-power Kummer cover of $X'_{\delta}$ (cf. Remark~\ref{remark: passing to p-power Kummer cover}), we may further assume that all $w_i$'s are fixed by $\Gamma$. Hence, we can apply Proposition~\ref{prop:unipotence_and_st} (note that the functors $\mathbf{N}_{\mathrm{dR}}$ and $D^{\dagger}_{\mathrm{rig}}$ coincide after inverting $t$) and conclude the proof.
\end{proof}

\vspace{0.3in}

%%% Section 9
\section{Around points of type V: general case}  \label{section: type V general}

In this section, we prove the local $p$-adic monodromy (Theorem~\ref{thm: local main thm}) for general type V points by reducing to the special case treated in \S \ref{section: type 5}. Since every type V point on $X_{\C_K}$ is $K'$-rational for some finite extension $K'/K$, it suffices to prove the folowing.

\begin{theorem}\label{thm: local main thm type V general}
Let $X$ be a smooth rigid analytic curve over $K$, viewed as an adic space over $\spa(K, \mathcal{O}_K)$. Let $\L$ be a de Rham local system on $X$. Let $\xi\in X$ be a $K$-rational type V point. Then, after replacing $K$ by a finite extension if necessary, there exists an open neighborhood $U$ of $\xi$ and a finite \'etale cover $U'\rightarrow U$ such that the pullback of $\mathbb{L}|_U$ to $U'$ is pointwise semistable.
\end{theorem}

A key ingredient in the proof is a spread-out argument (cf. Theorem~\ref{thm: spread-out thm}), which will also play an essential role in the proof of the main theorem in \S \ref{section: proof of the main theorem}.

\vspace{0.1in}
%% 9.1
\subsection{A spread-out theorem}\label{ss: spread-out thm}
\noindent
\vspace{0.1in}

\noindent The following spread-out theorem will be used in both \S \ref{ss: type 5} and in the proof of Theorem~\ref{thm:main_intro} in \S \ref{section: proof of the main theorem}.

\begin{theorem}\label{thm: spread-out thm}
Let $X$ be a smooth projective curve over $K$ and let $X^{\mathrm{ad}}$ denote the associated adic space over $\spa(K, \mathcal{O}_K)$. Let $\xi$ be a non-type I point on $X^{\mathrm{ad}}$ and let $U \subset X^{\mathrm{ad}}$ be an affinoid open neighborhood of $\xi$. Let $f\colon  V \ra U$ be a finite \'etale cover. Then, up to shrinking $U$, there exists a finite cover $\tilde{f}\colon Y\rightarrow X$ from another smooth projective curve $Y$ over $K$ such that $f\colon V\rightarrow U$ is precisely the pullback of $\tilde{f}^{\mathrm{ad}}\colon  Y^{\mathrm{ad}}\rightarrow X^{\mathrm{ad}}$ along $U\subset X^{\mathrm{ad}}$.
\end{theorem} 

 Theorem~\ref{thm: spread-out thm} relies on the following approximation result (Proposition~\ref{prop: p_adic_Runge}). This is closely related to a result of Raynaud \cite[Corollaire 3.5.2]{Raynaud} (in the setting of Berkovich curves), % 
which he refers to as the \emph{$p$-adic Runge's Theorem}, as it can be viewed as an analogue of Runge's theorem in the theory of complex Riemann surfaces. Although we only require the case of curves, we present the result in full generality.

\begin{proposition}\label{prop: p_adic_Runge}
Let $X$ be a smooth projective variety over $K$ and let $K(X)$ denote the field of rational functions on $X$. Let $\xi$ be a point in $X^{\mathrm{ad}}$ and let $U \subset X^{\mathrm{ad}}$ be a connected affinoid open neighborhood of $\xi$. Then, after shrinking $U$ if necessary, $K(X) \cap \Gamma (U, \mO_{X^{\mathrm{ad}}})$  is dense in $\Gamma (U, \mO_{X^{\mathrm{ad}}})$ with respect to the topology induced from the spectral norm on $U$. Here the intersection is taken in the fractional field of $\Gamma (U, \mO_{X^{\mathrm{ad}}})$. 
\end{proposition}

\begin{proof} The idea is to reduce to the case of the projective spaces $\P^n_K$. Let $n$ be the dimension of $X$ over $K$. Let us pick  a finite surjective cover  $f\colon  X \ra \P^n_{K}$ such that $f$ is smooth (thus finite \'etale) on an affinoid neighborhood containing $\xi$. To see such a cover exists, one can embed $X \subset \P_K^N$ as a closed subscheme where $N > n$, and consider a linear subspace $L \subset \P_K^{N}$ of complementary dimension (so $L$ is isomorphic to a copy of $\P_K^{N-n-1}$) such that $L$ is disjoint from $X$. The projection map from $\P_K^{N}$ that projects ``away from $L$'' determines a finite map $f\colon  X \ra \P_K^n$. We will arrange $L$ such that the singular locus does not contain $\xi$. To see that this is possible, let $\xi$ be a point of (cohomological) dimension $d_0$ and let $Z$ be a $d_0$-dimensional irreducible closed subscheme of $X$ that contains $\xi$ (i.e., such that $\xi \in Z^{\mathrm{ad}}$). We just need to arrange $f$ so that it is generically smooth on $Z$. Indeed, we pick a closed point $z \in Z$ and arrange the linear subspace $L$ so that it does not intersect the tangent space $T_z X$ of $X$ inside $\P_K^N$ (by comparing dimensions). The resulting finite map $f\colon  X \ra \P_K^n$ is then smooth at the point $z \in Z$, so it is generically smooth on $Z$, as desired. 

Now let $U$ be an affinoid open neighborhood of $\xi$ in $X^{\mathrm{ad}}$. By shrinking $U$, we may assume that $f$ is finite \'etale on $U$. Let $\eta = f(\xi)$ and choose a connected affinoid open neighborhood $W$ of $\eta$. We may assume that $W$ is contained the analytification of an affine open subscheme of $\P^n_{K}$. Let $V := f^{-1} (W)$. Then $V$ is affinoid as $f$ is finite. We will first prove the density result on $V$ (after shrinking $W$ around $\eta$ if necessary), and then pass to a smaller neighborhood of $\xi$ contained in the neighborhood $U$ that we started with. To this end, 
consider the following commutative diagram
\[ 
\begin{tikzcd}
K(\P_K^n) \cap \Gamma(W, \mO_{\P_K^{n,\mathrm{ad}}}) \arrow[r, "\alpha"] \arrow[d] & \Gamma(W,  \mO_{\P_K^{n, \mathrm{ad}}}) \arrow[d]    \\     K(X) \cap \Gamma(V, \mO_{X^{\mathrm{ad}}}) \arrow[r, "\beta"]  &  \Gamma (V, \mO_{X^{\mathrm{ad}}}). 
\end{tikzcd}
\] 
We caution the reader that $V$ might be disconnected, in which case $K(X) \cap \Gamma(V, \mO_{X^{\mathrm{ad}}})$ stands for the subring of analytic functions on $V$ obtained by restricting rational functions on $X$. We claim that, if  the map $\alpha$ has dense image, then $\beta$ also has dense image. 
To prove the claim, let us consider an affine open subscheme $W_0 = \spec A_0 \subset \P^n_K$ such that $W \subset W_0$, and let $U_0 = W_0 \times_{\P_K^n} X$ be the pullback scheme. We may assume that $W_0$ is small enough so that $U_0 = \spec B_0$ for a finite \'etale algebra $B_0$ over $A_0$ (and we identify $A_0$ as a sub-algebra of $B_0$). In particular, there exist elements $x_1, \ldots, x_m \in B_0$ which generate $B_0$ as a module over $A_0$. Let $B = \Gamma(V, \mO_{X^{\mathrm{ad}}})$ and $A = \Gamma(W, \mO_{\P_K^{n, \mathrm{ad}}})$. Then $B$ is generated by $x_1, \ldots, x_m$ as an $A$-module. By scaling, we may assume that $x_1, \ldots, x_m$ all have spectral norm $\le 1$ on $V$. Now, for any $b \in B$, we can write $b = \sum_{i=1}^m a_i x_i$ with $a_i \in A$. By assumption, for any $\epsilon > 0$, there exists $a'_i \in K(\P_K^n) \cap \Gamma(W, \mO_{\P_K^{n, \mathrm{ad}}})$ such that $|a_i - a'_i|_W < \epsilon$. Take \[b' = \sum_{i=1}^m a'_i x_i \in K(X) \cap \Gamma (V, \mO_{X^{\mathrm{ad}}}).\] Then $|b - b'|_V < \epsilon$ (notice that $|a_i-a'_i|_V\le |a_i-a'_i|_W$). This finishes the proof of the claim. By Lemma \ref{lemma:density_for_projective_spaces} below, we may further shrink the neighborhood $W$ of $\eta$ so that $\alpha$ has dense image. Therefore, for this choice of $W$, the map $\beta$ in the diagram above also has dense image. 

Next, let $U_0 \subset U \cap V$ be a rational subset of $V$ containing $\xi$, given by
\[U_0=\{ x \in V \: | \: |h_j(x)| \le |h(x)|, 1 \le j \le r\,; \: h(x) \ne 0 \};\]
namely, 
\[\Gamma(U_0, \mO_{X^{\mathrm{ad}}}) = B\langle\frac{h_1, \ldots, h_r}{h}\rangle\]
where $B =  \Gamma(V, \mO_{X^{\mathrm{ad}}})$ as above. Pick $a \in K^\times$ such that $|a|<|h(\xi)|$. Taking the intersection of $U_0$ with $\{x\in V \: | \: |h(x)|\ge |a| \}$, we may assume that $U_0$ has the form 
\[
U_0 = \left\{   x \in V \: | \: |h_j(x)| \le |h(x)|, 1 \le j \le r, |h(x)|\ge |a|\right\}.
\]
By the density result for $V$ proved above, there exist elements $q, q_1,\ldots, q_r \in   K(X) \cap \Gamma(V, \mO_{X^{\mathrm{ad}}})$ such that $|h-q|_{V} < |a|$ and $|h_j - q_j|_V < |a|$ for all $j = 1, ..., r$. Then we can identify $U_0$ as
\[
U_0 = \left\{x \in V \: | \: |q_j(x)| \le |q(x)|, 1 \le j \le r, |q(x)|\ge |a|\right\};
\]
namely,
\begin{equation}\label{eq:density_U}
\Gamma(U_0, \mO_{X^{\mathrm{ad}}}) = B \langle\frac{q_1, \ldots, q_r, a}{q}\rangle
\end{equation} 
Consider the inclusions of rings 
\[
 \Big(  K(X) \cap \Gamma(V, \mO_{X^{\mathrm{ad}}}) \Big) [\frac{1}{q}] \: \subset \:   K(X) \cap \Gamma(U_0, \mO_{X^{\mathrm{ad}}})  \: \subset \:   \Gamma(U_0, \mO_{X^{\mathrm{ad}}}). 
\]
By the description in (\ref{eq:density_U}), we know that $(K(X) \cap \Gamma(V, \mO_{X^{\mathrm{ad}}})) [\frac{1}{q}] $ is dense in $\Gamma(U_0, \mO_{X^{\mathrm{ad}}})$. Therefore, $K(X) \cap \Gamma(U_0, \mO_{X^{\mathrm{ad}}})$ is dense in $\Gamma(U_0, \mO_{X^{\mathrm{ad}}})$. Finally, by considering the connected component of $U_0$ that contains $\xi$, we arrive at the desired proposition. 
\end{proof}

The following lemma is used in the proof of the proposition above. 

\begin{lemma} \label{lemma:density_for_projective_spaces}
Let $U$ be a connected affinoid open neighborhood of some point $\xi \in \P^{n, \mathrm{ad}}_{K}$. Then, up to shrinking $U$, the intersection $K(\P_K^n) \cap \Gamma(U, \mO_{\P_K^{n, \mathrm{ad}}})$ is dense in $\Gamma(U, \mO_{\P_K^{n, \mathrm{ad}}})$ with respect to the spectral norm on $U$.
\end{lemma}

\begin{proof}
  We first place $U$ inside a standard $n$-dimensional unit polydisk $\mathbb{D}_K^n$ such that the coordinates $t_1, \ldots, t_n \in K(\P_K^n)$ are global rational functions on $\P_K^n$. Without loss of generality, we may assume that $U$ is a Laurent domain of the form \[U = \{z \in \mathbb{D}_K^n \: | \:| f_i(z)| \le a_i, |g_j(z)|\ge b_j\}\] for some $f_i, g_j \in \mO(\mathbb{D}_K^n)$ and $a_i, b_j\in p^{\Q}$. Now we may replace $f_i, g_j$ by $f_i', g_j'$ which are polynomials in $t_1, \ldots, t_n$ while keeping 
\[U = \{z \in \mathbb{D}_K^n \: | \:| f'_i(z)| \le a_i, |g'_j(z)|\ge b_j\}.\] In particular, these $f'_i$'s and $g'_j$'s are rational functions on $\P_K^n$, in which case, the claim becomes evident.
\end{proof}

To proceed, we also need the following perturbation lemma.

\begin{lemma}[$p$-adic perturbation of finite covers of affinoids]\label{lemma:p_adic_perturbation}
Let $Y$ be an affinoid smooth rigid analytic curve over $K$, viewed as an adic space. Let $\xi$ be a non-type I point on $Y$ and let $V = \spa (A, A^+)$ be an open neighborhood of $\xi$ that is a rational subset of $Y$. Let $V'$ be a finite cover of $V$ defined by a monic irreducible polynomial 
\[ F(Z) =  Z^d +  a_{d-1} Z^{d-1} + \cdots    +  a_0 \] where $a_i \in A$, which remains irreducible as a polynomial over $\mO_{Y, \xi}$. Then there exists sufficiently small $\epsilon \in \R_{> 0}$ such that the following holds: for any polynomial
\[ 
 G(Z) =  Z^d +  b_{d-1} Z^{d-1} + \cdots + b_0 
\]
with $b_i\in A$ such that $|a_i - b_i|_{V} < \epsilon$ for the spectral norm $|\cdot |_V$ on $V$, up to replacing $V$ by a smaller rational subset containing $\xi$ if necessary, $G(Z)$ defines the same finite cover $V' \ra V$.     
\end{lemma}

\begin{proof} It suffices to show that, as long as $\epsilon$ is sufficiently small, $F(Z)$ and $G(Z)$ define the same extension of the stalk $\mO_{Y, \xi}$. First suppose that $\xi$ has type I\!I, I\!I\!I, or I\!V. Then the natural map $\mO_{Y, \xi} \ra \mO_{Y, \xi}^{\wedge} = \sH (\xi)$ induces an equivalence between the categories of finite extensions of $\mO_{Y, \xi}$ and that of $\sH(\xi)$. This follows, for example, from considering the Berkovich stalk $\mO_{Y^{\mathrm{Ber}}, \xi}$ and the fact that it is \emph{quasicomplete} in the sense of \cite[Theorem 2.3.3, Proposition 2.4.1]{Berkovich_etale}, and using the fact that we have an identification  $\mO_{Y^{\mathrm{Ber}}, \xi}^{\wedge} = \mO_{Y, \xi}^{\wedge} = \sH (\xi)$  (see \S \ref{subsection: comparison_between_stalks}). Therefore, it suffices to show that for sufficiently small $\epsilon$, $F(Z)$ and $G(Z)$ define the same extension of $\sH(\xi)$. Note that the condition $|a_i - b_i|_{V} < \epsilon$ in particular implies that $|a_i - b_i|_{\xi} < \epsilon$. Thus, the claim follows from Krasner's lemma since $\sH(\xi)$ is a complete non-archimedean field.
% Actually insert the proof here, just in case referee asks.  

Now suppose that $\xi$ has type V and let $\xi_{\mathrm{I\!I}}$ denote the unique type I\!I point that specializes to $\xi$. Then $\xi_{\mathrm{I\!I}} \in V$. By Lemma~\ref{lemma: base_change_from_type2_to_type5}, the inclusion $\mO_{Y, \xi} \hookrightarrow \mO_{Y, \xi_{\mathrm{I\!I}}}$ induces an equivalence between finite extensions of $\mO_{Y, \xi}$ and that of $\mO_{Y, \xi_{\mathrm{I\!I}}}$. Given this claim, it suffices to check that $F(Z)$ and $G(Z)$ define the same extension over $\mO_{Y, \xi_{\mathrm{I\!I}}}$. Notice that the condition  $|a_i - b_i|_{V} < \epsilon$ implies that $|a_i - b_i|_{\xi_{\mathrm{I\!I}}} < \epsilon$. So we are reduced to the case of type I\!I points, and we are done.
\end{proof}

\begin{proof}[{Proof of Theorem~\ref{thm: spread-out thm}}]
By shrinking $U$ (and hence $V$) and passing to connected components of $V$, we may assume that \[\Gamma(V, \mO_{X^{\mathrm{ad}}}) = \Gamma(U, \mO_{X^{\mathrm{ad}}})[T]/F(T)\] for a monic polynomial $F(T) = \sum_{i=0}^l\beta_i T^i$ with coefficients $\beta_i \in \Gamma(U, \mO_{X^{\mathrm{ad}}})$ (in particular, $\beta_l=1$) such that $F(T)$ is irreducible in $\mO_{X^{\mathrm{ad}}, \xi}[T]$.  By Proposition~\ref{prop: p_adic_Runge}, upon further shrinking $U$ if necessary, for any $\epsilon > 0$, there exists $\sq \beta_i\in K(X)\cap \Gamma(U, \mO_{X^{\mathrm{ad}}})$ such that $|\beta_i-\sq \beta_i|_U<\epsilon$ for all $0\le i\le l-1$ (set $\sq \beta_l=1$). Here $|\cdot|_U$ stands for the spectral norm on $U$. Applying Lemma~\ref{lemma:p_adic_perturbation}, we deduce that for $\epsilon$ sufficiently small, up to further shrinking $U$, the polynomials $F(T)$ and $\sq F(T)=\sum_{i=0}^l \sq \beta_i T^i$ define the same finite \'etale cover $V\ra U$. Now, since $\sq F(T) \in K(X)[T]$, it defines a finite separable extension $L/K(X)$. Let $Y$ be the normalization of $X$ in $L$, so that we have $V \cong Y^{\mathrm{ad}} \times_{X^{\mathrm{ad}}} U$. 
\end{proof}

\vspace{0.1in}
%%% 7.6
\subsection{Open neighborhoods of type V points}\label{ss: nbhd type V}
\noindent
\vspace{0.1in}

\noindent Our proof of Theorem~\ref{thm: local main thm type V general} will be based on the existence of a nice toric chart on a neighborhood of $\xi$, which enables us to reduce to the case where $X$ is a suitable open subset of a disk (cf. \S \ref{ss: type 5}). To this end, we first introduce some convenient open neighborhoods of type V points. 

Let  $\delta \in \D_K$ be the type V point corresponding to the open disk $\D_K (0, < 1)$ and let $\delta_{\mathrm{I\!I}}$ be its associated type I\!I point, which is the Gauss point of the boundary circle. Let $X_{\sigma, f}$ be a neighborhood of $\delta$ as described in Notation \ref{notation:closed unit disk removing open unit disks}, where $\sigma \in p^\Q \cap (0, 1)$.  Let us introduce the following notation.

\begin{definition} \label{def:notation_of_nbhd_of_type_5_point}
Let $\mathcal{Y}$ be a smooth connected projective curve over $\mathcal{O}_K$. Let $Y$ be the adic generic fiber over $K$ and let $\mathcal{Y}_s$ be the special fiber over $k$. Let $\mathrm{sp}\colon Y\rightarrow \mathcal{Y}_s$ be the specialization map. Then the preimage of the generic point of $\mathcal{Y}_s$ under $\mathrm{sp}$ consists of a single type I\!I point $\xi_{\mathcal{Y}}$ on $Y$.
\begin{enumerate}
\item For every closed point $z\in \mathcal{Y}_s$, let $D_{z}$ denote the preimage $\mathrm{sp}^{-1}(z)$. Note that the interior $\mathring{D}_z$ of $D_z$ is isomorphic to an open unit disk\footnote{In fact, $\mathring{D}_z$ is a \emph{virtual open unit disk} (over $K$); namely, the disk does not necessarily contain a $K$-rational point. However, when $z$ is a $k$-point, the disk $\mathring{D}_z$ indeed contains a $K$-rational point, and hence it is a ``genuine'' open unit disk over $K$.} and $D_z$ is the disjoint union of $\mathring{D}_z$ with a type V point $\xi_z$.\footnote{Note that the closure of the type I\!I point $\xi_{\mathcal{Y}} \in Y$ precisely consists of $\xi_{\mathcal{Y}}$ and all of such $\xi_z$'s. }  We identify $\mathring{D}_z$ with an open unit disk with the coordinate $T_z$.\footnote{This amounts to choosing a center of the disk, which is  not unique.} 
\item For $\sigma\in (0,1)\cap p^{\mathbb{Q}}$, let $\mathring{D}_{z;\sigma}$ denote the half-open annulus $\{\sigma\le |T_z|<1\}$ inside $\mathring{D}_z$, and let \[D_{z;\sigma} \coloneqq \mathring{D}_{z;\sigma}\cup \{\xi_z\}.\]
\item Let $m\ge 1$. Let $z, z_1, \ldots, z_m$ be $m+1$ (ordered) distinct closed points in $Y_s$, and let $\sigma \in (0,1]\cap p^{\mathbb{Q}}$. We define $Y_{\sigma}(z, z_1, \ldots, z_m)$  
to be the open subset of $Y$ given by the conditions $|T_{z}|\ge \sigma$ and $|T_{z_i}| = 1$ for $i=1, \ldots, m$.\footnote{Roughly speaking, the condition ``$|T_{z}|\ge \sigma$'' means ``removing the open disk of radius $\sigma_i$'' inside $D_{z_i}$. (Here, if $\sigma_i=1$ for some $i$, the condition $|T_{z_i}|\ge 1$ simply means removing the entire $D_{z_i}$.)} When the points $z, z_1, ..., z_m$ are understood from the context, we shall simply write $Y_{\sigma}$ for $Y_{\sigma} (z, z_1, ..., z_m)$.
\end{enumerate}
\end{definition}

Now, let $X$ be a smooth rigid analytic curve over $K$ (viewed as an adic space) and let $\xi$ be a $K$-rational type V point on $X$. Let $\xi_{\mathrm{I\!I}}$ be the associated type I\!I point, which is  necessarily also $K$-rational.

\begin{lemma}\label{lemma: open nbhd type V--general}
\begin{enumerate}
\item[(1)]
After enlarging $K$ if necessary, there exists an open neighborhood $V$ of $\xi$ such that 
\begin{itemize}
\item $V = Y_{\sigma_0}(z,z_1, \ldots, z_m)$ for a smooth connected projective curve $\mathcal{Y}$ over $\mathcal{O}_K$, distinct closed points $z, z_1, \ldots, z_m$ on $\mathcal{Y}_s$, and a real number $\sigma_0\in (0,1)\cap p^{\mathbb{Q}}$ as in Definition~\ref{def:notation_of_nbhd_of_type_5_point};
\item $\xi_{\mathrm{I\!I}}=\xi_{\mathcal{Y}}$;
\item $\xi=\xi_z$ (in particular, $z$ is a $k$-point of $\mathcal{Y}_s$);\footnote{Since $z$ is a $k$-point of $Y_s$, the disk $\mathring{D}_z$ is a ``genuine'' open unit disk over $K$; namely, $\mathring{D}_z\cong \D_K(0,1^-)$, with coordinate $T_z$. Hence, $\mathring{D}_{z;\sigma}\cong A_K([\sigma,1); T_z)$ is a ``genuine'' half-open annulus over $K$ with radius $[\sigma,1)$.}
\end{itemize} 
Note that the special fiber $\mathcal{Y}_s$ of $\mathcal{Y}$ is uniquely determined by the imperfect residue field of the completed stalk $\mathscr{H}(\xi_{\mathrm{I\!I}})$.
\item[(2)] In the situation of (1), the collection
\[
\Psi=\left\{Y_{1}(z,z_1, \ldots, z_m, z'_1, \ldots, z'_l)\right\}
\]
form a basis of open neighborhoods of $\xi_{\mathrm{I\!I}}$ in $X$, where $z'_1, \ldots, z'_l$ are closed points of $\mathcal{Y}_s$ distinct from $z, z_1, \ldots, z_m$. Note that $Y_{1}(z,z_1, \ldots, z_m, z'_1, \ldots, z'_l)$ does not contain the type V point $\xi$.
\item[(3)] In the situation of (1), the collection
\[
\Psi'=\left\{Y_{\sigma}(z,z_1, \ldots, z_m, z'_1, \ldots, z'_l)\,;\,\,\sigma\in [\sigma_0,1)\cap p^{\Q}\right\}
\]
form a basis of open neighborhoods of $\xi$ in $X$, where $z'_1, \ldots, z'_l$ are closed points of $\mathcal{Y}_s$ distinct from $z, z_1, \ldots, z_m$.
\end{enumerate}
\end{lemma}

\begin{proof}
This follows from \cite[Proposition 2.2.1]{Berkovich_differential}.
\end{proof}

\vspace{0.1in}
%% 7.7
\subsection{The reduction step: a nice toric chart}\label{ss: type 5}
\noindent
\vspace{0.1in}

\noindent To prove Theorem~\ref{thm: local main thm type V general}, our strategy is to reduce the general situation to the case where $X= X_{\sigma}$ as in Theorem \ref{thm:main_type5_disk}.

Let $\xi$ be a $K$-rational type V point on $X$ and let $\xi_{\mathrm{I\!I}}$ be the associated type I\!I point. Up to enlarging $K$, we may assume that $\xi$ admits an open neighborhood of the form $V\cong Y_{\sigma_0}=Y_{\sigma_0} (z, z_1, ..., z_m)$ as in Definition~\ref{def:notation_of_nbhd_of_type_5_point}. By further enlarging $K$ and shrinking $V$ (that is, enlarging the set of closed points $\{z,z_1, \ldots, z_m\}$ and increasing $\sigma_0$) if necessary, there exists a toric chart
\[ h\colon  V =Y_{\sigma_0} \rightarrow \mathbb{D}_K\] such that the following holds.  
\begin{enumerate}
\item There exist $c_1, \ldots, c_n\in \overline{k}_K^{\times}$ and $\sigma \in(0,1)\cap p^{\mathbb{Q}}$, such that $h$ factors as  
\[ h\colon  V \rightarrow X_{\sigma, f} \hookrightarrow\mathbb{D}_K\] for some $f = \prod_{i=1}^n (T - [c_i])\in \mathcal{O}_{K_0}[T]$, where the first map is finite \'etale; 
\item $\delta_{\mathrm{I\!I}}=h(\xi_{\mathrm{I\!I}})$ is the Gaussian point of $\mathbb{D}_K$, and $\delta=h(\xi)$ is the type V point in $\mathbb{D}_K$ corresponding to the open unit disk centered at 0. Moreover, we have $h^{-1}(\delta_{\mathrm{I\!I}})=\{\xi_{\mathrm{I\!I}}\}$ and $h^{-1}(\delta)=\{\xi\}$.
\item $h^{-1}\left(A_K[\sigma,1)\right)=\mathring{D}_{z;\sigma_0}\cong A_K([\sigma_0,1);T_z)$; 
\item $h^{-1} \left( X_{1, f} \right) =  Y_1 \subset Y_{\sigma_0}$; 
\end{enumerate}

\noindent We remark that, by enlarging $K$ if necessary, the affinoid neighborhood 
\[ U \coloneqq  Y_1 = Y \minus   (D_{z} \sqcup \cdots \sqcup D_{z_m})
\]
of $\xi_{\mathrm{I\!I}}$ can be realized as the adic generic fiber of a smooth affine formal scheme $\mathfrak U$,  where the latter is the completion of $\mathcal Y$ along the locally closed subscheme $\mathcal{Y}_s - z, z_1, \ldots, z_m\}$. Likewise, $X_{1, f}$ is the adic generic fiber of a smooth affine formal scheme $\mathfrak D_{f}$, obtained as the completion of the $p$-adic formal scheme $\widehat \A^1_{\mO_K}$ along $\A^1_{k} - \{0,c_1, \ldots, c_n\}$. By shrinking $X_{1, f}$ (and hence also shrinking $V$) if necessary, we may further assume that $\mathfrak D_{f}$ admits a formal model $\mathfrak{D}^{\circ}_{f}$ over $W(k)$.

Using the toric chart $h$, we construct another toric chart that behaves well at $\xi_{\mathrm{I\!I}}$.

\begin{proposition} \label{prop:reduction_step_fet_map}
By enlarging $K$ and shrinking both $U$ and $V$ (that is, enlarging the set of closed points $\{z,z_1, \ldots, z_m\}$ and increasing $\sigma_0$) if necessary, there exist $c'_1, \ldots, c'_l \in \overline{k}_K^{\times}$, $\sigma' \in (0, 1] \cap p^{\Q}$, and a finite \'etale map 
\[ g\colon  V=Y_{\sigma_0}\ra  X_{\sigma', f'} \] 
where $f' = \prod_{i =1}^{l} (T - [c_i'])\in \mathcal{O}_{K_0}[T]$, such that 
\begin{enumerate}
\item[(1)] $g(\xi_{\mathrm{I\!I}})$ is the Gaussian point of $\mathbb{D}_K$, and $ g(\xi)$ is the type V point in $\mathbb{D}_K$ corresponding to the open unit disk centered at 0. Moreover, we have $g^{-1}(\delta_{\mathrm{I\!I}})=\{\xi_{\mathrm{I\!I}}\}$ and $g^{-1}(\delta)=\{\xi\}$;
\item[(2)] $g^{-1}\big(X_{1, f'}\big)=U= Y_1$;
\item[(3)] the restriction $g|_U\colon  U \ra  X_{1, f'}$ is the adic generic fiber of a finite \'etale map  
$\widehat g\colon  \mathfrak U \ra \mathfrak D_{f'}$ 
between smooth affine formal schemes. 
\end{enumerate}
\end{proposition} 

\begin{proof}
We will first construct such a map on a neighborhood $U$ of $\xi_{\mathrm{I\!I}}$ satisfying condition (3), then show that it extends to a neighborhood $V$ of $\xi$, after possibly enlarging the base field $K$.\\

\noindent \textbf{Step 1} (Factor $h$ on the integral models).

\vspace{0.1in}

\noindent Recall the finite \'etale map $h\colon  V=Y_{\sigma_0} \rightarrow X_{\sigma, f}$  which restricts to a finite \'etale map 
$
h\colon  U =Y_1 \rightarrow  X_{1, f}$ from the beginning of this subsection.  
By the theory of Raynaud \cite{Raynaud74}, up to shrinking $U$ and $X_{1, f}$ if necessary, there exists a finite map 
$\widehat{h}\colon  \mathfrak U \ra \mathfrak D_{f}$ 
that gives rise to $h$ upon passing to adic generic fibers (\emph{a priori}, we only have such a map up to taking admissible blowing-ups, but this can be avoided once we are allowed to shrink the source and target). Passing to special fibers, we obtain
\[\widehat{h}_s\colon  \mathfrak U_s = \mathcal{Y}_s \backslash \{z, z_1, \ldots, z_m\} \ra \mathfrak D_{f; s} = \A^1_{k} \backslash \{0, c_1, \ldots, c_n\},\] 
which extends to a finite map $\mathcal{Y}_s \ra \P^1_k$ between smooth projective curves over $k$. This map factors as 
\begin{equation}\label{eq:factorization h_s}
\mathcal{Y}_s \xrightarrow{\phi_s} \mathcal{Y}_s^{(N)} \xrightarrow{\:\psi_s \:} \P^1_k 
\end{equation}
where $\mathcal{Y}_s^{(N)} = \mathcal{Y}_s \times_{\spec k, \mathrm{Frob}^{N}} \spec k$ is the base change of $\mathcal{Y}_s$ along the $p^N$-th power Frobenius map on $k$, $\phi_s$ denotes the $p^N$-th power relative Frobenius map, and $\psi_s$ is a finite separable map. We would like to conclude that (\ref{eq:factorization h_s}) restricts to a factorization
$\mathfrak{U}_s\xrightarrow{\phi_s}\mathfrak{U}^{(N)}_s \xrightarrow{\:\psi_s \:}\mathfrak{D}_{f; s}$ of $\widehat{h}_s$. However, the open subscheme $\mathfrak U_s \subset \mathcal{Y}_s$ need not be stable under the $p^N$-th power Frobenius map. To remedy this issue, consider
\[
\mathfrak U'_s = (\phi_s)^{-1} (\mathfrak U_s^{(N)}) \cap \mathfrak U_s=\mathcal{Y}_s\backslash \left((\phi_s)^{-1} (\{z, z_1, \ldots, z_m\}) \cup \{z, z_1, \ldots, z_m\}\right).
\]
This corresponds to an open affine formal subscheme $\mathfrak U' \subset \mathfrak U$
and a rational subset $U' = \mathrm{sp}^{-1} (\mathfrak U_s')\subset U$ on the adic generic fiber. By further shrinking $X_{1, f}$ and shrinking $U'$ and $U$ accordingly if necessary, we can ensure that $\widehat h_s|_{\mathfrak{U}'_s}$ factors as 
\[
\mathfrak U_s' \xrightarrow{\phi_s} \mathfrak U_s^{(N)} \xrightarrow{\:\psi_s\:} \mathfrak D_{f; s},
\] 
where the second map $\psi_s$ is finite \'etale. \\ 

\noindent \textbf{Step 2} (Lift the finite \'etale map $\psi_s$).

\vspace{0.1in}

\noindent After untwisting Frobenius, we have a finite \'etale map  
$\psi_s\colon  \mathfrak{U}_s\rightarrow \mathfrak D_{f; s}^{(-N)}$, which we still denote by $\psi_s$ by a slight abuse of the notation. 
It lifts to a finite \'etale map
\[ \widehat \psi\colon  \mathfrak U \ra \mathfrak D_{f}^{(-N)}  \coloneqq  \mathfrak D^{\circ}_{f} \times_{\spf W(k), \mathrm{Frob}^{-N}} \spf \mO_K\]
and thus gives rise to a finite \'etale map 
$\psi\colon  U \ra X_{1, f}^{(-N)}$  
on the adic generic fibers. Here $X_{1, f}^{(-N)}$ denotes the adic generic fiber of $\mathfrak D_{f}^{(-N)}$; in particular, we have
\[ 
X_{1, f}^{(-N)}\cong X_{1, f'}
\]
with $f' := \prod_{i} (T - [c_i'])$, where $c_i' = \mathrm{Frob}^{-N}(c_i)$. 
Thus we obtain a finite \'etale map 
$
\widehat \psi\colon  \mathfrak U \ra \mathfrak D_{f'}=\mathfrak D_{f}^{(-N)}$ such that, on the adic generic fibers, it gives rise to a finite \'etale map 
\begin{equation} \label{eq:the_map_psi_over_U}
     \psi\colon  U \ra X_{1, f'}. 
\end{equation} 
This map sends $\xi_{\rm{I\!I}}$ to the Gaussian point $\delta_{\rm{I\!I}}$ of the disk.  \\

\noindent \textbf{Step 3} (Construction of the desired map).

\vspace{0.1in}

\noindent 
Next, we show that, up to enlarging $K$, shrinking $U, V$, and re-centering the unit disk $\D_K$ (and thus modifying the $c_i'$'s) if necessary, there exists a finite \'etale map
\[g \colon  V \ra  X_{\sigma', f'} \]
for some $\sigma'\in (0,1)\cap p^{\Q}$, which satisfies the requirement of the proposition. In particular, it sends the type V point $\xi$ to the type V point $\delta$ on $\mathbb{D}_K$ corresponding to the open unit disk centered at $0$.  We will construct this map $g$ from $\psi$ in (\ref{eq:the_map_psi_over_U}).\footnote{Ideally we would like to simply ``extend'' $\psi$ to such a finite \'etale map $g$. However, in the construction of $g$ we will invoke elementary deformation theory of smooth affine schemes in characteristic $p$. As a result, the restriction of $g$ to $U$ may differ from $\psi$ by an automorphism. } 
 
First, using the projective coordinate on $\P^1_K$, we identify $X_{1, f'}\subset \D_K\subset \P^1_K$ as an affinoid subspace. By Theorem~\ref{thm: spread-out thm}, up to shrinking $U$ and  $X_{1, f'}$, there exists a smooth projective curve $X‘$ over $K$ and a finite map $f'\colon  X' \ra \P^1_K$ such that the restriction of $f'$ to (the preimage of) $X_{1, f'}$ is given by $\psi$. Let us keep track of the identification 
\[ i: U \isom (f')^{-1} (X_{1, f'}) \subset X'
\] and denote the image by $U_{\tu{aux}} = i (U)$. 
Note that $f' \circ i = \psi$ by construction. 
In the rest of this proof, we write $X'_{\mathscr B}$ for the Berkovich space associated to $X'$ and write $X'_{\mathrm{ad}}$ for the associated adic space. Write $\xi_{\mathrm{I\!I}}' := i (\xi_{\mathrm{I\!I}}) \in U_{\tu{aux}}$, which is the image of $\xi_{\mathrm{I\!I}}$ in $ X'_{\mathrm{ad}}$. 

Next, after enlarging $K$ and shrinking $U$,  $V$ if necessary, we will construct an embedding $V \hookrightarrow X_{\mathrm{ad}}'$, which realizes $V$ as an affinoid open neighborhood inside $X'_{\mathrm{ad}}$, while keeping track of the effect on the formal model of $U$. 
To proceed, note that  by \cite[Proposition 2.2.1] {Berkovich_differential}, up to enlarging $K$ if necessary, the point $\xi_{\mathrm{I\!I}}'$ admits a Berkovich neighborhood in $X'_{\mathscr B}$ together with an isomorphism   
\[ \lambda:   Y'_{\mathscr B}  \minus \sqcup_{i \in I} B_i' \isom  V_{\mathscr B} \subset  X'_{\mathscr B},\] where $\mathcal Y'$ is a smooth connected projective curve over $\mO_K$ and $\{B_i'\}_{i \in I}$ is a finite collection of closed disks (in distinct residue disks) on the Berkovich generic fiber $ Y'_{\mathscr B}$ of $\mathcal Y'$. Write $D'_i \subset  Y'_{\mathscr B}$ for the residue disk containing $B'_i$ and let ${D'}_{i, \mathrm{ad}}$ denote the associated adic open disk. Up to shrinking $U$ if necessary, we may assume that  
\[ U_{\mathrm{aux}} = i(U) = \lambda \Big(  Y'_{\tu{ad}} \minus \sqcup_{i \in I} {D'}_{i, \mathrm{ad}} \Big). 
\]   
Let  $V' = \pi^{-1} (V_{\mathscr B})$ denote the pre-image of $V_{\mathscr B}$ under the natural projection $\pi\colon  X'_{\mathrm{ad}} \ra X'_{\mathscr B}$. Note that $V'$ contains all type V points on $X'_{\mathrm{ad}}$ that live in the closure of $\xi_{\mathrm{I\!I}}'$. By construction, we have 
$U_{\mathrm{aux}} \subset V' \subset X'_{\mathrm{ad}}.$ 
Let $\mathfrak U_{\mathrm{aux}}$ denote the smooth affine formal model of $ Y'_{\tu{ad}} \minus \sqcup_{i \in I} {D'}_{i, \mathrm{ad}}$ obtained from $\mathcal Y'$. We view it as a formal model of $U_{\mathrm{aux}}$ via the isomorphism $\lambda$. The isomorphism $i$ uniquely determines an isomorphism 
\[
\widehat i: \mathfrak U \isom \mathfrak U_{\mathrm{aux}}. 
\]
One of the difficulty in the construction of the embedding $V \hookrightarrow X_{\mathrm{ad}}'$ is that, $Y'_{\mathrm{ad}}$ is not necessarily isomorphic to $Y$. However, the special fibers $\mathcal{Y}_s$ and $\mathcal{Y}'_s$ are isomorphic, as they are determined by the imperfect residue field of the completed stalk $\mathscr H(\xi_{\mathrm{I\!I}})$. Thus, the map on special fibers induced by $\widehat i$ uniquely extends to an isomorphism 
$\iota: \mathcal{Y}_s \isom  \mathcal{Y}'_s.$  
Let $x_0 \in \mathcal{Y}_s$ be a closed point away from $\{z, z_1, \ldots, z_m\}$ and let $x_0' = \iota(x_0)\in \mathcal{Y}'_s$ be its image under $\iota$. Let $D_0 = \mathrm{sp}^{-1} (x_0) \subset Y$ be the residue disk of $x_0$ and similarly define $D_0'$. The map we shall construct will be defined on a subset of $Y \minus D_0$. For this, we want to remove $D_0$ from $Y$, but this might destroy finite \'etaleness over $X_{1, f'}$. In order to prevent this, we may need to shrink further. More precisely, after enlarging $K$ if necessary, write $b_0 := \psi_s (x_0) \in \mathfrak D_{f'; s}$ and let $D_{b_0} = \mathrm{sp}^{-1} (b_0) \subset X_{1, f'}$ be the residue disk of $b_0$. We shall replace $X_{1, f'}$ by $X_{1, f'} \minus D_{b_0}$ and replace $U$ by $U \minus \psi^{-1} (D_{b_0})$. This in turn leads to replacing $V$ by $V \minus \psi^{-1} (D_{b_0})$ and removing corresponding residue disks from $Y$ (using the isomorphism $\iota$). Note that the removed residue disks on $Y$ include $D_0$. From now on, we continue to denote the spaces by $U, V, U_{\mathrm{aux}}$ after this shrinking.

Let $\mathscr Y$ and $\mathscr Y'$ be the formal completions of $\mathcal Y$ and $\mathcal Y'$ along $\mathcal{Y}_s \minus \{x_0\}$ and $\mathcal{Y}'_s \minus \{x_0'\}$ respectively.  
Since the formal deformation of the smooth affine curve $\mathcal{Y}_s \minus \{x_0\}$ is unobstructed and in fact trivial, we may choose an isomorphism 
\[
\widehat{j}_0: \mathscr Y \isom \mathscr Y',
\]
which further induces   an isomorphism 
\[
 j_0: Y \minus D_0 \isom  Y' \minus D_0'
\]
of affinoid adic spaces on the generic fiber.  
By construction, $\widehat j_0$ restricts to an isomorphism $\widehat j_U: \mathfrak U \isom \mathfrak U'$, and we have $j_0 (U) = \lambda^{-1} (U_{\tu{aux}}) =  Y'_{\tu{ad}} \minus \sqcup_{i \in I} {D'}_{i, \mathrm{ad}}$, since $\widehat i$ and $\widehat j_U$ reduce to the same isomorphism on the special fiber. 
By further shrinking $V$ if necessary (i.e., increasing $\sigma_0$), we have $\lambda \circ j_0 (V) \subset V'  \subset X'_{\mathrm{ad}}$. Write $j := \lambda \circ j_0: V \ra    X'_{\mathrm{ad}}$ for the resulting embedding. 

Finally, after enlarging $K$, shrinking $V$ (i.e., enlarging $\sigma_0$), and choosing a different coordinate of $\D_K$ if necessary, we obtain the desired finite \'etale map 
$g \coloneqq  f' \circ j \colon  V \ra  X_{\sigma', f'}$ 
that sends $\xi$ to the type V point $\delta$ corresponding to the open unit disk in $\D_K$ centered at 0. This finishes the proof of the proposition. 
\end{proof}

\vspace{0.1in}
\subsection{Proof of $p$-adic monodromy for type V points}
\noindent
\vspace{0.1in}

\noindent Finally, we prove Theorem \ref{thm: local main thm type V general}.

\begin{proof}[Proof of Theorem \ref{thm: local main thm type V general}]
Let $\xi$ be the type V point as in the statement of the theorem and let $\xi_{\mathrm{I\!I}}$ be the associated type I\!I point. Recall that both $\xi$ and $\xi_{\mathrm{I\!I}}$ are $K$-rational. We want to show that, up to a finite extension of $K$, there exist an open neighborhood $W$ of $\xi$ and a finite \'etale cover $W'\rightarrow W$ such that the pullback of $\L|_W$ to $W'$ is pointwise semistable. By Proposition~\ref{prop: local main thm type II}, we may assume that $\L$ is already semistable in an open neighborhood of $\xi_{\mathrm{I\!I}}$.

Consider the map 
\[
g\colon  V \ra X_{\sigma', f'}
\]
given by Proposition~\ref{prop:reduction_step_fet_map}. In particular, we have $g^{-1}(\delta_{\mathrm{I\!I}})=\{\xi_{\mathrm{I\!I}}\}$ and $g^{-1}(\delta)=\{\xi\}$ where $\delta_{\mathrm{I\!I}}$ is the Gaussian point of $\D_K$ and $\delta$ is the type V point in $ X_{\sigma', f'}$ that  corresponds to the open unit disk in $\D_K$ centered at 0. Consider the $p$-adic local system $g_*\L$ on $ X_{\sigma', f'}$. By assumption (2) in Proposition~\ref{prop:reduction_step_fet_map}, $g$ restricts to $U\rightarrow  X_{1, f'}$ which admits a finite \'etale integral model $\mathfrak{U}\rightarrow \mathfrak{D}_{f'}$. Since $\L$ is assumed to be semistable in a neighborhood of $\xi_{\mathrm{I\!I}}$, we deduce that $g_*\L$ is semistable in a neighborhood of $\delta_{\mathrm{I\!I}}$ by applying \cite[Theorem 1.6]{GY}.

By Theorem \ref{thm:main_type5_disk}, up to enlarging $K$, up to further shrinking $X_{\sigma', f'}$ to a smaller neighborhood of $\delta$, there exists a finite \'etale cover $\psi : X_{\sigma', f'}'  \ra  X_{\sigma', f'}$ such that the local system $(g_*\L)|_{X_{\sigma', f'}}$ becomes pointwise semistable after pulling back along the map $\psi$. Let $W \coloneqq g^{-1}( X_{\sigma', f'})$ and $W' \coloneqq W\times_{ X_{\sigma', f'}}  X_{\sigma', f'}'$. They fit into a commutative diagram
\[
\begin{tikzcd}
    W' \arrow[r, "\psi'"] \arrow[d, "g'"] & W \arrow[d, "g"] \\ 
    X_{\sigma', f'}'  \arrow[r, "\psi"] &  X_{\sigma', f'}. 
\end{tikzcd}
\]
Then $(\psi')^* (g^* g_* \L)=(g')^* \psi^* (g_* \L)$ is pointwise semistable on $W'$. Note that for any classical point $x\in W$, the restriction $\L|_x$ is a direct summand of $(g^*g_*\L)|_x$ as $\Gal_{\kappa(x)}$-representations. It follows that $(\psi')^* \L$ is also pointwise semistable.
This finishes the proof of Theorem \ref{thm: local main thm type V general} (that is, the general case of Theorem~\ref{thm: local main thm} for type V points). 
\end{proof}

\begin{remark}
Combining the results established in \S \ref{ss: type 1}, \S \ref{ss: type 2},  \S \ref{section: type 3}, \S \ref{section: type 4}, \S \ref{section: type 5}, and \S \ref{section: type V general}, the proof of Theorem~\ref{thm: local main thm} is now complete.
\end{remark}

\vspace{0.1in}
%% 7.9
%
\section{Proof of the main theorem} \label{section: proof of the main theorem}

 Finally, it is time to finish the proof of the main theorem.

\begin{proof}[Proof of Theorem~\ref{thm:main_intro}]
Applying Theorem~\ref{thm: local main thm} and Theorem~\ref{thm: spread-out thm} to each (adic) point $\xi\in X_{\C_K}^{\mathrm{ad}}$, we obtain a triple $(L(\xi), U_{\xi}, f_{\xi}\colon Y_{\xi}\rightarrow X_{L(\xi)})$ where 
\begin{itemize}
\item $L(\xi)$ is a finite extension of $K$;
\item $U_{\xi}\subset X_{L(\xi)}$ is an open neighborhood of the image of $\xi$ in $X_{L(\xi)}$;
\item $f_{\xi}\colon Y_{\xi}\rightarrow X_{L(\xi)}$ is a finite cover between smooth projective curves over $L(\xi)$;
\end{itemize}
such that the pullback of $\mathbb{L}$ to $Y^{\mathrm{ad}}_{\xi}\times_{X_{L(\xi)}^{\mathrm{ad}}}U_{\xi}$ is pointwise semistable. Notice that $\{U_{\xi, \C_K}\}$ cover $X_{\C_K}^{\mathrm{ad}}$. By compactness, $X_{\C_K}^{\mathrm{ad}}$ is covered by a finite sub-collection $\{U_{\xi, \C_K}\}_{\xi\in \mathscr{C}}$. Pick a finite extension $L/K$ containing $L(\xi)$ for all $\xi\in \mathscr{C}$, so that $\{U_{\xi, L}\}_{\xi\in \mathscr{C}}$ cover $X_L^{\mathrm{ad}}$.

Now, consider the normalization $Y$ of the fiber product of all $Y_{\xi, L}$'s over $X_L$. By construction, the pullback of $\mathbb{L}$ to $Y$ is pointwise semistable. By enlarging $K$, we may further assume that $Y$ has semistable reduction. It then follows from the pointwise criterion (Theorem~\ref{thm:pointwise_criterion}) that the pullback of $\mathbb{L}$ to $Y$ is semistable.
\end{proof}

%%%%%%%%%%%%%%%%%%%%
%%%%%%%%%%%%%%%%%%%%

%%%%%%%%%%%%%%%%%%%%
%%%%%%%%%%%%%%%%%%%%

%%%%%%%%%%%%%%%%%%%%
%%%%%%%%%%%%%%%%%%%%

%%%%%%%%%%%%%%%%%%%%
%%%%%%%%%%%%%%%%%%%%

%%%%%%%%%%%%%%%%%%%%
%%%%%%%%%%%%%%%%%%%%

\vspace{0.1in}

\appendix  

\section{Notations for period rings} \label{appendix:notation}
In this appendix, we gather a list of notations of period rings that are commonly used in $p$-adic Hodge theory. As the table below suggests, for the period rings in the classical setting (cf. \S \ref{sec:classical_phi_Gamma}), our notation is in the same spirit as in the work of Kedlaya--Liu (\cite{Kedlaya-LiuI}, \cite{Kedlaya-LiuII}). We make one slight modification to their notation --- we write $\widetilde{\mathbf{C}}^{\dagger}_K$ and ${\mathbf{C}}^{\dagger}_K$ for the rings denoted by $\widetilde{\mathbf{C}}_K$ and ${\mathbf{C}}_K$ in \cite{Kedlaya-LiuI, Kedlaya-LiuII} to emphasize that these rings are obtained by taking the union over $r>0$, as in Definitions \ref{def:perfect_robba_ring_type_C_over_K} and \ref{definition:A^r_as_intersection}, corresponding to ``approaching the boundary of the annulus''.

In the table below, let $s, r$ be real numbers with $0<s\le r$ and $s,r\in p^{\mathbb{Q}}$. Put $r'=\frac{1}{r}\cdot\frac{p-1}{p}$ and $s'=\frac{1}{s}\cdot\frac{p-1}{p}$. 
Moreover, a dash ``---'' indicates that the notation there agrees with the notation used in this paper. We also ignore differences in font choices (for example, $\mathscr{R}$ versus $\mathcal{R}$).

\vskip 0.2in

\newpage 
%\begin{landscape}
 
\begin{center} 
\scalebox{0.85}{
   \addtolength{\leftskip} {-2cm} % increase (absolute) value if needed
  \addtolength{\rightskip}{-2cm}
{\renewcommand{\arraystretch}{1.6}%
\begin{tabular}{ |c|c|c|c|c|c |c|} 
 \hline
 Our choice %&  \cite{Kedlaya-slope-filtrations, Fontaine_90}   
 & \cite{Cherbonnier_Colmez}  &  \cite{Berger-differential} & \cite{Colmez08} & \cite{Kedlaya-LiuI, Kedlaya-LiuII}    \\ 
   \hline  
 $\widetilde{\mathbf{A}}_K$ & 
 %--- & 
 --- &
   
 --- &
 %$\widetilde{\mathbf{A}}_K   $ & 
 --- &
 %$\widetilde{\mathbf{A}}_K$ & 
 ---  
 %$\widetilde{\mathbf{A}}_K$  
 \\
 $\widetilde{\mathbf{A}}^r_K$ & 
% $\sq{\mathscr R}_K^{\mathrm{int}, r}$, $(\sq{\mathscr R}^{\mathrm{int}, r})^{H_K}$ & 
 $\widetilde{\mathbf{A}}^{\dagger}_{{r'}^{-}\!, K}, (\widetilde{\mathbf{A}}^{\dagger}_{{r'}^{-}})^{H_K}$ &
  
 $\widetilde{\mathbf{A}}^{\dagger, r'}_K,   \widetilde{\mathbf{A}}^{H_K}_{[r'; +\infty]}$ & $\widetilde{\mathbf{A}}^{(0,r]}_K$ & 
 --- 
 \\
 $\widetilde{\mathbf{A}}^{\dagger}_K$ &
%  $\sq{\mathscr R}_K^{\mathrm{int}}$ & 
  --- &
  --- &
% $\widetilde{\mathbf{A}}^{\dagger}_K$ & 
  --- &
%$\widetilde{\mathbf{A}}^{\dagger}_K$ & 
  ---  
%$\widetilde{\mathbf{A}}^{\dagger}_K$ 
 \\
 \hline
 $\widetilde{\mathbf{B}}_K$ &
% --- & 
 ---  &
 ---  &
% $\widetilde{\mathbf{B}}_K$ & 
 ---  &
%$\widetilde{\mathbf{B}}_K$ & 
 ---
%$\widetilde{\mathbf{B}}_K$ 
 \\
 $\widetilde{\mathbf{B}}^r_K$ &
% $\sq{\mathscr R}_K^{\mathrm{bd}, r}$, $(\sq{\mathscr R}^{\mathrm{bd}, r})^{H_K}$  & 
  $\widetilde{\mathbf{B}}^{\dagger}_{{r'}^{-}\!, K}, (\widetilde{\mathbf{B}}^{\dagger}_{{r'}^{-}})^{H_K}$ &
 $\widetilde{\mathbf{B}}^{\dagger, r'}_K, \widetilde{\mathbf{B}}^{H_K}_{[r'; +\infty]}$ & $\widetilde{\mathbf{B}}^{(0,r]}_K$ & 
 ---
 \\
 $\widetilde{\mathbf{B}}^{\dagger}_K$ & 
% $\sq{\mathscr R}_K^{\mathrm{bd}}$ & 
 ---  & 
 ---  & 
% $\widetilde{\mathbf{B}}^{\dagger}_K$ & 
 ---  &
 %$\widetilde{\mathbf{B}}^{\dagger}_K$ & 
 ---  
 %$\widetilde{\mathbf{B}}^{\dagger}_K$ 
 \\
 \hline
 $\widetilde{\mathbf{C}}^{[s,r]}_K$ & 
 % n/a & 
 n/a &
 $\widetilde{\mathbf{B}}_{[r' ;s']}^{H_K}$ &  
 n/a & 
 --- 
 \\
 $\widetilde{\mathbf{C}}^r_K$ &
 %  $\sq{\mathscr R}_K^{r}\;,  \: $  $(\sq{\mathscr R}^{r})^{H_K} \: \:$ & 
 n/a &
 $\widetilde{\mathbf{B}}^{\dagger, r'}_{\mathrm{rig}, K}, \widetilde{\mathbf{B}}^{H_K}_{[r'; +\infty[}$ & $\widetilde{\mathbf{B}}^{]0,r]}_K$ & 
 ---  
 \\
 $\widetilde{\mathbf{C}}^{\dagger}_K$ &
%   $\sq{\mathscr R}_K$ & 
 n/a &
 $\widetilde{\mathbf{B}}^{\dagger}_{\mathrm{rig}, K}$ & $\widetilde{\mathbf{B}}^{\dagger}_{\mathrm{rig}, K}$ & $\widetilde{\mathbf{C}}_K$ 
% \\   0 & 1 & 2 & 3 & 4 & 5
 \\
 \hline
 \hline
 $\mathbf{A}_K$ & 
%   --- & 
 --- &
 --- & 
 %$\mathbf{A}_K$ & 
 --- &
 %$\mathbf{A}_K$ & 
 --- 
 \\
 $\mathbf{A}^r_K$ & 
% ${\mathscr R}_K^{\mathrm{int}, r}$     & 
 $\mathbf{A}^{\dagger}_{{r'}^{-}\!, K}, (\mathbf{A}^{\dagger}_{{r'}^{-}})^{H_K}$   &
 $\mathbf{A}^{\dagger, r'}_K, \mathbf{A}^{H_K}_{[r'; +\infty]}$ & 
 $\mathbf{A}^{(0,r]}_K$ & 
 ---
 \\
 $\mathbf{A}^{\dagger}_K$ &
 %  ${\mathscr R}_K^{\mathrm{int}}$ & 
 --- &
 --- &
 %$\mathbf{A}^{\dagger}_K$ & 
 --- &
% $\mathbf{A}^{\dagger}_K$ & 
 --- 
 \\
 \hline
 $\mathbf{B}_K$ & 
%   --- & 
 --- &
 --- &
 % $\mathbf{B}_K$ & 
 --- &
 %$\mathbf{B}_K$ & 
 --- 
 \\
 $\mathbf{B}^r_K$ 
% &   ${\mathscr E}_K^{(0, r]}$, ${\mathscr R}_K^{\mathrm{bd}, r}$ %, $({\mathcal R}^{\mathrm{bd}, r})^{H_K}$  
   & 
 $\mathbf{B}^{\dagger}_{{r'}^{-}\!, K}, (\mathbf{B}^{\dagger}_{{r'}^{-}})^{H_K}$ &
 $\mathbf{B}^{\dagger, r'}_K,\mathbf{B}^{H_K}_{[r'; +\infty]}$ 
 & $\mathbf{B}^{(0,r]}_K$ & 
 ---  
 \\
 $\mathbf{B}^{\dagger}_K$ &
 %  ${\mathscr E}_K^{\dagger}$, ${\mathscr R}_K^{\mathrm{bd}}$  & 
 --- &
 --- &
 %$\mathbf{B}^{\dagger}_K$ & 
 --- & 
 %$\mathbf{B}^{\dagger}_K$ & 
 --- 
 \\
 \hline
 $\mathbf{C}^{[s,r]}_K$ & 
  % n/a & 
 n/a &
 $\mathbf{B}_{[r' ;s']}^{H_K}$ 
 & 
 n/a & 
 ---
 \\
 $\mathbf{C}^r_K$ & 
 %  ${\mathscr E}_K^{]0, r]}$, ${\mathscr R}_K^{r}$  & 
 n/a &
 $\mathbf{B}^{\dagger, r'}_{\mathrm{rig}, K},\mathbf{B}^{H_K}_{[r'; +\infty[}$ 
 & $\mathbf{B}^{]0,r]}_K$ 
 & 
 ---
 \\
 $\mathbf{C}^{\dagger}_K$ & 
%    ${\mathscr R}_K $  & 
 n/a &
 $\mathbf{B}^{\dagger}_{\mathrm{rig}, K}$ & $\mathbf{B}^{\dagger}_{\mathrm{rig}, K}$ & $\mathbf{C}_K$ 
% \\  0 & 1 & 2 & 3 & 4 & 5
 \\ 
 \hline
\end{tabular}}
} %end scalebox
\end{center}

%\end{landscape}
\vspace{0.2in}

Note that a similar table exists without the subscript $K$. Let us also point out that Kedlaya adopts a completely different set of notation in \cite{Kedlaya_p-adic_monodormy}. 

For relative period rings, we use
$\widetilde{\mathbf{A}}^r_X, \widetilde{\mathbf{A}}^{\dagger}_X, \widetilde{\mathbf{B}}^r_X, \widetilde{\mathbf{B}}^{\dagger}_X, \widetilde{\mathbf{C}}^{[s,r]}_X, \widetilde{\mathbf{C}}^r_X, \widetilde{\mathbf{C}}^{\dagger}_X$ 
to denote perfect period rings, as opposed to the notation $\widetilde{\mathbf{A}}^r_{\psi}, \widetilde{\mathbf{A}}^{\dagger}_{\psi}, \widetilde{\mathbf{B}}^r_{\psi}, \widetilde{\mathbf{B}}^{\dagger}_{\psi}, \widetilde{\mathbf{C}}^{[s,r]}_{\psi}, \widetilde{\mathbf{C}}^r_{\psi}, \widetilde{\mathbf{C}}_{\psi}$ used in \cite{Kedlaya-LiuII}. In fact, in both \cite{Kedlaya-LiuI} and \cite{Kedlaya-LiuII}, Kedlaya and Liu use $\widetilde{\mathbf{A}}^r_X, \widetilde{\mathbf{A}}^{\dagger}_X, \widetilde{\mathbf{B}}^r_X, \widetilde{\mathbf{B}}^{\dagger}_X, \widetilde{\mathbf{C}}^{[s,r]}_X, \widetilde{\mathbf{C}}^r_X, \widetilde{\mathbf{C}}_X$ to refer to sheaves of period rings.

\vspace{15mm}

\bibliographystyle{amsalpha}
\bibliography{library}

\vspace{15mm}

\begin{tabular}{l}
    H.D.\\
    Yau Mathematical Sciences Center \& Department of Mathematical Sciences, 
    Tsinghua University   \\
    Beijing, China\\
    \textit{E-mail address: }\texttt{hdiao@mail.tsinghua.edu.cn }\\
    \\
Y.S.M.\\
Beijing Institute of Mathematical Sciences and Applications\\
Beijing, China\\
\textit{E-mail address: }\texttt{ysmoon@bimsa.cn }\\
\\
Z.Y.\\
Department of Mathematics, 
University of California Santa Barbara\\
Santa Barbara, California, USA\\
\textit{E-mail address: }\texttt{yao@math.ucsb.edu }
\end{tabular}

\end{document}